%% file: kh_correspondence.tex
\documentclass[a4paper,10pt, english]{smfbook}
\usepackage{amsmath,amscd,amssymb}

\input{notation}
\input{new_theorem}

\author{Takuro Mochizuki}
\address{
Department of Mathematics, Kyoto University, Kyoto 606-8502, Japan}
\email{takuro@math.kyoto-u.ac.jp}

\title[Kobayashi-Hitchin correspondence]
 {Kobayashi-Hitchin correspondence
 for \\tame harmonic bundles\\ and an application}

\begin{document}
\frontmatter

\begin{abstract}
\input{abstract}
\end{abstract}


\thanks{This work was prepared with the partial support
 of Ministry of Education, Culture, Sports, Science and Technology.}

\maketitle

\tableofcontents

\mainmatter

\chapter{Introduction}
\input{1}

\chapter{Preliminary}
\label{chapter;05.10.1.2}
\input{2}

\chapter{Parabolic Higgs Bundle and Regular Filtered Higgs Bundle}
\label{chapter;05.9.8.1}
\input{3}

\chapter{An Ordinary Metric for a Parabolic Higgs Bundle}
\label{chapter;05.9.8.10}
\input{4}

\chapter[Associated Parabolic Higgs bundle]
 {Parabolic Higgs Bundle Associated to Tame Harmonic Bundle}
\label{chapter;05.10.1.4}
\input{5}

\chapter[Preliminary correspondence]
{Preliminary Correspondence and
 Bogomolov-Gieseker Inequality}
\label{chapter;05.10.1.5}
\input{6}

\chapter{Construction of a Frame}
\label{chapter;06.8.21.100}
\input{7}

\chapter{Some Convergence Results}
\label{chapter;06.8.21.101}
\input{8}

\chapter{Existence of Adapted Pluri-Harmonic Metric}
\label{chapter;06.8.21.102}
\input{9}

\chapter[The Deformation of Representations]
 {Torus Action and the Deformation of Representations}
\label{chapter;05.10.1.11}
\input{10}

\chapter{$G$-Harmonic Bundle (Appendix)}
\label{chapter;04.10.26.90}
\input{11}

\backmatter

\input{kh_correspondenceref}
\end{document}

%% file: notation.tex
\newcommand{\nbiga}{\mathcal{A}}
\newcommand{\nbigb}{\mathcal{B}}

\newcommand{\nbigd}{\mathcal{D}}
\newcommand{\nbige}{\mathcal{E}}
\newcommand{\nbigf}{\mathcal{F}}
\newcommand{\nbigg}{\mathcal{G}}

\newcommand{\nbigi}{\mathcal{I}}
\newcommand{\nbigj}{\mathcal{J}}
\newcommand{\nbigk}{\mathcal{K}}
\newcommand{\nbigl}{\mathcal{L}}

\newcommand{\nbign}{\mathcal{N}}
\newcommand{\nbigo}{\mathcal{O}}
\newcommand{\nbigp}{\mathcal{P}}

\newcommand{\nbigs}{\mathcal{S}}
\newcommand{\nbigt}{\mathcal{T}}
\newcommand{\nbigu}{\mathcal{U}}
\newcommand{\nbigv}{\mathcal{V}}
\newcommand{\nbigw}{\mathcal{W}}
\newcommand{\nbigx}{\mathcal{X}}

\newcommand{\nbigz}{\mathcal{Z}}

\newcommand{\proj}{\mathbb{P}}
\newcommand{\seisuu}{\mathbb{Z}}
\newcommand{\rnum}{{\boldsymbol Q}}

\newcommand{\cnum}{{\boldsymbol C}}
\newcommand{\real}{{\boldsymbol R}}
\newcommand{\hyperh}{\mathbb{H}}

\newcommand{\DD}{\mathbb{D}}
\newcommand{\EE}{\mathbb{E}}

\newcommand{\gbigd}{\mathfrak D}

\newcommand{\gbigx}{\mathfrak X}

\newcommand{\gminig}{\mathfrak g}
\newcommand{\gminih}{\mathfrak h}

\newcommand{\gminil}{\mathfrak l}
\newcommand{\gminim}{\mathfrak m}

\newcommand{\gminiu}{\mathfrak u}

\newcommand{\vecv}{{\boldsymbol v}}
\newcommand{\vecu}{{\boldsymbol u}}
\newcommand{\vecw}{{\boldsymbol w}}

\newcommand{\vecl}{{\boldsymbol l}}

\newcommand{\vecalpha}{{\boldsymbol \alpha}}
\newcommand{\veca}{{\boldsymbol a}}
\newcommand{\vecb}{{\boldsymbol b}}

\newcommand{\vecc}{{\boldsymbol c}}

\newcommand{\veck}{{\boldsymbol k}}
\newcommand{\vecm}{{\boldsymbol m}}

\newcommand{\vecF}{{\boldsymbol F}}
\newcommand{\vecE}{{\boldsymbol E}}

\newcommand{\lrarr}{\longrightarrow}

\newcommand{\pf}{{\bf Proof}\hspace{.1in}}

\def\Hom{\mathop{\rm Hom}\nolimits}
\def\End{\mathop{\rm End}\nolimits}

\def\Cok{\mathop{\rm Cok}\nolimits}

\def\Image{\mathop{\rm Im}\nolimits}

\def\Gr{\mathop{\rm Gr}\nolimits}

\def\GL{\mathop{\rm GL}\nolimits}

\def\rank{\mathop{\rm rank}\nolimits}
\def\Spec{\mathop{\rm Spec}\nolimits}

\def\Ker{\mathop{\rm Ker}\nolimits}
\def\ker{\mathop{\rm Ker}\nolimits}

\def\Gr{\mathop{\rm Gr}\nolimits}
\def\Sym{\mathop{\rm Sym}\nolimits}

\def\ad{\mathop{\rm ad}\nolimits}
\def\Res{\mathop{\rm Res}\nolimits}

\def\degpar{\mathop{\rm par\textrm{-}deg}\nolimits}
\def\parch{\mathop{\rm par\textrm{-}ch}\nolimits}
\def\ch{\mathop{\rm ch}\nolimits}
\def\parchern{\mathop{\rm par\textrm{-}c}\nolimits}

\def\tr{\mathop{\rm tr}\nolimits}

\def\vol{\mathop{\rm dvol}\nolimits}
\def\dvol{\mathop{\rm dvol}\nolimits}

\def\id{\mathop{\rm id}\nolimits}

\def\codim{\mathop{\rm codim}\nolimits}
\def\gap{\mathop{\rm gap}\nolimits}

\def\Irrep{\mathop{\rm Irrep}\nolimits}
\def\Irr{\mathop{\rm Irr}\nolimits}

\def\wt{\mathop{\rm wt}\nolimits}

\newcommand{\del}{\partial}
\newcommand{\delbar}{\overline{\del}}

\newcommand{\pardeg}{\degpar}

\newcommand{\nhom}{{\mathcal Hom}}

\newcommand{\itibar}{\underline{1}}
\newcommand{\nibar}{\underline{2}}

\newcommand{\harmonicbundle}{(E,\delbar_E,\theta,h)}

\newcommand{\barz}{\overline{z}}
\newcommand{\zbar}{\barz}
\newcommand{\zetabar}{\overline{\zeta}}
\newcommand{\baralpha}{\overline{\alpha}}
\newcommand{\alphabar}{\baralpha}

\newcommand{\sbar}{\overline{s}}

\newcommand{\prolong}[1]{{}^{\diamond}{#1}}

\newcommand{\prolongg}[2]{{}_{#1}{#2}}
\newcommand{\leftupdown}[2]{{}^{#1}_{#2}}

\newcommand{\Par}{{\mathcal Par}}

\newcommand{\lefttop}[1]{{}^{#1}}

\def\Pic{\mathop{\rm Pic}\nolimits}
\def\Tor{\mathop{\rm Tor}\nolimits}

\newcommand{\closedopen}[2]{[#1,#2[}
\newcommand{\openclosed}[2]{]#1,#2]}
\newcommand{\openopen}[2]{]#1,#2[}

\newcommand{\unitary}{\gminiu}
\newcommand{\gl}{\gminig\gminil}

\newcommand{\omegatilde}{\widetilde{\omega}}
\newcommand{\vecEhat}{\widehat{\vecE}}
\newcommand{\thetahat}{\widehat{\theta}}
\newcommand{\Ehat}{\widehat{E}}
\newcommand{\hhat}{\widehat{h}}

\newcommand{\rhotilde}{\widetilde{\rho}}
\newcommand{\gun}{G}
\newcommand{\xtilde}{\widetilde{x}}

\newcommand{\Etilde}{\widetilde{E}}
\newcommand{\vecEtilde}{\widetilde{\vecE}}
\newcommand{\thetatilde}{\widetilde{\theta}}
\newcommand{\rhobar}{\overline{\rho}}
\newcommand{\kappatilde}{\widetilde{\kappa}}
\newcommand{\qtilde}{\widetilde{q}}
\newcommand{\atilde}{\widetilde{a}}
\newcommand{\stilde}{\widetilde{s}}
\newcommand{\alphatilde}{\widetilde{\alpha}}
\newcommand{\ftilde}{\widetilde{f}}
\newcommand{\htilde}{\widetilde{h}}
\newcommand{\wbar}{\overline{w}}

\newcommand{\vecvtilde}{\widetilde{\vecv}}
\newcommand{\gtilde}{\widetilde{g}}

\newcommand{\bbar}{\overline{b}}
\newcommand{\gminimbar}{\overline{\gminim}}
\newcommand{\Fbar}{\overline{F}}
\newcommand{\epsilonbar}{\overline{\epsilon}}
\newcommand{\Wtilde}{\widetilde{W}}
\newcommand{\Ghat}{\widehat{G}}

\newcommand{\Gammatilde}{\widetilde{\Gamma}}
\newcommand{\Qtilde}{\widetilde{Q}}

\newcommand{\Phitilde}{\widetilde{\Phi}}
\newcommand{\Ftilde}{\widetilde{F}}
\newcommand{\Stilde}{\widetilde{S}}

\newcommand{\Xtilde}{\widetilde{X}}
\newcommand{\Dtilde}{\widetilde{D}}
\newcommand{\vtilde}{\widetilde{v}}
\newcommand{\Ptilde}{\widetilde{P}}
\newcommand{\Utilde}{\widetilde{U}}
\newcommand{\vecatilde}{\widetilde{\veca}}
\newcommand{\pitilde}{\widetilde{\pi}}
\newcommand{\nbigdstar}{\nbigd^{\star}}
\newcommand{\nbigutilde}{\widetilde{\nbigu}}
\newcommand{\zerotilde}{\widetilde{0}}
\newcommand{\ttilde}{\widetilde{t}}
\newcommand{\ctilde}{\widetilde{c}}
\newcommand{\nbigitilde}{\widetilde{\nbigi}}
\newcommand{\Ytilde}{\widetilde{Y}}
\newcommand{\Hbar}{\overline{H}}

%% file: new_theorem.tex
\newtheorem{thm}{Theorem}[chapter]
\newtheorem{cor}[thm]{Corollary}
\newtheorem{lem}[thm]{Lemma}
\newtheorem{prop}[thm]{Proposition}

\theoremstyle{definition}
\newtheorem{df}[thm]{Definition}

\newtheorem{condition}[thm]{Condition}
\newtheorem{assumption}[thm]{Assumption}

\theoremstyle{remark}
\newtheorem{rem}[thm]{Remark}

%% file: abstract.tex
We establish the correspondence between tame harmonic bundles
and $\mu_L$-polystable parabolic Higgs bundles
with trivial characteristic numbers.
We also show the Bogomolov-Gieseker type inequality
for $\mu_L$-stable parabolic Higgs bundles.

Then we show that any local system on a smooth quasi projective variety
can be deformed to a variation of polarized Hodge structure.
As a consequence, we can conclude that
some kind of discrete groups cannot be a split quotient
of the fundamental group of a smooth quasi projective variety.

%% file: 1.tex
\section{Background}

\subsection{Kobayashi-Hitchin correspondence}

We briefly recall some aspects of the so-called
Kobayashi-Hitchin correspondence.
(See the introduction of \cite{lubke-teleman}
for more detail.)
In 1960's, 
M. S. Narasimhan and C. S. Seshadri 
proved the correspondence
between irreducible flat unitary bundles and
stable vector bundles with degree $0$,
on a compact Riemann surface
(\cite{narasimhan-seshadri}).
Clearly, it was desired to extend their result
to the higher dimensional case and the non-flat case.

In early 1980's, S. Kobayashi introduced the Einstein-Hermitian condition
for holomorphic bundles on Kahler manifolds (\cite{koba3}, \cite{koba2}).
He and M. L\"{u}bke (\cite{lubke}) 
proved that the existence of Einstein-Hermitian metric implies
the polystability of the underlying holomorphic bundle.
S. K. Donaldson pioneered the way for the inverse problem
(\cite{don4} and \cite{don3}).
He attributed the problem to Kobayashi and N. Hitchin.
The definitive result was given by K. Uhlenbeck, S. T. Yau
and Donaldson (\cite{uy} and \cite{don}).
We also remark that V. Mehta and A. Ramanathan (\cite{mehta-ramanathan2})
proved the correspondence in the case where the Chern class is trivial,
i.e.,
the correspondence of flat unitary bundles
and stable vector bundles with trivial Chern classes.

On the other hand, it was quite fruitful to
consider the correspondences
for vector bundles with some additional structures like Higgs fields,
which was initiated by Hitchin (\cite{hitchin}).
He studied the Higgs bundles on a compact Riemann surface
and the moduli spaces.
His work has influenced various fields of mathematics.
It involves a lot of subjects and ideas,
and one of his results is the correspondence
of the stability and the existence of Hermitian-Einstein metrics
for Higgs bundles on a compact Riemann surface.

\subsection{A part of C. Simpson's work}
\label{subsection;04.10.26.10}

C. Simpson studied the Higgs bundles over
higher dimensional complex manifolds,
influenced by the work of Hitchin, but motivated by his own subject:
Variation of Polarized Hodge Structure.
He made great innovations in various areas of algebraic geometry.
Here, we recall just a part of his huge work.

Let $X$ be a smooth irreducible projective variety
over the complex number field,
and $E$ be an algebraic vector bundle on $X$.
Let $(E,\theta)$ be a Higgs bundle,
i.e.,
$\theta$ is a holomorphic section of $\End(E)\otimes\Omega^{1,0}_X$
satisfying $\theta^2=0$.
The ``stability'' and the ``Hermitian Einstein metric''
are naturally defined for Higgs bundles,
and Simpson proved that
there exists a Hermitian-Einstein metric 
of $(E,\theta)$
if and only if $(E,\theta)$ is polystable.
In the special case where the Chern class of the vector bundle
is trivial,
the Hermitian-Einstein metric gives the pluri-harmonic metric.
Together with the result of K. Corlette
who is also a great progenitor of the study of harmonic bundles
(\cite{corlette}),
Simpson obtained the Trinity on a smooth projective variety:
{\small
\begin{equation} \label{eq;04.10.14.1}
 \fbox{
\begin{tabular}{c}
 {\bf Algebraic Geometry}\\
polystable Higgs bundle\\
(trivial Chern class)
\end{tabular}
}
\leftrightarrow
\fbox{
 \begin{tabular}{c}
 {\bf Differential Geometry}\\
 harmonic bundle
 \end{tabular}
}
\leftrightarrow
\fbox{\begin{tabular}{c}
 {\bf Topology}\\
 semisimple \\ local system
 \end{tabular}
 }
\end{equation}
}

If $(E,\theta)$ is a stable Higgs bundle,
then $(E,\alpha\cdot\theta)$ is also a stable Higgs bundle.
Hence we obtain the family of stable Higgs bundles
$\bigl\{(E,\alpha\cdot\theta)\,\big|\,\alpha\in\cnum^{\ast}\bigr\}$.
Correspondingly, we obtain the family of 
flat bundles $\bigl\{L_{\alpha}\,\big|\,\alpha\in\cnum^{\ast}\bigr\}$.
Simpson showed that
we obtain the variation of polarized Hodge structure
as a limit $\lim_{\alpha\to 0}L_{\alpha}$.
In particular, it can be concluded that
any flat bundle can be deformed to a variation of polarized Hodge
structure.
As one of the applications,
he obtained the following remarkable result (\cite{s5}):
\begin{thm}[Simpson]
 \label{thm;04.10.26.10}
Let $\Gamma$ be a rigid discrete subgroup of a real algebraic group
which is not of Hodge type.
Then $\Gamma$ cannot be a split quotient of
the fundamental group of a smooth irreducible projective variety.
\hfill\qed
\end{thm}

There are classical known results on the rigidity of subgroups
of Lie groups.
The examples of rigid discrete subgroups
can be found in 4.7.1--4.7.4 in the 53 page of \cite{s5}.
The classification of real algebraic group of Hodge type
was done by Simpson. 
The examples of real algebraic group which is not of Hodge type
can be found in the 50 page of \cite{s5}.
As a corollary, he obtained the following.
\begin{cor}
$SL(n,\seisuu)$ $(n\geq 3)$ cannot be a split quotient
of the fundamental group of a smooth irreducible projective variety.
\hfill\qed
\end{cor}

\section{Main Purpose}

\subsection{Kobayashi-Hitchin correspondence for parabolic Higgs bundles}
\label{subsection;05.10.4.1}

It is an important and challenging problem
to generalize the correspondence (\ref{eq;04.10.14.1})
to the quasiprojective case from the projective case.
As for the correspondence of harmonic bundles and
semisimple local systems,
an excellent result was obtained by J. Jost and K. Zuo \cite{JZ2},
which says there exists a tame pluri-harmonic metric
on any semisimple local system over a quasiprojective variety.
The metric is called the Corlette-Jost-Zuo metric.

In this paper, we restrict ourselves to the correspondence
between Higgs bundles and harmonic bundles
on a quasiprojective variety $Y$.
More precisely,
we should consider not Higgs bundles on $Y$
but {\em parabolic} Higgs bundles on $(X,D)$,
where $(X,D)$ is a pair of a smooth irreducible projective variety
and a normal crossing divisor such that $Y=X-D$.
Such a generalization has been studied by several people.
In the non-Higgs case, J. Li \cite{li2} and
B. Steer-A. Wren \cite{steer-wren} established the correspondence.
In the Higgs case,
Simpson established the correspondence
in the one dimensional case \cite{s2},
and O. Biquard established it in the case
where $D$ is smooth \cite{b}.
\begin{rem}
Their results also include the correspondence in the case
where the characteristic numbers are non-trivial.
\hfill\qed
\end{rem}

For applications, however, it is desired that
the correspondence for parabolic Higgs bundles should be given
in the case where $D$ is not necessarily smooth,
which we would like to discuss in this paper.

We explain our result more precisely.
Let $X$  be a smooth irreducible
projective variety over the complex number field
provided an ample line bundle $L$.
Let $D$ be a simple normal crossing divisor of $X$.
The main purpose of this paper is to establish the correspondence
between tame harmonic bundles and $\mu_L$-parabolic Higgs bundles
whose characteristic numbers vanish.
(See Chapter \ref{chapter;05.9.8.1} for the meaning of the words.)

\begin{thm}
[Proposition  \ref{prop;05.9.8.30}--\ref{prop;05.8.29.500}, 
 and  Theorem \ref{thm;04.10.24.1}]
\label{thm;05.9.13.1} \label{thm;04.10.24.10}
Let $\bigl(\vecE_{\ast},\theta\bigr)$
be a regular filtered Higgs bundle on $(X,D)$,
and we put $E:=\vecE_{|X-D}$.
 It is $\mu_L$-polystable with
 trivial characteristic numbers,
 if and only if
 there exists a pluri-harmonic metric $h$ of $(E,\theta)$ on $X-D$
 which is adapted to the parabolic structure.
 Such a metric is unique up to an obvious ambiguity.
\hfill\qed
\end{thm}

\begin{rem}
Regular Higgs bundles and parabolic Higgs bundles are equivalent.
See Chapter {\rm \ref{chapter;05.9.8.1}}.
\hfill\qed
\end{rem}

\begin{rem}
More precisely on the existence result,
we can show the existence of
the adapted pluri-harmonic metric
for $\mu_L$-stable
reflexive saturated
regular filtered Higgs sheaf on $(X,D)$
with trivial characteristic numbers.
(See Sections \ref{section;05.8.23.1}--\ref{section;06.8.17.50}
for the definition.)
Then, due to our previous result in {\rm \cite{mochi2}},
it is a regular filtered Higgs bundle on $(X,D)$, in fact.
\hfill\qed
\end{rem}

We are mainly interested in the $\mu_L$-stable parabolic Higgs bundles
whose characteristic numbers vanish.
But we also obtain the following theorem
on more general $\mu_L$-stable parabolic Higgs bundles.
\begin{thm}
 [Theorem \ref{thm;05.7.30.30}]
\label{thm;05.9.13.2}
Let $X$ be a smooth irreducible
projective variety of an arbitrary dimension,
and $D$ be a simple normal crossing divisor.
Let $L$ be an ample line bundle on $X$.
Let $(\vecE_{\ast},\theta)$ be a 
$\mu_L$-stable regular Higgs bundle
in codimension two on $(X,D)$.
Then the following inequality holds:
\[
 \int_X\parch_{2,L}(\vecE_{\ast})
-\frac{\int_X\parchern_{1,L}^2(\vecE_{\ast})}{2\rank E}
\leq 0.
\]
Such an inequality is called
Bogomolov-Gieseker inequality.
\hfill\qed
\end{thm}

\subsection{Strategy for the proof of Bogomolov-Gieseker inequality}

We would like to explain our strategy for the proof of
the main theorems. 
First we describe an outline for
Bogomolov-Gieseker inequality (Theorem \ref{thm;05.9.13.2}),
which is much easier.
We have only to consider the case $\dim X=2$.
Essentially, it consists of the following two parts.
\begin{description}
\item[(1) The correspondence in the graded semisimple case]\mbox{{}}\\
 We establish the Kobayashi-Hitchin correspondence 
 for {\em graded semisimple} parabolic Higgs bundles.
 In particular, we obtain the Bogomolov-Gieseker inequality
 in this case.
\item[(2) Perturbation of the parabolic structure and taking the limit]\mbox{{}}\\
 Let $(\prolongg{\vecc}{E},\vecF,\theta)$ be a given $\vecc$-parabolic
$\mu_L$-stable Higgs bundle, which is not necessarily
 graded semisimple.
 For any small positive number $\epsilon$,
 we take a perturbation $\vecF^{(\epsilon)}$ of $\vecF$
 such that
 $(\prolongg{\vecc}{E},\vecF^{(\epsilon)},\theta)$
 is a graded semisimple $\mu_L$-stable parabolic Higgs bundle.
 Then the Bogomolov-Gieseker inequality holds
 for $(\prolongg{\vecc}{E},\vecF^{(\epsilon)},\theta)$.
 By taking a limit for $\epsilon\lrarr 0$,
 we obtain the Bogomolov-Gieseker inequality
 for the given $(\prolongg{\vecc}{E},\vecF,\theta)$.
\end{description}

Let us describe for more detail.

\noindent
{\bf (1)}\,\,
In \cite{s5}, 
Simpson constructed a Hermitian-Einstein metric
for Higgs bundle by the following process:
\begin{description}
\item[(i)] Take an appropriate initial metric.
\item[(ii)] Deform it along the heat equation.
\item[(iii)]Take a limit, and then we obtain the Hermitian-Einstein metric.
\end{description}
If the base space is compact, 
the steps (ii) and (iii) are the main issues,
and the step (i) is trivial.
Actually, Simpson also discussed the case where 
the base Kahler manifold is non-compact,
and he showed the existence of a Hermitian-Einstein metric
if we can take an initial metric 
whose curvatures satisfy some finiteness condition.
(See Section \ref{section;05.10.1.1}
for more precise statements.)
So, for a $\mu_L$-stable $\vecc$-parabolic Higgs bundle 
$(\prolongg{\vecc}{E},\vecF,\theta)$ on $(X,D)$,
where $X$ is a smooth projective surface
and $D$ is a simple normal crossing divisor,
ideally,
we would like to take an initial metric of
$E:=\prolongg{\vecc}{E}_{|X-D}$
adapted to the parabolic structure.
But, it is rather difficult, and the author is not sure
whether such a good metric can always be taken
for any parabolic Higgs bundles.
It seems one of the main obstacles 
to establish the Kobayashi-Hitchin correspondence
for parabolic Higgs bundles.

However, we can easily take such a good initial metric,
if we assume the vanishing of the nilpotent part of the residues
of the Higgs field on the graduation of the parabolic filtration.
Such a parabolic Higgs bundle will be called {\em graded semisimple}
in this paper.
We first establish the correspondence in this easy case.
(Proposition \ref{prop;05.7.30.10}).

\noindent
{\bf (2)}\,\,
Let $(\prolongg{\vecc}{E},\vecF,\theta)$ be a $\mu_L$-stable
$\vecc$-parabolic Higgs bundle
on $(X,D)$, where $\dim X=2$.
We take a perturbation of $\vecF^{(\epsilon)}$
as in Section \ref{section;05.7.30.15}.
In particular, $(\prolongg{\vecc}{E},\vecF^{(\epsilon)},\theta)$
is a $\mu_L$-stable graded semisimple $\vecc$-parabolic Higgs bundle,
and the following holds:
\[
 \parchern_1(\prolongg{\vecc}{E},\vecF)
=\parchern_1(\prolongg{\vecc}{E},\vecF^{(\epsilon)}),
\]
\[
 \left|
 \int_X \parch_{2}(\prolongg{\vecc}{E},\vecF)
-\int_X\parch_2(\prolongg{\vecc}{E},\vecF^{(\epsilon)})
\right|\leq C\cdot \epsilon.
\]
Then we obtain the Bogomolov-Gieseker inequality for
 $(\prolongg{\vecc}{E},\vecF^{(\epsilon)},\theta)$
by using the Hermitian-Einstein metric
obtained in (1).
By taking the limit $\epsilon\to 0$,
we obtain the desired inequality for the given
$(\prolongg{\vecc}{E},\vecF,\theta)$.

\subsection{Strategy for the proof of Kobayashi-Hitchin correspondence}

Let $X$ be a smooth projective surface,
and $D$ be a simple normal crossing divisor.
Let $L$ be an ample line bundle on $X$,
and $\omega$ be the Kahler form representing $c_1(L)$.
Roughly speaking,
the correspondence on $(X,D)$ as in Theorem \ref{thm;05.9.13.1}
can be divided into the following two parts:
\begin{itemize}
\item
 For a given tame harmonic bundle $(E,\delbar_E,\theta,h)$ on $X-D$,
 we obtain the $\mu_L$-polystable parabolic Higgs bundle
 $(\prolongg{\vecc}{E},\vecF,\theta)$
 with the trivial characteristic numbers.
\item
 On the converse,
 we obtain a pluri-harmonic metric of $(E,\delbar_E,\theta)$ on $X-D$
 for such $(\prolongg{\vecc}{E},\vecF,\theta)$.
\end{itemize}

As for the first issue,
most problem can be reduced to the one dimensional case,
which was established by Simpson \cite{s2}.
However, we have to show the vanishing of the characteristic numbers,
for which our study of the asymptotic behaviour of tame harmonic bundles
(\cite{mochi2}) is useful.

\vspace{.1in}

As for the second issue,
we use the perturbation method, again.
Namely,
let $(\prolongg{\vecc}{E},\vecF,\theta)$ be
a $\mu_L$-stable $\vecc$-parabolic Higgs bundle on $(X,D)$.
Take a perturbation $\vecF^{(\epsilon)}$ of the filtration
$\vecF$
for a small positive number $\epsilon$.
We also take metrics appropriate $\omega_{\epsilon}$ of $X-D$
such that $\lim_{\epsilon\to 0}\omega_{\epsilon}=\omega$,
and then we obtain Hermitian-Einstein metrics $h_{\epsilon}$
for the Higgs bundle $(E,\delbar_{E},\theta)$ on $X-D$
with respect to $\omega_{\epsilon}$, 
which is adapted to the parabolic structure $\vecF^{(\epsilon)}$.
Ideally, we would like to consider the limit
$\lim_{\epsilon\to 0}h_{\epsilon}$,
and we expect that the limit gives the Hermitian-Einstein metric $h$
for $(E,\delbar_E,\theta)$ with respect to $\omega$,
which is adapted to the given filtration $\vecF$.
Perhaps, it may be correct,
but it does not seem easy to show, in general.

We restrict ourselves to the simpler case
where the characteristic numbers of
$(\prolongg{\vecc}{E},\vecF,\theta)$ are trivial.
Under this assumption, we show such a convergence.
More precisely,
we show that there is a subsequence $\{\epsilon_i\}$
such that 
$\bigl\{(E,\delbar_{E},h_{\epsilon_i},\theta)\bigr\}$
converges to a harmonic bundle
$(E',\delbar_{E'},\theta',h')$ on $X-D$,
and we show that the given
$(\prolongg{\vecc}{E},\vecF,\theta)$
is isomorphic to the parabolic Higgs bundles 
obtained from $(E',\delbar_{E'},\theta',h')$.

\begin{rem}
We obtained a similar correspondence for
$\lambda$-connections in {\rm\cite{mochi5}}.
Although the argument is essentially same,
we need some additional argument
in the case of $\lambda$-connections.
\hfill\qed
\end{rem}

\section{Additional Results}

\subsection{Torus action and the deformation of a $G$-flat bundle}

Once Theorem \ref{thm;05.9.13.1} is established,
we can use some of the arguments for the applications given
in the projective case.
For example, we can deform any flat bundle to 
the one which comes from a variation of polarized Hodge structure.
We follow the well known framework given by Simpson
with a minor modification.
We briefly recall it,
and we will mention the problem that we have to care about
in the process.

Let $X$ be a smooth irreducible projective variety,
and $D$ be a simple normal crossing divisor
with the irreducible decomposition $D=\bigcup_{i\in S}D_i$.
Let $x$ be a point of $X-D$.
Let $\Gamma$ denote the fundamental group $\pi_1(X-D,x)$.
Any representation of $\Gamma$ can be 
deformed to a semisimple representation,
and hence we start with a semisimple one.

Let $(E,\nabla)$ be a flat bundle over $X-D$
such that the induced representation
$\rho:\Gamma\lrarr \GL(E_{|x})$ is semisimple.
Recall we can take a Corlette-Jost-Zuo metric of $(E,\nabla)$,
as mentioned in Subsection \ref{subsection;05.10.4.1}.
Hence we obtain a tame pure imaginary harmonic
bundle $(E,\delbar_E,\theta,h)$ on $X-D$,
and the induced $\mu_L$-polystable $\vecc$-parabolic Higgs bundle 
$(\prolongg{\vecc}{E},\vecF,\theta)$ on $(X,D)$,
where $\vecc$ denotes any element of $\real^S$.
We have the canonical decomposition
$(\prolongg{\vecc}{E},\vecF,\theta)=
 \bigoplus_{i}(\prolongg{\vecc}{E}_i,\vecF_i,\theta_i)^{\oplus\,m_i}$,
where each $(\prolongg{\vecc}{E}_i,\vecF_i,\theta_i)$ is $\mu_L$-stable.

Let us consider the family of $\vecc$-parabolic Higgs bundles
$\bigl(\prolongg{\vecc}{E},\vecF,t\cdot\theta\bigr)$
for $t\in\cnum^{\ast}$,
which are $\mu_L$-polystable.
Due to the standard Langton's trick \cite{langton},
we have the semistable $\vecc$-parabolic Higgs sheaves
$(\prolongg{\vecc}{\widetilde{E}_i},\widetilde{\vecF}_i,\widetilde{\theta_i})$
which are limits of
$(\prolongg{\vecc}{E}_i,\vecF_i,t\cdot\theta_i\bigr)$
in $t\to 0$.
On the other hand,
we can take a pluri-harmonic metric
$h_t$ of the Higgs bundle $(E,\delbar_{E},t\cdot \theta)$ on $X-D$
for each $t$,
which is adapted to the parabolic structure.
(Theorem \ref{thm;05.9.13.1}).
 Then we obtain the family of flat bundles
 $(E,\DD^{1}_t)$,
 and the associated family of the representations
$\bigl\{\rho_t:\Gamma\lrarr \GL(E_{|x})
 \,\big|\,t\in \cnum^{\ast}\bigr\}$.
Since $(E,\delbar_E,t\cdot\theta,h_t)$ is tame pure imaginary
in the case $t\in\real_{>0}$,
the representations $\rho_t$ are semisimple.
The family $\{\rho_t\,\big|\,t\in\cnum^{\ast}\}$ should be continuous
with respect to $t$,
and the limit $\lim_{t\to 0}\rho_t$ should exist, ideally.
We formulate the continuity of $\rho_t$ with respect to $t$
and the convergence of $\rho_t$ in $t\to 0$,
as follows.
Let $V$ be a $\cnum$-vector space such that $\rank(V)=\rank(E)$.
Let $h_V$ denote the metric of $V$,
and let $U(h_V)$ denote the unitary group for $h_V$.
We put $R(\Gamma,V):=Hom(\Gamma,\GL(V))$.
By the conjugate, $U(h_V)$ acts on the space $R(\Gamma,V)$.
Let $M(\Gamma,V,h_V)$ denote the usual quotient space.
Let $\pi_{\GL(V)}:R(\Gamma,V)\lrarr M(\Gamma,V,h_V)$ denote the projection.

By taking any isometry $(E_{|x},h_{t|x})\simeq (V,h_V)$,
we obtain the representation $\rho_t':\Gamma\lrarr \GL(V)$.
We put $\nbigp(t):=\pi_{\GL(V)}(\rho_t')$,
and we obtain the map $\nbigp:\cnum^{\ast}\lrarr M(\Gamma,V,h_V)$.
It is well defined.
Then, we obtain the following partial result.

\begin{prop}
 [Theorem \ref{thm;05.9.10.5}, Lemma \ref{lem;05.9.11.3}, 
 Proposition \ref{prop;04.10.23.101}]
 \label{prop;05.6.3.1} 
\begin{enumerate}
\item
 The induced map 
 $\nbigp$  is continuous.
\item
 $\nbigp\bigl(\{0<t\leq 1\}\bigr)$ is relatively compact in
 $M(\Gamma,V,h_V)$.
\item
 If each
 $\bigl(\prolongg{\vecc}{\widetilde{E}_i},\widetilde{\vecF}_i,\widetilde{\theta}_i\bigr)$
 is stable,
 then the limit $\lim_{t\to 0}\nbigp(t)$ exists,
 and the limit flat bundle underlies
 the variation of polarized Hodge structure.
 As a result, we can deform any flat bundle
 to a variation of polarized Hodge structure.
\hfill\qed
\end{enumerate}
\end{prop}

We would like to mention the point which we will care about.
For simplicity, we assume $(\prolongg{\vecc}{E},\vecF,\theta)$ is $\mu_L$-stable,
and $\bigl(\prolongg{\vecc}{E},\vecF,t\cdot\theta\bigr)$
converges to the $\mu_L$-stable parabolic Higgs bundle
$(\prolongg{\vecc}{\widetilde{E}},\widetilde{\vecF},\widetilde{\theta})$.
Let $\{t_i\}$ denote a sequence converging to $0$.
By taking an appropriate subsequence,
we may assume that the sequence
$\{(E,\delbar_E,h_{t_i},t_i\!\cdot\!\theta_i)\}$ converges
to a tame harmonic bundle $(E',\delbar_{E'},h',\theta')$
weakly in $L_2^p$ locally over $X-D$,
which is due to Uhlenbeck's compactness theorem and
the estimate for the Higgs fields.
Then we obtain the induced parabolic Higgs bundle
$(\prolongg{\vecc}{E'},\vecF',\theta')$.
We would like to show that
$(\prolongg{\vecc}{\widetilde{E}},\widetilde{\vecF},\widetilde{\theta})$
and $(\prolongg{\vecc}{E'},\vecF',\theta')$ are isomorphic.
Once we have known the existence
of a non-trivial map $G:\prolongg{\vecc}{E'}\lrarr \prolongg{\vecc}{\widetilde{E}}$ 
which is compatible with the parabolic structure and the Higgs field,
it is isomorphic due to the stability of
$(\prolongg{\vecc}{\widetilde{E}},\widetilde{\vecF},\widetilde{\theta})$.
Hence the existence of such $G$ is the main issue
for this argument.
We remark that the problem is rather obvious
if $D$ is empty.

\begin{rem}
Even if $(\prolongg{\vecc}{\widetilde{E}_i},\widetilde{\vecF}_i,\widetilde{\theta_i})$
are not $\mu_L$-stable,
the conclusion in the third claim of Proposition {\rm\ref{prop;05.6.3.1}}
should be true.
In fact, Simpson gave a detailed argument to show it,
 in the case where $D$ is empty  {\rm(\cite{s6}, \cite{s7})}.
 More strongly, he obtained the homeomorphism of
 the coarse moduli spaces of semistable flat bundles
 and semistable Higgs bundles.

 In this paper, we do not discuss the moduli spaces,
 and hence we omit to discuss the general case.
 Instead, we use an elementary inductive argument 
 on the rank of local systems,
 which is sufficient to obtain a deformation to
 a variation of polarized Hodge structure.
 However, it would be desirable to arrive at
 the thorough understanding
 as Simpson's work, in future.
\hfill\qed
\end{rem}

\begin{rem}
 For an application,
 we have to care about the relation between
 the deformation and the monodromy groups.
 We will discuss only a rough relation
 in Section {\rm\ref{section;05.10.4.2}}.
 More precise relation will be studied elsewhere.
\hfill\qed
\end{rem}

Once we can deform any local system on
a smooth quasiprojective variety
to a variation of polarized Hodge structure,
preserving some compatibility with the monodromy group,
we obtain the following corollary.
It is a natural generalization of Theorem \ref{thm;04.10.26.10}.

\begin{cor}
Let $\Gamma$ be a rigid discrete subgroup
of a real algebraic group which is not of Hodge type.
Then $\Gamma$ cannot be a split quotient
of the fundamental groups of any
smooth irreducible quasiprojective variety.
\hfill\qed
\end{cor}

\begin{rem}
Such a deformation of flat bundles on
a quasiprojective variety
was also discussed in {\rm\cite{jlz}} in a different way.
\hfill\qed
\end{rem}

\subsection{Tame pure imaginary pluri-harmonic reduction (Appendix)}

Let $G$ be a linear algebraic group defined over $\cnum$ or $\real$.
We will discuss a characterization of
reductive representations $\pi_1(X-D,x)\lrarr G$
via the existence of tame pure imaginary pluri-harmonic reduction.
Here a representation is called reductive, if
the Zariski closure of the image is reductive.
Such a kind of characterization was given
by Jost and Zuo (\cite{JZ2}) directly for $G$,
although their definition of reductivity looks different
from ours.
It is our purpose to explain that 
the problem can be reduced to the case $G=GL(n)$
by Tannakian consideration.
Some results are used in Chapter \ref{chapter;05.10.1.11}.

\section{Outline}

Chapter \ref{chapter;05.10.1.2} is an elementary preparation
for the discussion in the later chapters.
The reader can skip this chapter.
Chapter \ref{chapter;05.9.8.1} is preparation 
about parabolic Higgs bundles.
We discuss the perturbation of a given filtration
in Section \ref{section;05.7.30.15},
which is one of the keys in this paper.

In Chapter \ref{chapter;05.9.8.10},
an ordinary metric for parabolic Higgs bundle is given.
The construction is standard.
Our purpose is to establish the relation between
the parabolic characteristic numbers 
and some integrals,
in the case of  graded semisimple parabolic Higgs bundles.

In Chapter \ref{chapter;05.10.1.4},
we show the fundamental properties of the parabolic Higgs bundles
obtained from tame harmonic bundles.
Namely, we show the $\mu_L$-stability and 
the vanishing of the characteristic numbers.
In Chapter \ref{chapter;05.10.1.5},
we show the preliminary Kobayashi-Hitchin correspondence
for graded semisimple parabolic Higgs bundles.
Bogomolov-Gieseker inequality can be obtained as an easy corollary
of this preliminary correspondence and the perturbation argument
of the parabolic structure.

In Chapter \ref{chapter;06.8.21.100},
we construct a frame around the origin
for a tame harmonic bundle on a punctured disc.
It is a technical preparation to discuss
the convergence of a sequence of tame harmonic bundle.
Such a convergence is shown
in Chapter \ref{chapter;06.8.21.101}.
We also give a preparation for the existence theorem
of pluri-harmonic metric,
which is completed in Chapter \ref{chapter;06.8.21.102}.

Once the Kobayashi-Hitchin correspondence for tame harmonic bundles
is established,
we can apply Simpson's argument of the tours action,
and we can obtain some topological consequence
of quasiprojective varieties.
It is explained in Chapter \ref{chapter;05.10.1.11}.
Chapter \ref{chapter;04.10.26.90} is regarded as an appendix,
in which we recall something related to pluri-harmonic metrics
of $G$-flat bundles.

\section{Acknowledgement}

The author owes much thanks to C. Simpson.
The paper is a result of an effort to understand
his works, in particular, \cite{s1} and \cite{s5},
where the reader can find the framework and many ideas used in this paper.
The problem of Bogomolov-Gieseker inequality was passed to the author
from him a few years ago.
The author expresses his gratitude to the referee
for careful and patient reading, 
 pointing out some mistakes,
 and suggesting simplification.
The author thanks Y. Tsuchimoto  and A. Ishii
for their constant encouragement.
He is grateful to the colleagues of Department of Mathematics
at Kyoto University for their help.
He acknowledges the Max-Planck Institute for Mathematics
for their excellent hospitality and support,
where he revised this paper.

%% file: 2.tex
This chapter is a preparation for the later discussions.
We will often use the notation given in Sections
\ref{section;05.10.8.1}--\ref{section;05.10.1.1},
especially.

\section{Notation and Words}
\label{section;05.10.8.1}

We use the notation
$\seisuu$, $\rnum$,
$\real$ and $\cnum$
to denote the set of integers,
rational numbers,
real numbers
and complex numbers,
respectively.
For a real number $a$,
we put $\real_{>a}:=\{x\in\real\,|\,x>a\}$.
We use the notation
$\seisuu_{>a}$, $\seisuu_{\geq a}$,
$\rnum_{>a}$, etc. in a similar meaning.

For real numbers $a,b$,
we put as follows:
\[
\begin{array}{ll}
[a,b]:=\{x\in\real\,|\,a\leq x\leq b\}
 &
\closedopen{a}{b}:=\{x\in\real\,|\,a\leq x<b\} \\
\openclosed{a}{b}:=\{x\in\real\,|\,a<x\leq b\}
 &
\openopen{a}{b}:=\{x\in\real\,|\,a<x<b\}
\end{array}
 \]

The notation $\delta_{i,j}$ will be Kronecker's delta,
i.e., $\delta_{i,j}=1$ ($i=j$) and $\delta_{i,j}=0$ ($i\neq j$).

\vspace{.1in}

A normal crossing divisor $D$
of a complex manifold $X$ will be called {\em simple},
if each irreducible component is non-singular.
Let $D=\bigcup_{i\in S}D_i$ be the irreducible decomposition.
For elements $\veca\in \real^S$,
$a_i$ will denote the $i$-th component of $\veca$
$(i\in S)$.
The notation $\prolongg{\veca}{E}$
is often used to denote a vector bundle on $X$,
and we often put $E:=\prolongg{\veca}{E}_{|X-D}$.

\vspace{.1in}

Let $Y$ be a manifold, $E$ be a vector bundle on $Y$,
and $\{f_i\}$ be a sequence of sections of $E$.
We say $\{f_i\}$ converges to $f$ weakly in $L_l^p$ locally on $Y$,
if the restriction $\{f_{i\,|\,K}\}$ converges to $f_{|K}$
weakly in $L_l^p(K)$.

Let $\bigl\{(E^{(i)},\delbar^{(i)},\theta^{(i)})\bigr\}$ 
be a sequence of Higgs bundles on $Y$.
We say that the sequence
$\bigl\{(E^{(i)},\delbar^{(i)},\theta^{(i)})\bigr\}$ 
converges to
$(E^{(\infty)},\delbar^{(\infty)},\theta^{(\infty)})$
weakly in $L_2^p$ 
(resp. in $C^1$) locally on $Y$,
if there exist locally $L_2^p$-isomorphisms
(resp. $C^1$-isomorphisms)
$\Phi^{(i)}:E^{(i)}\lrarr E^{(\infty)}$ on $Y$
such that the sequences
$\{\Phi^{(i)}\bigl(\delbar^{(i)}\bigr)\}$ and
$\{\Phi^{(i)}(\theta^{(i)})\}$ weakly converge to
$\delbar^{(\infty)}$ and $\theta^{(\infty)}$ respectively
in $L_1^p$ (resp. $C^0$) locally on $Y$.

Let $E$ be a vector bundle on $Y$
with a hermitian metric $h$.
For an operator $F\in \End(E)\otimes\Omega^{p,q}_Y$,
we use the notation
$F^{\dagger}_h\in \End(E)\otimes\Omega^{q,p}_Y$
to denote the adjoint of $F$ with respect to $h$.
We often use $F^{\dagger}$,
if there are no risk of confusion.

\vspace{.1in}

Let $(S_i,\varphi_i)$ $(i=1,2,\ldots,\infty)$
be a pair of discrete subsets $S_i\subset\real$
and functions $\varphi_i:S_i\lrarr\seisuu_{>0}$.
We say that $\bigl\{(S_i,\varphi_i)\,\big|\,i=1,2,\ldots\bigr\}$ 
converges to $(S_{\infty},\varphi_{\infty})$,
if there exists $i_0$ for any $\epsilon>0$
such that 
(i) any $b\in S_{i}$ $(i>i_0)$
is contained in $\openopen{a-\epsilon}{a+\epsilon}$
for some $a\in S_{\infty}$,
(ii)
$\sum_{b\in S_i,|a-b|<\epsilon}\varphi_i(b)
   =\varphi_{\infty}(a)$
is satisfied.

\section[Review of Results of Simpson]
{Review of some Results of Simpson on 
 Kobayashi-Hitchin Correspondence}
\label{section;05.10.1.1}

\subsection{Analytic stability and the Hermitian-Einstein metric}
\label{subsection;05.10.5.1}

Let $Y$ be an $n$-dimensional complex manifold
which is not necessarily compact.
Let $\omega$ be a Kahler form of $Y$.
The adjoint for the multiplication of $\omega$
is denoted by $\Lambda_{\omega}$, or simply by $\Lambda$
if there are no confusion.
The Laplacian for $\omega$ is denoted by $\Delta_{\omega}$.

\begin{condition}
\label{condition;05.7.30.1}
\mbox{{}}
\begin{enumerate}
\item
 The volume of $Y$ with respect to $\omega$ is finite.
\item
 There exists an exhaustion function $\phi$ on $Y$
 such that 
$0\leq \sqrt{-1}\del\delbar\phi\leq C\cdot\omega$
 for some positive constant $C$.
 \item
 There exists an increasing function
 $\real_{\geq\,0}\lrarr\real_{\geq\,0}$
 such that $a(0)=0$ and $a(x)=x$ for $x\geq 1$,
 and the following holds:
 \begin{itemize}
\item
 Let $f$ be a positive bounded function on $Y$
 such that $\Delta_{\omega} f\leq B$ for some positive number $B$.
 Then $\sup_{Y}|f|\leq C(B)\cdot a\Bigl(\int_Yf\Bigr)$
 for some positive constant $C(B)$ depending on $B$.
 Moreover $\Delta_{\omega}f\leq 0$ implies
 $\Delta_{\omega} f=0$.
\hfill\qed
 \end{itemize}
\end{enumerate}
\end{condition}

Let $(E,\delbar_E,\theta)$ be a Higgs bundle on $Y$.
Let $h$ be a hermitian metric of $E$.
Then we have the $(1,0)$-operator $\del_{E}$
determined by 
$\delbar h(u,v)
=h\bigl(\delbar_Eu,v\bigr)+h\bigl(u,\del_{E}v\bigr)$.
We also have the adjoint $\theta^{\dagger}$.
If we emphasize the dependence on $h$,
we use the notation $\del_{E,h}$ and $\theta^{\dagger}_{h}$.
We obtain the connections
$D_h:=\delbar_E+\del_E$
and $\DD^1:=D_h+\theta+\theta^{\dagger}$.
The curvatures of $D_h$ and $\DD^1$
are denoted by $R(h)$ and $F(h)$ respectively.
When we emphasize the dependence on $\delbar_E$,
they are denoted by
$R(\delbar_E,h)$ and $F(\delbar_E,h)$.
We also use $R(E,h)$ and $F(E,h)$,
if we emphasize the bundle.

\begin{condition}
\label{condition;05.7.30.2}
$F(h)$ is bounded with respect to
$h$ and $\omega$.
\hfill\qed
\end{condition}

When Condition \ref{condition;05.7.30.2} is satisfied,
we put as follows:
\[
 \deg_{\omega}(E,h):=
 \frac{\sqrt{-1}}{2\pi}
 \int_Y\tr \bigl(F(h)\bigr)\cdot\omega^{n-1}
=\frac{\sqrt{-1}}{2\pi}
 \int_Y\tr \Lambda(F(h))
 \cdot\frac{\omega^n}{n}
\]
Note $\tr F(h)=\tr R(h)$.
Recall that a subsheaf $V\subset E$ is called
saturated if the quotient $E/V$ is torsion-free.
For any saturated Higgs subsheaf $V\subset E$,
there is a Zariski closed subset $Z$ of codimension two
such that $V_{|Y-Z}$ gives a subbundle of $E_{|Y-Z}$,
on which the metric $h_V$ of $V_{|Y-Z}$ is induced.
Let $\pi_V$ denote the orthogonal projection
of $E_{|Y-Z}$ onto $V_{|Y-Z}$.
Let $\tr_V$ denote the trace for endomorphisms of $V$.

\begin{prop}
 [\cite{s1} Lemma 3.2]
 \label{prop;05.9.8.3}
When the conditions {\rm\ref{condition;05.7.30.1}}
and {\rm\ref{condition;05.7.30.2}} are satisfied,
the integral
\[
 \deg_{\omega}(V,K):=\frac{\sqrt{-1}}{2\pi}
 \int_Y\tr_V\bigl(F(h_V)\bigr)\cdot\omega^{n-1}
\]
is well defined, and it takes the value in $\real\cup\{-\infty\}$.
The Chern-Weil formula holds as follows,
for some positive number $C$:
\[
 \deg_{\omega}(V,h_V)=
 \frac{\sqrt{-1}}{2\pi}
 \int_Y \tr\Bigl(
 \pi_V\circ\Lambda_{\omega} F(h)
 \Bigr)
\cdot\frac{\omega^n}{n}
-C\int_Y\bigl|D''\pi_V\bigr|_h^2\cdot\dvol_{\omega}.
\]
Here we put $D''=\delbar_E+\theta$.
In particular,
if the value $\deg_{\omega}(V,h_V)$ is finite,
$\delbar_E(\pi_V)$ and $[\theta,\pi_V]$ are $L^2$.
\hfill\qed
\end{prop}

For any $V\subset E$,
we put $\mu_{\omega}(V,h_V):=\deg_{\omega}(V,h_V)/\rank V$.

\begin{df}
[\cite{s1}]
A metrized Higgs bundle $(E,\delbar_E,\theta,h)$ is called analytic stable,
if the inequalities
$ \mu_{\omega}(V,h_V)<\mu_{\omega}(E,h)$
hold for any non-trivial Higgs saturated subsheaves
$(V,\theta_V)\subsetneq (E,\theta)$.
\hfill\qed
\end{df}

The following important theorem is crucial
for our argument.

\begin{prop}[Simpson] 
\label{prop;05.7.30.35}
 Let $(Y,\omega)$ be a Kahler manifold
 satisfying Condition {\rm\ref{condition;05.7.30.1}},
 and let $(E,\delbar_E,\theta,h_0)$ be 
 a metrized Higgs bundle satisfying
 Condition {\rm\ref{condition;05.7.30.2}}.
 If it is analytic stable,
 then there exists a hermitian metric $h=h_0\cdot s$
 satisfying the following conditions:
\begin{itemize}
\item
 $h$ and $h_0$ are mutually bounded.
\item
 $\det(h)=\det(h_0)$.
 In particular,
 we have $\tr F(h)=\tr F(h_0)$.
\item
 $D''(s)$ is $L^2$ with respect to $h_0$ and $\omega$.
\item
 It satisfies the Hermitian-Einstein condition
 $\Lambda_{\omega} F(h)^{\bot}=0$,
 where $F(h)^{\bot}$ denotes the trace free part
 of $F(h)$.
\item
 The following equalities hold:
\begin{equation}
 \label{eq;06.8.4.1}
 \int_{Y}\tr\Bigl(F(h)^2\Bigr)\cdot\omega^{n-2}
=\int_{Y}\tr\Bigl(F(h_0)^2\Bigr)\cdot\omega^{n-2},
\end{equation}
\begin{equation}
 \label{eq;06.8.4.2}
  \int_Y\tr\Bigl(F(h)^{\bot\,2}\Bigr)\cdot\omega^{n-2}
=\int_Y\tr\Bigl(F(h_0)^{\bot\,2}\Bigr)\cdot\omega^{n-2}.
\end{equation}
\end{itemize}
\end{prop}
\pf
Condition \ref{condition;05.7.30.2} implies
$\Lambda_{\omega}F(h)$ is bounded.
Applying Theorem 1 in \cite{s1},
we obtain the hermitian metric $h$
satisfying the first four conditions.
Due to Proposition 3.5 in \cite{s1},
we obtain the inequality 
$\int_Y\tr\bigl(F(h)^2\bigr)\cdot\omega^{n-2}
\leq
 \int_Y\tr\bigl(F(h_0)^2\bigr)\cdot\omega^{n-2}$.
Since we have assumed the boundedness of
$F(h_0)$,
we also obtain
$\int_Y\tr\bigl(F(h)^2\bigr)\cdot\omega^{n-2}
\geq
 \int_Y\tr\bigl(F(h_0)^2\bigr)\cdot\omega^{n-2}$
due to Lemma 7.4 in \cite{s1},
as mentioned in the remark just before the lemma.
Therefore, we obtain (\ref{eq;06.8.4.1}).
Since we have $\tr F(h_0)=\tr F(h)$,
we also obtain (\ref{eq;06.8.4.2}).
\hfill\qed

\subsection{Uniqueness}

The following proposition can be
proved by the methods in \cite{s1}.

\begin{prop}
 \label{prop;05.9.8.100}
Let $(Y,\omega)$ be a Kahler manifold
satisfying
Condition {\rm\ref{condition;05.7.30.1}},
and $(E,\delbar_E,\theta)$ be a Higgs bundle on $Y$.
Let $h_i$ $(i=1,2)$ be hermitian metrics of $E$
such that $\Lambda_{\omega}F(h_i)=0$.
We assume that $h_1$ and $h_2$ are mutually bounded.
Then the following holds:
\begin{itemize}
\item
 We have the decomposition of Higgs bundles
$(E,\theta)=\bigoplus (E_a,\theta_a)$
 which is orthogonal with respect to both of $h_i$.
\item 
 The restrictions of $h_i$ to $E_a$ are denoted by $h_{i,a}$.
 Then there exist positive numbers $b_a$ such that
 $h_{1,a}=b_a\cdot h_{2,a}$.
\end{itemize}
\end{prop}
\pf
We take the endomorphism
$s_1$ determined by $h_2=h_1\cdot s_1$.
Then we have the following inequality
due to Lemma 3.1 (d) in \cite{s1} on $X-D$:
\[
\Delta_{\omega}\log \tr\bigl(s_1\bigr)\leq 
 \bigl|\Lambda_{\omega}F(h_1)\bigr|
+\bigl|\Lambda_{\omega}F(h_2)\bigr|=0.
\]
Here we have used $\Lambda_{\omega}F(h_i)=0$.
Then we obtain $\Delta_{\omega}\tr\bigl(s_1\bigr)\leq 0$.
Since the function $\tr(s_1)$ is bounded on $Y$,
we obtain the harmonicity $\Delta_{\omega}\tr(s_1)=0$
due to Condition \ref{condition;05.7.30.1}.

We put $D''=\delbar+\theta$ and
$D':=\del_{E,h_1}+\theta^{\dagger}_{h_1}$,
where $\theta^{\dagger}_{h_1}$ denotes the adjoint of
$\theta$ with respect to the metric $h_{1}$.
Then we also have the following equality:
\[
 0=F\bigl(h_2\bigr)-F\bigl(h_1\bigr)
=D''\bigl(s_1^{-1}D's_1\bigr)
=-s_1^{-1}D''s_1\cdot s_1^{-1}\cdot D's_1
 +s_1^{-1}D''D's_1.
\]
Hence we obtain
$D''D's_1=D''s_1\cdot s_1^{-1}\cdot D's_1$.
As a result, we obtain the following equality:
\[
 \int \bigl|s_1^{-1/2}D''s_1\bigr|^2_{h_1}\dvol_{\omega}
=-\sqrt{-1}\int \Lambda_{\omega}\tr \bigl(D''D's_1\bigr)\dvol_{\omega}
=-\int \Delta_{\omega}\tr(s_1)\dvol_{\omega}
=0.
\]
Hence we obtain $D''s_1=0$,
i.e.,
$\delbar s_1=\bigl[\theta,s_1\bigr]=0$.
Since $s_1$ is self-adjoint with respect to $h_1$,
we obtain the flatness
$\bigl(\delbar+\del_{E,h_1}\bigr)s_1=0$.
Hence we obtain the decomposition
$E=\bigoplus_{a\in S} E_a$
such that $s_a=\bigoplus b_a\cdot \id_{E_a}$
for some positive constants $b_a$.
Let $\pi_{E_a}$ denote the orthogonal projection
onto $E_a$.
Then we have $\delbar \pi_{E_a}=0$.
Hence the decomposition $E=\bigoplus_{a\in S} E_a$ is holomorphic.
It is also compatible with the Higgs field.
Hence we obtain the decomposition as the Higgs bundles.
Then the claim of Proposition \ref{prop;05.9.8.100} is clear.
\hfill\qed

\begin{rem}
We have only to impose $\Lambda_{\omega}F(h_1)=\Lambda_{\omega}F(h_2)$
instead of $\Lambda_{\omega}F(h_i)=0$,
which can be shown by a minor refinement of the argument.
\hfill\qed
\end{rem}

\subsection{The one dimensional case}

In the one dimensional case,
Simpson established the Kobayashi-Hitchin correspondence for
parabolic Higgs bundle.
Here we recall only the special case.
(See Chapter \ref{chapter;05.9.8.1} for some definitions.)
\begin{prop}[Simpson]
 \label{prop;05.9.160}
Let $X$ be a smooth projective curve,
and $D$ be a divisor of $X$.
Let $\bigl(\vecE_{\ast},\theta\bigr)$ be a
filtered regular Higgs bundle on $(X,D)$.
We put $E=\prolongg{\vecc}{E}_{|X-D}$.
The following conditions are equivalent:
\begin{itemize}
\item
 $(\vecE_{\ast},\theta)$ is poly-stable
 with $\pardeg(\vecE_{\ast})=0$.
\item
 There exists a harmonic metric $h$ of
 $(E,\theta)$, which is adapted to the parabolic structure of
     $\vecE_{\ast}$.
\end{itemize}
Moreover,
such a metric is unique up to obvious ambiguity.
Namely, let $h_i$ $(i=1,2)$ be two harmonic metrics.
Then we have the decomposition of Higgs bundles
$(E,\theta)=\bigoplus (E_a,\theta_a)$ satisfying the following:
\begin{itemize}
\item
 The decomposition is orthogonal with respect to
 both of $h_i$.
\item
 The restrictions of $h_i$ to $E_a$ are denoted by $h_{i,a}$.
 Then there exist positive numbers $b_a$ 
 such that $h_{1,a}=b_{a}\cdot h_{2,a}$.
\end{itemize}
\end{prop}
\pf
See \cite{s2}.
We give only a remark on the uniqueness.
Let $(E,\delbar_E,\theta)$ be a Higgs bundle
on $X-D$, and $h_i$ $(i=1,2)$ be harmonic metrics on it.
Assume that the induced prolongments
$\prolongg{\vecc}{E}(h_i)$ are isomorphic.
(See Section \ref{section;05.9.8.110} for prolongment.)
Recall the norm estimate for tame harmonic bundles
in the one dimensional case (\cite{s2}),
which says that the harmonic metrics are determined
up to boundedness by the parabolic filtration and the weight filtration.
Hence we obtain the mutually boundedness of  $h_1$ and $h_2$.
Then the uniqueness follows from Proposition \ref{prop;05.9.8.100}.
\hfill\qed

\section{Weitzenb\"{o}ck Formula}

Let $(Y,\omega)$ be a Kahler manifold.
Let $h$ be a Hermitian-Einstein metric for
a Higgs bundle $(E,\delbar_E,\theta)$ on $Y$.
More strongly, we assume $\Lambda_{\omega}F(h)=0$.
The following lemma is a minor modification 
of Weitzenb\"{o}ck formula for harmonic bundles by Simpson (\cite{s2}).

\begin{lem}\label{lem;05.8.25.3}
Let $s$ be any holomorphic section of $E$ such that $\theta s=0$.
Then we have
$\Delta_{\omega}\log|s|_h^2\leq 0$,
where $\Delta_{\omega}$ denotes the Laplacian for $\omega$.
\end{lem}
\pf
We have
$\del\delbar|s|_h^2
=\del\bigl(s,\del_Es\bigr)
=(\del_Es,\del_Es)+(s,\delbar_E\del_Es)
=(\del_Es,\del_Es)+(s,R(h)s)$.
Then we obtain the following:
\[
 \del\delbar \log|s|_h^2
=\frac{\del\delbar |s|^2}{|s|^2}
-\frac{\del|s|^2\cdot\delbar|s|^2}{|s|^4}
=\frac{(s,R(h)s)}{|s|^2}
+\frac{(\del_E s,\del_E s)}{|s|^2}
-\frac{\del|s|^2\cdot\delbar|s|^2}{|s|^4}.
\]
We have
$R(h)=-(\theta^{\dagger}\theta+\theta\theta^{\dagger})+F(h)^{(1,1)}$,
where $F(h)^{(1,1)}$ denotes the $(1,1)$-part of $F(h)$.
Hence we have the following:
\begin{multline}
 \Lambda_{\omega}\bigl(s,R(h)s\bigr)
=\Lambda_{\omega}
 \Bigl(s,\,\,(-\theta\theta^{\dagger}-\theta^{\dagger}\theta)s\Bigr)
+\Lambda_{\omega}\bigl(s,F(h)^{(1,1)}s\bigr) \\
=-\Lambda_{\omega}\bigl(\theta^{\dagger}s,\theta^{\dagger}s\bigr)
-\Lambda_{\omega}\bigl(\theta s,\theta s\bigr)
+\Lambda_{\omega}\bigl(s,F(h)^{(1,1)}s\bigr) 
=-\Lambda_{\omega}(\theta^{\dagger}s,\,\theta^{\dagger}s).
\end{multline}
Here we have used
$\Lambda_{\omega} F(h)=\Lambda_{\omega}F(h)^{(1,1)}=0$.
Therefore we obtain the following:
\[
 -\sqrt{-1}\Lambda_{\omega}\bigl(s,R(h)s\bigr)
=\sqrt{-1}\Lambda_{\omega}(\theta^{\dagger}s,\theta^{\dagger}s)
=-\bigl|\theta^{\dagger}s\bigr|_h^2.
\]
On the other hand, we also have the following:
\[
 -\sqrt{-1}\Lambda_{\omega}\left(
 \frac{(\del s,\del s)}{|s|^2}
-\frac{\del |s|^2\delbar|s|^2}{|s|^4}
 \right)\leq 0.
\]
Hence we obtain 
$\Delta_{\omega}\log|s|^2\leq 0$.
\hfill\qed

\section{A Priori Estimate of Higgs Fields}

\label{section;05.8.27.1}

\subsection{On a disc}

We put $X(T):=\bigl\{z\in\cnum\,|\,|z|<T\bigr\}$
for any positive number $T$.
In the case $T=1$, $X(1)$ is denoted by $X$.
We will use the usual Euclidean metric 
$g=dz\cdot d\bar{z}$
and the induced measure $\dvol_g$.
Let $\Delta$ denote the Laplacian $-\del_z\delbar_z$.
By the standard theory of Dirichlet problem,
there exists a constant $C'$ such that the following holds:
\begin{itemize}
\item
We have the solution $\psi$
of the equation $\Delta \psi=\kappa$ 
such that $|\psi(P)|\leq C'\cdot\|\kappa\|_{L^2}$
for any $L^2$-function $\kappa$ 
and for any $P\in X$.
\end{itemize}

Let $(E,\delbar_E,\theta)$ be a Higgs bundle
on $X$ with a hermitian metric $h$.
We have the expression $\theta=f\cdot dz$.
We would like to estimate of the norm $\bigl|f\bigr|_h$
by the eigenvalues of $g$ and the $L^2$-norm
$\bigl\|F(h)\bigr\|_{L^2}:=
 \int_{X}|F(h)|_{h,g}^2\cdot\dvol_g$.

\begin{prop}
 \label{prop;05.8.25.2}
Let $t$ be any positive number such that $t<1$.
There exist constants $C$ and $C'$
such that the following inequality holds on $X(t)$:
\[
 |f|_h^2\leq C\cdot e^{10C'||F(h)||_{L^2}}.
\]
The constant $C'$ is as above,
and the constant $C$ depends only on $t$, the rank of $E$
and the eigenvalues of $f$.
\end{prop}
\pf
Let us begin with the following lemma,
which is just a minor modification of the fundamental inequality
in the theory of harmonic bundles.
\begin{lem}
 \label{lem;05.9.18.1}
We have the inequality:
\[
 \Delta\log|f|_h^2\leq 
-\frac{\bigl|[f,f^{\dagger}]\bigr|^2_h}{|f|_h^2}
+|F(h)|_{h,g}.
\]
\end{lem}
\pf
By a general formula,
we have the following inequality:
\[
 -\sqrt{-1}\Lambda\del\delbar\log|f|_h^2
\leq 
-\sqrt{-1}\Lambda\frac{\bigl(f,\,[R(h),f]\bigr)}{|f|_h^2}.
\]
We obtain the desired inequality
from $R(h)=F(h)-[\theta,\theta^{\dagger}]=
 F(h)-[f,f^{\dagger}]\cdot dz\cdot d\bar{z}$.
\hfill\qed

\vspace{.1in}

Let us take a function $A$ satisfying
$\Delta A=|F(h)|_h$ and $|A|\leq C'||F(h)||_{L^2}$.
Then we obtain the following:
\[
 \Delta\bigl(\log|f|_h^2-A \bigr)
=\Delta\log\bigl(|f|_h^2\cdot e^{-A}\bigr)
\leq 
 -\frac{\bigl|[f,f^{\dagger}]\bigr|^2_h}{|f|_h^2}.
\]

For any $Q\in X$,
let $\alpha_1(Q),\ldots,\alpha_{\rank(E)}(Q)$
denote the eigenvalues of $f_{|Q}$.
We put $\nu(Q):=\sum_{i=1}^{\rank(E)}|\alpha_i(Q)|^2$
and $\mu(Q):=|f_{|Q}|_h^2-\nu(Q)$.
It can be elementarily shown that
there exists a constant $C_1$ which depends only on 
the rank of $E$, such that 
$C_1\cdot\mu^2\leq \bigl|[f,f^{\dagger}]\bigr|_h^2$.
Hence, the following inequality holds:
\[
 \Delta\log\bigl(e^{-A}\cdot |f|_h^2\bigr)
\leq
 -C_1\cdot \frac{\mu^2}{|f|_h^2}.
\]
We also have a constant $C_2$
which depends only on the eigenvalues
of $f$, such that  $\nu\leq C_2$ holds.

Let $T$ be a number such that $0<T<1$,
and $\phi_T:X(T)\lrarr \real$ is given by the following:
\[
 \phi_T(z)=\frac{2T^2}{(T^2-|z|^2)^2}.
\]
Then we have
$\Delta\log\phi_T=-\phi_T$ and $\phi_T\geq 2$.
In particular, we have
$\nu \leq C_2\cdot\phi_{T}/2$.
The following lemma is clear.
\begin{lem}
\label{lem;05.8.25.1}
Either one of  $|f_{|Q}|_h^2\leq C_2\cdot\phi_T(Q)$
or $|f_{|Q}|_h^2\leq 2\mu(Q)$ holds
for any $Q\in X$.
\hfill\qed
\end{lem}

We take a constant $\widehat{C}_3>0$ satisfying 
$\widehat{C}_3>C_2$ and $\widehat{C}_3>4\cdot C_1^{-1}$,
and we put $C_3:=\widehat{C}_3\cdot e^{C'\|F(h)\|_{L^2}}$.
We put 
$S_T:=\bigl\{
 P\in X(T)\,\big|\, \bigl(e^{-A}\cdot|f|^2\bigr)(P)>C_3\cdot \phi_T(P)
 \bigr\}$.
For any point $P\in S_T$,
we have 
$ |f(P)|_h^2
 >C_3\cdot e^{A(P)}\cdot\phi_T(P)
 >C_2\cdot\phi_T(P)$.
Due to Lemma \ref{lem;05.8.25.1},
we obtain the following:
\[
  \Delta\log \bigl(e^{-A}\cdot|f|_h^2\bigr)(P)
\leq
 -\frac{C_1}{4}\cdot|f(P)|_h^2 
\leq
 -\frac{1}{C_3}\bigl(e^{-A}\cdot|f|_h^2\bigr)(P).
\]
On the other hand, we have the following:
\[
 \Delta\log(C_3\cdot\phi_T)=-\frac{1}{C_3}(C_3\cdot\phi_T).
\]
Moreover, it is easy to see $\del S_T\cap \{|z|=T\}=\emptyset$.
Hence, we obtain $S_T=\emptyset$ by a standard argument.
(See \cite{a}, \cite{s2} or the proof of
 Proposition 7.2 in \cite{mochi2}.)
Namely, we obtain the inequality 
$e^{-A}|f|_h^2 \leq 
 \widehat{C}_3\cdot e^{C'\|F(h)\|_{L^2}}\cdot\phi_T$
on $X(T)$.
Taking a limit for $T\to 1$,
we obtain 
$|f|_h^2\leq 
 e^{2C'||F(h)||_{L^2}}\cdot \widehat{C}_3\cdot (1-|z|^2)^{-1}$
on $X$.
Then the claim of Proposition \ref{prop;05.8.25.2}
follows.
\hfill\qed

\subsection{A Priori Estimate on a Multi-disc}

For a positive number $T$,
we put $Y(T):=\bigl\{(z_1,\ldots,z_n)\,\big|\,|z_i|<T\bigr\}$.
Let $g$ denote the metric $\sum dz_i\cdot d\zbar_i$ of
$Y(T)$.
Let $\omega$ be a Kahler form on $Y(T)$
such that there exists a constant $C>0$ such that
$C^{-1}\cdot\omega\leq  g\leq C\cdot \omega$.
Let $(E,\delbar_E,\theta)$ be a Higgs bundle
with a hermitian metric $h$,
which is Hermitian-Einstein
with respect to $\omega$.
For simplicity,
we restrict ourselves to the case $\Lambda_{\omega} F(h)=0$.
We assume
$||F(h)||_{L^2}<\infty$,
where $\|F(h)\|_{L^2}$ denotes
the $L^2$-norm of $F(h)$ with respect to $\omega$ and $h$.
We have the expression $\theta=\sum f_i\cdot dz_i$
for holomorphic sections $f_i\in \End(E)$ on $Y(T)$.

\begin{lem}
\label{lem;05.8.25.4}
Take $0<T_1<T$.
There exist some constants $C_1$ and $C_2$
such that the following inequality holds
for any $P\in Y(T_1)$:
\[
 \log|f_i|^2(P)
\leq
 C_1\cdot \bigl\|F(h)\bigr\|_{L^2}
+C_2.
\]
The constants $C_1$ and $C_2$ are good
in the sense that they depend only on
$T$, $T_1$, $\rank E$, the eigenvalues of $f_i$ $(i=1,2,\ldots,n)$
and the constant $C$.
\end{lem}
\pf
We take a positive number $T_2$ such that $T_1<T_2<T$.
The induced Higgs field and the metric of
$\End(E)$ are denoted by $\widetilde{\theta}$ and $\widetilde{h}$.
Then the metric $\widetilde{h}$ is a Hermitian-Einstein metric
of $\bigl(\End(E),\widetilde{\theta}\bigr)$
such that  $\Lambda_{\omega} F(\widetilde{h})=0$.
Because of $\widetilde{\theta}(f_i)=0$,
we have  the subharmonicity
$\Delta_{\omega} \log|f_i|_h^2\leq 0$
due to Lemma \ref{lem;05.8.25.3}.
For $P\in Y(T_1)$, we obtain the following inequality
(see Theorem 9.20 in \cite{GT}, for example):
\[
 \log|f_i|^2(P)
\leq
 C_3\cdot \int_{Y(T_2)^2}
 \log^+|f_i|^2\cdot \dvol_g.
\]
Here we put $\log^+(y):=\max\{0,\log y\}$,
and $C_3$ denotes a good constant.

The $(1,1)$-part of $F(h)$ is expressed as
$\sum F_{i,j}\cdot dz_i\cdot d\bar{z}_j$.
Due to Proposition \ref{prop;05.8.25.2},
there exist good constants $C_j$ $(j=4,5)$
such that the following inequality holds for any point
$(z_1,\ldots,z_n)\in Y(T_1)$:
\[
 \log|f_1|^2(z_1,\ldots,z_n)\leq
 C_4\cdot\left(
\int_{|w_1|\leq T_2}|F_{1,1}(w_1,z_2,\ldots,z_n)|^2\cdot 
\sqrt{-1}dw_1\wedge d\overline{w}_1
\right)^{1/2}\!\!\!\!\!+C_5.
\]
Then the claim of Lemma \ref{lem;05.8.25.4} follows.
\hfill\qed

\section[Norm Estimate in Two Dimensional Case]
 {Norm Estimate for Tame Harmonic Bundle
in Two Dimensional Case}

\subsection{Norm estimate}
\label{subsection;06.8.15.11}

We recall some results in \cite{mochi2}.
We use bold symbols like $\veca$ to denote a tuple,
and $a_i$ denotes the $i$-th component of $\veca$.
We say $\veca\leq \vecb$ for $\veca,\vecb\in\real^2$
if $a_i\leq b_i$.
We put $X:=\{(z_1,z_2)\in\cnum^2\,\big|\,|z_i|<1\}$,
$D_i:=\{z_i=0\}$ and $D:=D_1\cup D_2$.
Let $\harmonicbundle$ be a tame harmonic bundle
on $X-D$.
For each $\vecc=(c_1,c_2)\in\real^2$,
we obtain the locally free sheaf
$\prolongg{\vecc}{E}$ on $X$
with parabolic structure $\lefttop{i}F$ $(i=1,2)$,
as in Section \ref{section;05.9.8.110}.
We also obtain the Higgs field $\theta$
of $\prolongg{\vecc}{E}_{\ast}$.
The residue of $\theta$ induces
the endomorphism
$\Gr^F\Res_i(\theta)\in 
\End(\lefttop{i}\Gr^F(E_{|D_i}))$
whose eigenvalues are constant on $D_i$.
Thus, the nilpotent part $\nbign_i$ of
$\Gr^F\Res_i(\theta)$ is well defined.
It is shown that the conjugacy classes
of $\nbign_{i\,|\,P}$ are independent of $P\in D_i$.
Let $\lefttop{\itibar}W$ denote
the weight filtration of $\nbign_{1}$
on $\lefttop{1}\Gr^F(E_{|D_1})$.

We have two filtrations
$\lefttop{i}{F}$ $(i=1,2)$ on $\prolongg{\vecc}{E}_{|O}$.
We put 
$\lefttop{\nibar}\Gr^F_{\veca}
:=\lefttop{2}\Gr^F_{a_2}
 \lefttop{1}\Gr^F_{a_1}(\prolongg{\vecc}{E}_{|O})$.
The maps $\nbign_i$ induce
the endomorphisms of $\lefttop{\nibar}\Gr^F_{\veca}$
which are denoted by $\lefttop{\nibar}\nbign_i$.
Let $\lefttop{\nibar}W$ denote the weight filtration of
$\lefttop{\nibar}\nbign_1+\lefttop{\nibar}\nbign_2$.
We also have the filtration
induced by $\lefttop{\itibar}W$,
which is denoted by the same notation.
We can take a decomposition
$ \prolongg{\vecc}{E}=
 \bigoplus_{(\veca,\veck)\in\real^2\times\seisuu^2}
 U_{(\veca,\veck)}$
satisfying the following conditions:
\begin{itemize}
\item
$\lefttop{i}F_b(\prolongg{\vecc}{E}_{|D_i})=
 \bigoplus_{a_i\leq b}U_{\veca,\veck\,|\,D_i}$
and 
$\lefttop{1}F_{b_1}(\prolongg{\vecc}{E}_{|O})
 \cap \lefttop{2}F_{b_2}(\prolongg{\vecc}{E}_{|O})
=\bigoplus_{\veca\leq\vecb}U_{\veca,\veck|O}$
\item
We have
$\lefttop{\itibar}W_k\bigl(
 \lefttop{1}\Gr^F_b(\prolongg{\vecc}{E}_{|D_1})\bigr)
=\bigoplus_{a_1=b,k_1\leq k}U_{\veca,\veck\,|\,D_1}$
under the isomorphism
$\lefttop{1}\Gr^F_b(\prolongg{\vecc}{E}_{|D_1})
\simeq \bigoplus_{a_1=b}U_{\veca,\veck\,|\,D_1}$.
\item
We have
$\lefttop{\itibar}W_{k_1}\cap \lefttop{\nibar}W_{k_2}
\bigl(
 \lefttop{\nibar}\Gr^{F}_{\veca}(\prolongg{\vecc}{E}_{|O})
 \bigr)
=\bigoplus_{\vecl\leq \veck}U_{\veca,\vecl}$
under the isomorphism
$\lefttop{\nibar}\Gr^F_{\veca}(\prolongg{\vecc}{E}_{|O})
\simeq
\bigoplus_{\veck}U_{\veca,\vecl}$.
\end{itemize}

We take a holomorphic frame $\vecv=(v_1,\ldots,v_r)$
which is compatible with the decomposition,
i.e.,
for each $v_i$ we have
$(\veca(v_i),\veck(v_i))\in\real^2\times\seisuu^2$
such that $v_i\in U_{\veca(v_i),\veck(v_i)}$.
Let $\hhat_1$ be a hermitian metric of $E$
given as follows:
\[
 \hhat_1(v_i,v_j)=\delta_{i,j}
\cdot |z_1|^{-2a_1(v_i)}|z_2|^{-2a_2(v_i)}
 \bigl(-\log|z_1|\bigr)^{k_1(v_i)}
\bigl(-\log|z_2|\bigr)^{k_2(v_i)-k_1(v_i)}
\]
We put $Z:=\bigl\{(z_1,z_2)\,\big|\,|z_1|<|z_2|\bigr\}$.
\begin{lem}
\label{lem;06.8.15.2}
$h$ and $\hhat_1$ are mutually bounded on $Z$.
\hfill\qed
\end{lem}

\subsection{Some estimate for related metrics}
\label{subsection;06.8.15.10}

We put $\Xtilde:=\bigl\{(\zeta_1,\zeta_2)\,\big|\,|\zeta_i|<1\bigr\}$,
$\Dtilde_i=\{\zeta_i=0\}$
and $\Dtilde=\Dtilde_1\cup\Dtilde_2$.
Let $\pi:\Xtilde-\Dtilde\lrarr X-D$ denote the map
given by $\pi(\zeta_1,\zeta_2)=(\zeta_1\zeta_2,\zeta_2)$.
Then, we have $\pi^{-1}(Z)=\Xtilde-\Dtilde$.
Hence Lemma \ref{lem;06.8.15.2} is reworded
as $\pi^{\ast}h$ and $\pi^{\ast}\hhat_1$ are mutually bounded.

We give a preparation for later use.
We put $\Etilde:=\pi^{\ast}E$.
For $\veca=(a_1,a_2)\in\real^2$,
we put $\vecatilde:=(a_1,a_1+a_2)$.
Then, we put
$\pi^{\ast}U_{\veca,\veck}=:
 \Utilde_{\vecatilde,\veck}$.
We put $\vecvtilde:=\pi^{\ast}\vecv$.
We put $\atilde_1(\vtilde_i):=a_1(v_i)$,
$\atilde_2(\vtilde_i)=a_1(v_i)+a_2(v_i)$,
$k_j(\vtilde_i):=k_j(v_i)$.
Then, $\vtilde_i$ is a section of
$\Utilde_{\vecatilde(\vtilde_i),\veck(\vtilde_i)}$.

Let $\chi$ be a non-negative valued function 
on $\real$ such that
$\chi(t)=1$ $(t\leq 1/2)$ and $\chi(t)=0$
$(t\geq 2/3)$.
Let $\rho(\zeta):\cnum^{\ast}\lrarr \real$
be the function given by
$\rho(\zeta)=-\chi(|\zeta|)\cdot\log|\zeta|^2$.
Then, we will use the following metrics later
(Section \ref{section;06.8.17.100})
\[
 h_0(\vtilde_i,\vtilde_j):=
 \delta_{i,j}\cdot \prod_k|\zeta_k|^{-2a_k(v_i)}
\]
\[
  h_1(\vtilde_i,\vtilde_j):=
 h_0(\vtilde_i,\vtilde_j)\cdot
 \bigl(1+\rho(\zeta_1)+\rho(\zeta_2)\bigr)^{k_1(\vtilde_i)}
\cdot\bigl(1+ \rho(\zeta_2)\bigr)^{k_2(\vtilde_i)-k_1(\vtilde_i)}
\]
Then, $h_1$ and $\pi^{\ast}h$ are mutually bounded.
The curvature $R(h_0)$ is $0$.
Let $\omegatilde$ denote the Poincar\'e metric of $\Xtilde-\Dtilde$:
\[
 \omegatilde=\sum_{i=1,2}
 \frac{d\zeta_i\cdot d\zetabar_i}
 {|\zeta_i|^2(-\log|\zeta_i|^2)^2}
\]

\begin{lem}
\label{lem;06.8.15.46}
$R(h_1)$ and $\del_{h_1}-\del_{h_0}$
are bounded with respect to 
$(\omegatilde,h_i)$ $(i=0,1)$.
\end{lem}
\pf
$\del\log\bigl(1+\rho(\zeta_2)\bigr)$,
$\del\delbar\log\bigl(1+\rho(\zeta_2)\bigr)$,
$\del\log(1+\rho(\zeta_1)+\rho(\zeta_2))$
and $\del\delbar\log(1+\rho(\zeta_1)+\rho(\zeta_2))$
are bounded
with respect to $\omegatilde$.
Then, the boundedness of $R(h_1)$
and $\del_{h_1}-\del_{h_0}$
follow.
\hfill\qed

\section{Preliminary from Elementary Calculus}

Take $\epsilon>0$ and $N>1$.
In this section,
we use the following volume form $\dvol_{\epsilon,N}$
of a punctured disc $\Delta^{\ast}$:
\[
  \dvol_{\epsilon,N}:=
 \bigl(\epsilon^{N+2}\cdot|z|^{2\epsilon}+|z|^2\bigr)^{-1}
 \frac{\sqrt{-1}dz\wedge d\bar{z}}{|z|^2}
\]
Let $f$ be a function on a punctured disc $\Delta^{\ast}$
such that
$\|f\|_{L^2}^2:=
\int_{\Delta^{\ast}} |f|^2\cdot\dvol_{\epsilon,N}<\infty$.
We use the polar coordinate $z=r\cdot e^{\sqrt{-1}\theta}$.
For the decomposition $f=\sum f_n(r)\cdot e^{\sqrt{-1}n\theta}$,
we have $\|f\|_{L^2}^2=2\pi\sum_n \|f_n\|_{L^2}^2$,
where $\|f_n\|_{L^2}^2$ are given as follows:
\[
 \|f_n\|_{L^2}^2:=
\int_0^1|f_n(\rho)|^2\cdot
 \bigl(\epsilon^{N+2}\rho^{2\epsilon}+\rho^2\bigr)^{-1}
\frac{d\rho}{\rho}.
\]

\begin{prop}
 \label{prop;05.8.25.15}
Let $f$ be as above.
Then we have a function $v$ satisfying the following:
\[
 \delbar\del v=f\cdot \frac{d\bar{z}\wedge dz}{|z|^2},
\quad
 \bigl|v(z)\bigr|
\leq C\cdot \Bigl(|z|^{\epsilon}\epsilon^{(N-1)/2}+|z|^{1/2}\Bigr)
\cdot \|f\|_{L^2}.
\]
The constant $C$ can be independent of $\epsilon$, $N$ and $f$.
\end{prop}
\pf
We use the argument of S. Zucker in \cite{z}.
First let us consider the equation
$\delbar u=f\cdot d\bar{z}/\bar{z}$.
For the decomposition
$u=\sum u_n(\rho)\cdot e^{\sqrt{-1}n\theta}$,
it is equivalent to the following equations:
\[
 \frac{1}{2}\left(
 r\frac{\del}{\del r}u_n-n\cdot u_n
 \right)=f_n,\quad
(n\in\seisuu).
\]
We put as follows:
\[
 u_n:=\left\{
 \begin{array}{ll}
 2r^n\int_0^r\rho^{-n-1}f_n(\rho)\cdot d\rho
 & (n\leq 0),\\
 \mbox{{}}\\
 2r^n\int_A^r\rho^{-n-1}f_n(\rho)\cdot d\rho
 & (n>0).
\end{array}
 \right.
\]
Then $u=\sum u_n\cdot e^{\sqrt{-1}n\theta}$ satisfies
the equation $\delbar u=f\cdot d\bar{z}/\bar{z}$.

\begin{lem}\label{lem;05.8.25.10}
There exists $C_1>0$ such that
\[
 |u_n(r)|\leq C_1\cdot \|f_n\|_{L^2}
\cdot\left(
 \frac{\epsilon^{(N+2)/2}\cdot r^{\epsilon}}{|2\epsilon-2n|^{1/2}}
+\frac{r^{1/2}}{(1+|n|)^{1/2}}
 \right).
\]
The constant $C_1$ is independent of $n$, $\epsilon$, $N$ and $f$.
\end{lem}
\pf
In the case $n\leq 0$, we have the following:
\begin{multline}
  |u_n(r)|
\leq 
 2r^{n}\left(
 \int_0^r|f_n(\rho)|^2(\epsilon^{N+2}\rho^{2\epsilon}+\rho^2)^{-1}
 \frac{d\rho}{\rho}
 \right)^{1/2}
 \\
 \times
 \left(\int_0^r\rho^{-2n-1}(\epsilon^{N+2}\rho^{2\epsilon}+\rho^2)
 \cdot d\rho
 \right)^{1/2}
\end{multline}
We have the following:
\[
 \int_0^r\rho^{-2n-1}(\epsilon^{N+2}\rho^{2\epsilon}+\rho^2)d\rho
=\frac{\epsilon^{N+2}\cdot r^{2\epsilon-2n}}{2\epsilon-2n} 
+\frac{r^{-2n+2}}{-2n+2}.
\]
Hence we obtain the following:
\[
 |u_n(r)|\leq
 2\|f_n\|_{L^2}\cdot 
 \left(
 \frac{\epsilon^{(N+2)/2}\cdot r^{\epsilon}}{|2\epsilon-2n|^{1/2}}
+\frac{r}{|2-2n|^{1/2}}
 \right).
\]

In the case $n>0$, we also have the following:
\[
   |u_n(r)|\leq 
 2r^n\cdot
\|f_n\|_{L^2}
 \left|
 \int_A^r\rho^{-2n-1}(\epsilon^{N+2}\rho^{2\epsilon}+\rho^2)d\rho
\right|^{1/2}.
\]
We have the following:
\[
 \left|
\int_A^r\rho^{-2n-1}\epsilon^{N+2}\cdot\rho^{2\epsilon}\cdot d\rho
 \right|
\leq \frac{\epsilon^{N+2}}{|-2n+2\epsilon|}r^{-2n+2\epsilon}.
\]
We also have the following:
\[
 \int_A^r\rho^{-2n+1}d\rho
=\left\{
\begin{array}{ll}
 \log r-\log A & (n=1)\\
 \mbox{{}}\\
 (-2n+2)^{-1}(r^{-2n+2}-A^{-2n+2}) & (n\geq 2)
\end{array}
 \right.
\]
Therefore we obtain the following:
\[
 |u_n(r)|\leq
 C\cdot \|f_n\|_{L^2}
\left(
 \frac{\epsilon^{(N+2)/2}\cdot r^{\epsilon}}{|2\epsilon-2n|^{1/2}}
+\frac{r^{1/2}}{(1+|n|)^{1/2}}
 \right)
\]
Thus we are done.
\hfill\qed

\vspace{.1in}

Then let us consider the equation $\del v=u\cdot dz/z$.
For the decomposition $v=\sum v_n\cdot e^{\sqrt{-1}n\theta}$,
it is equivalent to the following equations:
\[
 \frac{1}{2}\left(
 r\frac{\del v_n}{\del r}+n\cdot v_n
 \right)=u_n,
\quad (n\in\seisuu).
\]
We put as follows:
\[
 v_n(r):=\left\{
 \begin{array}{ll}
 2r^{-n}\cdot\int_0^r\rho^{n-1}u_n(\rho)\cdot d\rho & (n\geq 0)\\
 \mbox{{}}\\
 2r^{-n}\cdot\int_A^r\rho^{n-1}u_n(\rho)\cdot d\rho & (n<0).
 \end{array}
 \right.
\]
Then we have $\del v=u\cdot dz/z$ for
$v:=\sum v_n\cdot e^{\sqrt{-1}n\theta}$.
From Lemma \ref{lem;05.8.25.10},
we obtain the following in the case $n>0$:
\begin{multline}
 |v_n(r)|\leq 2r^{-n}
 \int_0^r\rho^{n-1}
 \left(\frac{\epsilon^{(N+1)/2}\cdot\rho^{\epsilon}}{|2\epsilon-2n|^{1/2}}
+\frac{\rho^{1/2}}{(1+|n|)^{1/2}}
 \right)
 d\rho\cdot \|f_n\|_{L^2} \\
\leq
 C_2\cdot \|f_n\|_{L^2}
\cdot
 \left(
\frac{\epsilon^{(N+2)/2}}{|2\epsilon-2n|^{1/2}}
\frac{r^{\epsilon}}{|n+\epsilon|}
+\frac{1}{(1+|n|)^{1/2}}\frac{r^{1/2}}{n+1/2}
 \right).
\end{multline}
We have a similar estimate in the case $n<0$.
Hence we obtain the following:
\[
 |v(z)|
\leq
 \sum_{n}|v_n(r)|
\leq 
 C_4\cdot(\epsilon^{(N-1)/2}r^{\epsilon}+r^{1/2})\cdot 
\|f\|_{L^2}.
\]
Thus the proof of Proposition \ref{prop;05.8.25.15} is finished.
\hfill\qed

\section{Reflexive Sheaf}
\label{section;06.8.27.1}

We recall some general facts about 
reflexive sheaves.
See \cite{har2} and \cite{mehta-ramanathan1}
for some more properties of reflexive sheaves.
Let $X$ be a complex manifold.
Recall that a coherent $\nbigo_X$-module $\nbige$
is called reflexive,
if $\nbige$ is isomorphic to the double dual 
$\nbige^{\lor\lor}:=
\nhom\bigl(\nhom(\nbige,\nbigo_X),\nbigo_X\bigr)$
of $\nbige$.
Recall we can take a resolution {\em locally} on $X$
(Lemma 3.1 of \cite{mehta-ramanathan1}):
\begin{equation}
\label{eq;06.8.12.30}
 0\lrarr \nbige\lrarr \nbigv_0\lrarr \nbigv_1\lrarr 0
\end{equation}
Here $\nbigv_0$ is locally free
and $\nbigv_1$ is torsion-free.
The following Hartogs type theorem is well known.
\begin{lem}
Let $Z$ be a closed subset of $X$
whose codimension is larger than $2$.
Let $f$ be a section of 
a reflexive sheaf $\nbige$ on $X\setminus Z$.
Then $f$ is naturally extended to
the section of $\nbige$ over $X$.
\end{lem}
\pf
We have only to check the claim locally.
Let us take a resolution (\ref{eq;06.8.12.30}),
and then $f$ induces the section of $\widetilde{f}$
of $\nbigv_0$ on $X-Z$.
Due to the Hartogs' theorem,
$\ftilde$ can be extended to the section 
on $X$.
Since it is mapped to $0$ in $\nbigv_1$,
we obtain the section of $\nbige$
on $X$.
\hfill\qed

\vspace{.1in}
The converse is also true.
\begin{lem}
\label{lem;06.8.12.35}
Let $\nbigf$ be a torsion-free coherent sheaf on $X$
such that any section $f$ of $\nbigf$ on $U-Z$
is extended to the section on $U$,
where $U$ denotes an open subset 
and $Z$ denotes a closed subset with $\codim Z\geq 2$.
Then $\nbigf$ is reflexive.
\end{lem}
\pf
We have the inclusion $\iota:\nbigf\lrarr\nbigf^{\lor\lor}$,
which is isomorphic outside of the subset $Z_0\subset X$
with $\codim(Z_0)\geq 2$.
Then, we obtain the surjectivity of $\iota$
from the given property of $\nbigf$,
and thus $\iota$ is isomorphic.
\hfill\qed

\begin{lem}
If $\nbige$ is reflexive,
$\nbige\otimes\nbigo_{D}$ is torsion-free
for a divisor $D$.
\end{lem}
\pf
Take a resolution as in (\ref{eq;06.8.12.30}).
Because of $\Tor^1(\nbigv_1,\nbigo_{D})=0$,
we obtain the injection
$\nbige\otimes\nbigo_{D}\lrarr \nbigv_0\otimes\nbigo_{D}$,
and hence $\nbige\otimes\nbigo_{D}$ is torsion-free.
\hfill\qed

\begin{lem}
\label{lem;06.8.24.1}
If $\nbige$ is a reflexive sheaf,
$\nhom(\nbigf,\nbige)$ is also reflexive
for any coherent sheaf $\nbigf$.
\end{lem}
\pf
Let us check the condition in Lemma \ref{lem;06.8.12.35}.
Let $U$ be a small open subset,
on which we have a resolution
$\nbigv_{-1}\stackrel{a}\lrarr
 \nbigo_U^{\oplus r}
 \stackrel{b}{\lrarr} \nbigf\lrarr 0$ on $U$.
Let $f$ be a homomorphism
$\nbigf\lrarr\nbige$ on $U\setminus Z$,
where $\codim Z\geq 2$.
The morphism $\nbigo_U^{\oplus r}\lrarr \nbige$ is naturally induced
on $U\setminus Z$,
which is naturally extended to the morphism 
$\varphi:\nbigo_U^{\oplus r}\lrarr \nbige$ on $U$
by the Hartogs property.
Since $\varphi\circ a$ is $0$,
$\varphi$ induces the extension of $f$.
\hfill\qed

\section{Moduli Spaces of Representations}
\label{subsection;05.9.12.1}

Let $\Gamma$ be a finitely presented group,
and $V$ be a finite dimensional vector space over $\cnum$.
For $a,f\in\GL(V)$,
we put $\ad(a)(f):=a\circ f\circ a^{-1}$.
The space of homomorphisms $R(\Gamma,V):=Hom(\Gamma,\GL(V))$
is naturally an affine variety over $\cnum$.
We regard it as a Hausdorff topological space with the usual topology,
not the Zariski topology.
We have the natural action of $\GL(V)$
on $R(\Gamma,V)$ given by $\ad$.
Let $h_V$ be a hermitian metric of $V$,
and let $U(h_V)$ denote the unitary group of $V$
with respect to $h_V$.
The usual quotient space $R(\Gamma,V)\big/U(h_V)$
is denoted by $M(\Gamma,V,h_V)$.
Let $\pi_{\GL(V)}$ denote the projection
$R(\Gamma,V)\lrarr M(\Gamma,V,h_V)$.

More generally,
we consider the moduli space of representations
to a complex reductive subgroup $G$ of $\GL(V)$.
We put $R(\Gamma,G):=\Hom(\Gamma,G)$,
which we regard as a Hausdorff topological space
with the usual topology.
It is the closed subspace of $R(\Gamma,V)$.

Let $K$ be a maximal compact subgroup of $G$.
Assume that the hermitian metric $h_V$ of $V$ is $K$-invariant.
We put 
$N_G\bigl(h_V\bigr):=
 \bigl\{u\in U(h_V)\,\big|\,\ad(u)(G)=G \bigr\}$
which is compact.
We have the natural adjoint action of $N_G\bigl(h_V\bigr)$ on $G$,
which induces the action on $R(\Gamma,G)$.
The usual quotient space is denoted by $M(\Gamma,G,h_V)$.
Let $\pi_G$ denote the projection
$R(\Gamma,G)\lrarr M(\Gamma,G,h_V)$.
We have the naturally defined map
$\Phi:M(\Gamma,G,h_V)\lrarr M(\Gamma,V,h_V)$.
The map $\Phi$ is clearly proper
in the sense that the inverse image of any compact subset via $\Phi$
is also compact.

A representation $\rho\in R(\Gamma,G)$ is called Zariski dense,
if the image of $\rho$ is Zariski dense in $G$.
Let $\nbigu$ be the subset of $R(\Gamma,G)$,
which consists of Zariski dense representations.
Then the restriction of $\Phi$ to $\nbigu$ is injective.

\vspace{.1in}

Let $\rho$ and $\rho'$ be elements of $R(\Gamma,G)$.
We say that $\rho$ and $\rho'$ are isomorphic in $G$,
if there is an element $g\in G$ such that
$\ad(g)\circ \rho=\rho'$.
We say $\rho'$ is a deformation of $\rho$ in $G$,
if there is a continuous family of representations
$\rho_t:[0,1]\times\Gamma\lrarr G$
such that $\rho_0=\rho$ and $\rho_1=\rho'$.
We say $\rho'$ is a deformation of $\rho$ in $G$
modulo $N_G(h_V)$,
if there is an element $u\in N_G(h_V)$ such that
$\rho$ can be deformed to $\ad(u)\circ\rho'$
in $G$.
The two notions are different
if $N_G(h_V)$ is not connected,
in general.
We also remark that 
$\rho$ can be deformed to $\rho'$ in $G$ modulo $N_G(h_V)$,
if and only if
$\pi_G(\rho)$ and $\pi_G(\rho')$
are contained in the same connected component
of $M(\Gamma,G,h_V)$.

\vspace{.1in}

We recall some deformation invariance from \cite{s5}.
A representation $\rho\in R(\Gamma,G)$ is called rigid,
if the orbit $G\cdot\rho$ is open in $R(\Gamma,G)$.
\begin{lem}
 \label{lem;05.9.10.1}
Let $\rho\in R(\Gamma,G)$ be  a rigid and Zariski dense representation.
Then any deformation $\rho'$ of $\rho$ in $G$
is isomorphic to $\rho$ in $G$.
\end{lem}
\pf
If $\rho$ is Zariski dense,
then $G\cdot\rho$ is closed in $R(\Gamma,G)$.
Hence it is a connected component.
\hfill\qed

%% file: 3.tex
We recall the notion of parabolic structure,
and then we give some detail about the characteristic numbers
for parabolic sheaves.
In Section \ref{section;05.7.30.15},
a perturbation of the filtration is given,
which will be useful in our later argument.

\section{Parabolic Higgs Bundle}
\label{section;05.8.23.1}

\subsection{$\vecc$-Parabolic Higgs sheaf}

Let us recall the notion of parabolic structure
and the Chern characteristic numbers of parabolic bundles
following 
\cite{li2}, \cite{my}, \cite{s1}, \cite{s2},
\cite{steer-wren}
and \cite{y}.
Our convention is slightly different from theirs.

Let $X$ be a connected complex manifold
and $D$ be a simple normal crossing divisor
with the irreducible decomposition $D=\bigcup_{i\in S} D_i$.
Let $\vecc=(c_i\,\big|\,i\in S)$ be an element of $\real^S$.
Let $\nbige$ be a torsion-free coherent $\nbigo_X$-module.
Let us consider
a collection of the increasing filtrations
$\lefttop{i}\nbigf$ $(i\in S)$
indexed by $\openclosed{c_i-1}{c_i}$ such that
$\lefttop{i}\nbigf_{a}(\nbige)\supset \nbige(-D_i)$
for any $a\in\openclosed{c_i-1}{c_i}$.
We put
$\lefttop{i}\Gr^{\nbigf}_{a}\nbige
 :=\lefttop{i}\nbigf_{a}(\nbige)
 \big/\lefttop{i}\nbigf_{<a}(\nbige)$.
We assume that
the sets
$\bigl\{a\,\big|\,\lefttop{i}\Gr^{\nbigf}_{a}\nbige\neq 0\bigr\}$
are finite for any $i$.
Such tuples of filtrations are called 
the $\vecc$-parabolic structure of $\nbige$ at $D$,
and the tuple $\bigl(\nbige,\{\lefttop{i}\nbigf\,|\,i\in S\}\bigr)$
is called a $\vecc$-parabolic sheaf on $(X,D)$.
We will sometimes omit to denote $\vecc$.
We say $\bigl(\nbige,\{\lefttop{i}\nbigf\,|\,i\in S\}\bigr)$
is reflexive, if $\nbige$ is reflexive.
(See \cite{har2} and \cite{mehta-ramanathan1}
for reflexive sheaves.
See also Section \ref{section;06.8.27.1}.)
\begin{df}
For a reflexive $\vecc$-parabolic sheaf
 $\bigl(\nbige,\{\lefttop{i}\nbigf\,|\,i\in S\}\bigr)$,
we say that the parabolic structure is saturated,
if $\nbige/\lefttop{i}\nbigf_a$ are torsion-free
$\nbigo_{D_i}$-modules
for any $i$ and $a$.
\hfill\qed
\end{df}

We remark that each $\lefttop{i}\nbigf_a$
are also reflexive.
To see it, let us see the inclusion
$\lefttop{i}\nbigf_a\lrarr\lefttop{i}\nbigf_a^{\lor\lor}$.
Since $\nbige$ is reflexive,
the inclusion $\lefttop{i}\nbigf_a\lrarr \nbige$
is extended to the injection
$\lefttop{i}\nbigf_a^{\lor\lor}\lrarr \nbige$.
(See the proof of Lemma \ref{lem;06.8.24.1}.)
Hence we obtain the inclusion
$\lefttop{i}\nbigf_a^{\lor\lor}/\lefttop{i}\nbigf_a
\lrarr \nbige/\lefttop{i}\nbigf_a$.
The codimension of the support of
$\lefttop{i}\nbigf_a^{\lor\lor}/\lefttop{i}\nbigf_a$
is larger than $2$,
and $\nbige/\lefttop{i}\nbigf_a$ is
torsion-free as an $\nbigo_{D_i}$-module.
Hence we obtain
$\lefttop{i}\nbigf_a^{\lor\lor}/\lefttop{i}\nbigf_a=0$

\vspace{.1in}

We will use the notation ${\nbige}_{\ast}$
instead of $\bigl(\nbige,\{\lefttop{i}\nbigf\}\bigr)$
for simplicity.
When we emphasize $\vecc$,
we will often use the notation $\prolongg{\vecc}{\nbige}$
and $\prolongg{\vecc}{\nbige}_{\ast}$
instead of $\nbige$ and $\nbige_{\ast}$.
In the case $\vecc=(0,\ldots,0)$,
the notation $\prolong{\nbige}_{\ast}$ is used.
We will also use the following notation.
\begin{equation}
 \label{eq;05.9.14.1}
 \Par\bigl(\nbige_{\ast},i\bigr):=
 \bigl\{
 a\,\big|\,
 \lefttop{i}\Gr^{\nbigf}_a(\nbige)\neq 0
 \bigr\},
\quad
 \Par'\bigl(\nbige_{\ast},i\bigr):=
 \Par\bigl(\nbige_{\ast},i\bigr)
\cup\{c_i,c_i-1\},
\end{equation}
\begin{equation}
 \label{eq;05.9.14.2}
 \gap (\nbige_{\ast},i):=
 \min\bigl\{
 |a-b|\,\big|\,
 a,b\in\Par'(\nbige_{\ast},i),\,\,
 a\neq b
 \bigr\},
\quad
 \gap(\nbige_{\ast}):=
 \min_{i\in S}
 \gap\bigl(\nbige_{\ast},i\bigr).
\end{equation}

Let us recall a Higgs field (\cite{y}) of a $\vecc$-parabolic sheaf
on $(X,D)$.
A holomorphic homomorphism
$\theta:\nbige\lrarr \nbige\otimes\Omega^{1,0}_X(\log D)$
is called a Higgs field of $\nbige_{\ast}$,
if the following holds:
\begin{itemize}
\item 
The naturally defined composite
$\theta^2=\theta\wedge\theta:
 \nbige\lrarr\nbige\otimes\Omega^{2,0}_X(\log D)$
 vanishes.
\item
$\theta\bigl(\lefttop{i}\nbigf_{a}\bigr)
\subset \lefttop{i}\nbigf_a\otimes\Omega^{1,0}_X(\log D)$
\end{itemize}
Such a tuple $(\nbige_{\ast},\theta)$ is called
a $\vecc$-parabolic Higgs sheaf on $(X,D)$.

A $\vecc$-parabolic Higgs sheaf $(\nbige_{\ast},\theta)$ on $(X,D)$
is called reflexive and saturated,
if the underlying $\vecc$-parabolic sheaf is reflexive and saturated.
A morphism between $\vecc$-parabolic Higgs sheaves
is defined to be a morphism of the underlying sheaf
which is compatible with the parabolic structures
and the Higgs fields.

\begin{lem}
\label{lem;06.8.5.40}
Let $(\nbige_{\ast},\theta)$ be any $\vecc$-parabolic Higgs sheaf
on $(X,D)$.
Then there exists the reflexive saturated
parabolic Higgs sheaf
$(\nbige'_{\ast},\theta')$,
such that we have the morphism
$(\nbige_{\ast},\theta)\lrarr (\nbige_{\ast}',\theta')$
which is isomorphic in codimension one,
i.e. isomorphic outside of the subset with codimension two.
Such $(\nbige'_{\ast},\theta')$ is unique up to the canonical
isomorphism.
\end{lem}
\pf
Let $\nbige'$ denote the double dual of $\nbige$.
We have the canonical morphism
$\nbige\lrarr\nbige'$ which is isomorphic
outside of the subset $Z$ of codimension two.
Let $\lefttop{i}\nbigf^1_a$ denote the subsheaf of $\nbige'$
which consists of the sections $f$ of $\nbige'$
such that $f_{|X-Z}\in \lefttop{i}\nbigf_a$.
Such a subsheaf is coherent (\cite{siu-trautmann}).
We have $\nbige'(-D_i)\subset \lefttop{i}\nbigf^1_a$
for any $a\in \openclosed{c_i-1}{c_i}$.
We have the natural surjection
$\pi_{i,a}:\nbige'\lrarr\nbige'/\lefttop{i}\nbigf^1_a$,
and the target is the $\nbigo_{D_i}$-module.
Let $T_{i,a}$ denote the torsion part of
$\nbige'/\lefttop{i}\nbigf^1_a$ as an $\nbigo_{D_i}$-module,
and we put $\lefttop{i}\nbigf'_a:=\pi_{i,a}^{-1}(T_{i,a})$.
Then, it is easy to see that
$\bigl\{\lefttop{i}\nbigf'\,\big|\,i\in S\bigr\}$
gives the saturated $\vecc$-parabolic structure of
$\nbige'$.
The Higgs field $\theta$ naturally induces
the morphism
$\nbige\lrarr\nbige'\otimes\Omega_X^{1,0}(\log D)$.
Due to the reflexivity of $\nbige'$,
we obtain 
$\theta':\nbige'\lrarr\nbige'\otimes\Omega_X^{1,0}(\log D)$
satisfying $\theta^2=0$.
It is easy to check
$\theta(\lefttop{i}\nbigf'_a)\subset 
 \lefttop{i}\nbigf'_a\otimes\Omega^{1,0}_X(\log D)$.
The uniqueness is clear.
\hfill\qed

\vspace{.1in}

For a $\vecc$-parabolic Higgs sheaves $(\nbige_{i\,\ast},\theta_i)$
$(i=1,2)$ on $(X,D)$,
we obtain the sheaf of the morphisms
$\nhom\bigl((\nbige_{1\,\ast},\theta_1),(\nbige_{2\,\ast},\theta_2)\bigr)$.
\begin{lem}
\label{lem;06.8.22.1}
If $(\nbige_{2\,\ast},\theta_2)$ is reflexive and saturated,
$\nhom\bigl((\nbige_{1\,\ast},\theta_1),(\nbige_{2\,\ast},\theta_2)\bigr)$
is reflexive.
\end{lem}
\pf
We have only to check
the condition in Lemma \ref{lem;06.8.12.35}.
Let $f$ be a section of
$\nhom\bigl((\nbige_{1\,\ast},\theta_1),
 (\nbige_{2\,\ast},\theta_2) \bigr)$
on $U\setminus Z$,
where $U$ denotes an open subset
and $Z$ denotes a closed subset with $\codim(Z)\geq 2$.
Since $\nbige_2$ is reflexive,
it is extended to the homomorphism
$\widetilde{f}:\nbige_{1}\lrarr\nbige_2$ on $U$,
which is compatible with $\theta_i$.
We have the induced map
$\varphi:\lefttop{i}\nbigf(\nbige_1)\lrarr
 \nbige_2/\lefttop{i}\nbigf(\nbige_2)$.
The codimension of the support of $\Image(\varphi)$
is larger than $2$,
and $\nbige_2/\lefttop{i}\nbigf(\nbige_2)$
is a torsion-free $\nbigo_{D_i}$-module.
Hence, we obtain $\Image(\varphi)=0$,
i.e.,
$\widetilde{f}$ preserves the filtration.
\hfill\qed

\vspace{.1in}
Assume $X$ is projective.
Let $Y$ be a sufficiently ample  and generic hypersurface of $X$.
We put $D_Y:=D\cap Y$, which is assumed to be
a simple normal crossing divisor of $Y$.
Let $(\nbige_{i\,\ast|Y},\theta_{iY})$ denote the induced
parabolic Higgs sheaf on $(Y,D_Y)$ by $(\nbige_{i\,\ast},\theta_i)$.
If $\nbige_{i\,\ast}$ is reflexive and saturated,
so is $\nbige_{i\,\ast|Y}$.
(See Corollary 3.1.1 of \cite{mehta-ramanathan1}.)

\begin{lem}
\label{lem;06.8.22.2}
Assume $\dim X\geq 2$ and 
that $\nbige_{2\,\ast}$ is saturated and reflexive.
For any morphism $f:(\nbige_{1\,\ast|Y},\theta_{1\,Y})
 \lrarr (\nbige_{2\,\ast|Y},\theta_{2\,Y})$,
we have $F:(\nbige_{1\,\ast},\theta_1)\lrarr (\nbige_{2\,\ast},\theta_2)$
which induces $f$.
\end{lem}
\pf
Let $\theta_{i|Y}:\nbige_{i\ast|Y}
 \lrarr \nbige_{i\ast|Y}\otimes\Omega^{1,0}_X(\log D)_{|Y}$ 
denote the restriction of $\theta_i$ to $Y$.
We have the induced morphism
$G:f\circ\theta_{1|Y}-\theta_{2|Y}\circ f:
 \nbige_{1\ast|Y}\lrarr
 \nbige_{2\ast|Y}\otimes\Omega^{1,0}_{X}(\log D)_{|Y}$.
Because of 
$f\circ \theta_{1\,Y}-\theta_{2\,Y}\circ f=0$
in $\nhom\bigl(\nbige_{1\,\ast|Y},\nbige_{2\,\ast|Y}\bigr)
 \otimes\Omega^{1,0}_Y(\log D_Y)$,
$G$ induces the map
$\nbige_{1\,\ast|Y}\lrarr \nbige_{2\,\ast|Y}\otimes\nbigo(-Y)_{|Y}$.
We regard it as the section of
$\nbigj:=\nhom\bigl(\nbige_{1\,\ast},\nbige_{2\,\ast}\bigr)
 \otimes\nbigo(-Y)_{|Y}$.
Since $\nbigg:=\nhom\bigl(\nbige_{1\,\ast},\nbige_{2\,\ast}\bigr)$ is
reflexive,
we have $H^i\bigl(X,\nbigg\otimes\nbigo(-Y)\bigr)=0$
 $(i=0,1)$, if $Y$ is sufficiently ample.
(See the proof of Proposition 3.2 in \cite{mehta-ramanathan1}.)
Hence, we have $H^0(Y,\nbigj)=0$,
i.e., $G=0$.
Then, the claim of the lemma follows from
Generalized Enriques Severi Lemma 
(Proposition 3.2 in \cite{mehta-ramanathan1})
and Lemma \ref{lem;06.8.22.1}.
\hfill\qed

\begin{rem}
We also have the parallel notion of $\vecc$-parabolic sheaves
on smooth varieties with simple normal crossing
divisors over a field $k$.
\hfill\qed
\end{rem}

\begin{rem}
\label{rem;06.8.15.101}
Sometimes,
it will be convenient to consider filtrations
$\lefttop{i}\nbigf$
such that
$S(\lefttop{i}\nbigf)=
 \bigl\{a\in\real\,\big|\,\lefttop{i}\Gr^{\nbigf}_a(\nbige)\neq 0\bigr\}$
is not contained in an interval $\openclosed{c_i-1}{c_i}$
for some $c_i$.
In that case, we will call
$\{\lefttop{i}\nbigf\,|\,i\in S\}$
a generalized parabolic structure.
Higgs field is also defined as in the standard case,
i.e.,
a holomorphic map
$\theta:\nbige\lrarr\nbige\otimes\Omega^{1,0}_X(\log D)$
such that $\theta^2=0$
and $\theta\bigl(\lefttop{i}\nbigf_a\bigr)
\subset\lefttop{i}\nbigf_a\otimes\Omega_X^{1,0}(\log D)$.
\hfill\qed
\end{rem}

\subsection{The parabolic first Chern class and   the degree}
\label{subsection;06.8.4.10}

For a $\vecc$-parabolic sheaf $\nbige_{\ast}$ on $(X,D)$,
we put as follows:
\[
 \wt(\nbige_{\ast},i):=
 \sum_{a\in\openclosed{c_i-1}{c_i}}
 a\cdot\rank_{D_i} \lefttop{i}\Gr^{\nbigf}_a(\nbige).
\]
Here $\rank_{D_i}\lefttop{i}\Gr^{\nbigf}_a(\nbige)$
denotes the rank as an $\nbigo_{D_i}$-module.
In the following, we will often denote it
by $\rank\lefttop{i}\Gr^{\nbigf}_a(\nbige)$,
if there are no risk of confusion.
The parabolic first Chern class of $\nbige_{\ast}$
is defined as follows:
\[
 \parchern_1(\nbige_{\ast}):=
 c_1(\nbige)-\sum_{i\in S}\wt (\nbige_{\ast},i)\cdot [D_i]
\in H^2(X,\real).
\]
Here $[D_i]$ denotes the cohomology class given by $D_i$.
If $X$ is an $n$-dimensional compact Kahler manifold with a Kahler form $\omega$,
we put as follows:
\[
 \pardeg_{\omega}(\nbige_{\ast}):=
\int_X \parchern_1(\nbige_{\ast})\cdot\omega^{n-1},
\quad
 \mu_{\omega}(\nbige_{\ast}):=
 \frac{\pardeg_{\omega}(\nbige_{\ast})}{\rank \nbige}.
\]
If $\omega$ is the first Chern class of an ample line bundle $L$,
we also use the notation
$\pardeg_L(\nbige_{\ast})$
and $\mu_L(\nbige_{\ast})$.

\begin{lem}
\label{lem;06.8.24.2}
Let $\nbige^{(i)}_{\ast}$ $(i=1,2)$ be 
$\vecc$-parabolic sheaves
on $(X,D)$,
and let $f:\nbige^{(1)}_{\ast}\lrarr \nbige^{(2)}_{\ast}$
be a morphism which is generically isomorphic.
Then, we have 
$\mu(\nbige^{(1)}_{\ast})\leq \mu(\nbige^{(2)}_{\ast})$.
If the equality occurs,
$f$ is isomorphic in codimension one.
\end{lem}
\pf
By considering the restriction to
a generic complete intersection curve,
we have only to discuss the case $\dim X=1$.
Let $P$ be any point of $D$.
We put
$F_a^{(i)}:=
\Image\bigl(\lefttop{P}\nbigf_a(\nbige^{(i)})\bigr)_{|P}
\lrarr \nbige^{(i)}_{|P}$ for
$a\in \openclosed{c(P)-1}{c(P)}$,
which give the filtration $F^{(i)}$
of $\nbige^{(i)}_{|P}$.
We have the induced map
$f_{|P}:\nbige^{(1)}_{|P}\lrarr \nbige^{(2)}_{|P}$
which preserve the filtrations.
We put $I:=\Image(f_{|P})$,
$K:=\Ker(f_{|P})$ and $C:=\Cok(f_{|P})$.
Let $F(K)$ (resp. $F^{(1)}(I)$)
denote the induced filtration on $K$  (resp. $I$)
by $F^{(1)}$.
Let $F(C)$ (resp. $F^{(2)}(I)$)
denote the induced filtration on $C$ (resp. $I$)
by $F^{(2)}$.
We put as follows:
\[
 w(K):=\sum a\cdot \Gr^{F}_a(K),
\quad
 w^{(i)}(I):=\sum a\cdot \Gr^{F^{(i)}}_a(I),
\quad
 w(C):=\sum a\cdot \Gr^{F}_a(C)
\]
Then, we have
$-w^{(1)}(I)\leq -w^{(2)}(I)$
and $-w(K)<-w(C)+r_0$,
where $r_0=\rank K=\rank C$.
It is easy to obtain the claims of the lemma
from these relations.
\hfill\qed

\begin{rem}
For the parabolic first Chern class on 
algebraic varieties,
we have only to replace the cohomology group 
and the integral
by the Chow group and 
the degree of the $0$-cycles.
\hfill\qed
\end{rem}

\subsection{$\mu_L$-Stability}

Let $X$ be a smooth projective variety with an ample line bundle $L$
over a field $k$,
and $D$ be a simple normal crossing divisor of $X$.
The $\mu_L$-stability of 
$\vecc$-parabolic Higgs sheaves is defined as usual.
Namely, a $\vecc$-parabolic Higgs sheaf $\bigl(\nbige_{\ast},\theta\bigr)$
is called $\mu_L$-stable, if the inequality 
$\pardeg_L(\nbige'_{\ast})<\pardeg_L(\nbige_{\ast})$ holds
for any saturated non-trivial subsheaf $\nbige'\subsetneq \nbige$
such that $\theta(\nbige')\subset \nbige'\otimes\Omega^{1,0}(\log D)$.
(Recall a subsheaf $\nbige'\subset\nbige$ is called saturated,
if $\nbige/\nbige'$ is torsion-free.)
Here the parabolic structure of $\nbige'_{\ast}$
is the naturally induced one from 
the parabolic structure of $\nbige_{\ast}$.
Similarly, $\mu_L$-semistability and $\mu_L$-polystability
are also defined in a standard manner.

\vspace{.1in}

Let $\bigl(\nbige^{(i)}_{\ast},\theta^{(i)}\bigr)$ $(i=1,2)$
be $\mu_L$-semistable $\vecc$-parabolic Higgs sheaves
such that $\mu_L(\nbige^{(1)}_{\ast})=\mu_L(\nbige^{(2)}_{\ast})$.
Let $f:(\nbige^{(1)}_{\ast},\theta^{(1)})
 \lrarr (\nbige^{(2)}_{\ast},\theta^{(2)})$ be a non-trivial morphism.
Let $(\nbigk_{\ast},\theta_{\nbigk})$ 
denote the kernel 
of $f$ with the naturally induced parabolic structure
and the Higgs field.
Let $\nbigi$ denote the image of $f$,
and $\nbigitilde$ denote the saturated subsheaf of $\nbige^{(2)}$
generated by $\nbigi$.
The parabolic structures of $\nbige^{(1)}_{\ast}$ 
and $\nbige^{(2)}_{\ast}$
induce the parabolic structures of
$\nbigi$ and $\nbigitilde$, respectively.
We denote the induced parabolic sheaves by
$(\nbigi_{\ast},\theta_{\nbigi})$
and $(\nbigitilde_{\ast},\theta_{\nbigitilde})$.

\begin{lem}
\label{lem;06.8.4.3}
$(\nbigk_{\ast},\theta_{\nbigk})$,
$(\nbigi_{\ast},\theta_{\nbigi})$
and $(\nbigitilde_{\ast},\theta_{\nbigitilde})$
are also $\mu_L$-semistable
such that
$\mu_L(\nbigk_{\ast})=\mu_L(\nbigi_{\ast})
=\mu_L(\nbigitilde_{\ast})
=\mu_L(\nbige^{(i)}_{\ast})$.
Moreover,
 $\nbigi_{\ast}$ and $\nbigitilde_{\ast}$
are isomorphic in codimension one.
\end{lem}
\pf
Using Lemma \ref{lem;06.8.24.2} and $\mu_L$-semistability
of $(\nbige^{(i)}_{\ast},\theta^{(i)})$,
we have 
$ \mu(\nbige^{(1)}_{\ast})
\leq \mu(\nbigi_{\ast})
\leq \mu(\nbigitilde_{\ast})
\leq \mu(\nbige^{(2)}_{\ast})$.
Since the equalities hold,
the claim of the lemma follows.
\hfill\qed

\begin{lem}
\label{lem;06.8.22.20}
Let $\bigl(\nbige^{(i)}_{\ast},\theta^{(i)}\bigr)$ $(i=1,2)$
be $\mu_L$-semistable reflexive saturated
parabolic Higgs sheaves
such that $\mu_L(\nbige^{(1)}_{\ast})=\mu_L(\nbige^{(2)}_{\ast})$.
Assume either one of the following:
\begin{enumerate}
\item
One of $(\nbige^{(i)}_{\ast},\theta^{(i)})$ is $\mu_L$-stable,
and $\rank(\nbige^{(1)})=\rank(\nbige^{(2)})$ holds.
\item
Both of $(\nbige^{(i)}_{\ast},\theta^{(i)})$ are $\mu_L$-stable.
\end{enumerate}
If there is a non-trivial map 
$f:(\nbige^{(1)}_{\ast},\theta^{(1)})
 \lrarr
 (\nbige^{(2)}_{\ast},\theta^{(2)})$,
then $f$ is isomorphic.
\end{lem}
\pf
If $(\nbige^{(1)}_{\ast},\theta^{(1)})$ is $\mu_L$-stable,
the kernel of $f$ is trivial due to Lemma \ref{lem;06.8.4.3}.
If $(\nbige^{(2)}_{\ast},\theta^{(2)})$ is $\mu_L$-stable,
the image of $f$ and $\nbige^{(2)}$ are 
same at the generic point of $X$.
Thus, we obtain that $f$ is generically isomorphic
in any case.
Then, we obtain that $f$ is isomorphic
in codimension one, due to Lemma \ref{lem;06.8.24.2}.
Since both of $\nbige^{(i)}_{\ast}$
are reflexive and saturated,
we obtain that $f$ is isomorphic.
\hfill\qed

\begin{cor}
 \label{cor;05.9.14.5}
Let $(\nbige_{\ast},\theta)$ be a $\mu_L$-polystable 
reflexive saturated Higgs sheaf.
Then we have the unique decomposition:
\[
 (\nbige_{\ast},\theta)
=\bigoplus_{j} (\nbige^{(j)}_{\ast},\theta^{(j)})\otimes\cnum^{m(j)}.
\]
Here, $(\nbige^{(j)}_{\ast},\theta^{(j)})$ are $\mu_L$-stable
with $\mu_L(\nbige^{(j)}_{\ast})=\mu(\nbige_{\ast})$,
and they are mutually non-isomorphic.
It is called the canonical decomposition in the rest of the paper.
\hfill\qed
\end{cor}

\subsection{$\vecc$-Parabolic Higgs bundle in codimension $k$}
\label{subsection;05.10.10.10}

We will often use the notation $\prolongg{\vecc}{E}$ instead
of $\nbige$.
We put as follows, for each $i\in S$:
\[
 \lefttop{i}F_a\bigl(\prolongg{\vecc}{E}_{|D_i}\bigr):=
\Image\Bigl(
 \lefttop{i}\nbigf_a(\prolongg{\vecc}{E})_{|D_i}
\lrarr
 \prolongg{\vecc}{E}_{|D_i}
 \Bigr).
\]
The tuple
$\bigl(\lefttop{i}\nbigf\,\big|\,i\in S\bigr)$
can clearly be reconstructed from 
the tuple of the filtrations
$\vecF:=\bigl(\lefttop{i}F\,\big|\,i\in S\bigr)$.
Hence we will often consider
$\bigl(\prolongg{\vecc}{E},\vecF\bigr)$ instead of
$\bigl(\prolongg{\vecc}{E},\{\lefttop{i}\nbigf\,|\,i\in S\}\bigr)$,
when $\prolongg{\vecc}{E}$ is locally free.

\begin{df}
Let $\prolongg{\vecc}{E}_{\ast}=(\prolongg{\vecc}{E},\vecF)$ be a
 $\vecc$-parabolic sheaf
such that $\prolongg{\vecc}{E}$ is locally free.
If the following conditions are satisfied,
$\prolongg{\vecc}{E}_{\ast}$ is called a $\vecc$-parabolic bundle.
\begin{itemize}
\item
 Each $\lefttop{i}F$ of $\prolongg{\vecc}{E}_{|D_i}$ is
 the filtration in the category of vector bundles on $D_{i}$.
 Namely,
 $\lefttop{i}\Gr_a^F(\prolongg{\vecc}{E}_{|D_i})=\lefttop{i}F_{a}\big/\lefttop{i}F_{<a}$
 are locally free $\nbigo_{D_i}$-modules.
\item
 The tuple of the filtrations $\vecF$ is compatible
 in the sense of Definition {\rm 4.37} in {\rm\cite{mochi2}}.
 (In this case, the decompositions are trivial.)
\end{itemize}
We remark that the second condition is trivial
in the case $\dim X=2$.
\hfill\qed
\end{df}

The notion of $\vecc$-parabolic Higgs bundle is too restrictive
in the case $\dim X>2$.
Hence we will also use the following notion in the case $k=2$.
\begin{df}
Let $\prolongg{\vecc}{E}_{\ast}$ be a $\vecc$-parabolic sheaf
on $(X,D)$.
It is called a $\vecc$-parabolic Higgs bundle
in codimension $k$,
if the following condition is satisfied:
\begin{itemize}
\item
There is a Zariski closed subset $Z\subset D$ 
with $\codim_X(Z)>k$
such that the restriction of $\prolongg{\vecc}{E}_{\ast}$ 
to $(X-Z,D-Z)$ is a $\vecc$-parabolic bundle.
\hfill\qed
\end{itemize}
\end{df}

It is easy to observe that
a reflexive saturated $\vecc$-parabolic Higgs sheaf
is a $\vecc$-parabolic Higgs bundle in codimension two.

\subsection{The characteristic number for $\vecc$-parabolic bundle in
   codimension two}
\label{subsection;05.10.6.1}

For any $\vecc$-parabolic bundle $\prolongg{\vecc}{E}_{\ast}$
in codimension two,
the parabolic second Chern character 
$\parch_2(\prolongg{\vecc}{E}_{\ast})\in H^4(X,\real)$
is defined as follows:
\begin{multline}
 \parch_2(\prolongg{\vecc}{E}_{\ast}):=
 \ch_2(\prolongg{\vecc}{E})
-\!\sum_{\substack{\\
 i\in S \\
 a\in\Par(\prolongg{\vecc}{E}_{\ast},i)}}
 \!
 a\cdot\iota_{i\,\ast}\Bigl(
 c_1\bigl(
 \lefttop{i}\Gr^{F}_a(\prolongg{\vecc}{E})\bigr)
 \Bigr)  \\
+\frac{1}{2}
 \sum_{\substack{
 i\in S \\ a\in \Par(\prolongg{\vecc}{E}_{\ast},i)}}
 a^2\cdot\rank\left(
 \lefttop{i}\Gr^F_a\bigl(\prolongg{\vecc}{E}\bigr)
 \right)\cdot[D_i]^2
 \\
+\frac{1}{2}\sum_{\substack{(i,j)\in S^2\\ i\neq j}}
 \sum_{\substack{
 P\in\Irr(D_i\cap D_j)\\
 (a_i,a_j)\in \Par(\prolongg{\vecc}{E}_{\ast},P) }}
 a_i\cdot a_j\cdot
 \rank
 \lefttop{P}\Gr^F_{(a_i,a_j)}(\prolongg{\vecc}{E})
 \cdot[P].
\end{multline}
Let us explain some of the notation:
\begin{itemize}
\item
 $\ch_2(\prolongg{\vecc}{E})$ denotes the second Chern character 
 of $\prolongg{\vecc}{E}$.
\item
 $\iota_i$ denotes the closed immersion $D_i\lrarr X$,
 and $\iota_{i\,\ast}:H^2(D_i)\lrarr H^4(X)$
 denotes the associated Gysin map.
\item
 $\Irr(D_i\cap D_j)$ denotes the set of the irreducible
 components of $D_i\cap D_j$.
\item
 Let $P$ be an element of $\Irr(D_i\cap D_j)$.
 The generic point of the component is also denoted by $P$.
 We put
 $\lefttop{P}F_{(a,b)}:=\lefttop{i}F_{a\,|\,P}
 \cap \lefttop{j}F_{b\,|\,P}$
and $\lefttop{P}\Gr^F_{\veca}:=
 \lefttop{P}F_{\veca}\big/\sum_{\veca'\lneq\veca}\lefttop{P}F_{\veca'}$.
 Then $\rank\lefttop{P}\Gr^F_{\veca}$ denotes the rank
 of $\lefttop{P}\Gr^F_{\veca}$
 as an $\nbigo_{P}$-module.
\item
 We put
 $\Par(\prolongg{\vecc}{E}_{\ast},P):=
 \bigl\{
 \veca\,\big|\,
 \lefttop{P}\Gr^F_{\veca}(\prolongg{\vecc}{E})\neq 0
 \bigr\}$.
\item
 $[D_i]\in H^2(X,\real)$ and $[P]\in H^4(X,\real)$
 denote the cohomology classes
 given by $D_i$ and $P$ respectively.
\end{itemize}

If $X$ is an $n$-dimensional compact Kahler manifold
with a Kahler form $\omega$,
we put as follows:
\[
\parch_{2,\omega}(\prolongg{\vecc}{E}_{\ast}):=
\parch_2(\prolongg{\vecc}{E}_{\ast})
\cdot\omega^{n-2},
\quad
 \parchern_{1,\omega}^2(\prolongg{\vecc}{E}_{\ast}):=
 \parchern_1(\prolongg{\vecc}{E}_{\ast})^2\cdot\omega^{n-2}.
\]
If $\omega$ is the first Chern class of an ample line bundle $L$,
we use the notation
$\parchern^2_{1,L}(\prolongg{\vecc}{E}_{\ast})$
and $\parch_{2,L}(\prolongg{\vecc}{E}_{\ast})$.
In the case $\dim X=2$,
we have the obvious equalities
$\parchern^2_{1,L}(\prolongg{\vecc}{E}_{\ast})
=\parchern^2_1(\prolongg{\vecc}{E}_{\ast})$
and
$\parch_{2,L}(\prolongg{\vecc}{E}_{\ast})
=\parch_{2}(\prolongg{\vecc}{E}_{\ast})$.

\begin{df}
Let $X$ be a smooth projective variety
with an ample line bundle $L$,
and let $D$ be a simple normal crossing divisor.
Let $(\prolongg{\vecc}{E}_{\ast},\theta)$
be a $\mu_L$-polystable
reflexive saturated
$\vecc$-parabolic Higgs sheaf on $(X,D)$.
We say that $(\prolongg{\vecc}{E}_{\ast},\theta)$
has trivial characteristic numbers,
if any stable component 
$(\prolongg{\vecc}{E}'_{\ast},\theta')$
of $(\prolongg{\vecc}{E}_{\ast},\theta)$
satisfies
$\pardeg_L(\prolongg{\vecc}{E'}_{\ast})
=\int_X\parch_{2,L}(\prolongg{\vecc}{E'}_{\ast})
=0$
\hfill\qed
\end{df}

\section{Filtered Sheaf}
\label{section;06.8.17.50}

\subsection{Definitions}
\label{subsection;05.9.29.1}

We recall the notion of filtered sheaf
by following \cite{s2}.
Let $X$ be a complex manifold, 
and $D$ be a simple normal crossing divisor
with the irreducible decomposition
$D=\bigcup_{i\in S}D_i$.
For $\veca\in\real^S$,
$a_i$ denotes the $i$-th component of $\veca$
for $i\in S$.
A filtered sheaf
on $(X,D)$ is defined to be a tuple
$\vecE_{\ast}=\bigl(
\vecE,\bigl\{\prolongg{\vecc}{E}
 \,\big|\,\vecc\in\real^S\bigr\}
 \bigr)$ 
as follows:
\begin{itemize}
\item
 $\vecE$ is a quasi coherent $\nbigo_X$-module.
 We put $E:=\vecE_{|X-D}$.
\item
 $\prolongg{\vecc}{E}$ is a coherent $\nbigo_X$-submodule
 of $\vecE$ for each $\vecc\in\real^S$
 such that
 $\prolongg{\vecc}{E}_{|X-D}=E$.
\item
 In the case $\veca\leq \vecb$,
 we have $\prolongg{\veca}{E}\subset\prolongg{\vecb}{E}$,
 where $\veca\leq \vecb$ means
 $a_i\leq b_i$ for all $i\in S$.
 We also have $\bigcup_{\veca\in\real^S}\prolongg{\veca}{E}=\vecE$.
\item
 We have 
 $\prolongg{\veca'}{E}=
 \prolongg{\veca}{E}\otimes\nbigo_X(-\sum n_j\cdot D_j)$
as submodules of $\vecE$,
where $\veca'=\veca-(n_j\,\big|\,j\in S)$
for some integers $n_j$.
\item
For each $\vecc\in\real^S$,
the filtration $\lefttop{i}\nbigf$ of $\prolongg{\vecc}{E}$
indexed by $\openclosed{c_i-1}{c_i}$
is given as follows:
\[
 \lefttop{i}\nbigf_d(\prolongg{\vecc}{E})
:=\bigcup_{\substack{a_i\leq d\\ \veca\leq \vecc}}
 \prolongg{\veca}{E}.
\]
Then the tuple $\bigl(\prolongg{\vecc}{E},
 \{\lefttop{i}\nbigf\,|\,i\in S\}
 \bigr)$ is a $\vecc$-parabolic sheaf,
i.e.,
the sets
$\bigl\{a\in\openclosed{c_i-1}{c_i}\big|\,
 \lefttop{i}\Gr^{\nbigf}_a(\prolongg{\vecc}{E})\neq 0\bigr\}$
are finite.
\end{itemize}

\begin{rem}
By definition,
we obtain the $\vecc$-parabolic sheaf $\prolongg{\vecc}{E}_{\ast}$
obtained from filtered sheaf $\vecE_{\ast}$
for any $\vecc\in\real^S$,
which is called the $\vecc$-truncation of $\vecE_{\ast}$.
On the other hand,
a filtered sheaf $\vecE_{\ast}$ can be reconstructed
from any $\vecc$-parabolic sheaf $\prolongg{\vecc}{E}_{\ast}$.
So we can identify them.
\hfill\qed
\end{rem}

\begin{df}
A filtered sheaf $\vecE_{\ast}$ is called
reflexive and saturated,
if any $\vecc$-truncations are reflexive
and saturated.

A filtered sheaf $\vecE_{\ast}$
is called a filtered bundle in codimension $k$,
if any $\vecc$-truncations are $\vecc$-parabolic Higgs bundle
in codimension $k$.
\hfill\qed
\end{df}

\begin{rem}
In the definition,
``any $\vecc$'' can be replaced with ``some $\vecc$''.
\hfill\qed
\end{rem}

A Higgs field of $\vecE_{\ast}$ is defined to be 
a holomorphic homomorphism
$\theta:\vecE\lrarr\vecE\otimes\Omega^{1,0}(\log D)$
satisfying
$\theta(\prolongg{\vecc}{E})\subset
 \prolongg{\vecc}{E}\otimes\Omega^{1,0}_X(\log D)$.

\vspace{.1in}

Let $\vecE^{(i)}_{\ast}$ $(i=1,2)$ be a filtered bundle on $(X,D)$.
We put as follows:
\[
 \vecEtilde:=Hom(\vecE^{(1)},\vecE^{(2)}),
\quad\quad
 \prolongg{\veca}{\Etilde}:=
 \bigl\{ f\in \vecEtilde\,\big|\,
 f\bigl(\prolongg{\vecc}{E^{(1)}}\bigr)\subset
 \prolongg{\vecc+\veca}{E^{(2)}},\,\forall \vecc
 \bigr\}.
\]
\[
 \vecEhat:=\vecE^{(1)}\otimes\vecE^{(2)},
\quad\quad
 \prolongg{\veca}{\Ehat}:=
 \sum_{\veca_1+\veca_2\leq \veca}
 \prolongg{\veca_1}{E^{(1)}}
\otimes
 \prolongg{\veca_2}{E^{(2)}}.
\]
Then $\bigl(\vecEtilde,\{\prolongg{\veca}{\Etilde}\}\bigr)$
and $\bigl(\vecEhat,\{\prolongg{\veca}{\Ehat}\}\bigr)$ are also filtered
bundles.
They are denoted by
$Hom\bigl(\vecE^{(1)}_{\ast},\vecE^{(2)}_{\ast}\bigr)$
and $\vecE^{(1)}_{\ast}\otimes \vecE^{(2)}_{\ast}$.

Let $(\vecE_{\ast},\theta)$ be a regular filtered Higgs bundle.
Let $a$ and $b$ be non-negative integers.
Applying the above construction,
we obtain the parabolic structures and the Higgs fields
on $T^{a,b}(\vecE):=Hom\bigl(\vecE^{\otimes\,a},\vecE^{\otimes\,b}\bigr)$.
We denote it by $(T^{a,b}\vecE_{\ast},\theta)$.

\subsection{The characteristic numbers of filtered bundles
 in codimension two}

Let $X$ be a smooth projective variety
with an ample line bundle $L$,
and let $D$ be a simple normal crossing divisor.
Let $\vecE_{\ast}$ be a filtered bundle in codimension two
on $(X,D)$.

\begin{lem} 
 \label{lem;05.7.28.200}
For any $\vecc,\vecc'\in\real^S$,
we have
$\parchern_1(\prolongg{\vecc}{E}_{\ast})
=\parchern_1(\prolongg{\vecc'}{E}_{\ast})$ in $H^2(X,\real)$.
\end{lem}
\pf
The $j$-th components of $\vecc$ and $\vecc'$
are denoted by $c_j$ and $c_j'$ for any $j\in S$.
Take an element $i\in S$.
We have only to consider the case
$c_j=c_j'$ $(j\neq i)$.
We may also assume 
$c_i'\in \Par\bigl(\vecE_{\ast},i\bigr)$ and $c_i<c_i'$.
Moreover it can be assumed that $c_i$ is sufficiently close to $c'_i$.
Then we have the following exact sequence
of $\nbigo_X$-modules:
\[
 0\lrarr \prolongg{\vecc}{E}\lrarr
 \prolongg{\vecc'}{E}\lrarr
 \lefttop{i}\Gr^F_{c'_i}\bigl(\prolongg{\vecc'}{E}_{|D_i}\bigr)
 \lrarr 0.
\]
We put $c:=c_i'-1$.
Then we have the following:
\begin{equation} \label{eq;05.7.28.150}
 \lefttop{i}\Gr^F_{c}(\prolongg{\vecc}{E})\otimes\nbigo(D_i)
\simeq
 \lefttop{i}\Gr^F_{c'_i}(\prolongg{\vecc'}{E}),
\quad\quad
 \lefttop{i}\Gr^F_{a}(\prolongg{\vecc}{E})
\simeq
 \lefttop{i}\Gr^F_{a}(\prolongg{\vecc'}{E}),\,\,
(c<a<c_i').
\end{equation}
Therefore we have
$\wt(\prolongg{\vecc}{E}_{\ast},i)
=\wt(\prolongg{\vecc'}{E}_{\ast},i)
-\rank \lefttop{i}\Gr^F_c(\prolongg{\vecc}{E})$.
On the other hand, we have 
$ c_1\bigl(\prolongg{\vecc'}{E}\bigr)
=c_1\bigl(\prolongg{\vecc}{E}\bigr)
+c_1\bigl(\iota_{\ast}\lefttop{i}\Gr^F_{c'}(\prolongg{\vecc'}{E})\bigr)$.
There is a closed subset $W\subsetneq D_i$ such that
$\lefttop{i}\Gr^F_{c'}(\prolongg{\vecc'}{E})_{|D_i-W}$
is isomorphic to a direct sum of $\nbigo_{D_i-W}$.
We remark that $H^2(X,\real)\simeq H^2(X\setminus W,\real)$,
because the codimension of $W$ in $X$ is larger than two.
Then it is easy to check
$c_1\bigl(\iota_{\ast}\lefttop{i}\Gr^F_{c'}(\prolongg{\vecc'}{E})\bigr)
=\rank \lefttop{i}\Gr^F_{c}(\prolongg{\vecc}{E})\cdot [D_i]$.
Then the claim of the lemma immediately follows.
\hfill\qed

\begin{cor}\label{cor;05.8.23.2}
For any $\vecc,\vecc'\in \real^S$,
we have the following:
\[
 \pardeg_{L}(\prolongg{\vecc}{E}_{\ast})
=\pardeg_{L}(\prolongg{\vecc'}{E}_{\ast}),
\quad\quad
 \int_X\parchern_{1,L}^2(\prolongg{\vecc}{E}_{\ast})
=\int_X\parchern_{1,L}^2(\prolongg{\vecc'}{E}_{\ast}).
\]
In particular,
the characteristic numbers
$\pardeg_L(\vecE_{\ast}):=\pardeg_L(\prolongg{\vecc}{E}_{\ast})$
and $\int_X\parchern_{1,L}^2(\vecE_{\ast}):=
 \int_X\parchern_{1,L}^2(\prolongg{\vecc}{E}_{\ast})$
are well defined.
\hfill\qed
\end{cor}

\begin{rem}
The $\mu_L$-stability of a regular filtered Higgs bundle
is defined,
which is equivalent to the stability of
any $\vecc$-truncation.
Due to Corollary {\rm\ref{cor;05.8.23.2}},
it is independent of a choice of $\vecc$.
\hfill\qed
\end{rem}

\begin{prop}
\label{prop;05.7.28.100}
For any $\vecc,\vecc'\in\real^S$,
we have the following:
\[
 \int_X\parch_{2,L}(\prolongg{\vecc}{E}_{\ast})
=\int_X\parch_{2,L}(\prolongg{\vecc'}{E}_{\ast}).
\]
In particular,
$\int_X\parch_{2,L}(\vecE_{\ast}):=
 \int_X\parch_{2,L}(\prolongg{\vecc}{E}_{\ast})$
is well defined.
\end{prop}
\pf
We have only to consider the case $\dim X=2$.
We use the following lemma.
\begin{lem} 
Let $Y$ be a smooth projective surface,
and $D$ be a smooth divisor of $Y$.
Let $\nbigf$ be an $\nbigo_D$-coherent module.
Then we have the following:
\[
 \int_X\ch_2(\iota_{\ast}\nbigf)
=\deg_{D}\nbigf-\frac{1}{2}\rank_D(\nbigf)\cdot (D,D).
\]
\end{lem}
\pf
By considering the blow up of $D\times\{0\}$
in $Y\times\cnum$ as in \cite{fulton},
we can reduce the problem in the case
$Y$ is a projective space bundle over $D$.
We can also reduce the problem to the case
$\nbigf$ is a locally free sheaf on $D$.
Then, in particular, we may assume that
there is a locally free sheaf $\widetilde{\nbigf}$
such that $\widetilde{\nbigf}_{|D}=\nbigf$.
In the case,
we have
the $K$-theoretic equality
$\iota_{\ast}\nbigf=\widetilde{\nbigf}\cdot
 \bigl(\nbigo-\nbigo(-D)\bigr)$.
Therefore we have the following:
\[
 \ch(\iota_{\ast}\nbigf)=\ch(\widetilde{\nbigf})
\cdot\bigl(D-D^2/2\bigr)
=\rank \widetilde{\nbigf}\cdot D
+\left(
 -\frac{1}{2}\rank\widetilde{\nbigf}\cdot D^2
+c_1(\widetilde{\nbigf})\cdot D
 \right).
\]
Then the claim of the lemma is clear.
\hfill\qed

\vspace{.1in}

Let us return to the proof of Lemma \ref{prop;05.7.28.100}.
We use the notation in the proof of Lemma \ref{lem;05.7.28.200}.
We have the following equalities:
\begin{multline}
 \int_X\ch_{2}(\prolongg{\vecc'}{E})
=\int_X\ch_{2}(\prolongg{\vecc}{E})
+\deg_{D_i}(\lefttop{i}\Gr^F_{c'_i}(\prolongg{\vecc'}{E}))
-\frac{1}{2}\rank\lefttop{i}\Gr^F_{c'_i}(\prolongg{\vecc'}{E})
\cdot D_i^2 \\
=\int_X\ch_{2}(\prolongg{\vecc}{E})
+\deg_{D_i}(\lefttop{i}\Gr^F_c(\prolongg{\vecc}{E}))
+\frac{1}{2}\rank\lefttop{i}\Gr^F_c(\prolongg{\vecc}{E})
 \cdot D_i^2.
\end{multline}
Here we have used (\ref{eq;05.7.28.150}).
We also have the following:
\[
  c'_i\cdot \deg_{D_i}(\lefttop{i}\Gr^F_{c'_i}
 (\prolongg{\vecc'}{E}))
=(c+1)\cdot\Bigl(
 \deg_{D_i}(\lefttop{i}\Gr^F_{c}(\prolongg{\vecc}{E}))
+\rank \lefttop{i}\Gr^F_{c}(\prolongg{\vecc}{E})\cdot D_i^2
 \Bigr).
\]
We remark the isomorphism
$\lefttop{P}\Gr^F_{(c'_i,a)}(\prolongg{\vecc'}{E})
\simeq
 \lefttop{P}\Gr^F_{(c,a)}(\prolongg{\vecc}{E})$
and the following exact sequence:
\[
 0\lrarr\lefttop{j}\Gr^F_{a}(\prolongg{\vecc}{E})
\lrarr \lefttop{j}\Gr^F_{a}(\prolongg{\vecc'}{E})
\lrarr \bigoplus_{P\in D_i\cap D_j}
 \lefttop{P}\Gr^F_{(c'_i,a)}(\prolongg{\vecc'}{E})
\lrarr 0.
\]
Hence we obtain the following equality:
\[
 a\cdot\deg_{D_j}
 \bigl(\lefttop{j}\Gr^F_{a}(\prolongg{\vecc'}{E})\bigr)
=a\cdot\deg_{D_j}
 \bigl(\lefttop{j}\Gr^F_a(\prolongg{\vecc}{E})\bigr)
+a\cdot\sum_{P\in D_i\cap D_j}
 \rank \lefttop{P}\Gr^F_{(c,a)}(\prolongg{\vecc}{E}).
\]

We have the following equalities:
\begin{equation}
 \frac{1}{2}c_i^{\prime\,2}\cdot\rank\lefttop{i}\Gr^F_{c'_i}
 (\prolongg{\vecc'}{E})\cdot D_i^2
=\frac{1}{2}c^2\rank\lefttop{i}\Gr^F_{c}(\prolongg{\vecc}{E})
 \cdot D_i^2
+\left(c+\frac{1}{2}\right)\cdot
 \rank \lefttop{i}\Gr^F_{c'_i}(\prolongg{\vecc'}{E})\cdot D_i^2.
\end{equation}
\begin{equation}
 c'_i\cdot a\cdot\rank \lefttop{P}\Gr^{F}_{(c_i',a)}
 \bigl(\prolongg{\vecc'}{E}\bigr)
=c\cdot a\cdot\rank\lefttop{P}\Gr^F_{(c,a)}(\prolongg{\vecc}{E})
+a\cdot\rank \lefttop{P}\Gr^F_{(c,a)}(\prolongg{\vecc}{E}).
\end{equation}

Then we obtain the following:
\begin{multline}
 \int_X\parch_{2,L}(\prolongg{\vecc'}{E}_{\ast})
-\int_X\parch_{2,L}(\prolongg{\vecc}{E}_{\ast}) 
=
\deg_{D_i}(\lefttop{i}\Gr^F_{c}(\prolongg{\vecc}{E}))
+\frac{1}{2}\rank\lefttop{i}\Gr^F_{c}(\prolongg{\vecc}{E})
\cdot D_i^2 \\
-\deg_{D_i}(\lefttop{i}\Gr^F_c(\prolongg{\vecc}{E}))
-(c+1)\rank\lefttop{i}\Gr^F_{c}(\prolongg{\vecc}{E})D_i^2 
-\sum_{j\neq i}\sum_{P\in D_i\cap D_j}
 \sum_a a\cdot\rank \lefttop{P}\Gr^F_{(c,a)}(\prolongg{\vecc}{E})\\
+\left(c+\frac{1}{2}\right)
  \rank \lefttop{i}\Gr^F_{c}(\prolongg{\vecc}{E})D_i^2
+\sum_{j\neq i}
 \sum_{P\in D_i\cap D_j}
 \sum_{a}a\cdot\rank \lefttop{P}\Gr^F_{(c,a)}(\prolongg{\vecc}{E})
=0.
\end{multline}
Thus we are done.
\hfill\qed

\begin{df}
Let $(\vecE_{\ast},\theta)$
be a $\mu_L$-polystable reflexive saturated
regular filtered Higgs sheaf on $(X,D)$.
We say that $(\vecE_{\ast},\theta)$
has trivial characteristic numbers,
if any stable component 
$(\vecE'_{\ast},\theta')$
of $(\vecE_{\ast},\theta)$
satisfies
$\pardeg(\vecE'_{\ast})
=\int_X\parch_2(\vecE'_{\ast})
=0$.
\hfill\qed
\end{df}

\section{Perturbation of Parabolic Structure}
\label{section;05.7.30.15}

Let $X$ be a smooth projective {\em surface}
with an ample line bundle $L$,
and $D$ be a simple normal crossing divisor
with the irreducible decomposition $D=\bigcup_{i\in S} D_i$.
Let $(\prolongg{\vecc}{E},\vecF,\theta)$ be a $\vecc$-parabolic
Higgs bundle over $(X,D)$.
Due to the projectivity of $D_i$,
the eigenvalues of
$\Res_i(\theta)\in \End\bigl(\prolongg{\vecc}{E}_{|D_i}\bigr)$ are constant.
Hence we obtain the generalized eigen decomposition
with respect to $\Res_i(\theta)$:
\[
 \lefttop{i}\Gr^F_a\bigl(\prolongg{\vecc}{E}_{|D_i}\bigr)
=\bigoplus_{\alpha\in\cnum}
 \lefttop{i}\Gr^{F,\EE}_{(a,\alpha)}
\bigl(\prolongg{\vecc}{E}_{|D_i}\bigr).
\]
Let $\nbign_i$ denote the nilpotent part
of the induced endomorphism $\Gr^F\Res_i(\theta)$ on
$\lefttop{i}\Gr^F_{a}(\prolongg{\vecc}{E}_{|D_i})$.
\begin{df}
 \label{df;05.8.30.1}
The $\vecc$-parabolic Higgs bundle
$(\prolongg{\vecc}{E},\vecF,\theta)$ is called 
graded semisimple,
if $\nbign_i$ are $0$ for any $i\in S$.
\hfill\qed
\end{df}

For simplicity,
we assume $c_i\not\in \Par\bigl(\prolongg{\vecc}{E}_{\ast},i\bigr)$ for
any $i$,
where $\vecc=(c_i\,|\,i\in S)$.
\begin{prop}
Let $\epsilon$ be any positive number
satisfying
$\epsilon\cdot 100\rank(E)\leq
\gap(\prolongg{\vecc}{E},\vecF)$.
There exists a $\vecc$-parabolic structure
$\vecF^{(\epsilon)}=\bigl( \lefttop{i}F^{(\epsilon)}\,\big|\,i\in S\bigr)$
such that the following holds:
\begin{itemize}
\item
 $(\prolongg{\vecc}{E},\vecF^{(\epsilon)})$ is 
 a graded semisimple $\vecc$-parabolic Higgs bundle.
\item
 We have
 $\wt(\prolongg{\vecc}{E},\vecF^{(\epsilon)},i)
=\wt(\prolongg{\vecc}{E},\vecF,i)$.
 (See Subsection {\rm\ref{subsection;06.8.4.10}}
 for $\wt$.)
 In particular,
 we have 
 $\parchern_1(\prolongg{\vecc}{E},\vecF^{(\epsilon)})
 =\parchern_1(\prolongg{\vecc}{E},\vecF)$.
\item
 There is a constant $C$, which is independent of $\epsilon$,
 such that the following holds:
\[
 \left|
 \int_X\parch_{2}(\prolongg{\vecc}{E},\vecF^{(\epsilon)})
-\int_X\parch_{2}(\prolongg{\vecc}{E},\vecF)
 \right|
\leq C\cdot\epsilon,
\]
\item
$\gap(\prolongg{\vecc}{E},\vecF^{(\epsilon)})=\epsilon$.
\end{itemize}
Such $(\prolongg{\vecc}{E},\vecF^{(\epsilon)},\theta)$
is called an $\epsilon$-perturbation of
$(\prolongg{\vecc}{E},\vecF,\theta)$.
\end{prop}
\pf
To take a refinement of the filtration $\lefttop{i}F$,
we see the weight filtration induced on $\lefttop{i}\Gr^F$.
Let $\eta$ be a generic point of $D_i$.
We have the weight filtration
$W_{\eta}$ of the nilpotent map
$\nbign_{i,\eta}$ on
$\lefttop{i}\Gr^F\bigl(\prolongg{\vecc}{E}_{|D_i}\bigr)_{\eta}$,
which is indexed by $\seisuu$.
Then we can extend it to the filtration $W$ of
$\lefttop{i}\Gr^F\bigl(\prolongg{\vecc}{E}_{|D_i}\bigr)$
in the category of vector bundles on $D_i$
due to $\dim D_i=1$.
By our construction,
$\nbign_i(W_k)\subset W_{k-2}$
and $\dim \Gr^W_k=\dim \Gr^{W}_{-k}$.
The endomorphism
$\Res_i(\theta)$ preserves the filtration $W$
on $\lefttop{i}\Gr^{F}(\prolongg{\vecc}{E}_{|D_i})$,
and the nilpotent part of the induced endomorphisms
on $\Gr^W\lefttop{i}\Gr^F(\prolongg{\vecc}{E}_{|D_i})$
are trivial.

\vspace{.1in}

Let us take the refinement of the filtration $\lefttop{i}F$.
For any $a\in\openclosed{c_i-1}{c_i}$,
we have the surjection
$ \pi_a:\lefttop{i}F_a(\prolongg{\vecc}{E}_{|D_i})
\lrarr
 \lefttop{i}\Gr^F_a(\prolongg{\vecc}{E}_{|D_i})$.
We put 
$\lefttop{i}\widetilde{F}_{a,k}:=\pi_a^{-1}(W_k)$.
We use the lexicographic order on
$\openclosed{c_i-1}{c_i}\times\seisuu$.
Thus we obtain the increasing filtration
$\lefttop{i}\widetilde{F}$ indexed by
$\openclosed{c_i-1}{c_i}\times\seisuu$.
The set 
$\widetilde{S}_i:=
\bigl\{(a,k)\in \openclosed{c_i-1}{c_i}\times\seisuu
 \,\big|\,
 \lefttop{i}\Gr^{\widetilde{F}}_{(a,k)}\neq 0\bigr\}$
is finite.

Let 
$\varphi_i:\widetilde{S}_i\lrarr \openclosed{c_i-1}{c_i}$
be the increasing map
given by
$\varphi_i(a,k):=a+k\epsilon$.
We put as follows:
\[
 \lefttop{i}F^{(\epsilon)}_b=
\bigcup_{\varphi_i(a,k)\leq b}
 \lefttop{i}\widetilde{F}_{(a,k)}
\]
Thus we obtain the $\vecc$-parabolic structure
$\vecF^{(\epsilon)}=\bigl(\lefttop{i}F^{(\epsilon)}\,\big|\,
 i\in S\bigr)$.

Let $P$ be any point of $D_i$.
Take a holomorphic coordinate neighbourhood
$(U_P,z_1,z_2)$ around $P$
such that $U_P\cap D_i=\{z_1=0\}$.
Then we have the expression
$\theta=f_1(z_1,z_2)\cdot dz_1/z_1+f_2(z_1,z_2)\cdot dz_2$.
Then, $f_j(0,z_2)$  $(j=1,2)$ preserve
the filtration $\lefttop{i}F^{(\epsilon)}$.
Therefore, it is easy to see
that $\bigl(\prolongg{\vecc}{E},\vecF^{(\epsilon)},\theta\bigr)$
is $\vecc$-parabolic Higgs bundle on $(X,D)$.
By our construction, 
it has the desired property.
\hfill\qed

\vspace{.1in}

The following proposition is standard.
\begin{prop}
Assume that $\bigl(\prolongg{\vecc}{E},\vecF,\theta\bigr)$
is $\mu_L$-stable.
If $\epsilon$ is sufficiently small,
then the $\epsilon$-perturbation
$\bigl(\prolongg{\vecc}{E},\vecF^{(\epsilon)},\theta\bigr)$ is also
$\mu_L$-stable.
\end{prop}
\pf
Let $\prolongg{\vecc}{\widehat{E}}\subset \prolongg{\vecc}{E}$
be a saturated subsheaf
such that $\theta\bigl(\prolongg{\vecc}{\widehat{E}}\bigr)
\subset \prolongg{\vecc}{\widehat{E}}\otimes\Omega^{1,0}(\log D)$.
Let $\widehat{\vecF}$ and $\widehat{\vecF}^{(\epsilon)}$
be the tuples of the filtrations of
$\prolongg{\vecc}{\widehat{E}}$ induced by $\vecF$
and $\vecF^{(\epsilon)}$ respectively.
There is a constant $C$, 
which is independent of choices of 
$\prolongg{\vecc}{\widehat{E}}$ and small $\epsilon>0$,
such that
$ \bigl|
 \mu_{L}(\prolongg{\vecc}{\widehat{E}},\widehat{\vecF})
-\mu_{L}(\prolongg{\vecc}{\widehat{E}},\widehat{\vecF}^{(\epsilon)})
 \bigr|
\leq C\cdot\epsilon$.
Therefore,
we have only to show the existence of a positive number 
$\eta$ satisfying the inequalities
$\mu_{L}(\prolongg{\vecc}{\widehat{E}},\vecF)+\eta
<\mu_{L}(\prolongg{\vecc}{E},\vecF)$,
for any saturated Higgs subsheaf
$0\neq \prolongg{\vecc}{\widehat{E}}\subsetneq
 \prolongg{\vecc}{E}$
under the $\mu_L$-stability 
of $\bigl(\prolongg{\vecc}{E},\vecF,\theta\bigr)$.
It is standard, so we give only a brief outline.
Due to a lemma of A. Grothendieck
(see Lemma 2.5 in \cite{grothendieck})
we know the boundedness of the family $\nbigg(A)$
of saturated Higgs subsheaves $\prolongg{\vecc}{\widehat{E}}
\subsetneq \prolongg{\vecc}{E}$ such that 
$\deg_{L}(\prolongg{\vecc}{\widehat{E}})\geq -A$
for any fixed number $A$.

Let us consider the case where $A$ is sufficiently large.
Then $\mu_{L}(\prolongg{\vecc}{\widehat{E}}_{\ast})$
is sufficiently small
for any $\prolongg{\vecc}{\widehat{E}}\not\in \nbigg(A)$.
On the other hand,
since the family $\nbigg(A)$ is bounded,
the function $\mu_{L}$ on $\nbigg(A)$ 
have the maximum,
which is strictly smaller than 
$\mu_{L}(\prolongg{\vecc}{E}_{\ast})$
due to the $\mu_L$-stability.
Thus we are done.
\hfill\qed

\section{Mehta-Ramanathan Type Theorem}

\subsection{Statement}
\label{subsection;06.8.12.20}

We discuss the Mehta-Ramanathan type theorem
for parabolic Higgs sheaves.
Let $X$ be an $n$-dimensional smooth projective variety over 
$\cnum$ with an ample line bundle $L$.
For simplicity, we assume the characteristic number of $k$ is $0$.
Let $D$ be a simple normal crossing divisor of $X$.

\begin{prop}
\label{prop;06.8.12.15}
Let $(V_{\ast},\theta)$ be a parabolic Higgs sheaf over $(X,D)$.
It is $\mu_L$-(semi)stable,
if and only if $(V_{\ast},\theta)_{|Y}$
is $\mu_L$-(semi)stable,
where $Y$ denotes a complete intersection
of sufficiently ample generic hypersurfaces.
\end{prop}

We closely follow the arguments
of V. Mehta, A. Ramanathan 
(\cite{mehta-ramanathan1}, \cite{mehta-ramanathan2})
and Simpson (\cite{s5}).
See the papers for more detail.

\subsection{$\nbigw$-operator}

In the following, let $k$ denote a field of characteristic $0$.
Let $\nbigx$ be a smooth projective variety over $k$,
with an ample line bundle $L$.
Let $\nbigd$ be a simple normal crossing divisor of $\nbigx$.
Let $\nbigw$ be a vector bundle on $\nbigx$.
A $\nbigw$-valued operator of a parabolic sheaf $V_{\ast}$
on $(\nbigx,\nbigd)$
is defined to be a morphism
$\eta:V_{\ast}\lrarr V_{\ast}\otimes \nbigw$.
A $\nbigw$-subobject of $(V_{\ast},\eta)$
is a saturated subsheaf $F\subset V$
such that $\eta(F)\subset F\otimes \nbigw$.
We endow $F$ with the induced parabolic structure.
A parabolic sheaf with a $\nbigw$-valued operator
$(V_{\ast},\eta)$ is defined to be $\mu_L$-semistable
if and only if
$\mu_L(F_{\ast})\leq \mu_L(V_{\ast})$
holds
for any $\nbigw$-subobject
$F_{\ast}\subset V_{\ast}$.
The $\mu_L$-stability is also defined similarly.

In general, we have the $\nbigw$-subobjects
$F_{\ast}\subset V_{\ast}$
with the properties:
(i) $\mu_L(G_{\ast})\leq \mu_L(F_{\ast})$
 for any $\nbigw$-subobject $G_{\ast}$,
(ii) if $\mu_L(G_{\ast})=\mu_L(F_{\ast})$,
 we have $\rank(G)\leq \rank(F)$.
Such $F_{\ast}$ is uniquely determined,
and called the $\beta$-$\nbigw$-subobject
of $(V_{\ast},\eta)$.
By a similar argument,
we obtain the Harder-Narasimhan filtration.

\subsection{Weil's Lemma}

In general, for a given projective variety $\nbigx$
with a normal crossing divisor
$\nbigd=\bigcup_{j\in S}\nbigd_j$,
a pair of a line bundle $\nbigl$ on $\nbigx$
and a tuple $\veca=(a_j\,|\,j\in S)\in\real^S$
is called a parabolic line bundle on $(\nbigx,\nbigd)$.
We can regard them as the $\veca$-parabolic sheaf
on $(\nbigx,\nbigd)$ in an obvious manner.
Let $\Pic(\nbigx,\nbigd)$ denote the set of
parabolic line bundles on $(\nbigx,\nbigd)$.

\vspace{.1in}
Let us return to the setting in Subsection 
\ref{subsection;06.8.12.20}.
For simplicity, we assume $H^i(X,L^m)=0$
for any $m\geq 1$ and $i>0$.
We put $S_m:=H^0(X,L^m)$ for $m\in\seisuu_{\geq\,1}$.
For $\vecm=(m_1,\ldots,m_{n-1})\in\seisuu_{\geq\,1}^{n-1}$,
we put $S_{\vecm}:=\prod_{i=1}^t S_{m_i}$.
Let $Z_{\vecm}$ denote the correspondence variety,
i.e.,
$Z_{\vecm}=\bigl\{(x,s_1,\ldots,s_{n-1})\in X\times S_{\vecm},\,\big|\,
 s_i(x)=0,\,\,1\leq i\leq n-1 \bigr\}$.
The natural morphisms
$Z_{\vecm}\lrarr S_{\vecm}$
and $Z_{\vecm}\lrarr X$ are denoted by $q_{\vecm}$
and $p_{\vecm}$, respectively.
We put $Z^D_{\vecm}:=Z_{\vecm}\times_{X}D$
and $Z^{D_j}_{\vecm}:=Z_{\vecm}\times_XD_j$.
Recall that $Z^{D_j}_{\vecm}$ are irreducible,
because $Z^{D_j}_{\vecm}$ is a vector bundle
over $D_j$.
Let $K_{\vecm}$ denote the function field
of $S_{\vecm}$.
We put
$Y_{\vecm}:=
 Z_{\vecm}\times_{S_{\vecm}}K_{\vecm}$,
$Y^{D_j}_{\vecm}:=Z^{D_j}_{\vecm}\times_{S_{\vecm}}K_{\vecm}$
and $Y^D_{\vecm}:=Z^D_{\vecm}\times_{ S_{\vecm}}K_{\vecm}$.
The irreducible decomposition
of $Z^D_{\vecm}\times_{ S_{\vecm}}K_{\vecm}$
is given by $\bigcup_j Z^{D_j}_{\vecm}\times_{S_{\vecm}}K_{\vecm}$.
Recall the following result of Mehta and Ramanathan,
by whom such a type of lemma is called
Weil's Lemma.

\begin{lem}
Assume $n\geq 2$.
For $\vecm=(m_1,\ldots,m_{n-1})$
with each $m_i\geq 3$,
the natural map
$\Pic(X,D)\lrarr \Pic(Y_{\vecm},Y^D_{\vecm})$
is bijective.
\end{lem}
\pf
Since we have the natural correspondence between
the irreducible components
of $D$ and $Y^D_{\vecm}$,
the claim is obviously reduced
to Proposition 2.1 of \cite{mehta-ramanathan1}.
\hfil\qed

\subsection{A family of degenerating curves}
\label{subsection;06.8.12.1}

As in \cite{mehta-ramanathan1},
we fix a sequence of integers
$(\alpha_1,\ldots,\alpha_{n-1})$
with $\alpha_i\geq 2$.
We put $\alpha:=\prod\alpha_i$.
For a positive integer $m$,
let $(m)$ denote $(\alpha_1^{m},\ldots,\alpha_{n-1}^m)$.
Let $V_{\ast}$ be a coherent parabolic sheaf on $(X,D)$.
For each $m$,
we can take an open subset $U_m\subset S_{(m)}$
such that
(i) $q_{(m)}^{-1}(s)$ are smooth $(s\in U_m)$,
(ii) $q^{-1}_{(m)}(s)$ intersects with the smooth part 
of $D$ transversally,
(iii) $V_{\ast}$ is a parabolic bundle
on an appropriate neighbourhood of each $q^{-1}_{(m)}(s)\subset X$.
In the following, we will shrink $U_m$, if necessary.
In Section 5 of \cite{mehta-ramanathan1},
Mehta and Ramanathan constructed a family of
degenerating curves.
Take integers $l>m>0$.
Let $A$ be a discrete valuation ring over $k$
with the quotient field $K$.
Then there exists a curve $C$ over $\Spec A$
with a morphism $\varphi:C\lrarr X\times\Spec A$
over $\Spec A$
with the properties:
(i) $C$ is smooth,
(ii) the generic fiber $C_K$ gives
 a sufficiently general $K$-valued point
 in $U_l$,
(iii) the special fiber $C_k$ is reduced
 with smooth irreducible components
 $C_k^i$ $(i=1,\ldots,\alpha^{l-m})$
 which are
 sufficiently general $k$-valued points in $U_m$.
 We use the notation $D_C$
 to denote $C\times_X D$.
 We also use the notation
 $D_{j,C}$, $D_{j,C_K}$
 and $D_{j,C_{k}^i}$ in similar meanings.
 Then, we obtain the parabolic bundle
 $\varphi^{\ast}(V_{\ast})$ on $(C,D_C)$,
 which is denoted by $V_{\ast|C}$.
 The restriction to $C_K$ and $C_k^i$
 are denoted similarly.
 Let $W_{\ast}$ be a parabolic subsheaf
 of $V_{\ast\,|C_K}$.
 Recall that 
 $W$ can be extended to
 the subsheaf $\Wtilde\subset V_{|C}$,
 flat over $\Spec A$  with the properties:
 (i) $\Wtilde$ is a vector bundle over $C$,
 (ii) $\Wtilde_{|C_k^i}\lrarr V_{|C_k^i}$ are injective.
 (See Section 4 of \cite{mehta-ramanathan1}.)
In particular, we have 
$\deg_L(\det(\Wtilde_{|C_K}))=
 \sum \deg_L(\det(\Wtilde_{|C_{k}^i}))$.
 We have the induced parabolic structure
 of $\Wtilde_{|C_k^i}$ 
 as the subsheaf of $V_{\ast|C_k^i}$,
 for which we have
 $\wt(W_{l\ast},D_{j,C_K})\geq
 \wt(\Wtilde_{|C_k^i\ast},D_{j,C_k^i})$
 for each $D_j$.
 Therefore, we obtain
 $\mu_L(\Wtilde_{\ast|C_K})\leq
  \sum_i \mu_L(\Wtilde_{|C_{k}^i,\ast})$.
If the equality occurs,
we have $\wt(W_{l\ast},D_{j,C_K})=
 \wt(\Wtilde_{|C_k^i\ast},D_{j,C_k^i})$
for any $i$ and $j$,
and $\Wtilde_{\ast}$ with the induced parabolic structure
is the parabolic bundle.

\subsection{The arguments of Mehta and Ramanathan}

Let $\nbigw$ be a vector bundle on $X$.
Let $(V_{\ast},\eta)$ be a 
parabolic sheaf with a $\nbigw$-operator on $(X,D)$.

\begin{lem}
\label{lem;06.8.13.1}
$(V_{\ast},\eta)$ is $\mu_L$-semistable,
if and only if
there exists a positive integer $m_0$
such that $(V_{\ast},\eta)_{|Y_{(m)}}$
is also $\mu_L$-semistable
for any $m\geq m_0$.
\end{lem}
\pf
We have only to show
the ``only if'' part.
We reproduce the argument in \cite{mehta-ramanathan1}.
First,
assume $(V_{\ast},\eta)_{|Y_{(m)}}$
is $\mu_L$-semistable for some $m$,
and we show that
$(V_{\ast},\eta)_{|Y_{(l)}}$
is $\mu_L$-semistable for any $l> m$.
We take a family of degenerating curves $C$
as in Subsection \ref{subsection;06.8.12.1}.
We have the $\beta$-$\nbigw$-subobject
$W_{l,\ast}\subset V_{\ast|C_K}$.
We extend it to $\Wtilde\subset V_{|C}$.
Note that it is naturally the $\nbigw$-subobject.
Since we have
$\mu_L(W_{l\,\ast})\leq \sum_i \mu_L(\Wtilde_{|C_k^i\,\ast})$
and $\mu_L(V_{\ast|C_K})=\sum_i \mu_L(V_{\ast|C_k^i})$,
we obtain 
$\mu_L(W_{l\,\ast})\leq \mu_L(V_{\ast|C_K})$.
Thus, we obtain the semistability of $V_{\ast|Y_{(l)}}$.

We will show that $V_{\ast}$ is not semistable
if $V_{\ast|Y_{(m)}}$
are not semistable for any $m$.
By shrinking $U_m$ appropriately,
we may have $\nbigw$-subobjects
$W_{m\,\ast}$
of $p_{(m)}^{\ast}V_{\ast|q_{(m)}^{-1}U_m}$
such that 
$W_{m\,\ast|\,q_{(m)}^{-1}(s)}$
is the $\beta$-$\nbigw$-subobject of
$(V_{\ast},\eta)_{|q_{(m)}^{-1}(s)}$
for any $s\in U_m$.
The restriction $W_{m\,\ast|Y_{(m)}}$ is
the $\beta$-$\nbigw$-subobject
of $(V_{\ast},\eta)_{|Y_{(m)}}$.
We have the parabolic line bundle
$\nbigl_{m\,\ast}\in \Pic(X,D)$
corresponding to 
$\det(W_{m,\ast})_{|Y_{(m)}}\in\Pic(Y_{(m)},Y_{(m)}^{D})$.

We put $\beta_m:=\mu_L(W_{m,\ast|Y_{(m)}})$.
For $l>m$,
we obtain $\beta_{l}\leq\alpha^{l-m}\cdot\beta_m$
by using a family of degenerating curves.
Since we have
$\beta_m=\alpha^{m}\cdot \mu_L(\nbigl_{m\,\ast})/\rank(W_m)$,
we obtain
$\mu_L(\nbigl_{l\,\ast})/\rank W_l
 \leq \mu_L(\nbigl_{m\,\ast})/\rank W_m$.
On the other hand,
we have $\beta_m\geq \alpha^{m}\mu_L(V_{\ast})$,
the sequence $\{\mu_L(\nbigl_{m\,\ast})\}$ is bounded.
Since $\{\wt(\nbigl_m,D_j)\}$  is finite,
we may take a subsequence $Q\subset \{m\}$
such that $\deg_L(\nbigl_m)$, $\wt(\nbigl_m,D_j)$
and $\rank(W_m)$ 
are independent of the choice of $m\in Q$.

Let us show that $\nbigl_{m}$ $(m\in Q)$
are isomorphic.
Take $l>m$  in $Q$.
We take a family of degenerating curves
as above.
We extend $W_{l|C_K}$ to $\Wtilde$
on $C$.
From $\beta_l=\alpha^{l-m}\beta_m$,
$\beta_l=\mu_L(W_{l\,\ast})
\leq \sum \mu_L(\Wtilde_{|C_k^i\,\ast})$
and $\mu_L(\Wtilde_{|C_k^i\,\ast})\leq \beta_m$,
we obtain
$\mu_L(\Wtilde_{|C_k^i\,\ast})=\beta_m$,
and thus $\Wtilde_{|C_k^i\,\ast}$ 
are $\beta$-$\nbigw$-subobjects of $V_{\ast|C_k^i}$.
In particular,
$\mu_L(\det(\Wtilde_{|C_k^i\,\ast}))
=\mu_L(\nbigl_{l\,\ast|C_k^i})$.
We also obtain 
$\mu_L(W_{l\,\ast})=\sum \mu_L(\Wtilde_{|C_k^i\,\ast})$,
and hence
$\wt(\Wtilde_{|C_k^i\,\ast},D_{j,C_k^i})
=\wt(W_{l\,\ast},D_{j,C_K})
=\wt(\nbigl_{l\,\ast},D_j)$.
Hence we obtain
$\deg_L\bigl(\det(\Wtilde_{|C_k^i})\bigr)
=\deg_L(\nbigl_{l|C_k^i})$,
and thus
$\nbigl_{l|C}\simeq \det(\Wtilde)$.
Since the parabolic weights are also same,
we have
$\det(\Wtilde)_{\ast}\simeq \nbigl_{l\,\ast|C}$.
Since $C_k^i$ are sufficiently general in $U_m$,
we obtain 
$\nbigl_{l\,\ast|Y_{(m)}}\simeq \nbigl_{m\,\ast|Y_{(m)}}$,
and hence $\nbigl_{l\,\ast}$ and $\nbigl_{m\,\ast}$
are isomorphic.
Now, let $\nbigl_{\ast}$ denote $\nbigl_{l\,\ast}$
$(l\in Q)$.

Let us show the existence of
a $\nbigw$-subsheaf $\Wtilde$ of $V$,
such that
$\Wtilde_{|q_{(m)}^{-1}(s)}
=W_{m|q_{(m)}^{-1}(s)}$
for a sufficiently large $m$.
Such $\Wtilde$ will contradict with
the semistability of $(V_{\ast},\eta)$.
Let $U$ be an open subset of $X$
on which $V$ is a vector bundle.
We may assume that 
$\codim(X-U)\geq 2$.
We put $r=\rank(W_m)$ for $m\in Q$.
Let $G$ denote the bundle 
of Grassmann varieties on $U$,
whose fiber over $q\in U$
consists of the subspaces of $V_{|q}$
with rank $r$.
We have the natural embedding of
$G$ into the projectivization of
$\bigwedge^r V_{|U}$.
Let $\Ghat\subset \bigwedge^r V_{|U}$
denote the cone over $G$.

Let $F$ denote the double dual of $\bigwedge^rV$.
We have the naturally induced saturated parabolic
structure of $F$.
Let $\nhom(\nbigl_{\ast},F_{\ast})$
denote the sheaf of homomorphisms
from $\nbigl_{\ast}$ to $F_{\ast}$,
which is reflexive.
We put 
$H:=H^0\bigl(X,\nhom(\nbigl_{\ast},F_{\ast})\bigr)$.
For any $\phi\in H$,
we put 
$\Sigma(\phi):=\{x\in U\,|\,\phi(x)\in \Sigma\}$.
Since $\{\Sigma(\phi)\,|\,\phi\in H\}$ is bounded family,
we have $q_{(m)}^{-1}(s)\not\subset\Sigma(\phi)$
for a sufficiently large $m$ and $s\in U_m$,
unless $\Sigma(\phi)\neq U$.
On the other hand,
there exists a non-trivial morphism
$\phi\in H$ such that
$q_{(m)}^{-1}(s)\subset\Sigma(\phi)$
for such $m$ and $s$,
due to the above consideration and
 General Enriques-Severi Lemma
(Proposition 3.2 \cite{mehta-ramanathan1}).
Hence, we obtain $\Sigma(\phi)=U$
for such $\phi$.
The image of $\phi$ naturally induces
the saturated subsheaf $\Wtilde\subset V$.
If $m$ is sufficiently large,
we also obtain $\eta(\Wtilde)\subset\Wtilde\otimes\nbigw$.
To see it,
we recall the boundedness of the family $\nbigs$
of the saturated subsheaves $F$
of $V$ such that $\deg(F)\geq C$,
for some fixed $C$
(Lemma 2.5 in \cite{grothendieck}).
So we can take a large $m$ 
such that $\eta(F)\subset F\otimes\nbigw$
$(F\in\nbigs)$
if and only if
$\eta(F_{|q^{-1}(s)})\subset F_{|q^{-1}(s)}\otimes\nbigw$
for a sufficiently general $s\in U_m$.
Thus we are done.
\hfill\qed

\begin{lem}
\label{lem;06.8.12.10}
$(V_{\ast},\eta)$ is $\mu_L$-stable,
if and only if
there exists a positive integer $m_0$
such that $(V_{\ast},\eta)_{|Y_{(m)}}$
is also $\mu_L$-stable
for any $m\geq m_0$.
\end{lem}
\pf
We reproduce the argument in
\cite{mehta-ramanathan2}.
First, let us see 
$V_{|q_{(m)}^{-1}(s)}$ is simple for sufficiently large $m$
if $(V_{\ast},\eta)$ is $\mu_L$-stable.
To show it, we have only to consider the case
$V_{\ast}$ is reflexive and saturated.
Let $\nhom((V_{\ast},\eta),(V_{\ast},\eta))$ be the sheaf
of endomorphisms of $V$
which preserves the parabolic structure
and commutes with $\eta$.
Then, it is easy to check
$\nhom((V_{\ast},\eta),(V_{\ast},\eta))$ is reflexive
by using Lemma \ref{lem;06.8.12.35},
and hence the claim is shown
by applying General Enriques-Severi Lemma.

Let us recall the notion of {\em socle}
of semistable objects,
which is the direct sum of stable subobjects
(See \cite{mehta-ramanathan2}.
 Recall we have assumed the characteristic of $k$ is $0$.)
Assume that $(V_{\ast},\eta)_{|Y_{(m)}}$ is stable
for some $m$.
Then, it can be shown that
$(V_{\ast},\eta)_{|Y_{(l)}}$ is also stable
for any $l>m$
by using a family of degenerating curves
and the socle of $(V_{\ast},\eta)_{|Y_{(l)}}$,
instead of $\beta$-$\nbigw$-subobjects.
So we assume that $(V_{\ast},\eta)_{|Y_{(m)}}$
is not stable for any $m$,
and we will show that $(V_{\ast},\eta)$ is not $\mu_L$-stable.

Let $N$ be sufficiently large.
By shrinking $U_m$ appropriately for $m\geq N$,
we may assume 
(i) $V_{|q_{(m)}^{-1}(s)}$ is simple and semistable
for any $s\in U_m$,
(ii) the socle of $V_{|Y_{(m)}}$
is extended to 
$W_{m\,\ast}\subset p_{(m)}^{\ast}V_{\ast|q_{(m)}^{-1}(U_m)}$,
(iii)
$W_{m\,\ast|q_{(m)}^{-1}(s)}$
is the socle of $V_{|q_{(m)}^{-1}(s)}$
for any $s\in U_m$.
We have the parabolic line bundle
$\nbigl_{m\,\ast}$ on $(X,D)$
corresponding to
$\det(W_{m\,\ast|Y_{(m)}})$ on
$(Y_{(m)},Y^D_{(m)})$.
We have $\mu_L(\nbigl_{m,\ast})=\rank(W_m)\cdot\mu_L(V_{\ast})$.
Hence, we can take a subsequence $Q\subset\{m\}$
such that 
$\rank W_m$,
$\wt(\nbigl_{m\,\ast},D_i)$
and $\deg(\nbigl_m)$ are independent of $m\in Q$.
We put $r:=\rank W_m$ for $m\in Q$.

Let $G_m$ denote the bundle of Grassmann varieties
on $q_{(m)}^{-1}(U_m)$,
whose fiber over $Q\in q_{(m)}^{-1}(U_m)$
consists of the subspace of 
$p_{(m)}^{\ast}\bigl(V\bigr)_{|Q}$
with rank $r$.
We have the natural embedding of
$G_m$ into the projectivization of
$p_{(m)}^{\ast}\bigl(\bigwedge^r V\bigr)_{|q_{(m)}^{-1}(U_m)}$.
Let $\Ghat_m$ denote the cone over $G_m$.

Take $m_0\in Q$,
and let $E$ denote the set of 
$\nbigl_{\ast}\in \Pic(X,D)$
with
$\mu_L(\nbigl_{\ast})=r\cdot\mu_L(V_{\ast})$
such that there exists
$\phi:\nbigl_{\ast|Y_{(m_0)}}\lrarr
 \bigwedge^rV_{\ast|Y_{(m_0)}}$
with $\phi(L_{|Y_{(m_0)}})\subset \Ghat_{m_0}$.
By the same argument as the proof of
Lemma 2.7--2.8 of \cite{mehta-ramanathan2},
it can be shown that $E$ is finite.

Let us show that 
$\nbigl_l\in E$ for any $l\in Q$ with $l> m_0$.
Let $C$ be a family of degenerating curves.
We extend $W_{l|C_K}$
to $\Wtilde\subset V_{\ast}$.
We have the inequalities
$\mu_L(W_{l\,\ast})\leq
\sum \mu_L(\Wtilde_{|C_{k}^i\,\ast})$,
$\mu_L(\Wtilde_{|C_k^i\,\ast})\leq\alpha^{m} \mu_L(V_{\ast})$
and the equality $\mu_L(W_{l\,\ast})=\alpha^l\mu_L(V_{\ast})$.
Thus, the inequalities are actually equalities.
Hence, we have
$\wt(\det(\Wtilde)_{|C_k^i\,\ast},C_{D_{k}^i})
=\wt(\nbigl_{l\,\ast},C_{D_k}^i)$
and 
$\mu_L(\det(\Wtilde)_{|C_k^i,\ast})
=\mu_L(\nbigl_{l\,\ast|C_k^i})$.
Therefore, we obtain
$\nbigl_{l\,\ast|C}
\simeq
 \det(\Wtilde)_{\ast}$.
In particular,
$\nbigl_{l\,\ast|C_k^i}
\simeq
 \det(\Wtilde_{|C_k^i})_{\ast}$.
Since $C_k^i$ are sufficiently general,
we obtain $\nbigl_{l\,\ast}\in E$.

Then, we can take a subsequence $Q'\subset Q$
such that 
$\nbigl_{m\,\ast}$ are isomorphic
$(m\in Q')$.
The rest of the argument is same as the last part
of the proof of Lemma \ref{lem;06.8.13.1}.
\hfill\qed

\subsection{End of Proof of Proposition \ref{prop;06.8.12.15}}

We have only to show the ``only if'' part.
We reproduce the argument in \cite{s5}.
Assume the $\mu_L$-stability of $(V_{\ast},\theta)$.
Let $Y=Y_1\cap\cdots \cap Y_{t}$ be a generic complete
intersection,
where $\deg_L(Y_i)$ are appropriately large numbers.
We put $Y^{(i)}:=Y_1\cap\cdots\cap Y_i$
and $Y^{(0)}:=X$.
We also put $D^{(i)}:=D\cap Y^{(i)}$
and $D^{(0)}=D$.
We put $C_1:=\prod_{i=1}^t\bigl(
  \deg_L(Y_i)/\int_X c_1(L)^n\bigr)$.
We put $\nbigw^{(i)}:=\Omega_{Y^{(i)}}(\log D^{(i)})_{|Y}$.
Let $\theta_Y^{(i)}$ denote the induced
$\nbigw^{(i)}$-operation of $V_{\ast|Y}$.
We may assume that
$(V_{\ast|Y},\theta_Y^{(0)})$ is $\mu_L$-stable
due to Lemma \ref{lem;06.8.12.10}.
By applying the Mehta-Ramanathan type theorem
to the Harder-Narasimhan filtration of $V_{\ast}$,
we may have a constant $B$
such that
(i) it is independent of the choice of $Y_i$ and 
a sufficiently large $\pardeg_L(Y_i)$,
(ii) $\pardeg_L(F_{\ast})\leq B\cdot C_1$
for any $F_{\ast}\subset V_{\ast|Y}$.
We show that
$(V_{\ast},\theta^{(i)})$ are 
$\mu_L$-stable by an induction.

Assume that the claim holds for $i-1$.
Let $F_{\ast}$ be a $\nbigw^{(i)}$-object
of $V_{\ast|Y}$ such that
$\mu_L(F_{\ast})\geq 
 \mu_L(V_{\ast|Y})=\mu_L(V_{\ast})\cdot C_1$,
and we will derive the contradiction.
We put $G:=V/F$,
which is provided with the induced parabolic structure.
Then, we have the induced map
$\theta:F_{\ast}\lrarr G_{\ast}(-Y_{i})$.
Let $H$ denote the kernel.
Let $N$ denote the saturated subsheaf of
$G(-Y_{i})$ generated by $F/H$,
provided with the induced parabolic structure.
We have $\mu\bigl((F/H)_{\ast}\bigr)\leq \mu(N_{\ast})$.
Let $J\subset E_{\ast}(-Y_i)$ denote the pull back
of $N$ via $E(-Y_i)\lrarr G(-Y_i)$
with the induced parabolic structure.
We obtain the following:
\begin{multline}
B\cdot C_1\geq \pardeg_L(J(Y)_{\ast})
\geq
 \pardeg_L(F_{\ast})
+\pardeg_L\bigl(N_{\ast}(Y_i)\bigr)\\
\geq 2\pardeg_L(F_{\ast})
-\pardeg_L(H_{\ast})
+\rank(F/H)\cdot\deg_L(\nbigo(Y_i)_{|Y})\\
\geq \bigl(2\rank(F)\cdot\mu(V_{\ast})-B\bigr)\cdot C_1
 +\rank(F/H)\cdot\deg_L(\nbigo(Y_i)_{|Y})
\end{multline}
If $\deg_L(Y_i)$ is sufficiently large,
$\deg_L(\nbigo(Y_i)_{|Y})$ is much larger than
$C_1$.
Hence $\rank(F/H)$ must be $0$,
and hence $F$ is actually
a $\nbigw^{(i-1)}$-subobject,
which contradicts with the $\mu_L$-semistability
of $(V_{\ast|Y},\theta^{(i-1)})$.
Thus the induction can proceed.
\hfill\qed

\section{Adapted Metric}
\label{section;05.9.8.110}

We recall a `typical' example of filtered sheaf.
Let $E$ be a holomorphic vector bundle on $X-D$.
If we are given a hermitian metric $h$ of $E$,
we obtain the $\nbigo_X$-module $\prolongg{\vecc}{E}(h)$
for any $\vecc\in\real^S$,
as is explained in the following.
Let us take hermitian metrics $h_i$ of $\nbigo(D_i)$.
Let $\sigma_i:\nbigo\lrarr\nbigo(D_i)$ denote the canonical section.
We denote the norm of $\sigma_i$ with respect to $h_i$
by $|\sigma_i|_{h_i}$.
For any open set $U\subset X$,
we put as follows:
\[
 \Gamma\bigl(U,\prolongg{\vecc}{E}(h)\bigr):=
 \Bigl\{
 f\in \Gamma(U\setminus D,E)\,\Big|\,
 |f|_h=O\Bigl(\prod|\sigma_i|_{h_i}^{-c_i-\epsilon}\Bigr)
\,\,\forall \epsilon>0
 \Bigr\}.
\]
Thus we obtain the $\nbigo_X$-module
$\prolongg{\vecc}{E}(h)$.
We also put $\vecE(h):=\bigcup_{\vecc}\prolongg{\vecc}{E}(h)$.
\begin{rem}
In general,
$\prolongg{\vecc}{E}(h)$ are not coherent.
\hfill\qed
\end{rem}

\begin{df}
Let $\widetilde{\vecE}_{\ast}$ be a filtered vector bundle.
We put $E:=\widetilde{E}=\widetilde{\vecE}_{|X-D}$.
A hermitian metric $h$ of $E$ is called adapted to 
the parabolic structure of $\widetilde{\vecE}_{\ast}$,
if the isomorphism $E\simeq \widetilde{E}$
is extended to the isomorphisms
$\prolongg{\vecc}{E}(h)\simeq
 \prolongg{\vecc}{\widetilde{E}}$ for any $\vecc\in\real^S$.
\hfill\qed
\end{df}

\section{Convergence}
\label{section;05.10.1.3}

We give the definition of convergence of
a sequence of parabolic Higgs bundles.
Although we need such a notion only in the case 
where the base complex manifold is a curve,
the definition is given generally.
Let $X$ be a complex manifold,
and $D=\bigcup_{j\in S}D_j$ be a simple normal crossing divisor of $X$.
Let $p$ be a number which is sufficiently larger than $\dim X$.
Let $b$ be any positive integer.

\begin{df}
\label{df;05.8.29.300}
Let $\bigl(E^{(i)},\delbar^{(i)},\vecF^{(i)},\theta^{(i)}\bigr)$
$(i=1,2,\ldots)$
be a sequence of $\vecc$-parabolic Higgs bundles on $(X,D)$.
We say that the sequence
$\bigl\{(E^{(i)},\delbar^{(i)},\vecF^{(i)},\theta^{(i)})\bigr\}$
weakly converges to 
$\bigl(E^{(\infty)},\delbar^{(\infty)},
 \vecF^{(\infty)},\theta^{(\infty)}\bigr)$
in $L_b^p$ on $X$,
if there exist locally $L_b^p$-isomorphisms
$\Phi^{(i)}:E^{(i)}\lrarr E^{(\infty)}$ on $X$
satisfying the following conditions:
\begin{itemize}
\item
 The sequence
 $\{ \Phi^{(i)}(\delbar^{(i)})-\delbar^{(\infty)}\}$
 converges  to $0$ weakly in $L_{b-1}^p$ locally on $X$.
\item
 The sequence
$\{ \Phi^{(i)}(\theta^{(i)})-\theta^{(\infty)}\}$ converges  to $0$ 
 weakly in $L_{b-1}^p$ locally on $X$,
 as sections of $\End(E^{(\infty)})\otimes \Omega^{1,0}(\log D)$.
\item
 For simplicity, we assume that $\Phi^{(i)}$ are $C^{\infty}$ around  $D$.
\item
 The sequence
 $\{\Phi^{(i)}\bigl(\lefttop{j}F^{(i)}\bigr)\}$ converges to
 $\lefttop{j}F^{(\infty)}$ in an obvious sense.
More precisely,
for any $\delta>0$, $j\in S$ and $a\in \openclosed{c_j-1}{c_j}$,
 there exists $m_0$ such that
$\rank \lefttop{j}F^{(\infty)}_a=\rank\lefttop{j}F^{(i)}_{a+\delta}$
and that
 $\lefttop{j}F^{(\infty)}_a$
and 
 $\Phi^{(i)}\bigl(\lefttop{j}F^{(i)}_{a+\delta}\bigr)$
are sufficiently close in the Grassmaniann varieties,
 for any $i>m_0$.
\hfill\qed
\end{itemize}
\end{df}

The following lemma is standard.
\begin{lem}
\label{lem;05.10.7.100}
Let $X$ be a smooth projective variety,
and $D$ be a simple normal crossing divisor of $X$.
Assume that a sequence of $\vecc$-parabolic Higgs bundles
$\bigl\{(E^{(i)},\delbar^{(i)},\vecF^{(i)},\theta^{(i)})\bigr\}$ on $(X,D)$
converges to
$(E^{(\infty)},\delbar^{(\infty)},\vecF^{(\infty)},\theta^{(\infty)})$
weakly in $L_b^p$ on $X$.
Assume that there exist non-zero holomorphic sections $s^{(i)}$
of $(E^{(i)},\delbar^{(i)})$
such that
 $\theta^{(i)}(s^{(i)})=0$
and that
 $s^{(i)}_{|P}\in \lefttop{j}F_0\bigl(E^{(i)}_{|P}\bigr)$
 for any $P\in D_j$ and $j\in S$.

Then there exists a non-zero holomorphic section $s^{(\infty)}$
of $\bigl(E^{(\infty)},\delbar^{(\infty)}\bigr)$
such that 
 $\theta^{(\infty)}(s^{(\infty)})=0$
and that
 $s^{(\infty)}_{|P}\in \lefttop{j}F_0\bigl(E^{(\infty)}_{|P}\bigr)$
 for any $P\in D$ and $j\in S$.
\end{lem}
\pf
Let us take a $C^{\infty}$-metric $\widetilde{h}$
of $E^{(\infty)}$ on $X$.
We put $t^{(i)}:=\Phi^{(i)}(s^{(i)})$.
Since $p$ is large,
we remark that $\Phi^{(i)}$ are $C^0$.
Hence we have
$\max_{P\in X}|t^{(i)}(P)|_{\widetilde{h}}$.
We may assume 
$\max_{P\in X}|t^{(i)}(P)|_{\widetilde{h}}=1$.

We have
$\Phi^{(i)}(\delbar^{(i)})=\delbar^{(\infty)}+a_i$,
and hence $\delbar^{(\infty)}t^{(i)}=-a_i(t^{(i)})$.
Due to $|t^{(i)}|\leq 1$ and $a_i\lrarr 0$ weakly in $L_{b-1}^p$,
the $L_b^p$-norm of $t^{(i)}$ are bounded.
Hence we can take an appropriate subsequence
$\{t^{(i)}\,\big|\,i\in I\}$ which weakly converges to $s^{(\infty)}$
in $L_b^p$ on $X$.
In particular,
$\{t^{(i)}\}$ converges to a section $s^{(\infty)}$ in $C^0$.
Due to $\max_{P}|s^{(\infty)}(P)|_{\widetilde{h}}=1$,
the section $s^{(\infty)}$ is non-trivial.
We also have
$\delbar^{(\infty)}s^{(\infty)}=0$  in $L_{b-1}^p$,
and hence $s^{(\infty)}$ is a non-trivial holomorphic section of
$(E^{(\infty)},\delbar^{(\infty)})$.
It is easy to see that $s^{(\infty)}$ has the desired property.
\hfill\qed

\begin{cor}
 \label{cor;05.8.29.301}
Let $(X,D)$ be as in Lemma {\rm\ref{lem;05.10.7.100}}.
Assume that
a sequence of $\vecc$-parabolic Higgs bundles
 $\bigl\{(E^{(i)},\delbar^{(i)},\vecF^{(i)},\theta^{(i)})\bigr\}$
on $(X,D)$ weakly converges to both 
$(E,\delbar_E,\vecF,\theta)$ and $(E',\delbar_{E'},\vecF',\theta')$
in $L_b^p$ on $X$.
Then there exists a non-trivial holomorphic map
$f:(E,\delbar_{E})\lrarr (E',\delbar_{E'})$ on $X$
which is compatible with the parabolic structures
and the Higgs fields.
\hfill\qed
\end{cor}

%% file: 4.tex
In this chapter, we would like to explain
about an ordinary metric for parabolic Higgs bundles,
which is a metric adapted to the parabolic structure.
Such a metric has been standard
in the study of parabolic bundles
(for example, see  \cite{biquard2},
 \cite{li-na} and \cite{li2}).
It is our purpose to see that
it gives a rather good metric when
the parabolic Higgs bundle is {\em graded semisimple}.
(If parabolic Higgs bundle is not graded semisimple,
 we need more complicated metric
 as discussed in \cite{b} and \cite{s2}.)
After giving estimates around
the intersection and the smooth part
of the divisor
in Sections \ref{section;05.7.31.3}
and \ref{section;05.7.31.5},
we see some properties of an ordinary metric
in Section \ref{section;05.8.30.100}.

\section{Around the Intersection $D_i\cap D_j$}
\label{section;05.7.31.3}

\subsection{Construction of a metric}
\label{subsection;06.8.6.1}

We put $X:=\{(z_1,z_2)\in\cnum^2\,\big|\,|z_i|<1\}$,
$D_i:=\{z_i=0\}$ and $D=D_1\cup D_2$.
Take a positive number $\epsilon$,
and let $\omega_{\epsilon}$ denote the following metric,
for some positive number $N$:
\[
 \sum \bigl(
 \epsilon^{N+2}\cdot |z_i|^{2\epsilon}
 +|z_i|^2\bigr)
 \cdot \frac{dz_i\cdot d\bar{z}_i}{|z_i|^2}.
\]

Let $(\prolongg{\vecc}{E}_{\ast},\theta)$ be 
a $\vecc$-parabolic Higgs bundle on $(X,D)$.
We put $E:=\prolongg{\vecc}{E}_{|X-D}$.
We take a positive number $\epsilon$
such that
$10\epsilon<\gap(\prolongg{\vecc}{E}_{\ast})$.
We have the description:
\[
 \theta=f_1\cdot \frac{dz_1}{z_1}+f_2\cdot \frac{dz_2}{z_2},
\quad
 f_i\in End(\prolongg{\vecc}{E}).
\]
We have $\Res_i(\theta)=f_{i\,|\,D_i}$.

\begin{assumption}\mbox{{}}
\begin{itemize}
\item
 The eigenvalues of $\Res_i(\theta)$ are constant.
 The sets of the eigenvalues of $\Res_i(\theta)$
 are denoted by $S_i$.
\item
 We have the decomposition:
\[
 \prolongg{\vecc}{E}=
 \bigoplus_{\vecalpha\in S_1\times S_2} 
 \prolongg{\vecc}{E_{\vecalpha}}
\quad\mbox{\rm such that }\quad
f_i(\prolongg{\vecc}{E_{\vecalpha}})
 \subset \prolongg{\vecc}{E_{\vecalpha}}.
\]
 There are some positive constants
 $C$ and $\eta$ such that
 any eigenvalue $\beta$ of $f_{i\,|\,E_{\vecalpha}}$ 
 satisfies 
 $|\beta-\alpha_i|\leq C\cdot|z_i|^{\eta}$
 for $\vecalpha=(\alpha_1,\alpha_2)$.
\hfill\qed
\end{itemize}
\end{assumption}

\begin{rem}
\label{rem;06.8.15.102}
The first condition is satisfied, 
when we are given
a projective surface $X'$ with a simple normal crossing divisor $D'$
and a $\vecc$-parabolic Higgs bundle $(\prolongg{\vecc'}{E'}_{\ast},\theta')$
on $(X',D')$, such that $(X,D)\subset (X',D')$
and $(\prolongg{\vecc}{E}_{\ast},\theta)
=(\prolongg{\vecc'}{E'}_{\ast},\theta')_{|X}$.
The second condition is also satisfied,
if we replace $X$ with a smaller open subset around the origin $O=(0,0)$.
\hfill\qed
\end{rem}

In the following,
we replace $X$ with a smaller open subset
containing $O$ without mentioning,
if it is necessary.
Let us take a holomorphic decomposition
$ \prolongg{\vecc}{E_{\vecalpha}}
=\bigoplus_{\veca\in\real^2} U_{\vecalpha,\veca}$
satisfying the following conditions,
where $b_i$ denotes the $i$-th component of $\vecb$:
\[
 \bigoplus_{\vecb\leq \veca}
 U_{\vecalpha,\vecb\,|\,O}
=\lefttop{1}F_{a_1\,|\,O}\cap \lefttop{2}F_{a_2\,|\,O}\cap
 \prolongg{\vecc}{E_{\vecalpha\,|\,O}},
\quad\quad
 \bigoplus_{b_i\leq a}U_{\vecalpha,\vecb\,|\,D_i}
=\prolongg{\vecc}{E_{\vecalpha\,|\,D_i}}\cap \lefttop{i}F_{a}.
\]
We take a holomorphic frame $\vecv=(v_1,\ldots,v_r)$
compatible with the decomposition,
i.e.,
we have 
$(\veca(v_j),\vecalpha(v_j))\in\real^2\times\cnum^2$
for each $v_j$
such that $v_j\in U_{\vecalpha(v_j),\veca(v_j)}$.
Let $h_0'$ be the hermitian metric  of $\prolongg{\vecc}{E}$
for which $\vecv$ is orthonormal.
Let $h_0$ be the hermitian metric of $E$
such that
$h_0(v_i,v_j)=h_0'(v_i,v_j)\cdot |z_1|^{-2a_1(v_i)}\cdot
|z_2|^{-2a_2(v_i)}$,
where $a_j(v_i)$ denotes the $j$-th component of $\veca(v_i)$.
We put as follows:
\[
 A=A_1+A_2,\quad
 A_i=\bigoplus \left(-a_i\frac{dz_i}{z_i}\right)\cdot
 \id_{U_{\vecalpha,\veca}}.
\]
Then, we have $\del_{h_0}=\del_{h_0'}+A$.
We also have $R(h_0)=R(h_0')=0$.

\subsection{Estimate of $F(h_0)$
 in the graded semisimple case}

\begin{prop}
\label{prop;06.8.15.100}
\mbox{{}} 
If $(\prolongg{\vecc}{E}_{\ast},\theta)$ is graded semisimple
in the sense of Definition {\rm\ref{df;05.8.30.1}},
then $F(h_0)$ is bounded with respect to $\omega_{\epsilon}$ and $h_0$.
\end{prop}
\pf
Since we have
$F(h_0)=R(h_0)+[\theta,\theta^{\dagger}]
 +\del_{h_0}\theta+\delbar\theta^{\dagger}$,
we have only to estimate
$[\theta,\theta^{\dagger}]$,
$\del_{h_0}\theta$ and $\delbar\theta^{\dagger}$.
We have the natural decompositions
$ f_i=\bigoplus f_{i\,\vecalpha}$ for $i=1,2$,
where 
$f_{i\,\vecalpha}\in \End(\prolongg{\vecc}{E}_{\vecalpha})$.
Since the decomposition of
$E=\bigoplus E_{\vecalpha}$ is orthogonal with respect to $h_0$,
the adjoint $f_i^{\dagger}$ of $f_i$ with respect to $h_0$ 
preserves the decomposition.
Hence we have the decomposition
$f_i^{\dagger}=\bigoplus f^{\dagger}_{i\,\vecalpha}$,
and $f^{\dagger}_{i\,\vecalpha}$ is adjoint of $f_{i\,\vecalpha}$
with respect to $h_{0|U_{\veca,\vecalpha}}$.

\vspace{.1in}

Let us show that
$\bigl[\theta,\theta^{\dagger}\bigr]$ is bounded
with respect to $h_0$ and $\omega_{\epsilon}$.
We put 
$N_{i}:=f_{i}-\bigoplus_{\vecalpha}
\alpha_i\cdot \id_{\prolongg{\vecc}{E_{\vecalpha}}}$
for $i=1,2$,
and then 
we have
$ \bigl[f_i,f_j^{\dagger}\bigr]
=\bigoplus_{\vecalpha}
 \bigl[N_{i},N^{\dagger}_{j}\bigr]$.
Since $(\prolongg{\vecc}{E}_{\ast},\theta)$ is graded semisimple,
we have 
$N_{1\,|\,D_1}\bigl(\lefttop{1}F_a\bigr)
\subset \lefttop{1}F_{<a}$.
We also have 
$N_{1\,|\,D_2}\bigl(\lefttop{2}F_a\bigr)
\subset
 \lefttop{2}F_{a}$.
Hence, we obtain
$\bigl|N_{1}\bigr|_{h_0}\leq C\cdot |z_1|^{2\epsilon}$
for some positive constant $C$.
Similarly we can obtain the estimate
$\bigl|N_{2}\bigr|_{h_0}\leq C\cdot |z_2|^{2\epsilon}$.
Thus we obtain the boundedness of
$[\theta,\theta^{\dagger}_{h_0}]$
with respect to $h_0$ and $\omega_{\epsilon}$.

\vspace{.1in}

Let us see the estimate of $\del_{h_0}\theta$.
We have the following,
where $\alpha_1$ denotes the first component of $\vecalpha$:
\[
 \del_{h_0}
 \left(f_1\cdot\frac{dz_1}{z_1}\right)
=\del_{h_0}\left(
\sum_{\vecalpha}\alpha_1\cdot\id_{E_{\vecalpha}}
 \cdot \frac{dz_1}{z_1}\right)
+
\del_{h_0'}\left(N_{1}\frac{dz_1}{z_1}\right)
+\left[A_2,\,\,N_{1}\frac{dz_1}{z_1}\right].
\]
The first term is $0$.
We put $\Omega:=dz_1\wedge dz_2/z_1\cdot z_2$.
Let us see the second term
$\del_{h_0'}N_{1}\cdot dz_1/z_1=:
 G_0\cdot \Omega$.
Then, $G_0$ is a $C^{\infty}$-section
of $\End(E)$ satisfying
$G_{0\,|\,D_1}(\lefttop{1}F_a)\subset \lefttop{1}F_{<a}$
and 
$G_{0\,|\,D_2}=0$.
Let us see the third term
$[A_2,N_{1}]\cdot dz_2/z_2
=:G_{1}\cdot \Omega$.
Then, $G_1$ is a $C^{\infty}$-section
of $\End(E)$ such that
$G_{1\,|\,D_i}\bigl(\lefttop{i}F_a\bigr)
\subset\lefttop{i}F_{<a}$.
Hence, the second and the third terms are
bounded.
Thus we obtain the boundedness of
$\del_{h_0}\theta$.
Since $\delbar\theta^{\dagger}_{h_0}$ is adjoint
of $\del_{h_0}\theta$ with respect to $h_0$,
it is also bounded.
Thus the proof of Proposition \ref{prop;06.8.15.100}
is finished.
\hfill\qed

\vspace{.1in}

\section{Around a Smooth Point of the Divisor}
\label{section;05.7.31.5}

\subsection{Setting}
\label{subsection;05.8.30.10}

Let $Y$ be a complex curve,
and $L$ be a line bundle on $Y$.
Let $\nbigu$ be a neighbourhood of $Y$ in $L$.
The projection $L\lrarr Y$ induces
$\pi:\nbigu\lrarr Y$.
Let $\sigma$ denote the canonical section
of $\pi^{\ast}L$.
Let $|\cdot|$ be a hermitian metric of $\pi^{\ast}L$.
Thus, we obtain the function
$\bigl|\sigma\bigr|:\nbigu\lrarr\real$.
Let $J_0$ denote the complex structure of $\nbigu$
as the open subset of $L$,
and let $J$ be any other integrable complex structure
such that $J-J_0=O(|\sigma|)$.
We regard $\nbigu$ as a complex manifold
via the complex structure $J$.
The $(0,1)$-operator $\delbar$ is induced by $J$.

Let $(E,\delbar_E)$ be a holomorphic vector bundle
on $\nbigu$.
We put $E_Y:=E_{|Y}$,
and let $F$ be a filtration of $E_Y$ in the category
of holomorphic vector bundles indexed by $\real$.
For later use,
we consider the case where
$S(F)=\bigl\{a\,\big|\,\Gr^F_a(E)\neq 0\bigr\}$
is not necessarily contained in 
an interval $\openclosed{c-1}{c}$ of the length $1$.
Thus $E_{\ast}=(E,F)$ is a parabolic bundle
in a slightly generalized sense
(Remark \ref{rem;06.8.15.101}).
We put $\gap(F):=\max\{|a-b|\neq 0\,|\,a,b\in S(F)\}$.
Let $\epsilon$ be a positive number
such that $10\epsilon<\gap(F)$.
Let $\omega$ be a Kahler form of $\nbigu$.
Take a small positive number $C$
and a large real number $N$.
Then, we put
$\omega_{\epsilon}:=\omega+C\cdot \epsilon^N
\sqrt{-1}\del\delbar |\sigma|^{2\epsilon}$,
which gives a Kahler form of
$\nbigu\setminus Y$.

Let $\theta$ be a Higgs field of $E_{\ast}$
in the sense of Remark \ref{rem;06.8.15.101}.
We put $f:=\Res(\theta)\in End(E_Y)$.

\begin{assumption}
The eigenvalues of $f$ are assumed to be constant on $Y$.
(See Remark {\rm\ref{rem;06.8.15.102}.})
\hfill\qed
\end{assumption}

\subsection{Construction of a metric}
\label{subsection;06.8.6.2}

We construct a hermitian metric of $E_{|\nbigu-Y}$
adapted to the filtration,
by following \cite{li2} and \cite{li-na} essentially.
(See also \cite{biquard2}.)
We have the generalized eigen decomposition
$E_Y=\bigoplus_{\alpha\in\cnum}
 \Gr^{\EE}_{\alpha}(E_Y)$ with respect to $f$.
We also have the generalized eigen decomposition
$\Gr^F_a(E_Y)=\bigoplus_{\alpha}\Gr^{F,\EE}_{(a,\alpha)}(E_Y)$
of $\Gr_a^F(E_Y)$ with respect to $\Gr^F(f)$.
Then we put
$\Ehat_{Y,u}:=\Gr^{F,\EE}_{u}(E_Y)$
for $u\in \real\times\cnum$,
and 
$ \Ehat_{Y}:=\bigoplus \Ehat_{Y,u}$.

Let $h'_0$ be a $C^{\infty}$-metric of $E$ on $\nbigu$.
The holomorphic structure of $E$
and the metric $h'_0$ induces the unitary connection
$\nabla_{0}$ of $E$ on $\nbigu$.
We put $h_Y:=h'_{0|Y}$.
We assume that the decomposition
 $E_Y=\bigoplus \Gr^{\EE}_{\alpha}(E_Y)$
is orthogonal with respect to $h_Y$.
The holomorphic structure of $E_Y$
and the metric $h_Y$ induce
the unitary connection $\nabla_{E_Y}$
of $E_Y$.
Thus the connection $\nabla_{\pi^{\ast}E_Y}$ is induced
on $\pi^{\ast}E_Y$.
Then, we can take a $C^{\infty}$-isometry
$\Phi:\pi^{\ast}E_Y\lrarr E$ such that
$\nabla_{0}\circ\Phi-\Phi\circ\pi^{\ast}\nabla_{E_Y}
=O(|\sigma|)$ with respect to $\omega$, as in \cite{li2}.
To see it,
we take any isometry $\Phi'$
such that $\Phi'_{|Y}$ is the identity.
We identify $E$ and $\pi^{\ast}E$
via $\Phi'$ for a while.
Let $\gminiu(E)$ be the bundle
of anti-hermitian endomorphisms of $E$.
We have the section
$A=\nabla_0-\nabla_{\pi^{\ast}E_Y}$
of $\gminiu(E)\otimes\Omega^1_{\nbigu}$.
We can take a $C^{\infty}$-section
$B$ of $\gminiu(E)$ such that
$B=O(|\sigma|)$
and $\nabla_{\pi^{\ast}E_Y}B-A=O(|\sigma|)$,
which can be easily checked
by using the partition of unity on $Y$.
Then we obtain
$g^{-1}\circ\nabla_{\pi^{\ast}E_Y}g-
 \nabla_0=O(|\sigma|)$
for $g=\exp(B)$,
which implies the  existence
of an appropriate isometry $\Phi$.
We identify $E$ and $\pi^{\ast}E_Y$
via such a $\Phi$
as $C^{\infty}$-bundles.

The metric $h_Y$ induces the orthogonal decomposition 
$\Gr^{\EE}_{\alpha}(E_Y)=
 \bigoplus_{a\in\real} \nbigg_{(a,\alpha)}$
such that 
$\bigoplus_{a\leq b}\nbigg_{(a,\alpha)}
=F_b\Gr^{\EE}_{\alpha}(E)$.
We have the natural $C^{\infty}$-isomorphism
$\nbigg_u\simeq \Ehat_{Y,u}$,
and thus $E_Y\simeq \Ehat_Y$.
We identify them as $C^{\infty}$-bundles
via the isomorphism.
Let $h_{Y,u}$ denote the restriction of $h_Y$
to $\nbigg_{u}$ for $u\in\real\times\cnum$.
We put $E_u:=\pi^{\ast}\nbigg_u$,
and thus $E=\bigoplus E_u$
and $h_0'=\pi^{\ast}h_Y=\bigoplus \pi^{\ast}h_{Y,u}$.
We put as follows:
\begin{equation}
 \label{eq;06.8.16.2}
 h_0:=\bigoplus \pi^{\ast}h_{Y,(a,\alpha)}\cdot |\sigma|^{-2a}.
\end{equation}

\subsection{Estimate of $R(h_0)$}

We put $\Gamma:=\bigoplus a\cdot \id_{E_{a,\alpha}}$.
\begin{lem}
 \label{lem;05.9.15.5}
$R(h_0,\delbar_E)$ is bounded with respect to
$\omega_{\epsilon}$ and $h_0$.
More strongly,
we have the following estimate,
with respect to $h_0$ and $\omega_{\epsilon}$:
\begin{equation}
\label{eq;06.8.16.21}
 R(h_0,\delbar_E)
=\bigoplus_{u\in\real\times\cnum}
 \pi^{\ast} R(h_{Y,u},\delbar_{\Ehat_{Y,u}})
+\Gamma\cdot\delbar\del\log|\sigma|^{-2}
+O\bigl(|\sigma|^{\epsilon}\bigr).
\end{equation}
\end{lem}
\pf
Let $\delbar_1$ denote the $(0,1)$-part of
$\pi^{\ast}\nabla_{\Ehat_Y}$.
Let $T$ denote the $(0,1)$-part of
$\nabla_0-\pi^{\ast}\nabla_{E_Y}$.
We put $S=\delbar_{E_Y}-\delbar_{\Ehat_Y}$.
We put $Q=T+\pi^{\ast}S$.
Then, we have 
$\delbar_E=\delbar_1+Q$.
We have $S(F_{a})\subset F_{<a}\otimes\Omega^{0,1}_Y$,
and $T_{|Y}=0$ in 
 $\bigl(\End(E)\otimes\Omega^{1}_{\nbigu}\bigr)_{|Y}$.
Hence, we have $Q=O(|\sigma|^{4\epsilon})$.
The operator $\del_{1,h_0}$ is determined by the condition
$\delbar h_0(u,v)=h_0(\delbar_1u,v)+h_0(u,\del_{1,h_0}v)$
for smooth sections $u$ and $v$ of $E$.
Similarly, we obtain the operator
$\del_{1,h_0'}$.

Let $Q^{\dagger}_{h_0}$ denote the adjoint of
$Q$ with respect to $h_0$,
and then 
$\del_{E,h_0}=\del_{1,h_0}-Q^{\dagger}_{h_0}$.
Hence we obtain
$R(\delbar_E,h_0)
=\bigl[\delbar_1,\del_{1,h_0}\bigr]
-\delbar_1Q^{\dagger}_{h_0}
+\del_{1,h_0}Q
-\bigl[Q,Q^{\dagger}_{h_0}\bigr]$.
Since $Q$ and $Q^{\dagger}_{h_0}$
are $O(|\sigma|^{4\epsilon})$ with respect to
$\omega_{\epsilon}$ and $h_0$,
so is $[Q,Q^{\dagger}_{h_0}]$.
We have 
$\del_{1,h_0}Q=
 \del_{1,h_0'}Q
+\del\log|\sigma|^{-2}[\Gamma,Q]$.
Since $Q$ is sufficiently small,
the second term is $O(|\sigma|^{2\epsilon})$
with respect to $\omega_{\epsilon}$ and $h_0$.
Since $T_{|Y}$ is $0$
in $\bigl(\End(E)\otimes\Omega^{1}_{\nbigu}\bigr)_{|Y}$,
we have 
$\del_{1,h_0'}T=O(|\sigma|^{2\epsilon})$
with respect to $\omega_{\epsilon}$ and $h_0$.
Since $\bigl(\del_{1,h_0'}S\bigr)_{|Y}(F_{a})
\subset \bigl(
F_{<a}\otimes\Omega^{1,0}(\log Y)\otimes\Omega^{0,1}
\bigr)_{|Y}$,
we have $\del_{1,h_0'}S=O(|\sigma|^{2\epsilon})$
with respect to $h_0$ and $\omega_{\epsilon}$.
Thus, $\del_{1,h_0}Q$
and the adjoint 
$\delbar_1 Q^{\dagger}_{h_0}$
are also $O(|\sigma|^{2\epsilon})$
with respect to $\omega_{\epsilon}$ and $h_0$.
We have
$\bigl[\delbar_1,\del_{1,h_0}\bigr]=
 \bigl[\delbar_1,\del_{1,h_0'}\bigr]
+\Gamma\cdot \delbar\del\log|\sigma|^{-2}$.
Since we have
$\delbar_1+\del_{1,h_0'}=\nabla_{\pi^{\ast}\Ehat_Y}$
by our construction,
we obtain
$\bigl[\delbar_1,\del_{1,h_0'}\bigr]
=\pi^{\ast}R(h_Y,\delbar_{\Ehat_Y})
+[\delbar_1,\delbar_1]+[\del_{1,h_0'},\del_{1,h_0'}]
=\pi^{\ast}R(h_Y,\delbar_{\Ehat_Y})
 +O(|\sigma|^{2\epsilon})$
with respect to $\omega_{\epsilon}$ and $h_0$.
Thus Lemma \ref{lem;05.9.15.5}
is proved.
\hfill\qed

\begin{cor}
\label{cor;06.8.16.41}
We have the following estimate
with respect to $\omega_{\epsilon}$:
\[
 \tr R(h_0,\delbar_E)
=\sum_{(a,\alpha)}
 \pi^{\ast}\tr R(h_{Y,(a,\alpha)},\delbar_{E_{Y,(a,\alpha)}})
+\sum a\cdot\rank \Gr^{F}_a(E)\cdot \delbar\del\log|\sigma|^{-2}
+O(1)
\]
\hfill\qed
\end{cor}

\subsection{Estimate of $F(h_0)$ in the graded semisimple case}

\begin{prop}
\label{prop;05.8.30.50}
If $(E_{\ast},\theta)$ is graded semisimple,
$F(h_0)$ is bounded
with respect to $\omega_{\epsilon}$ and $h_0$.
\end{prop}
\pf
We put $\rho_0:=\bigoplus \alpha\cdot\id_{E_{(a,\alpha)}}$
and $\rhobar_0:=\bigoplus \alphabar\cdot\id_{E_{(a,\alpha)}}$.
Let $P$ be any point of $Y$.
Let $(U,z_1,z_2)$ be a holomorphic coordinate neighbourhood
of $\bigl(\nbigu,J\bigr)$ around $P$
such that $U\cap Y=\{z_1=0\}$.
We are given the Higgs field:
\[
 \theta=f_1\cdot\frac{dz_1}{z_1}+f_2\cdot dz_2.
\]

Since $f_{2|Y}$ preserves the filtration $F$,
$f_2$ is bounded with respect to $h_0$.
It is easy to see $[\rho_0,f_2]_{|Y}=0$.
Hence $[\rho_0,f_2]$ is $O(|\sigma|^{2\epsilon})$
with respect to $h_0$.
We put $f_1'=f_1-\rho_0$.
Due to the graded semisimplicity of $(E_{\ast},\theta)$,
we have $f'_{1|Y}\bigl(F_{a}\bigr)\subset F_{<a}$.
Hence $f_1'$ is $O(|\sigma|^{2\epsilon})$ with respect to $h_0$.
Then it is easy to check
the boundedness of $[\theta,\theta^{\dagger}]$
with respect to $\omega_{\epsilon}$ and $h_0$,
by a direct calculation.

We have the following:
\[
 \del_{E,h_0}(f_1)\cdot \frac{dz_1}{z_1}
=\del_{1,h_0'}(f_1')\cdot\frac{dz_1}{z_1}
+\bigl[\Gamma,f_1'\bigr]
 \cdot\del\log|\sigma|^2\cdot \frac{dz_1}{z_1}
-\bigl[Q^{\dagger}_{h_0},f_1\bigr]\cdot\frac{dz_1}{z_1}
\]
Then, 
$\del_{1,h_0'}f_1'=A\cdot dz_2\cdot dz_1/z_1$ is 
$C^{\infty}$-$(2,0)$-form of $\End(E)$,
and $A_{|Y}(F_{a})\subset F_{<a}$.
Hence the first term is $O(|\sigma|^{2\epsilon})$
with respect to $\omega_{\epsilon}$ and $h_0$.
Similarly, the same estimate holds for the second term.
Since $Q^{\dagger}_{h_0}=O(|\sigma|^{2\epsilon})$,
the third term is $O(|\sigma|^{\epsilon})$.

We have
$ \del_{E,h_0}f_2\cdot dz_2
=\del_{1,h_0'}f_2\cdot dz_2
+[\Gamma,f_2]\cdot\del\log|\sigma|^2\cdot dz_2
-[Q^{\dagger}_{h_0},f_2]\cdot dz_2$.
Since the first term is $C^{\infty}$-$2$-form
of $\End(E)$,
it is $O(|\sigma|^{2\epsilon})$
with respect to $\omega_{\epsilon}$ and $h_0$.
We have $[\Gamma,f_2](F_a)\subset F_{<a}$,
the same estimate holds for the second term.
Since $Q^{\dagger}_{h_0}$ is $O(|\sigma|^{2\epsilon})$,
the third term is $O(|\sigma|^{\epsilon})$
with respect to $\omega_{\epsilon}$ and $h_0$.
Then Proposition \ref{prop;05.8.30.50} is proved.
\hfill\qed

\subsection{Preliminary for the calculation of the integral}

Let $\hhat_Y=\bigoplus \hhat_{Y,u}$
be a hermitian metric of $E_Y$
for which $\bigoplus \Ehat_{Y,u}$ is orthogonal.
We put $\hhat:=\pi^{\ast}\hhat_Y$.
We put $A:=\del_{E,h_0}-\del_{E,\hhat}$.

\begin{lem}
We have the following estimates
with respect to $\omega_{\epsilon}$:
\begin{equation}
\label{eq;06.8.16.20}
\tr A=\sum a\cdot\rank\Gr^F_a(E)\cdot\del\log|\sigma|^{-2}
+O(1)
\end{equation}
\begin{multline}
\label{eq;06.8.16.30}
\tr\bigl(A\cdot R(h_0)\bigr)
=\sum \pi^{\ast}\tr R(\Ehat_{Y,u},h_{Y,u})\cdot
 a\cdot \del\log|\sigma|^{-2} \\
+\sum \rank \Ehat_{Y,u}\cdot a^2\cdot\delbar\del\log|\sigma|^{-2}
 \del\log|\sigma|^{-2}
+\tr \Bigl(
 Q^{\dagger}_{\hhat}\cdot\bigl[\Gamma\cdot\del\log|\sigma|^{-2},\,
 Q\bigr]\Bigr)
+O(1)
\end{multline}
\begin{multline}
\label{eq;06.8.16.32}
 \tr \bigl(A\cdot  R(\hhat)\bigr)=
\sum_{u}
 \pi^{\ast} \tr R(\Ehat_{Y,u},\hhat_u) \cdot
 a\cdot \del\log|\sigma|^{-2} \\
-\tr\bigl(\Gamma\cdot\del\log|\sigma|^{-2}
 \bigl[Q,Q^{\dagger}_{\hhat}\bigr] \bigr)
+O(1)
\end{multline}
Here, $u=(a,\alpha)$.
\end{lem}
\pf
We have
$\del_{E,h_0}=\del_{1,h_0'}-Q^{\dagger}_{h_0}
+\Gamma\cdot\del\log|\sigma|^{-2}$
and $\del_{E,\hhat}=\del_{1,\hhat}-Q^{\dagger}_{\hhat}$.
We put
$P=\del_{1,h_0'}-\del_{1,\hhat}$,
which is a $C^{\infty}$-section of
$\bigoplus \End(E_u)\otimes\Omega^{1,0}$.
Thus, we have
$A=P+Q^{\dagger}_{\hhat}-Q^{\dagger}_{h_0}
 +\Gamma\cdot\del\log|\sigma|^{-2}$.
Since $Q^{\dagger}_{\hhat}$ and $Q^{\dagger}_{h_0}$
are bounded with respect to $(\omega_{\epsilon},\hhat)$,
we obtain (\ref{eq;06.8.16.20}).

Let us show (\ref{eq;06.8.16.30}).
Since $P+Q^{\dagger}_{h_0}$ is bounded
with respect to $h_0$ and $\omega_{\epsilon}$,
we have the boundedness of
$\tr\bigl((P+Q^{\dagger}_{h_0})\cdot R(h_0)\bigr)$
with respect to $\omega_{\epsilon}$.
From (\ref{eq;06.8.16.21}),
we obtain the following:
\begin{multline}
 \tr\bigl(
 \Gamma\cdot\del\log|\sigma|^{-2}
\cdot R(h_0)
 \bigr)
=\sum_{a,\alpha} \pi^{\ast}\tr R(\Ehat_{Y,a,\alpha},h_{Y,a,\alpha})
 \cdot a\cdot\del\log|\sigma|^{-2} \\
+\sum_{a,\alpha}
 \rank\Ehat_{Y,a,\alpha}\cdot a^{2}\cdot\delbar\del\log|\sigma|^{-2}
 \cdot\del \log|\sigma|^{-2}
+O(1).
\end{multline}
Let us see $\tr\bigl(Q_{\hhat}^{\dagger}\cdot R(h_0)\bigr)$.
We decompose it as follows:
\begin{equation}
\label{eq;06.8.16.25}
\tr\bigl(Q_{\hhat}^{\dagger}\cdot [\delbar_1,\del_{1,h_0}] \bigr)
-\tr\bigl(Q_{\hhat}^{\dagger}\cdot \delbar_1Q^{\dagger}_{h_0} \bigr)
+\tr\bigl(Q_{\hhat}^{\dagger}\cdot \del_{1,h_0}Q \bigr)
-\tr\bigl(Q_{\hhat}^{\dagger} \cdot[Q,Q^{\dagger}_{h_0}] \bigr)
\end{equation}
Since $[\delbar_1,\del_{1,h_0}]$ is bounded 
with respect to $(\omega_{\epsilon},\hhat)$,
we obtain the boundedness of the first term.
Recall
$Q^{\dagger}_{h_0}=
 (\pi^{\ast}S)^{\dagger}_{h_0}+T^{\dagger}_{h_0}$.
Because of $T_{|Y}=0$ in $(\End(E)\otimes \Omega^{0,1})_{|Y}$
and $\del_{1,h_0}T=\del_{1,h_0'}T+
 \bigl[\Gamma\cdot\del\log|\sigma|^{-2},T\bigr]$,
we have $\del_{1,h_0}T_{|Y}=0$ in 
$\bigl(\End(E)\otimes\Omega^{1,0}(\log Y)\otimes\Omega^{0,1}\bigr)_{|Y}$.
Because of
$\delbar_1 T_{h_0}^{\dagger}=
 \bigl( \del_{1,h_0}T\bigr)^{\dagger}_{h_0}$,
it is easy to obtain
$\delbar_1 T_{h_0}^{\dagger}
=O(|\sigma|^{2\epsilon})$ with respect to
$(\hhat,\omega_{\epsilon})$.
We also have 
$T_{h_0}^{\dagger}=O(|\sigma|^{2\epsilon})$ with respect to
$(\hhat,\omega_{\epsilon})$.
Since
$\pi^{\ast}S$
is a  section of
$\bigoplus_{a>a'} Hom(E_{a,\alpha},E_{a',\alpha'})\otimes\Omega^{0,1}$,
we have
$\pi^{\ast}S^{\dagger}_{h_0}=O(|\sigma|^{2\epsilon})$
with respect to $(\hhat,\omega_{\epsilon})$.
Hence, $Q^{\dagger}_{h_0}$ 
and $[Q^{\dagger}_{h_0},Q]$
are $O(|\sigma|^{2\epsilon})$
with respect to $(\omega_{\epsilon},\hhat)$.
Therefore, the fourth term in (\ref{eq;06.8.16.25})
is bounded.
Because of
$\delbar_1\pi^{\ast}S^{\dagger}_{h_0}
=\bigl( \del_{1,h_0'}\pi^{\ast}S\bigr)^{\dagger}_{h_0}
+\bigl( [\Gamma\cdot\del\log|\sigma|^{-2},\pi^{\ast}S]
 \bigr)^{\dagger}_{h_0}$,
it is easy to obtain
$\delbar_1\pi^{\ast}S^{\dagger}_{h_0}
=O(|\sigma|^{2\epsilon})$
with respect to $(\omega_{\epsilon},\hhat)$.
Together with the estimate of $\delbar_1 T^{\dagger}_{h_0}$
above,
we obtain the boundedness of
$\delbar_1Q^{\dagger}_{h_0}$
with respect to $(\omega_{\epsilon},\hhat)$.
Hence, we obtain the boundedness
of the second term in (\ref{eq;06.8.16.25}).
We have
$\del_{1,h_0}Q=\del_{1,h_0'}Q+[\Gamma\cdot\del\log|\sigma|^{-2},Q]$,
and $\del_{1,h_0'}Q$ is bounded
with respect to $(\omega_{\epsilon},\hhat)$.
Therefore, the third term is
$O(1)+\tr\bigl(Q^{\dagger}_{h_0}[\Gamma\cdot\del\log|\sigma|^{-2},Q]\bigr)$.
Thus we obtain (\ref{eq;06.8.16.30}).

Let us show (\ref{eq;06.8.16.32}).
Since $P$, $Q^{\dagger}_{\hhat}$ and $Q^{\dagger}_{h_0}$
are bounded with respect to $(\omega_{\epsilon},\hhat)$,
we have $\tr\bigl(
(P+Q^{\dagger}_{\hhat}-Q^{\dagger}_{h_0})R(\hhat)\bigr)=O(1)$
with respect to $\omega_{\epsilon}$.
We have
$R(\hhat)=[\delbar_1,\del_{1,\hhat}]
-\delbar_1 Q^{\dagger}_{\hhat}
+\del_{1,\hhat}Q
-[Q^{\dagger}_{\hhat},Q]$.
Because of
$\del_{1,\hhat}T=O(|\sigma|^{2\epsilon})$
with respect to $(\omega_{\epsilon},\hhat)$
and 
$\del_{1,\hhat}\pi^{\ast}S\in
 \bigoplus_{a>a'}\Hom(E_{a,\alpha},E_{a',\alpha'})
 \otimes\Omega^{2}$,
we have 
$\tr\bigl(\Gamma\cdot \del\log|\sigma|^{-2}\cdot
 \del_{1,\hhat}Q\bigr)=O(|\sigma|^{2\epsilon})$
with respect to $\omega_{\epsilon}$.
By a similar reason,
$\tr\bigl(\Gamma\cdot\del\log|\sigma|^{-2}
 \delbar_1Q^{\dagger}_{\hhat}\bigr)=O(|\sigma|^{2\epsilon})$.
Since we have
$[\delbar_1,\del_{1,\hhat}]
=\pi^{\ast}R(\Ehat,\hhat_Y)+O(|\sigma|^{2\epsilon})$
with respect to $(\hhat,\omega_{\epsilon})$,
we obtain (\ref{eq;06.8.16.32}).
\hfill\qed

\begin{cor}
\label{cor;06.8.16.40}
We have the following estimates with respect to $\omega_{\epsilon}$:
\begin{multline}
 \tr\bigl(A\cdot R(h_0)+A\cdot R(\hhat)\bigr)
=\sum\pi^{\ast}\bigl(
 \tr R(\Ehat_{Y,u},h_{Y,u})
+\tr R(\Ehat_{Y,u},\hhat_{Y,u})
\bigr)
 \cdot a\cdot\del\log|\sigma|^{-2}\\
+ \sum a^2\cdot \rank \Ehat_{Y,u}\cdot\del\log|\sigma|^{-2}
 \cdot\del\delbar\log|\sigma|^{-2}
\end{multline}
Here, $u=(a,\alpha)$.
\hfill\qed
\end{cor}

\subsection{Estimate of a related metric}
\label{subsection;06.8.15.20}

For later use (Section \ref{section;06.8.17.100}),
we consider a related metric
in the case where one more filtration $W$ is given on
$\Gr_{(a,\alpha)}^{F,\EE}(E)$ indexed by $\seisuu$.
We put $\Etilde_{u,k}:=\Gr^W_k\Gr^{F,\EE}_{u}(E_Y)$
for $(u,k)\in(\real\times\cnum)\times \seisuu$
and $\Etilde_Y:=\bigoplus \Etilde_{u,k}$.
We put $\Ftilde_{(a,k)}\Gr^{\EE}_{\alpha}(E)
 :=\pi_a^{-1}(W_k)$,
where $\pi_a$ denotes the projection
$F_a\Gr^{\EE}_{\alpha}(E)\lrarr 
 \Gr_{(a,\alpha)}^{F,\EE}(E)$.

The metric $h_Y$ induces
the orthogonal decomposition
$\Gr^{\EE}_{\alpha}(E)=
 \bigoplus_{(a,k)\in\real\times\cnum}
 \nbigg_{a,\alpha,k}$
such that
$\Ftilde_{(b,l)}\Gr^{\EE}_{\alpha}(E)
=\bigoplus_{(a,k)\leq (b,l)}
 \nbigg_{a,\alpha,k}$.
We have the natural $C^{\infty}$-isomorphism
$\nbigg_{u,k}\simeq \Gr^W_k\Gr^{F,\EE}_{u}(E_Y)$
for $(u,k)\in(\real\times\cnum)\times\seisuu$.
Thus, we obtain the $C^{\infty}$-identification
of $E_Y$ and $\Etilde_Y$.
Let $h_{Y,u,k}$ denote the restriction of $h_Y$
to $\nbigg_{u,k}$.

Via identification $\Phi:\pi^{\ast}E_Y\simeq E$,
we obtain the $C^{\infty}$-decomposition
$E=\bigoplus E_{a,\alpha,k}$.
Then, we put as follows:
\[
 h_1:=\bigoplus_{a,\alpha,k}
 \pi^{\ast}h_{Y,a,\alpha,k}\cdot
 |\sigma|^{-2a}\cdot\bigl(-\log|\sigma|^2\bigr)^{k}.
\]
The metrics $h_0$ and $h_1$ are mutually bounded
up to log order,
i.e., $C^{-1}\cdot h_0\cdot (-\log|\sigma|)^{-N}
\leq h_1\leq C\cdot h_0\cdot (-\log|\sigma|)^N$
for some constants $C$ and $N$.

For appropriate constants $C_1$,
we put 
$\omegatilde:=\omega+C_1\cdot \del\delbar\log(-\log|\sigma|^2)$,
which gives the Poincar\'e like metric on $\nbigu\setminus Y$.

\begin{lem}
\label{lem;06.8.15.45}
$R(h_1)$ is bounded with respect to
$\omegatilde$ and $h_i$ $(i=0,1)$.
The difference $\del_{E,h_1}-\del_{E,h_0}$
is bounded with respect to 
$\omegatilde$ and $h_0$.
\end{lem}
\pf
Under the identification $E_Y=\Etilde_Y$,
we put $\Stilde=\delbar_{E_Y}-\delbar_{\Etilde_Y}$.
We put $S':=\Stilde-S$.
As before, we have
$\delbar_E=\delbar_2+\Qtilde$ and
$\Qtilde=T+\pi^{\ast}\Stilde$.
We also have
$\delbar_1=\delbar_2+\pi^{\ast}S'$.
Because of $T_{|Y}=0$ in 
$\bigl(\End(E)\otimes\Omega^1_{\nbigu}\bigr)_{|Y}$,
$T$ and $T^{\dagger}_{h_1}$
are $O(|\sigma|^{2\epsilon})$
with respect to $(h_i,\omegatilde)$ $(i=0,1)$.
Because of
$\Stilde(\Ftilde_{(a,k)})\subset 
 \Ftilde_{<(a,k)}\otimes\Omega^{0,1}_Y$,
$\Stilde$ and $\Stilde^{\dagger}_{h_1}$
are  $O\bigl((-\log|\sigma|)^{-1/2}\bigr)$
with respect to $(h_1,\omegatilde)$.
We also obtain 
$\Stilde=O(1)$
and $\Stilde^{\dagger}_{h_1}=O\bigl((-\log|\sigma|^{-1/2})\bigr)$
with respect to $(h_0,\omegatilde)$.
In particular, $\Qtilde$ and $\Qtilde^{\dagger}_{h_1}$
are bounded with respect to $(h_i,\omegatilde)$ $(i=0,1)$.

We put $\nbigk:=\bigoplus k/2\cdot \id_{E_{u,k}}$.
Then, we obtain the following:
\begin{multline}
\del_{E,h_1}=\del_{2,h_1}-\Qtilde^{\dagger}_{h_1}
=\del_{2,h_0}+\nbigk\cdot\del\log(-\log|\sigma|^{2})
 -\Qtilde^{\dagger}_{h_1}\\
=\del_{1,h_0}+(\pi^{\ast}S')^{\dagger}_{h_0}
+\nbigk\cdot\del\log(-\log|\sigma|^2)
-\Qtilde^{\dagger}_{h_1}\\
=
 \del_{E,h_0}+Q^{\dagger}_{h_0}
+\bigl(\pi^{\ast}S'\bigr)^{\dagger}_{h_0}
+\nbigk\cdot\del\log(-\log|\sigma|^2)
-\Qtilde^{\dagger}_{h_1}.
\end{multline}
It is easy to see that $\pi^{\ast}S'$ and
$(\pi^{\ast}S')^{\dagger}_{h_0}$
are bounded with respect to $h_0$.
Thus, we obtain the boundedness of
$\del_{E,h_1}-\del_{E,h_0}$
with respect to $(\omegatilde,h_0)$.

We decompose $R(h_1)$ as follows:
\begin{equation}
\label{eq;06.8.15.2}
R(h_1)
=\bigl[\delbar_2,\del_{2,h_1}\bigr]
+\del_{2,h_1}\Qtilde
-\delbar_2\Qtilde^{\dagger}_{h_1}
-\bigl[\Qtilde,\Qtilde^{\dagger}_{h_1}\bigr]
\end{equation}
We decompose the second term as follows:
\begin{multline}
\label{eq;06.8.15.1}
\bigl[\del_{2,h_1},\Qtilde\bigr]
=\bigl[\nbigk\cdot\del\log(-\log|\sigma|^2),\Qtilde\bigr] \\
+\bigl[\del_{2,h_0'}+\Gamma\cdot\del\log|\sigma|^{-2},
 T\bigr]
+\bigl[\del_{2,h_0'}+\Gamma\cdot\del\log|\sigma|^{-2},
 \Stilde\bigr]
\end{multline}
Since $\del\log(-\log|\sigma|^2)$ is bounded with respect
to $\omegatilde$,
we have the boundedness of
$\nbigk\cdot\del\log(-\log|\sigma|^2)$
with respect to $(\omegatilde,h_i)$ $(i=0,1)$.
Hence, the first term in (\ref{eq;06.8.15.1}) is bounded.
The adjoint with respect to $h_1$
also satisfies the same estimate.

We have $T=O(|\sigma|^{3\epsilon})$
with respect to $(\omegatilde,h_i)$ $(i=0,1)$
and $[\del_{2,h_0'},T]_{|Y}=0$
in $\bigl(\End(E)\otimes\Omega^{1,0}(\log D)\otimes\Omega^{0,1}
\bigr)_{|Y}$.
Hence 
$\bigl[\Gamma\cdot\del\log|\sigma|^2,T\bigr]$
and 
$[\del_{2,h_0'},T]$
are $O(|\sigma|^{3\epsilon})$ with respect to $(\omegatilde,h_i)$
$(i=0,1)$.
Their adjoints with respect to $h_1$
are also $O(|\sigma|^{2\epsilon})$
with respect to $(\omegatilde,h_i)$.
Therefore, we obtain the boundedness
of the second term in (\ref{eq;06.8.15.1})
and the adjoint.

Let $\Stilde=A\cdot d\zbar_1+B\cdot d\zbar_2$ be the expression
for a local coordinate $(U,z_1,z_2)$
such that $z_1^{-1}(0)=Y\cap U$.
Then, we have $A_{|Y}=0$
and $B_{|Y}(\Ftilde_{(a,k)})\subset \Ftilde_{<(a,k)}$.
We have $[\Gamma,B]_{|Y}(F_{a})\subset F_{<a}$.
Thus $[\Gamma\cdot\del\log|\sigma|^{-2},\Stilde]$,
and the adjoint with respect to $h_1$
are $O(|\sigma|^{2\epsilon})$ with respect to
$(\omegatilde,h_i)$ $(i=0,1)$.
We have
$\bigl[\del_{2,h_0'},\,A\cdot d\zbar_1\bigr]_{|Y}=0$
in $\bigl(\End(E)\otimes
 \Omega^{1,0}(\log Y)\otimes\Omega^{0,1}\bigr)_{|Y}$.
For the expression
$\bigl[\del_{2,h_0'},Bd\zbar_2\bigr]
=\bigl(
C_1\cdot dz_1/z_1+C_2dz_2\bigr)\cdot d\zbar_2$,
we have $C_{1|Y}=0$ and 
$C_{2|Y}(\Ftilde_{(a,k)})\subset \Ftilde_{<(a,k)}$.
Hence, 
$\bigl[\del_{2,h_0'},\Stilde\bigr]$
and the adjoint with respect to $h_1$
are bounded with respect to both of $(\omegatilde,h_i)$ $(i=0,1)$.
Therefore, we obtain the boundedness
of the third term in (\ref{eq;06.8.15.1})
and the adjoint.
Thus we obtain the boundedness of the 
second and third terms in (\ref{eq;06.8.15.2}).

We have
$\bigl[\delbar_2,\del_{2,h_1}\bigr]=
 [\delbar_2,\del_{2,h_0'}]
+\delbar\del\log|\sigma|^{-2}\cdot\Gamma
+\delbar\del\log(-\log|\sigma|^2)\cdot\nbigk$
which is bounded with respect to $(\omegatilde,h_i)$ 
$(i=0,1)$.
Thus we obtain the boundedness of $R(h_1)$.
\hfill\qed

\section{Global Ordinary Metric}
\label{section;05.8.30.100}

\subsection{Decomposition and metric of a base space}
\label{section;05.8.23.5}

Let $X$ be a smooth projective surface,
and $D$ be a simple normal crossing divisor
with the irreducible decomposition $D=\bigcup_{i\in S}D_i$.
We also assume that $D$ is ample.
Let $L$ be an ample line bundle on $X$,
and $\omega$ be a Kahler form which represents $c_1(L)$.
For any point $P\in D_i\cap D_j$,
we take a holomorphic coordinate
$(U_P,z_i,z_j)$ around $P$
such that $U_P\cap D_k=\{z_k=0\}$ $(k=i,j)$
and $U_P\simeq \Delta^2$ by the coordinate.
Let us take a hermitian metric $g_i$ of $\nbigo(D_i)$
and the canonical section $\nbigo\lrarr\nbigo(D_i)$
is denoted by $\sigma_i$.
We may assume
$|\sigma_k|^2_{g_k}=|z_k|^2$ $(k=i,j)$
on $U_P$ for $P\in D_i\cap D_j$.

Let us take a hermitian metric $g$ of the tangent bundle $TX$
such that $g=dz_i\cdot d\bar{z}_i+dz_j\cdot d\bar{z}_j$ on $U_P$.
It is not necessarily same as $\omega$.
The metric $g$ induces the exponential map
$\exp:TX\lrarr X$.
Let $N_{D_i}X$ denote the normal bundle of $D_i$ in $X$.
We can take a sufficiently small neighbourhood 
$U_i'$ of $D_i$ in $N_{D_i}X$
such that the restriction of $\exp_{|U_i'}$
gives the diffeomorphism of $U_i'$
and the neighbourhood $U_i$ of $D_i$ in $X$.
We may assume $U_i\cap U_j=\coprod_{P\in D_i\cap D_j}U_P$.

Let $p_i$ denote the diffeomorphism
$\exp_{|U_i'}:U_i'\lrarr U_i$.
Let $\pi_i$ denote the natural projection
$U_i'\lrarr D_i$.
Via the diffeomorphism $p_i$,
we also have the $C^{\infty}$-map $U_i\lrarr D_i$,
which is also denoted by $\pi_i$.
On $U_P$, $\pi_i$ is same as the natural projection 
$(z_i,z_j)\longmapsto z_j$.
Via $p_i$, we have two complex structures
$J_{U_i}$ and $J_{U_i'}$ on $U_i$.
Due to our choice of the hermitian metric $g$,
$p_i$ preserves the holomorphic structure
(i.e., $J_{U_i'}-J_{U_i}=0$)
on $U_P$.
The derivative of $p_i$ gives the isomorphism of
the complex bundles
$T(N_{D_i}(X))_{|D_i}\simeq TD_i\oplus N_{D_i}X\simeq TX_{|D_i}$
on $D_i$.
Hence we have $J_{U_i}-J_{U_i'}=O(|\sigma|)$.

\vspace{.1in}

Let $\epsilon$ be any number such that $0<\epsilon<1/2$.
Let us fix a real number $N$,
which is sufficiently large, say $N>10$.
We put as follows, for some positive number $C>0$:
\[
 \omega_{\epsilon}:= \omega+
\sum_i C\cdot \epsilon^N\cdot
       \sqrt{-1}\del\delbar |\sigma_i|_{g_i}^{2\epsilon}.
\]

\begin{prop}\label{prop;05.7.30.3}
If $C$ is sufficiently small,
then $\omega_{\epsilon}$ are Kahler metrics of $X-D$
for any $0<\epsilon<1/2$.
\end{prop}
\pf
We put $\phi_i:=|\sigma_i|_{g_i}^2$.
We have
$ \sqrt{-1}\cdot\del\delbar \phi_i^{\epsilon}
=\sqrt{-1}\cdot\epsilon^2\cdot\phi_i^{\epsilon}\cdot
 \del\log\phi_i\cdot\delbar \log\phi_i
+\sqrt{-1}\cdot\epsilon\cdot\phi_i^{\epsilon}\cdot
 \del\delbar\log\phi_i$.
Hence the claim of Proposition \ref{prop;05.7.30.3}
immediately follows from the next lemma.
\hfill\qed

\begin{lem}
We put $f_t(\epsilon):=\epsilon^{l}\cdot t^{2\epsilon}$
for $0<\epsilon\leq 1/2$ and for $l\geq 1$.
The following inequality holds:
\begin{equation} \label{eq;05.7.30.4} 
 f_t(\epsilon)\leq
 \left(
 \frac{l}{-\log t^2}
 \right)^{l}\cdot e^{-l}
\quad\quad
(0<t<e^{-l})
\end{equation}
\begin{equation} \label{eq;05.7.30.5}
 f_t(\epsilon)\leq \left(\frac{1}{2}\right)^{l}\cdot t
\quad\quad
(t\geq e^{-l})
\end{equation}
\end{lem}
\pf
We have
$f_t'(\epsilon)=\epsilon^{l-1}t^{2\epsilon}\cdot 
 \bigl(l+\epsilon\log t^2\bigr)$.
If $t<e^{-l}$, we have $\epsilon_0:=l\times(-\log t^2)^{-1}<1/2$
and $f'_t(\epsilon_0)=0$.
Hence $f_t$ takes the maximum at $\epsilon=\epsilon_0$,
and we obtain (\ref{eq;05.7.30.4}).
If $t\geq e^{-1}$,
we have $f'_t(\epsilon)> 0$ for any $0<\epsilon<1/2$,
and thus $f_t(\epsilon)$ takes the maximum at $\epsilon=1/2$.
Thus we obtain (\ref{eq;05.7.30.5}).
\hfill\qed

\vspace{.1in}
The Kahler forms $\omega_{\epsilon}$ behave well
around any point of $D$ in the following sense,
which is clear from the construction.
\begin{lem}
Let $P$ be any point of $D_i\cap D_j$.
Then there exist positive constants $C_i$ $(i=1,2)$
such that the following holds  on $U_P$,
for any $0<\epsilon<1/2$:
\[
 C_1\cdot \omega_{\epsilon}
\leq
\sqrt{-1}\cdot\epsilon^{N+2}\cdot
 \left(
 \frac{dz_i\cdot d\bar{z}_i}{|z_i|^{2-2\epsilon}}
+\frac{dz_j\cdot d\bar{z}_j}{|z_j|^{2-2\epsilon}}
 \right)
+\sqrt{-1}\bigl(
dz_i\cdot d\bar{z}_i
+dz_j\cdot d\bar{z}_j\bigr)
\leq
 C_2\cdot\omega_{\epsilon}.
\]
Let $Q$ be any point of $D_i^{\circ}$,
and $(U,w_1,w_2)$ be a holomorphic coordinate around $Q$
such that $U\cap D_i=\{w_1=0\}$.
Then there exist positive constants $C_i$ $(i=1,2)$
such that the following holds for any $0<\epsilon<1/2$ on $U$:
\[
 C_1\cdot \omega_{\epsilon}
\leq
\sqrt{-1}\cdot\epsilon^{N+2}\cdot
 \left(
 \frac{dw_1\cdot d\bar{w}_1}{|w_1|^{2-2\epsilon}}
 \right)
+\sqrt{-1}\bigl(
  dw_1\cdot d\bar{w}_1
+dw_2\cdot d\bar{w}_2\bigr)
\leq
 C_2\cdot\omega_{\epsilon}.
\]
\hfill\qed
\end{lem}

\begin{lem}[Simpson \cite{s1}, Li \cite{li2}]
Let us consider the case $\epsilon=1/m$ for 
some positive integer $m$.
Then the metric $\omega_{\epsilon}$ satisfies
Condition {\rm\ref{condition;05.7.30.1}}.
\end{lem}
\pf
We use the argument of Simpson in \cite{s1}.
The first condition is easy to check.
Since we have assumed that $D$ is ample,
we can take a $C^{\infty}$-metric $|\cdot|$ of $\nbigo(D)$
with the non-negative curvature.
We put $\phi:=-\log |\sigma|$,
where $\sigma$ denote the canonical section.
Then $\sqrt{-1}\del\delbar\phi$ is a non-negative 
$C^{\infty}$-two form,
and it is easy to check that the second condition is satisfied.

To check the condition 3,
we give the following remark.
Let $P$ be a point of $D_i\cap D_j$.
For simplicity, let us consider the case $(i,j)=(1,2)$.
We put $V_P:=\bigl\{(\zeta_1,\zeta_2)\,\big|\,|\zeta_i|<1\bigr\}$.
Let us take the ramified covering
$\varphi:V_P\lrarr U_P$
given by $(\zeta_1,\zeta_2)\longmapsto (\zeta_1^m,\zeta_2^m)$.
Then it is easy to check that
$\omegatilde=\varphi^{-1}\omega_{\epsilon}$ naturally gives
the $C^{\infty}$-Kahler form on $V_P$.
If $f$ is a bounded positive function on $U_P\setminus D$
satisfying $\Delta_{\omega_{\epsilon}}(f)\leq B$
for some constant $B$,
we obtain
$\Delta_{\widetilde{\omega}}\bigl(\varphi^{\ast}f\bigr)
 \leq B$ on $V_P-\varphi^{-1}(D\cap U_P)$.
Since $\widetilde{\omega}$ is $C^{\infty}$ on $V_P$,
we may apply the argument of Proposition 2.2 in \cite{s1}.
Hence $\Delta_{\widetilde{\omega}}\bigl(\varphi^{\ast}f\bigr)\leq B$
holds weakly on $V_P$.
Then we can apply the arguments of Proposition 2.1 in \cite{s1},
and we obtain an appropriate estimate for the sup norm of $f$.
By a similar argument,
we obtain such an estimate around any smooth points of $D$.
Thus we are done.
\hfill\qed

\subsection{A construction of an ordinary metric of the bundle}
\label{subsection;05.8.29.450}

Let $(\prolongg{\vecc}{E}_{\ast},\theta)$ be
a $\vecc$-parabolic Higgs bundle on $(X,D)$.
In the following,
we shrink the open sets $U_i$
without mentioning, if it is necessary.
We put $D_i^{\circ}:=D_i\setminus \bigcup_{j\neq i}D_j$.

On each $D_i$,
we have the generalized eigen decomposition
with respect to $\Res_i(\theta)$:
\begin{equation}
\label{eq;06.8.16.1}
\prolongg{\vecc}{E}_{|D_i}
=\bigoplus_{\alpha}\lefttop{i}\Gr^{\EE}_{\alpha}
 (\prolongg{\vecc}{E}_{|D_i})
\end{equation}
For each point $P\in D_i\cap D_j$,
we may assume that there is a decomposition
$\prolongg{\vecc}{E}_{|U_P}=
 \bigoplus \lefttop{P}U_{\veca,\vecalpha}$
as in Section \ref{section;05.7.31.3}.
Let $\lefttop{P}\vecv$ be a holomorphic frame
compatible with the decomposition.
We take a $C^{\infty}$-hermitian metric $\hhat_0$
of $\prolongg{\vecc}{E}$
such that 
$\lefttop{P}\vecv$ is an orthonormal frame on $U_P$
and that the decomposition
(\ref{eq;06.8.16.1}) is orthogonal.
Then, we can take a $C^{\infty}$-isomorphism
$\lefttop{i}\Phi:\pi_{i}^{\ast}
\bigl(
 \prolongg{\vecc}{E}_{|D_i}\bigr)
 \simeq \prolongg{\vecc}{E}$ on $U_i\setminus D$
such that 
(i) the restriction of $\lefttop{i}\Phi$ to $D_i$ is the identity,
(ii) the restriction of $\lefttop{i}\Phi$ to $U_P$
 is given by the frames 
 $\lefttop{P}\vecv$ and $\pi_i^{\ast}\bigl(
 \lefttop{P}\vecv_{|U_P\cap D_i}\bigr)$.
(\cite{li2}. See also 
the explanation in Subsection \ref{subsection;06.8.6.2}.)
We also obtain the orthogonal decompositions
$\Gr_{\alpha}^{\EE}\bigl(
 \prolongg{\vecc}{E}_{|D_i}\bigr)
=\bigoplus_{a\in\real} \lefttop{i}\nbigg_{(a,\alpha)}$
with respect to $\hhat_0$
such that 
$\lefttop{i}F_b\Gr^{\EE}_{\alpha}\bigl(
 \prolongg{\vecc}{E}_{|D_i} \bigr)
=\bigoplus_{a\leq b}\nbigg_{(a,\alpha)}$.
They induce the $C^{\infty}$-decompositions
$\prolongg{\vecc}{E}_{|U_i}
=\bigoplus \leftupdown{i}{\vecc}{E}_{(a,\alpha)}$.

\vspace{.1in}

We can take a hermitian metric $h_0$ of $E$ on $X-D$,
which is as in Subsection \ref{subsection;06.8.6.1} on $U_P$,
and as in Subsection \ref{subsection;06.8.6.2} on 
$U_i\setminus\bigcup U_P$.
More precisely,
we take a hermitian metric $h_{D_i}$
of $\prolongg{\vecc}{E}_{|D_i^{\circ}}$
such that
(i) the decomposition
$\prolongg{\vecc}{E}_{|D_i^{\circ}}
=\bigoplus \lefttop{i}\nbigg_{u|D_i^{\circ}}$
 is orthogonal,
(ii) $h_{D_i}(\lefttop{P}v_k,\lefttop{P}v_l)=
 \delta_{k,l}\cdot |z_j|^{-2a_j(\lefttop{P}v_k)}$
for each $P\in D_i\cap D_j$ $(j\neq i)$.
Let $h_{D_i,u}$ denote the restriction
of $h_{D_i}$ to $\lefttop{i}\nbigg_{u|D_i^{\circ}}$.
Then, $h_0$ is given by (\ref{eq;06.8.16.2})
on $U_i\setminus D$.
We have 
$h_0(\lefttop{P}v_k,\lefttop{P}v_l)=
 \delta_{k,l}\cdot|z_i|^{-2a_i(\lefttop{P}v_k)}
 \cdot |z_j|^{-2a_j(\lefttop{P}v_k)}$
on $U_P\setminus D$ for $P\in D_i\cap D_j$.
Thus, we obtain the metric of $E$ on
$\bigcup_i U_i\setminus D$.
We extend it to the metric of $E$ on $X-D$.
Such a metric $h_0$ is called an ordinary metric,
in this paper.
The following lemma immediately
follows from Proposition \ref{prop;06.8.15.100}
and Proposition \ref{prop;05.8.30.50}.
\begin{lem}
If $(\prolongg{\vecc}{E}_{\ast},\theta)$ is graded semisimple,
then 
$F(h_0)$ is bounded with respect to $h_0$ and $\omega_{\epsilon}$.
\hfill\qed
\end{lem}

\subsection{Calculation of the integrals}

\begin{lem}
 \label{lem;05.9.8.60}
\[
\left(\frac{\sqrt{-1}}{2\pi}\right)^2
 \int_{X-D}\bigl(\tr R(h_0)\bigr)^2
=\int_X\parchern_1^2(\prolongg{\vecc}{E}_{\ast}).
\]
\end{lem}
\pf
We have 
$ \bigl(\tr R(h_{0})\bigr)^2
=\bigl(\tr R(\hhat_0)\bigr)^2
+\tr R(\hhat_0)\cdot \delbar\tr A
+\tr R(h_0)\cdot \delbar\tr A$.
We have the following equality:
\[
 \left(\frac{\sqrt{-1}}{2\pi}\right)^2
\int_{X-D}\bigl(\tr R(\hhat_0)\bigr)^2
=\int_X c_1(\prolongg{\vecc}{E})^2.
\]
Due to (\ref{eq;06.8.16.20}),
we obtain the following:
\begin{multline}
 \left(\frac{\sqrt{-1}}{2\pi}\right)^2
\int_{X-D}\tr R(\hhat_0)\cdot \delbar\tr A
=\sum_i\frac{\sqrt{-1}}{2\pi}
\int_{D_i}\tr R(\hhat_{0|D_i},\prolongg{\vecc}{E}_{|D_i})\cdot 
  (-\wt(\prolongg{\vecc}{E}_{\ast},i)) \\
=\sum_i -\wt(\prolongg{\vecc}{E}_{\ast},i)
\cdot \deg_{D_i}(\prolongg{\vecc}{E}_{|D_i})
=-\sum_i\wt(\prolongg{\vecc}{E}_{\ast},i)
 \int_Xc_1(\prolongg{\vecc}{E})\cdot [D_i]
\end{multline}
We put $\Ehat_{D_i,u}:=\lefttop{i}\Gr^{F,\EE}_u(E_{|D_i})$,
which is naturally isomorphic to $\lefttop{i}\nbigg_u$
as $C^{\infty}$-bundles.
Hence the metric $h_{D_i,u}$ on $\Ehat_{D_i,u}$
is induced (Subsection \ref{subsection;05.8.29.450}).
Then, we obtain the following,
using Corollary  \ref{cor;06.8.16.41}:
\begin{multline}
  \left(\frac{\sqrt{-1}}{2\pi}\right)^2
\int_{X-D}\tr R(h_0)\cdot \delbar\tr A
=-\sum_i
 \wt(\prolongg{\vecc}{E}_{\ast},i)
 \sum_{u}
 \frac{\sqrt{-1}}{2\pi}
 \int_{D_i}
  \tr R(h_{D_i,u},\Ehat_{D_i,u}) \\
+\sum_i
 \wt(\prolongg{\vecc}{E}_{\ast},i)^2\cdot
 \frac{\sqrt{-1}}{2\pi}
 \int_{D_i}\delbar\del\log|\sigma_i|^{2}
\end{multline}
We have the naturally induced parabolic
structure of $\prolongg{\vecc}{E}_{|D_i}$
at $D_i\cap \bigcup_{j\neq i}D_j$.
Then we have the following equality:
\begin{multline}
\sum_{u}
 \frac{\sqrt{-1}}{2\pi}
 \int_{D_i}
 \tr R(h_{D_i,u},\Ehat_{D_i,u})
=\pardeg_{D_i}(\prolongg{\vecc}{E}_{|D_i\,\ast}) \\
=\deg_{D_i}(\prolongg{\vecc}{E}_{|D_i})
-\sum_{j\neq i}\wt(\prolongg{\vecc}{E}_{\ast},j)
\cdot \int_X[D_i]\cdot [D_j].
\end{multline}
We also have
$\frac{\sqrt{-1}}{2\pi}\int_{D_i}\delbar\del\log|\sigma_i|^2 
=\int_X [D_i]^2$.
Thus we obtain the following:
\begin{multline}
  \left(\frac{\sqrt{-1}}{2\pi}\right)^2
\int_{X-D}\tr R(h_0)\cdot \delbar\tr A
=-\sum_i\wt(\prolongg{\vecc}{E}_{\ast},i)
 \int_Xc_1(\prolongg{\vecc}{E})\cdot [D_i] \\
+\sum_i\sum_{j\neq i}
 \wt(\prolongg{\vecc}{E}_{\ast},i)
\cdot \wt(\prolongg{\vecc}{E}_{\ast},j)
 \int_X[D_i]\cdot [D_j]  
+\sum_i 
\wt(\prolongg{\vecc}{E}_{\ast},i)^2
\cdot \int_X[D_i]^2 \\
=-\sum_i\wt(\prolongg{\vecc}{E}_{\ast},i)
 \int_Xc_1(\prolongg{\vecc}{E})\cdot [D_i] 
+\sum_i\sum_{j}
 \wt(\prolongg{\vecc}{E}_{\ast},i)
\cdot \wt(\prolongg{\vecc}{E}_{\ast},j)
 \int_X[D_i]\cdot [D_j].
\end{multline}
Then the claim of the lemma follows.
\hfill\qed

\begin{cor}
 \label{cor;05.10.6.2}
\[
 \left(\frac{\sqrt{-1}}{2\pi}
 \right)^2\int_{X-D}
 \bigl(\tr F(h_{0})\bigr)^2
=\int_X\parchern_1^2(\prolongg{\vecc}{E}_{\ast}).
\]
\end{cor}
\pf
It follows from 
$\bigl(\tr F(h_0)\bigr)^2=\bigl(\tr R(h_0)\bigr)^2$
and the previous lemma.
\hfill\qed

\begin{prop}
 \label{prop;05.8.1.11}
If $(\vecE_{\ast},\theta)$ is graded semisimple,
then the following equality holds:
\[
 \left(
 \frac{\sqrt{-1}}{2\pi}
 \right)^2
\int_{X-D}\tr\Bigl(F(h_{0})^2\Bigr)
=2\int_X\parch_{2}(\prolongg{\vecc}{E}_{\ast}).
\]
\end{prop}
\pf
We have only to show the following two equalities:
\begin{equation} \label{eq;05.8.1.21}
 \int_{X-D}\tr\Bigl(F(h_{0})^2\Bigr)
=\int_{X-D}\tr\Bigl(R(h_0)^2\Bigr).
\end{equation}
\begin{equation} \label{eq;05.8.1.22}
\left(\frac{\sqrt{-1}}{2\pi}\right)^2\int_{X-D}\tr\Bigl(R(h_0)^2\Bigr)
=2\int_{X}\parch_{2}(\prolongg{\vecc}{E}_{\ast}).
\end{equation}

Let us show (\ref{eq;05.8.1.21}).
By a direct calculation or the classical Chern-Simons theory,
we obtain the following equality:
\begin{multline}
 \tr\bigl(F(h_0)^2\bigr)
=\tr\bigl( R(h_0)^2 \bigr)
+2\tr\Bigl(
 \bigl(\del_{h_0}\theta+\delbar\theta^{\dagger}_{h_0}\bigr)\cdot R(h_0)
 \Bigr) \\
+d\Bigl(
 \tr\bigl((\theta+\theta^{\dagger}_{h_0})\cdot
 (\del_{h_0}\theta+\delbar\theta^{\dagger}_{h_0})\bigr) 
+(2/3)\cdot\tr\bigl((\theta+\theta^{\dagger}_{h_0})^3\bigr)
 \Bigr).
\end{multline}
Since
$R(h_0)$, $\del_{h_0}\theta$
and $\delbar\theta^{\dagger}_{h_0}$ are
a $(1,1)$-form, a $(2,0)$-form and a $(0,2)$-form
respectively,
we obtain the vanishing of the second term
in the right hand side.
It is easy to obtain
$\tr\bigl((\theta+\theta^{\dagger}_{h_0})^3\bigr)=0$
from $\theta^2=\theta_{h_0}^{\dagger\,2}=0$.

We put
$Y_i(\delta):=
 \Bigl\{x\in X\,\Big|\,
 |\sigma_i(x)|=\min_j |\sigma_j(x)|=\delta\Bigr\}$
and $Y(\delta):=\bigcup_{i}Y_{i}(\delta)$.
From the estimate in Sections 
\ref{section;05.7.31.3}--\ref{section;05.7.31.5},
$\tr\bigl(\theta\cdot\delbar\theta^{\dagger}_{h_0}\bigr)$
and $\tr\bigl(\theta_{h_0}^{\dagger}\cdot\del_{h_0}\theta\bigr)$
are bounded with respect to $\omega_{\epsilon'}$
for some $0<\epsilon'<\epsilon$.
Hence, we obtain the following:
\[
 \lim_{\delta\to 0}
   \int_{Y(\delta)}\tr(\theta\cdot\delbar\theta^{\dagger}_{h_0})=
\lim_{\delta\to 0}
 \int_{Y(\delta)}
\tr(\theta^{\dagger}_{h_0}\cdot\del_{h_0}\theta)
=0
\]
Then, we obtain the formula (\ref{eq;05.8.1.21}).

\vspace{.1in}

Finally, let us see (\ref{eq;05.8.1.22}).
We put $A:=\del_{h_0}-\del_{\hhat_0}$.
Then we have the following:
\[
 \tr\bigl(R(h_0)^2\bigr)
=\tr\bigl(R(\hhat_0)^2\bigr)
+d\tr\bigl(A\cdot R(h_0)+A\cdot R(\hhat_0)\bigr).
\]
The contribution of the first term is as follows:
\[
 \left(\frac{\sqrt{-1}}{2\pi}\right)^2
\int_{X-D}\tr\bigl(R(\hhat_0)^2\bigr)
=2\ch_2(\prolongg{\vecc}{E}).
\]
From Corollary \ref{cor;06.8.16.40},
we obtain the following:
\begin{multline}
 \left(\frac{\sqrt{-1}}{2\pi}\right)^2
 \int_{X-D}d\tr\Bigl(A\cdot R(h_0)+A\cdot R(\hhat_0)\Bigr)=
-\sum_{i,a,\alpha}a\cdot\deg_{D_i}
 \bigl(\lefttop{i}\Gr^{F,\EE}_{a,\alpha}(\prolongg{\vecc}{E}_{|D_i})\bigr) \\
-\sum_{i,a,\alpha}a\cdot\pardeg_{D_i}
\bigl(
 \lefttop{i}\Gr^{F,\EE}_{a,\alpha}
 (\prolongg{\vecc}{E}_{|D_i})_{\ast}\bigr)
+\sum_{i,a,\alpha}a^2
 \rank\lefttop{i}\Gr^{F,\EE}_{a,\alpha}\bigl(\prolongg{\vecc}{E}_{|D_i}\bigr)
\int_X[D_i]^2.
\end{multline}
Here $\lefttop{i}\Gr^{F,\EE}_{a,\alpha}\bigl(
 \prolongg{\vecc}{E}_{|D_i}\bigr)_{\ast} $
is the parabolic bundle on
$\bigl(D_i,D_i\cap\bigcup_{j\neq i}D_j\bigr)$
with the canonically induced parabolic structure.
We have the following:
\begin{multline}
\sum_{\alpha}
 \pardeg_{D_i}\Bigl(
 \lefttop{i}\Gr^{F,\EE}_{a,\alpha}\bigl(
 \prolongg{\vecc}{E}_{|D_i}\bigr)_{\ast}
 \Bigr)
=\pardeg_{D_i}\Bigl(
 \lefttop{i}\Gr^{F}_{a}\bigl(
 \prolongg{\vecc}{E}_{|D_i} \bigr)_{\ast}
 \Bigr) \\
=\deg_{D_i}\bigl(\lefttop{i}\Gr^F_a\bigl(\prolongg{\vecc}{E}_{|D_i}\bigr)
 \bigr)
-\sum_{\substack{j\neq i,\\ P\in D_i\cap D_j}}
 \sum_{\substack{\veca\in\Par(\prolongg{\vecc}{E},P)\\ a_i=a}}
 a_j\cdot\rank\left(
 \lefttop{P}\Gr^F_{\veca}(\prolongg{\vecc}{E}_{|O})
 \right).
\end{multline}
Then (\ref{eq;05.8.1.22}) immediately follows.
\hfill\qed

\subsection{The degree of subsheaves}

Let $V$ be a saturated coherent $\nbigo_{X-D}$-submodule of $E$.
Let $\pi_V$ denote the orthogonal projection of $E$
onto $V$ with respect to $h_{0}$,
which is defined outside the Zariski closed subset
of codimension two.
It is also the orthogonal projection with respect to $h_{0}$.
Let $h_{V}$ be the metric of $V$
induced by $h_0$.
The following lemmas are the special case
of the results of J. Li \cite{li2}.
\begin{lem}
\label{lem;05.8.30.200}
$\delbar\pi_V$ is $L^2$ with respect to $h_0$ and $\omega_{\epsilon}$
if and only if there exists a coherent subsheaf
$\prolongg{\vecc}{V}\subset \prolongg{\vecc}{E}$
such that $\prolongg{\vecc}{V}_{|X-D}=V$.
\hfill\qed
\end{lem}

\begin{lem}
 \label{lem;05.9.8.2}
$\deg_{\omega_{\epsilon}}(V,h_V)
=\pardeg_{\omega}(\prolongg{\vecc}{V}_{\ast})$
holds.
\end{lem}
\pf
We give just an outline of a proof.
By considering the exterior product of $E$ and $V$,
we may assume that $\rank V=1$.
We may assume that $L$ is very ample.
Let $C$ be a smooth divisor of $X$
with $\nbigo(C)\simeq L$
such that $\prolongg{\vecc}{V}$  is locally free on a neighbourhood of $C$,
that $C$ intersects with the smooth part of $D$
transversally,
and that $\lefttop{i}F$ is a filtration in the category
of the vector bundles on $D_i$
around $C\cap D_i$.
We can take a smooth $(1,1)$-form $\tau$
whose support is contained in a sufficiently small
neighbourhood of $C$,
such that $\tau$ and $\omega$ represents
the same cohomology class.
Then we have
$\int\tr R(h_V)\cdot\omega=\int \tr R(h_V)\cdot \tau$.
It can be checked 
$\frac{\sqrt{-1}}{2\pi}\int \tr R(h_V)\cdot\tau
=\pardeg_{\omega}(\prolongg{\vecc}{V}_{\ast})$
by an elementary argument.
\hfill\qed

%% file: 5.tex
In this chapter,
we show the fundamental property of
the parabolic Higgs bundles associated to
tame harmonic bundles,
such as $\mu_L$-polystability and
the vanishing of characteristic numbers.
We also see the uniqueness of the adapted pluri-harmonic metric
for parabolic Higgs bundles.
These results give the half of Theorem \ref{thm;05.9.13.1}.

\section{Polystability and Uniqueness}

Let $X$ be a smooth irreducible projective variety,
and $D$ be a simple normal crossing divisor
with the irreducible decomposition
$D=\bigcup_{i\in S}D_i$.
Let $L$ be any ample line bundle of $X$.

\begin{prop}
 \label{prop;05.9.8.30}
\mbox{{}}
Let $\harmonicbundle$ be a tame harmonic bundle on $X-D$,
and let $(\prolongg{\vecc}{E}_{\ast},\theta)$
denote the associated $\vecc$-parabolic Higgs bundle
for any $\vecc\in\real^S$.
(See Section {\rm\ref{section;05.9.8.110}}.)
\begin{itemize}
\item
$(\prolongg{\vecc}{E}_{\ast},\theta)$ is $\mu_L$-polystable,
and $\pardeg_{L}(\prolongg{\vecc}{E}_{\ast})=0$.
\item
Let $(\prolongg{\vecc}{E}_{\ast},\theta)=
 \bigoplus_i (\prolongg{\vecc}{E_i}_{\ast},\theta_i)\otimes \cnum^{p(i)}$
be the canonical decomposition
(Corollary {\rm\ref{cor;05.9.14.5}}).
Then we have the orthogonal decomposition $h=\bigoplus_i h_i\otimes g_i$.
Here $h_i$ are pluri-harmonic metrics
for $(E_i,\delbar_{E_i},\theta_i)$,
and $g_i$ are hermitian metrics of $\cnum^{p(i)}$.
\end{itemize}
\end{prop}
\pf
The equality  $\pardeg_{L}(\prolongg{\vecc}{E}_{\ast})=0$
can be easily reduced to the curve case
(Proposition \ref{prop;05.9.160}).
It also follows from the curve case
that $(\prolongg{\vecc}{E}_{\ast},\theta)$ is $\mu_L$-semistable.

Let us show $(\prolongg{\vecc}{E}_{\ast},\theta)$
is $\mu_L$-polystable.
Let $(\prolongg{\vecc}{V}_{\ast},\theta_V)$ be a non-trivial 
saturated Higgs subsheaf of
$(\prolongg{\vecc}{E}_{\ast},\theta)$ such that 
$\mu_L(\prolongg{\vecc}{V}_{\ast})=
 \mu_L(\prolongg{\vecc}{E}_{\ast})=0$
and $\rank(V)<\rank(E)$.
Recall that we have the closed subset $Z\subset X$
such that $\prolongg{\vecc}{V}_{|X-Z}$ is the subbundle
of $\prolongg{\vecc}{E}_{|X-Z}$.
The codimension of $Z$ is larger than $2$.
We have the orthogonal projection $\pi_V:E\lrarr V$
on the open set $X-(Z\cup D)$.
Let $C\subset X$ be a smooth curve 
such that (i) $C$ intersects with the smooth part of $D$
transversally, (ii) $C\cap Z=\emptyset$.
Let $\theta_C$ denote the induced Higgs field
of $E_{|C\setminus D}$.
Due to the result in the curve case,
we obtain that $\pi_{V|C}$ is holomorphic
and that $\theta_C$ and $\pi_{V|C}$ commute.
Then, we obtain that $\pi_{V|X-(D\cup Z)}$
is holomorphic
and that $[\pi_V,\theta]=0$.
Since the codimension of $Z$ is larger than two,
$\pi_V$ naturally gives the holomorphic map
$E\lrarr E$ on $X-D$,
which is also denoted by $\pi_V$.
It is easy to see $\pi_V^2=\pi_V$,
and that the restriction of $\pi_V$ to $V$ is the identity.
Hence we obtain the decomposition
$E=V\oplus V'$, where we put $V'=\ker\pi_V$.
We can conclude that  $V$ and $V'$ are vector subbundles of $E$,
and the decomposition is orthogonal with respect to the metric $h$.
Since we have $[\pi_V,\theta]=0$,
the decomposition is also compatible with the Higgs field.
Hence we obtain the decomposition 
of $(E,\delbar_E,\theta,h)$ into
$(V,\delbar_V,\theta_V,h_V)\oplus
 (V',\delbar_{V'},\theta_{V'},h_{V'})$
as harmonic bundles.
Then it is easy that
$(\prolongg{\vecc}{E}_{\ast},\theta)$ is also 
decomposed into
$(\prolongg{\vecc}{V}_{\ast},\theta_V)
\oplus
 (\prolongg{\vecc}{V'}_{\ast},\theta_{V'})$.
Since both of
$(\prolongg{\vecc}{V}_{\ast},\theta_V)$ and 
$ (\prolongg{\vecc}{V'}_{\ast},\theta_{V'})$ 
are obtained from tame harmonic bundles,
they are $\mu_L$-semistable.
And we have $\rank(V)<\rank (E)$ and $\rank(V')<\rank (E)$.
Hence the $\mu_L$-polystability of $(\prolongg{\vecc}{E},\theta)$
can be shown by an easy induction on the rank.

From the argument above,
the second claim is also clear.
\hfill\qed

\begin{prop}
 \label{prop;05.9.8.300}
Let $(\prolongg{\vecc}{E}_{\ast},\theta)$ be
a $\vecc$-parabolic Higgs bundle on $(X,D)$.
We put $E:=\prolongg{\vecc}{E}_{|X-D}$.
Assume that we have pluri-harmonic metrics $h_i$ of $(E,\delbar_E,\theta)$
$(i=1,2)$, which are adapted to the parabolic structures.
Then we have the decomposition of Higgs bundles
$(E,\theta)=\bigoplus_a (E_a,\theta_a)$ satisfying the following conditions:
\begin{itemize}
\item
 The decomposition is orthogonal with respect to
 both of $h_i$.
 The restrictions of $h_i$ to $E_a$ are denoted by $h_{i,a}$.
\item
 There exist positive numbers $b_a$ such that
 $h_{1,a}=b_a\cdot h_{2,a}$.
\end{itemize}
We remark that the decomposition $(E,\theta)=\bigoplus (E_a,\theta_a)$ 
induces the decomposition of the $\vecc$-parabolic Higgs bundles
$ \bigl(\prolongg{\vecc}{E}_{\ast},\theta\bigr)
=\bigoplus (\prolongg{\vecc}{E_a}_{\ast},\theta_a)$.
\end{prop}
\pf
Recall the norm estimate for tame harmonic bundles (\cite{mochi2})
which says that the harmonic metrics are determined
up to boundedness by the parabolic filtration and the weight filtration.
Hence we obtain the mutually boundedness of $h_1$ and $h_2$.
Then the uniqueness follows from Proposition \ref{prop;05.9.8.100}.
(The Kahler metric of $X-D$ is given by the restriction of a Kahler
metric of  $X$. It satisfies Condition \ref{condition;05.7.30.1},
according to \cite{s1}.)
\hfill\qed

\section{Vanishing of Characteristic Numbers}
\label{section;06.8.17.100}

\begin{prop}
 \label{prop;05.8.29.500}
Let $\harmonicbundle$ be a tame harmonic bundle on $X-D$,
and $(\prolongg{\vecc}{E}_{\ast},\theta)$ be the induced $\vecc$-parabolic
Higgs bundle.
Then we have the vanishing of
the following characteristic numbers:
\[
\int_X \parch_{2,L}(\prolongg{\vecc}{E}_{\ast})=0,
\quad\quad
\int_X\parchern_{1,L}^2(\prolongg{\vecc}{E}_{\ast})=0.
\]
\end{prop}
\pf
We may and will assume $\dim X=2$.
Let $h_0$ be an ordinary metric for 
the parabolic Higgs bundle $(\prolongg{\vecc}{E}_{\ast},\theta)$.
We have only to show
$\int \tr \bigl(R(h_0)^2\bigr)
=\int \tr\bigl(R(h_0)\bigr)^2=0$.
Let $I(D)$ denote the set of the intersection points of $D$.
Let $\pi:\Xtilde\lrarr X$ be a blow up at $I(D)$.
We put $\Dtilde:=\pi^{-1}(D)$.
Let $\Dtilde_i$ denote the proper transform of $D_i$,
and let $\Dtilde_P$ denote the exceptional curve $\pi^{-1}(P)$.
We put $\widetilde{S}:=S\cup I(D)$.
Then, we have 
$\Dtilde=\bigcup_{i\in \Stilde}\Dtilde_i$.
We take neighbourhoods
$\Utilde_i$ of $\Dtilde_i$
with retractions $\pitilde_i:\Utilde_i\lrarr \Dtilde_i$
for $i\in \Stilde$,
as in Subsection \ref{section;05.8.23.5}.

We put $\Etilde:=\pi^{-1}(\prolongg{\vecc}{E})$.
On $\Etilde_{|\Dtilde_i}$ $(i\in S)$,
we have the naturally induced filtration $\lefttop{i}F$.
For any intersection point $P\in D_i\cap D_j$,
we have the isomorphism
$\Etilde_{|\Dtilde_P}
 \simeq
 \prolongg{\vecc}{E}_{|P}\otimes\nbigo_{\Dtilde_P}$.
We have the filtrations
$\lefttop{i}F$ and $\lefttop{j}F$
on $\prolongg{\vecc}{E}_{|P}$.
We take a decomposition
$ \prolongg{\vecc}{E}_{|P}
 =\bigoplus U_{\veca}$
such that
$\lefttop{1}F_{a_1}\cap \lefttop{2}F_{a_2}
=\bigoplus_{\vecb\leq\veca}U_{\vecb}$.
Then, we put
$\lefttop{P}F_b
 (\prolongg{\vecc}{E}_{|P})
 :=\bigoplus_{b_1+b_2\leq b}U_{\vecb}$,
which gives the filtration of 
 $\prolongg{\vecc}{E}_{|P}$.
The induced filtration on 
 $\Etilde_{|\Dtilde_P}$
is also denoted by $\lefttop{P}F$.
The tuple of filtrations
$\bigl(\lefttop{i}F\,\big|\,i\in\Stilde\bigr)$
is denoted by $\vecF$.

We put $\thetatilde:=\pi^{-1}\theta$.
Then $(\Etilde,\vecF,\thetatilde)$ is a generalized
parabolic Higgs bundle in the sense of
Remark \ref{rem;06.8.15.101}.
The residue $\Res_i\thetatilde$ preserves
the filtration $\lefttop{i}F$.
On each $i\in\Stilde$,
the residue $\Res_i\thetatilde$ induces
the endomorphism of $\lefttop{i}\Gr^F_a(\Etilde)$.
The eigenvalues are constant,
and hence the nilpotent part $\nbign_i$
is well defined.
The conjugacy classes of
$\nbign_{i\,|\,P}$ are independent of
the choice of $P\in \Dtilde_i$
(\cite{mochi2}).
Thus, we obtain the weight filtration 
$\lefttop{i}W$ on $\lefttop{i}\Gr^F_a(\Etilde)$.
We put $\Ftilde_{(a,k)}:=\pi_a^{-1}(\lefttop{i}W_k)$,
where $\pi_a$ denotes the projection
$\lefttop{i}F_a\lrarr \lefttop{i}\Gr^F_a(\Etilde)$.

Let $\Ptilde_i$ denote the intersection point
of $\Dtilde_i$ and $\Dtilde_P$ 
for $i\in S$ and $P\in I(D)$.
Around $\Ptilde_i$,
we have the holomorphic frame 
$\lefttop{P\,i}\vecvtilde$,
as in Subsection \ref{subsection;06.8.15.10}.
Namely, we take a holomorphic frame $\lefttop{P\,i}\vecv$
around $P$ as in Subsection \ref{subsection;06.8.15.11},
($D_i$ plays the role of $D_1$, there)
and we put 
$\lefttop{P\,i}\vecvtilde:=
 \pi^{-1}(\lefttop{P,i}\vecv)$ around $\Ptilde_i$.
We take a hermitian metric $\hhat_1$ of $\Etilde$
such that $\lefttop{P,i}\vecvtilde$ around $\Ptilde_i$
are orthonormal with respect to $\hhat_1$.
By using it, we take $C^{\infty}$-isomorphisms
$\Phitilde_i:\pitilde_i^{\ast}\Etilde_{|\Dtilde_i}
\simeq \Etilde_{|\Utilde_i}$ on $\Utilde_i$ $(i\in \Stilde)$
as in Subsection \ref{subsection;05.8.29.450}.
Then, we can take a hermitian metric $h_1$
which is as in Subsection \ref{subsection;06.8.15.10}
for the frame $\lefttop{P\,i}\vecv$ around $\Ptilde_i$,
and as in Subsection \ref{subsection;06.8.15.20}
around $\Dtilde_i$.

\begin{lem}
We have
$\int \tr\bigl(R(h_0)^2\bigr)=\int \tr\bigl(R(h_1)^2\bigr)$
and $\int\tr\bigl(R(h_0)\bigr)^2=\int\tr\bigl(R(h_1)\bigr)^2$.
\end{lem}
\pf
Let $\omegatilde$ denote a Poincar\'e like metric on $\Xtilde-\Dtilde$.
Let $\htilde_0$ be an ordinary metric for $(\Etilde,\vecF,\thetatilde)$
as constructed in Subsection \ref{subsection;05.8.29.450}.
Then, $\pi^{\ast}h_0$ and $\htilde_0$
are mutually bounded.
Both of $\pi^{\ast}R(h_0)$ and $R(\htilde_0)$
are bounded with respect to $\htilde_0$
and $\omegatilde$.

Let us see that
$A_0=\del_{\Etilde,\pi^{-1}h_0}-\del_{\Etilde,\htilde_0}$
is bounded with respect to $\omegatilde$ and $\htilde_0$.
Let us recall the description of $A_0$ around $\Ptilde_i$.
We take a holomorphic coordinate neighbourhood
$(U,z_1,z_2)$ such that 
$\{z_1\cdot z_2=0\}=U\cap(\Dtilde_i\cup\Dtilde_P)$.
We put $D'_j:=\{z_j=0\}$.
We have two holomorphic decomposition
$\Etilde_{|U}=\bigoplus U_{\veca}=\bigoplus \Utilde_{\veca}$
such that
$\lefttop{j}F_{b}=
 \bigoplus_{a_j\leq b} U_{\veca|D'_j}
=\bigoplus_{a_j\leq b}\Utilde_{\veca|D_j'}$,
where $a_j$ denotes the $j$-th component of $\veca$.
We put 
$\Gamma_j=\bigoplus a_j\cdot \id_{U_{\veca}}$
and $\Gammatilde_j=\bigoplus a_j\cdot \id_{\Utilde_{\veca}}$.
We have
$\del_{\Etilde,\pi^{-1}(h_0)}=
   \nbigd_1-\sum \Gamma_j\cdot dz_j/z_j$
and
$\del_{\Etilde,\htilde_0}=
   \nbigd_2-\sum \Gammatilde_j\cdot dz_j/z_j$,
where $\nbigd_1$ (resp. $\nbigd_2$) is the $(1,0)$-operator
of $\Etilde_{|U}$,
preserving the decomposition
$\Etilde_{|U}=\bigoplus U_{\veca}$
(resp. $\Etilde_U=\bigoplus \Utilde_{\veca}$).
Then, we have
$A_0=\sum (\Gammatilde_j-\Gamma_j)dz_j/z_j
 +(\nbigd_1-\nbigd_2)$.
Because of
$(\Gammatilde_j-\Gamma_j)_{|D_j'}\lefttop{j}F_a
\subset \lefttop{j}F_{<a}$
and $(\nbigd_1-\nbigd_2)_{|D_j'}\lefttop{j}F_a
 \subset\lefttop{j}F_a$,
we obtain the boundedness of $A_0$
with respect to $\pi^{-1}h_0$ and $\omegatilde$,
around $\Ptilde_i$.

Let us recall the description of $A_0$
around $Q\in \Dtilde_i$ $(i\in \Stilde)$.
Let $(U,z_1,z_2)$ be a holomorphic coordinate
around $Q$ such that $z_1^{-1}(0)=U\cap \Dtilde_i$.
We have two $C^{\infty}$-decomposition
$\Etilde_{|U}=\bigoplus E_a=\bigoplus \Etilde_a$
such that 
$\lefttop{i}F_b=\bigoplus_{a\leq b}E_{a|\Dtilde_i}
 =\bigoplus_{a\leq b}\Etilde_{a|\Dtilde_i}$.
We put $\Gamma:=\bigoplus a\cdot \id_{E_a}$
and $\Gammatilde:=\bigoplus a\cdot \id_{\Etilde_a}$.
We have a description
$\del_{\pi^{-1}h_0}=\del_{1,\pi^{-1}h_0'}
 -\Gamma dz_1/z_1+O(1)$,
where $O(1)$ denotes the bounded one form
with respect to $h_0$ and $\omegatilde$,
and $\del_{1,\pi^{-1}h_0'}$ is operator
on $\Etilde_{|U}$ (not on $\Etilde_{U\setminus \Dtilde_i}$)
such that $\del_{1,\pi^{-1}h_0'|\Dtilde_i}$ preserves 
the filtration $\lefttop{i}F$.
(See the proof of Lemma \ref{lem;05.9.15.5}.)
Similarly, we have 
$\del_{\htilde_0}=\del_{1,\htilde_0'}
 -\Gammatilde dz_1/z_1+O(1)$.
For the expression
$\del_{1,\pi^{-1}h_0'}-\del_{1,\htilde_0'}
=B_1\cdot dz_1/z_1+B_2\cdot dz_2$,
we have $B_{1|\Dtilde_i}=0$
and $B_{2|\Dtilde_i}(\lefttop{i}F_{a})
   \subset \lefttop{i}F_{a}$.
We also have
$(\Gammatilde-\Gamma)_{|\Dtilde_i}
 \lefttop{i}F_a\subset \lefttop{i}F_{<a}$.
Thus, we obtain the boundedness of
$\del_{\pi^{-1}h_0}-\del_{\htilde_0}$.
Now,
it is easy to obtain
$\int \tr\bigl(R(h_0)^2\bigr)=\int \tr\bigl(R(\htilde_0)^2\bigr)$
and
$\int\tr\bigl(R(h_0)\bigr)^2=\int\tr\bigl(R(\htilde_0)\bigr)^2$.

Due to the lemmas \ref{lem;06.8.15.46},
\ref{lem;05.9.15.5}
and \ref{lem;06.8.15.45},
$R(\htilde_0)$, $R(h_1)$
and $A_0:=\del_{h_1}-\del_{\htilde_0}$
are bounded with respect to $(\htilde_0,\omegatilde)$.
Hence,
$\tr\bigl(A_0\bigr)$,
$\tr\bigl(R(\htilde_0)\bigr)$,
$\tr\bigl( R(\htilde_0)\cdot A_0\bigr)$,
$\tr\bigl(R(h_1)\bigr)$
and $\tr\bigl(R(h_1)\cdot A_0\bigr)$
are bounded with respect to $\omegatilde$.
Then, it is easy to show 
$\int \tr \bigl(R(\htilde_0)^2\bigr)
=\int \tr\bigl(R(h_1)^2\bigr)$
and 
$\int \tr\bigl(R(\htilde_0)\bigr)^2
=\int \tr\bigl(R(h_1)\bigr)^2$.
\hfill\qed

\vspace{.1in}

Due to the norm estimate (Lemma \ref{lem;06.8.15.2}),
$\htilde:=\pi^{\ast}h$ and $h_1$
are mutually bounded.
We also have that $R(\htilde)$ is bounded
with respect to $\htilde$ and $\omegatilde$.
Let $s$ denote the self-adjoint
endomorphism of $\pi^{-1}(E)$ with respect to 
$\htilde$ and $h_1$,
determined by $\htilde=h_1\cdot s$.
We have
$\del_{\htilde}-\del_{h_1}=
 s^{-1}\del_{h_1}s$ and
$\delbar \bigl(s^{-1}\del_{h_1}s\bigr)
=R(\htilde)-R(h_1)$,
which is bounded with respect to $h_1$ and
$\omegatilde$.

Let us show the following equality
for any test function $\chi$ on $\Xtilde-\Dtilde$:
\begin{equation}
\label{eq;06.8.15.50}
 \int\bigl(s^{-1}\del_{h_1}(\chi\cdot s),\,
 \del_{h_1}(\chi\cdot s)\bigr)_{h_1}\cdot\omegatilde
=\int\bigl(\chi\cdot \delbar(s^{-1}\del_{h_1}s),\,
 \chi\cdot s\bigr)\cdot\omegatilde
+\int \del\chi\cdot\delbar\chi\cdot \tr(s)\cdot\omegatilde.
\end{equation}
We have the following:
\begin{multline}
 \int\bigl(s^{-1}\del_{h_1}(\chi\cdot s),\,
 \del_{h_1}(\chi\cdot s)\bigr)_{h_1}
=\int\Bigl(
 \delbar\Bigl(
 s^{-1}\cdot\del_{h_1}(\chi\cdot s)\Bigr),\,\,
 \chi\cdot s
 \Bigr)_{h_1}  \\
=\int\bigl(\delbar\del\chi,\,\,\,\chi\cdot s\bigr)_{h_1}
+\int\Bigl(
 \chi\cdot\delbar\bigl(s^{-1}\del_{h_1}s_0\bigr),\,\,
 \chi\cdot s
 \Bigr)_{h_1}
+\int\bigl(\delbar\chi\cdot s^{-1}\del_{h_1}s,\,\,
 \chi\cdot s\bigr)_{h_1}.
\end{multline}
Moreover, we have the following:
\begin{multline}
 \bigl(\delbar\del\chi,\chi\cdot s\bigr)_{h_1}
+\bigl(\delbar\chi\wedge s^{-1}\del_{h_1}s,\,\,
 \chi\cdot s\bigr)_{h_1}
=\tr\bigl(\delbar\del\chi\cdot \chi\cdot s\bigr)
+\tr\bigl(\delbar\chi\cdot\del_{h_1}s\cdot \chi\bigr) \\
=-\del\Bigl(\tr\bigl(\delbar\chi\cdot\chi\cdot s\bigr)\Bigr)
-\tr\bigl(\delbar\chi\del\chi\cdot s\bigr).
\end{multline}
Thus we obtain (\ref{eq;06.8.15.50})

\begin{lem}
 \label{lem;05.9.8.50}
$s^{-1}\del_{h_1}s$ is $L^2$
with respect to $\widetilde{\omega}$ and $\htilde$.
\end{lem}
\pf
Let $\rho$ be a non-negative valued function
on $\real$
satisfying $\rho(t)=1$ for $t\leq 1/2$
and $\rho(t)=0$ for $t\geq 2/3$.
Take hermitian metrics $g_i$ of 
the line bundles $\nbigo(\Dtilde_i)$ $(i\in \Stilde)$.
Let $\sigma_i$ denote the canonical section
of $\nbigo(\Dtilde_i)$,
and $|\sigma_i|$ denote the norm function of $\sigma_i$
with respect to  $g_i$.
We may assume $|\sigma_i|<1$.
We put
$\chi_N:=
 \prod_{i\in\Stilde}\rho\bigl(-N^{-1}\log|\sigma_i|^2\bigr)$.
Then, $\del\chi_N$ is bounded with respect to $\omegatilde$,
independently of $N$.
By using (\ref{eq;06.8.15.50}),
we obtain 
$\int \bigl|s^{-1}\del_{h_1}(\chi_Ns)\bigr|^2_{h_1}
 \dvol_{\omegatilde}
 <C$ for some constant $C$,
and thus we obtain the claim of the lemma.
\hfill\qed

\vspace{.1in}

We put $A_1:=s^{-1}\del_{h_1}s$,
which is $L^2$ with respect to $\omegatilde$ and $h_1$.
We have $R(h_1)=R(\htilde)-\delbar A_1$.
Since we have $\tr R(\htilde)=\tr F(\htilde)=0$,
we have
$\tr \bigl(R(h_1)\bigr)^2
=-d\bigl(\tr R(h_1)\cdot \tr A_1\bigr)$.
Since $R(h_1)$ is bounded with respect to $\omegatilde$
and $h_1$,
we obtain that
$\tr R(h_1)\cdot \tr A_1$
is $L^2$ with respect to $\omegatilde$.
We also know that 
$d\bigl(\tr R(h_1)\cdot\tr(A_1)
 \bigr)$ is integrable.
Then we obtain the vanishing,
due to Lemma 5.2 in \cite{s1}:
\[
 \int \bigl(\tr R(h_1)\bigr)^2
=
 \int d\bigl( \tr R(h_1)
 \cdot \tr A\bigr)=0.
\]
(Note that $\widetilde{\omega}$ satisfies
the condition of the lemma.)
Thus, we obtain
$\int \parchern_1(\prolongg{\vecc}{E}_{\ast})^2=
 (\frac{\sqrt{-1}}{2\pi})^{2}\int \bigl(\tr R(h_0)\bigr)^2
=(\frac{\sqrt{-1}}{2\pi})^2
\int \bigl(\tr R(h_1)\bigr)^2=0$.

Because of $R(\htilde)=-[\theta,\theta_{\htilde}^{\dagger}]$
and $\theta^2=0$,
we easily obtain $\tr\bigl(R(\htilde)^2\bigr)=0$.
Thus we obtain the following:
\[
 \tr\bigl(R(h_1)^2\bigr)
+\delbar
 \Bigl(
 \tr\bigl(A_1\cdot R(h_1)\bigr)
+\tr\bigl(A_1\cdot R(\htilde)\bigr)
 \Bigr)=0
\]
From the boundedness of $R(h_1)$ and $R(\htilde)$
with respect to $\omegatilde$ and $h_1$,
we obtain that $\tr(A_1\cdot R(h_1))$
and $\tr(A_1\cdot R(\htilde))$ are $L^2$
with respect to $\omegatilde$.
Thus we obtain the vanishing,
by using Lemma 5.2 in \cite{s1} again:
\[
 \int \delbar\Bigl(
 \tr\bigl(A_1\cdot R(h_1)\bigr)+\tr\bigl(A_1\cdot R(\htilde)\bigr)
\Bigr)=0.
\]
Thus, we obtain
$\int_X \parch_{2}(\prolongg{\vecc}{E}_{\ast})=
 (\frac{\sqrt{-1}}{2\pi})^{2}\int \tr\bigl(R(h_0)^2\bigr)
=(\frac{\sqrt{-1}}{2\pi})^{2}\int \tr\bigl(R(h_1)^2\bigr)=0$.
\hfill\qed

%% file: 6.tex
In this chapter,
we show the existence of the adapted pluri-harmonic metric
for {\em graded semisimple} 
$\mu_L$-stable parabolic Higgs bundles
on a surface (Proposition \ref{prop;05.7.30.10}).
We will use it together with the perturbation of the parabolic structure
(Section \ref{section;05.7.30.15})
to derive more interesting results.
One of the immediate consequences 
is Bogomolov-Gieseker inequality (Theorem \ref{thm;05.7.30.30}).

\section{Graded Semisimple Parabolic Higgs Bundles on Surface}


We show an existence of Hermitian-Einstein metric
for $\mu_L$-stable parabolic Higgs bundle on a surface
under the graded semisimplicity assumption,
which makes the problem much easier.
Later, we will discuss such existence theorem
for parabolic Higgs bundle with trivial characteristic numbers
in the case where 
the graded semisimplicity is not assumed.

\begin{prop} 
 \label{prop;05.7.30.10}
Let $X$ be a smooth irreducible projective surface 
with an ample line bundle $L$,
and $D$ be a simple normal crossing divisor.
Let $\omega$ be a Kahler form of $X$,
which represents $c_1(L)$.
Let $(\prolongg{\vecc}{E}_{\ast},\theta)$ be a
$\vecc$-parabolic Higgs bundle on $(X,D)$,
which is $\mu_{L}$-stable and
{\em graded semisimple}.
Let us take a positive number $\epsilon$ satisfying the following:
\begin{itemize}
\item
$10\epsilon<\gap(\prolongg{\vecc}{E}_{\ast})$,
and $\epsilon=m^{-1}$ for some positive integer $m$.
\end{itemize}
We take a Kahler form $\omega_{\epsilon}$ of $X-D$,
as in Subsection {\rm\ref{section;05.8.23.5}}.
We put $E=\prolongg{\vecc}{E}_{|X-D}$,
and the restriction of $\theta$ to $X-D$ is denoted 
by the same notation.
Then there exists a hermitian metric $h$ of $E$
satisfying the following conditions:
\begin{itemize}
\item
 Hermitian-Einstein condition 
 $\Lambda_{\omega_{\epsilon}} F(h)=a\cdot \id_{E}$
 for some constant $a$ determined by the following equation:
\begin{equation}
 \label{eq;05.9.17.1}
 a\cdot\frac{\sqrt{-1}}{2\pi}
 \frac{\rank E}{2}
 \int_{X-D}\omega_{\epsilon}^2
=a\cdot\frac{\sqrt{-1}}{2\pi}
 \frac{\rank(E)}{2}
 \int_X\omega^2
=\pardeg_{\omega}(\prolongg{\vecc}{E}_{\ast}).
\end{equation}
\item
 $h$ is adapted to the parabolic structure of
 $\prolongg{\vecc}{E}_{\ast}$.
\item
 $\deg_{\omega_{\epsilon}}(E,h)
 =\pardeg_{\omega}(\prolongg{\vecc}{E}_{\ast})$.
\item
 We have the following equalities:
\[
 \int_X 2\parch_{2}(\prolongg{\vecc}{E}_{\ast})
=\left(\frac{\sqrt{-1}}{2\pi}\right)^2\int_{X-D}
 \tr\Bigl(F(h)^2\Bigr),
\]
\[
  \int_X \parchern_1^2(\prolongg{\vecc}{E}_{\ast})
=\left(\frac{\sqrt{-1}}{2\pi}\right)^2\int_{X-D}
 \tr\bigl(F(h)\bigr)^2.
\]
\end{itemize}
\end{prop}
\pf
Let us take an ordinary metric $h_0$
for the parabolic bundle $(\prolongg{\vecc}{E}_{\ast},\theta)$
as in Section \ref{section;05.8.30.100}.
Note we have
$\Lambda_{\omega_{\epsilon}}\tr R(h_0)
=\Lambda_{\omega_{\epsilon}}\tr F(h_0)$.
We put $\gamma_i:=\wt(\prolongg{\vecc}{E}_{\ast},i)$.

Let us see the induced metric $\det(h_0)$ of $\det(E)$.
Due to our construction,
$\det(h_0)$ is of the form 
$\tau\cdot |z_i|^{-2 \gamma_i}\cdot |z_j|^{-2\gamma_j}$
around $P\in D_i\cap D_j$,
where $\tau$ denotes a positive $C^{\infty}$-metric
of $\det\bigl(\prolongg{\vecc}{E}\bigr)_{|U_P}$.
If $P$ is a smooth point of $D_i$.
then the metric $\det(h_0)$ is of the form
$\tau\cdot |\sigma_i|_{g_i}^{-2\gamma_i}$,
where $\tau$ and $\gamma_i$ are as above.
Therefore,
$\tr R(h_0)=R\bigl(\det(h_0)\bigr)$
is $C^{\infty}$ on $X$.
If $a$ is determined by (\ref{eq;05.9.17.1}),
we have 
$\int_{X-D}\bigl(
 \tr \Lambda_{\omega_{\epsilon}}F(h_0)
-\rank(E)\cdot a\bigr)
\cdot \omega_{\epsilon}^2=0$.
Recall $\epsilon=m^{-1}$ for some positive integer $m$.
Then the following lemma can be shown 
by a consideration of orbifolds.
\begin{lem}
We can take a bounded $C^{\infty}$-function $g$ on $X-D$
satisfying the conditions
{\rm(i)}
$\Delta_{\omega_{\epsilon}}g
=\sqrt{-1}\Lambda_{\omega_{\epsilon}}\tr\bigl(F(h_0)\bigr)
   -\sqrt{-1}\rank(E)\cdot a$,
where $a$ is determined by the equation
{\rm (\ref{eq;05.9.17.1})},
{\rm(ii)}
$\del g$, $\delbar g$, and $\del\delbar g$ are
bounded with respect to $\omega_{\epsilon}$.
\hfill\qed
\end{lem}

We put $g':=g/\rank E$ and
$h_{in}:=h_0\cdot \exp(-g')$.
We remark that the adjoints $\theta$
for $h_0$ and $h_{in}$ are same.
We also remark that
$\del_{h_{in}}-\del_{h_0}$ and
$R(h_{in})-R(h_0)$ are
just multiplications 
$-\del g'\cdot \id_{E}$ and $\del\delbar g'\cdot \id_{E}$
respectively,
which are bounded with respect to $\omega_{\epsilon}$.

\begin{lem}
\label{lem;06.8.17.1}
The metric $h_{in}$ satisfies the following conditions:
\begin{itemize}
\item
 $h_{in}$ is adapted to the parabolic structure of
 $\prolongg{\vecc}{E}_{\ast}$.
\item
 $F(h_{in})$ is bounded with respect to $h_{in}$ and $\omega_{\epsilon}$.
\item
 Let $V$ be any saturated coherent subsheaf of $E$,
 and let $\pi_V$ denote the orthogonal projection of $E$
 onto $V$.
  Then $\delbar\pi_V$ is $L^2$ with respect to
 $h_{in}$ and $\omega_{\epsilon}$,
 if and only if 
 there exists a saturated coherent subsheaf $\prolongg{\vecc}{V}$
 of $\prolongg{\vecc}{E}$ such that 
 $\prolongg{\vecc}{V}_{|X-D}=V$.
 Moreover we have
 $\pardeg_{\omega}(\prolongg{\vecc}{V}_{\ast})
 =\deg_{\omega_{\epsilon}}(V,h_{in,V})$,
 where $h_{in,V}$ denotes the metric of $V$ induced by $h_{in}$.
\item
 $\tr\Lambda_{\omega_{\epsilon}}F(h_{in})\
  =\rank(E)\cdot a$ 
 for the constant $a$
 determined by the equation
 {\rm(\ref{eq;05.9.17.1})}.
\item The following equalities hold:
\[
 \left(\frac{\sqrt{-1}}{2\pi}\right)^2
 \int_{X-D}\tr \Bigl(F(h_{in})^2\Bigr)
=\int_X2\parch_{2}(\prolongg{\vecc}{E}_{\ast}),
\]
\[
 \left(\frac{\sqrt{-1}}{2\pi}\right)^2
 \int_{X-D}\tr\Bigl(F(h_{in})\Bigr)^2
=\int_X\parchern_{1}^2(\prolongg{\vecc}{E}_{\ast}).
\]
\end{itemize}
Due to the third condition,
 $(E,h_{in},\theta)$ is analytic stable with respect to $\omega_{\epsilon}$,
 if and only if $(\prolongg{\vecc}{E}_{\ast},\theta)$ is $\mu_L$-stable.
\end{lem}
\pf
Since $g'$ is bounded and since $h_0$ is adapted to
the parabolic structure,
$h_{in}$ is also adapted to the parabolic structure.
We have $F(h_{in})=F(h_0)+\del\delbar g'\cdot \id_{E}$.
Hence the boundedness of $F(h_{in})$ with respect to 
$\omega_{\epsilon}$ and $h_0$ follows from 
those of $F(h_0)$ and $\del\delbar g'$.

For any saturated subsheaf $V\subset E$,
the orthogonal decomposition $\pi^{h_0}_V$
and $\pi^{h_{in}}_V$ are same.
Hence $\delbar \pi_V^{h_{in}}$ is $L^2$,
if and only if 
there exists a coherent subsheaf
$\prolongg{\vecc}{V}\subset\prolongg{\vecc}{E}$
such that $\prolongg{\vecc}{V}_{|X-D}=V$,
by Lemma \ref{lem;05.8.30.200}.
Let $h_{0,V}$ and $h_{in,V}$ denote the metrics of $V$
induced by $h_0$ and $h_{in}$.
We have
$\tr F(h_{in,V})=\tr F(h_{0,V})+\rank(V)\cdot \del\delbar g'$.
Then we obtain
$\deg_{\omega_{\epsilon}}(V,h_{0,V})
=\deg_{\omega_{\epsilon}}(V,h_{in,V})$
from the boundedness of $\del\delbar g'$ and $\del g'$
with respect to $\omega_{\epsilon}$.
Therefore the third condition is satisfied.
The fourth condition is satisfied by our construction.
The fifth condition is also checked
by using the boundedness of
$F(h_{in})$, $F(h_0)$, $\delbar\del g'$  and $\delbar g'$.
\hfill\qed

\vspace{.1in}

Now Proposition \ref{prop;05.7.30.10}
follows from Lemma \ref{lem;06.8.17.1}
and Proposition \ref{prop;05.7.30.35}.
\hfill\qed

\section{Bogomolov-Gieseker Inequality}

We have an immediate and standard corollary 
of Proposition \ref{prop;05.7.30.10},
as in \cite{s1}.
\begin{cor}
\label{cor;05.7.30.20}
Let $X$ be a smooth irreducible projective surface 
with an ample line bundle $L$,
and let $D$ be a simple normal crossing divisor of $X$.
Let $(\prolongg{\vecc}{E}_{\ast},\theta)$ be
a $\mu_L$-stable $\vecc$-parabolic graded semisimple Higgs bundle
on $(X,D)$.
Then we have the following inequality:
\[
 \int_X\parch_{2}(\prolongg{\vecc}{E}_{\ast})
-\frac{\int_X\parchern_1^2(\prolongg{\vecc}{E}_{\ast})}{2\rank E}
\leq 0.
\]
\end{cor}
\pf
Let $h$ be the metric of $E$ as in Proposition \ref{prop;05.7.30.10}.
Then we have the following:
\[
 \int_X\parch_{2}(\prolongg{\vecc}{E}_{\ast})
-\frac{\int_X\parchern_1^2(\prolongg{\vecc}{E}_{\ast})}{2\rank E}
=\left(\frac{\sqrt{-1}}{2\pi} \right)^2
 \int_{X-D}
 \tr\Bigl(F(h)^{\bot\,2}\Bigr).
\]
Then the claim follows from
$\tr\Bigl(F(h)^{\bot\,2}\Bigr)\geq 0$.
(See the pages 878--879 in \cite{s1}.)
\hfill\qed

\vspace{.1in}

By using the perturbation of the parabolic structure,
we can remove the assumption of graded semisimplicity.
We can also remove the assumption $\dim X=2$
by using Mehta-Ramanathan type theorem.
\begin{thm}[Bogomolov-Gieseker inequality]
 \label{thm;05.7.30.30}
Let $X$ be a smooth irreducible
projective variety of an arbitrary dimension
with an ample line bundle $L$,
and let $D$ be a simple normal crossing divisor.
Let $(\vecE_{\ast},\theta)$ be a 
$\mu_L$-stable regular Higgs bundle
in codimension two on $(X,D)$.
Then the following inequality holds:
\[
 \int_X\parch_{2,L}(\vecE_{\ast})
-\frac{\int_X\parchern_{1,L}^2(\vecE_{\ast})}{2\rank E}
\leq 0.
\]
(See Subsection {\rm\ref{subsection;05.10.6.1}}
for the characteristic numbers.)
\end{thm}
\pf
Due to the Mehta-Ramanathan type theorem
(Proposition \ref{prop;06.8.12.15}),
the problem can be reduced to the case
where $X$ is a surface.
Take a real number
$c_i\not\in \Par(\vecE_{\ast},i)$ for each $i$,
and let us consider the $\vecc$-truncation
$(\prolongg{\vecc}{E}_{\ast},\theta)$.
Let $\vecF$ denote the induced $\vecc$-parabolic structure
of $\prolongg{\vecc}{E}$.
Let $\epsilon$ be any sufficiently small positive number,
and let us take an $\epsilon$-perturbation $\vecF^{(\epsilon)}$
of $\vecF$ as in Section \ref{section;05.7.30.15}.
Since $(\prolongg{\vecc}{E},\vecF^{(\epsilon)},\theta)$
is $\mu_{L}$-stable and graded semisimple,
we obtain the following inequality due to
Corollary \ref{cor;05.7.30.20}:
\[
 \int_X\parch_{2}(\prolongg{\vecc}{E},\vecF^{(\epsilon)})
-\frac{\int_X\parchern_1^2(\prolongg{\vecc}{E},\vecF^{(\epsilon)})}
 {2\rank E}
\leq 0.
\]
By taking the limit in $\epsilon\to 0$,
we obtain the desired inequality.
\hfill\qed

\begin{cor}
\label{cor;06.8.5.30}
Let $X$ be a smooth irreducible projective surface
with an ample line bundle $L$,
and let $D$ be a simple normal crossing divisor.
Let $(\vecE_{\ast},\theta)$ be a $\mu_L$-stable 
parabolic Higgs bundle on $(X,D)$.
Assume $\int_X\parch_2(\vecE_{\ast})=\pardeg_{L}(\vecE_{\ast})=0$.
Then we have $\parchern_1(\vecE_{\ast})=0$.
\end{cor}
\pf
$\pardeg_{L}(\vecE_{\ast})=0$ implies
$\int_X \parchern_1(\vecE_{\ast})\cdot c_1(L)=0$.
Due to the Hodge index theorem,
it implies 
$-\int\parchern_1^2(\vecE_{\ast})\geq 0$,
and if the equality holds then $\parchern_1(\vecE_{\ast})=0$.
On the other hand,
we have the following inequality,
due to Theorem \ref{thm;05.7.30.30}:
\[
 -\frac{\int_X\parchern_{1}^2(\vecE_{\ast})}{2\rank E}
\leq
 -\int_X\parch_2(\vecE_{\ast})=0.
\]
Thus the claim follows.
\hfill\qed

%% file: 7.tex
We put
$X(T):=\{z\in\cnum\,|\,|z|<T\}$
and $X^{\ast}(T):=X(T)-\{O\}$,
where $O$ denotes the origin.
In the case $T=1$,
we omit to denote $T$.
Let $\harmonicbundle$ be a tame harmonic bundle
on $X^{\ast}$.
Recall that the coefficients $a_j(z)$
of $P(z,t):=\det(t-f_0(z))=\sum a_j(z)\cdot t^j$ 
are holomorphic on $X$,
where $f_0\in \End(E)$ is given by $\theta=f_0\cdot dz/z$.
The set of the solutions of the polynomial $P(0,t)$ is denoted by $S_0$.

\begin{assumption}
\label{assumption;05.9.2.4}
We assume the following:
\begin{enumerate}
\item
We have the decomposition
$E=\bigoplus_{\alpha\in S_0} E_{\alpha}$,
such that
$f_0(E_{\alpha})\subset E_{\alpha}$.
In particular, we have the decomposition
$f_0=\bigoplus f_{0\,\alpha}$.
\item
There exist some positive numbers 
$T_0<1$, $C_0$ and $\epsilon_0$
such that 
$|\beta-\alpha|<C_0\cdot|z(Q)|^{\epsilon_0}$
holds for any eigenvalue $\beta$ of $f_{\alpha\,|\,Q}$
 $(Q\in X^{\ast}(T_0))$.
\item
We put
 $\xi:=\sum_{\alpha\in S_0}
 \rank(E_{\alpha})\cdot |\alpha|^2$.
We assume $\xi<K_0$
for a given constant $K_0$.
\hfill\qed
\end{enumerate}
\end{assumption}

\begin{rem}
\label{rem;06.8.20.1}
The conditions {\rm 1} and {\rm 2}
are always satisfied,
if we replace $X$ by a smaller open set.
Moreover, it is controlled by the behaviour of
the eigenvalues of $f_0$.
\hfill\qed
\end{rem}

We obtain the parabolic Higgs bundle
$\bigl(\prolongg{a}{E}_{\ast},\theta\bigr)$
for $a\in\real$ from $(E,\delbar_E,h)$,
where $\prolongg{a}{E}$ is 
as in Section \ref{section;05.9.8.110}
(\cite{s2}).
In the case $a=0$,
we use the notation $\prolong{E}$.
Thus we have the parabolic filtration $F$ of $\prolongg{a}{E}_{|O}$
and the sets
$\Par(\prolongg{a}{E}):=
 \bigl\{b\,\big|\,\Gr^F_b(\prolongg{a}{E}_{|O})\neq 0\bigr\}$.
For any $b\in\Par(\prolongg{a}{E})$,
we put $\gminim(b):=\dim\Gr^F_b(\prolongg{a}{E}_{|O})$.
Recall $\det(\prolongg{a}{E})\simeq \prolongg{\atilde}{\det(E)}$,
where $\atilde$ is given as follows:
\[
 \widetilde{a}:=\sum_{b\in\Par(\prolongg{a}{E})}\gminim(b)\cdot b.
\]

Let $U_0$ be a finite subset of $\openopen{a-1}{a}$,
and let $\eta_0$ be a sufficiently small positive numbers
such that $U_0\subset\openopen{a-1+10\cdot\eta_0}{a-10\cdot\eta_0}$
and $|b-c|>10\cdot\eta_0$ for any distinct elements $b,c\in U_0$.
We make an additional assumption.
\begin{assumption}
 \label{assumption;06.8.10.1}
For any $c\in\Par(\prolongg{a}{E})$,
 there exists $b\in U_0$ such that
 $|c-b|<\eta_0$.
\hfill\qed
\end{assumption}
We put
$\nbigp(b):=\bigl\{c\in\Par(\prolongg{a}{E})\,
 \big|\,|c-b|<\eta_0\bigr\}$.
We obtain the decomposition
$\Par(\prolongg{a}{E})=\coprod_{b\in U_0}\nbigp(b)$.
We put $\bbar:=\max\nbigp(b)$.

\vspace{.1in}

In the following of this chapter,
we say that a constant $C$ is good,
if it depends only on $T_0$, $C_0$, $\epsilon_0$,
$K_0$, $\eta_0$ and $r:=\rank(E)$.
We say a constant $C(B)$ is good
if it depends also on additional data $B$.

\begin{prop} 
 \label{prop;05.8.28.50}
Let $(E,\delbar_E,\theta,h)$ be a tame harmonic bundle
on $X^{\ast}$ satisfying
the assumptions {\rm\ref{assumption;05.9.2.4}}
and {\rm\ref{assumption;06.8.10.1}}.
\begin{itemize}
\item
There exist holomorphic sections $F_1,\ldots,F_r$
of $\prolongg{a}{E}$ on $X(\gamma_0)$
with the numbers $b_1,\ldots,b_r\in U_0$,
such that
$|F_i|_{h}
 \leq C_{10}\cdot|z|^{-\bbar_i}\cdot (-\log|z|^2)^N$ holds.
Here $\gamma_0$, $C_{10}$ and $N$
are good constants.
 We have $\#\{b_i=b\}=\#\nbigp(b)$.
\item
 $C_{11}^{-1}\cdot|z|^{-\atilde}
 \leq \bigl|\bigwedge_{i=1}^rF_i\bigr|_h\leq C_{11}\cdot |z|^{-\atilde}$
 holds for a good constant $C_{11}$.
 In particular, $F_1,\ldots,F_r$ give the frame of $\prolongg{a}{E}$.
\item
On any compact subset $H\subset X^{\ast}(\gamma_0)$,
we have $\|F_{i|H}\|_{L_1^p,h}\leq C_{12}(H,p)$,
where $p$  is an arbitrarily large number.
\end{itemize}
\end{prop}

We will prove the proposition
in the rest of this chapter.

\section{A Priori Estimate of Higgs Field on a Punctured Disc}

Let $(E,\delbar_E,\theta,h)$ be a tame harmonic bundle
on $X^{\ast}$
as in Proposition \ref{prop;05.8.28.50}.
We know that the curvature $R(h)$ of $\delbar_E+\del_{E}$
is bounded with respect to $h$ and the Poincar\'{e} metric 
$\gtilde=|z|^{-2}(-\log|z|)^{-2}dz\cdot d\zbar$
on $X^{\ast}(T)$ for $T<1$.
(\cite{s2}. See also \cite{mochi2}).
We would like to show
that the estimate is uniform,
when we vary the set $S_0$ boundedly.

\begin{prop}
\label{prop;06.8.5.5}
$\bigl| R(h)\bigr|_{h,\gtilde}\leq K_{10}$
holds on $X^{\ast}(T_1)$
for some good constants $T_1$ and $K_{10}$.
\end{prop}
\pf
In the following argument,
$K_i$, $\epsilon_i$ and $T_i$ will denote 
good constants,
and $\Delta$ denotes the Laplacian
$-\del_z\del_{\zbar}$.
Let $\nbigl$ be a line bundle $\nbigo_{X^{\ast}}\cdot e$
with the Higgs field $\theta_{\nbigl}$ and the metric $h_{\nbigl}$
given by $\theta(e)=e\cdot \beta\cdot dz/z$
$(\beta\in\cnum)$ and $h_{\nbigl}(e,e)=1$.
Since we have only to consider
$(E,\delbar_E,\theta,h)\otimes
 (\nbigl,\theta_{\nbigl},h_{\nbigl})$,
we may and will assume $0<K_1<\xi<K_2$.

By an elementary argument,
we can take a decomposition
$S_0=\coprod_{i=1}^{k_0}S^{(1)}_i$
with the following property:
\begin{itemize}
\item
 $|\alpha_j-\alpha_k|\leq 1$ 
 for any $\alpha_j,\alpha_k\in S_i^{(1)}$.
\item
$|\alpha_j-\alpha_k|>\rank(E)^{-1}$
for $\alpha_j\in S_i^{(1)}$ and 
 $\alpha_k\in S_0-S_i^{(1)}$.
\end{itemize}
We put $\nbigs(1):=\{1,\ldots,k_0\}\subset\seisuu_{>0}$.
Inductively on $n$,
we take a subset
$\nbigs(n)\subset\seisuu_{>0}^n$
and a decomposition
$S_0=\coprod_{I\in\nbigs(n)}S^{(n)}_I$
as follows.
Assume $\nbigs(n)$ and $S^{(n)}_I$ $(I\in\nbigs(n))$
are already given.
We can take a decomposition
$S^{(n)}_I=\coprod_{i=1}^{k(I)} S_{I,i}^{(n+1)}$
with the following property:
\begin{itemize}
\item
$|\alpha_j-\alpha_k|\leq (n+1)^{-1}$
 for $\alpha_j,\alpha_k\in S_{I,i}^{(n+1)}$,
\item
$|\alpha_j-\alpha_k|>(n+1)^{-1}\cdot \rank(E)^{-1}$
 for $\alpha_j\in S^{(n+1)}_{I,i}$
 and $\alpha_k\in S^{(n)}_I-S^{(n+1)}_{I,i}$.
\end{itemize}
Then we put 
$\nbigs(n+1):=\bigl\{(I,i)\,\big|\,I\in
\nbigs(n),\,i=1,\ldots,k(I)\bigr\}$
and $S^{(n+1)}_{(I,i)}:=S^{(n+1)}_{I,i}$,
where $(I,i)\in \seisuu^{n+1}_{>0}$ 
denotes the element naturally determined by $I$ and $i$.

We have the lexicographic order on $\seisuu_{>0}^n$,
which induces the order on $\nbigs(n)$.
Take a total order $\leq_1$ on $S_0$,
which satisfies the following condition for any $n$:
\begin{itemize}
\item
 Let $\alpha\in S^{(n)}_I$ and $\beta\in S^{(n)}_J$.
 If $I<J$ in $\nbigs(n)$,
 we have $\alpha\leq_1 \beta$.
\end{itemize}
We put $F_{\alpha}E:=\bigoplus_{\beta\leq_1 \alpha}E_{\beta}$
and $F_{<\alpha}E:=\bigoplus_{\beta<_1\alpha}E_{\beta}$.
Let $E_{\alpha}'$ denote the orthogonal complement
of $F_{<\alpha}(E)$ in $F_{\alpha}(E)$.
We put 
$\rho:=\bigoplus_{\alpha\in S_0}\alpha\cdot\id_{E_{\alpha}}$
and $\rho':=\bigoplus_{\alpha\in S_0}\alpha\cdot \id_{E_{\alpha}'}$.
We have $|\rho'|_h^2=\xi$.
The following lemma is shown 
in the proof of Simpson's Main estimate.
(See \cite{s2} and the proof of Proposition 7.2 of
\cite{mochi2}.)

\begin{lem}
\label{lem;06.8.7.15}
$|f_0-\rho'|_h\leq K_{11}\cdot\bigl(-\log|z|\bigr)^{-1}$
holds on $X^{\ast}(T_1)$.
\hfill\qed
\end{lem}

For $J\in \nbigs(n)$,
we put 
$E^{(n)}_{J}:=\bigoplus_{\alpha\in S^{(n)}_{J}} E_{\alpha}$
and 
$E^{\prime(n)}_{J}:=\bigoplus_{\alpha\in S^{(n)}_{J}}E'_{\alpha}$.
We have the natural decomposition
$\End(E)=\bigoplus_{J_1,J_2\in\nbigs(n)}
 Hom(E^{\prime(n)}_{J_1},E^{\prime(n)}_{J_2})$.
For $I\in\nbigs(n-1)$ and $A\in \End(E)$,
let $A_{n,I,i,j}$ denote the 
 $Hom(E^{\prime(n)}_{I,i},E^{\prime(n)}_{I,j})$-component.

\begin{lem}
\label{lem;06.8.7.10}
We have
$\bigl|
\bigl[\rho^{\prime\,\dagger},f_0\bigr]_{n,I,i,j}
\bigr|_h\leq K_{30}\cdot (-\log|z|)^{2}$
for $i\neq j$
on $X^{\ast}(T_2)$.
\end{lem}
\pf
We put
$ \kappa:=\bigoplus_{i=1}^{k(I)}
 i\cdot \id_{E^{(n)}_{I,i}}$ and
$ \kappa':=\bigoplus_{i=1}^{k(I)}
 i\cdot \id_{E^{\prime\,(n)}_{I,i}}$.
We also put 
$q:=\kappa-\kappa'\in 
 \nbigt:=\bigoplus_{J_1>J_2}
 Hom(E^{\prime\,(n)}_{J_1},E^{\prime\,(n)}_{J_2})$.
First, we give some estimate of $q$.

Let $\varphi:X^{\ast}\lrarr X^{\ast}$
denote the map given by $\varphi(z)=z^n$.
We remark that $\varphi^{\ast}\harmonicbundle$
satisfies Assumption \ref{assumption;05.9.2.4}
independently of $n$,
if we replace $C_0$ with a larger good constant.
We put $\htilde:=\varphi^{\ast}h$.
We put $\ftilde_0:=n\cdot \varphi^{-1}(f_0)$,
i.e.,
$\varphi^{-1}\theta=\ftilde_0\cdot dz/z$.
Let $\ftilde_0^{\dagger}$ denote the adjoint of
$\ftilde_0$ with respect to $\htilde$.
We also put
$\rhotilde^{\prime}:=n\cdot \varphi^{-1}(\rho')$.

Let $F_{\ftilde_0}$ denote the endomorphism
of $\varphi^{-1}\nbigt$ induced by the adjoint of $\ftilde_0$,
i.e.,
$F_{\ftilde_0}(x)=[\ftilde_0,x]$.
Let $\pi_{\nbigt}$ denote the orthogonal projection
of $\varphi^{-1}\End(E)$ onto $\varphi^{-1}\nbigt$.
The composite of the adjoint of $\ftilde_0^{\dagger}$
and $\pi_{\nbigt}$ induces the endomorphism 
$G_{\ftilde_0^{\dagger}}$ of $\varphi^{-1}\nbigt$.

\begin{lem}
\footnote{The author thanks the referee
 for simplification of the proof}
\label{lem;06.8.7.5}
$F_{\ftilde_0}$ and 
$G_{\ftilde_0^{\dagger}}$ are 
invertible,
and the norms of their inverses 
are dominated 
by a good constant.
\end{lem}
\pf
Let $H$ denote the endomorphism of $\varphi^{-1}\nbigt$
induced by the adjoint of $\rhotilde^{\prime}$,
and we put $H_1:=G_{f_0^{\dagger}}-H$.
For any $\alpha\in S^{(n)}_I$ and $\beta\in S^{(n)}_J$ $(I\neq J)$,
we have $n\cdot |\alpha-\beta|> \rank(E)^{-1}$.
Hence the norm of $H^{-1}$
is dominated by a good constant.
From
$|\ftilde_0-\rhotilde^{\prime}|_{\htilde}\leq 
 K_{31}\cdot (-\log|z|)^{-1}$,
the norm of $H_1$ is dominated by a good constant.
It is easy to see that 
$H_2:=H^{-1}\circ H_1$ is nilpotent.
Then the claim for
$G_{f_0^{\dagger}}^{-1}=
 H^{-1}\circ\bigl(\id+\sum (-1)^jH_2^j\bigr)$
can be easily checked.
The claim for $F_{f_0}^{-1}$ can be checked similarly.
\hfill\qed

\vspace{.1in}

We put $\kappatilde:=\varphi^{-1}\kappa$,
$\kappatilde':=\varphi^{-1}\kappa'$
and $\qtilde:=\varphi^{-1}q$.
We have
$0=\bigl[\ftilde_0,\kappatilde\bigr]
=\bigl[\ftilde_0-\rhotilde^{\prime},\kappatilde'\bigr]
+\bigl[\ftilde_0,\qtilde\bigr]$.
Due to Lemma \ref{lem;06.8.7.15}
and Lemma \ref{lem;06.8.7.5},
we obtain the estimate
$\bigl|\qtilde\bigr|_{\htilde}\leq K_{32}(-\log|z|)^{-1}$.
From $[\kappatilde,\ftilde_0^{\dagger}]
=[\kappatilde-\kappatilde',\ftilde_0^{\dagger}]
+[\kappatilde',\ftilde_0^{\dagger}]$,
we obtain
$\bigl|[\kappatilde,\ftilde_0^{\dagger}]\bigr|_{\htilde}^2
\geq 
 \bigl|\pi_S\big([\kappatilde-\kappatilde',
 \ftilde_0^{\dagger}])\bigr|_{\htilde}^2
 =\bigl|G_{\ftilde_0^{\dagger}}(\qtilde)\bigr|_{\htilde}^2$.
Hence, we obtain
$|\qtilde|_{\htilde}^2
\leq
K_{33}\bigl|[\kappatilde,\ftilde_0^{\dagger}]\bigr|_{\htilde}$.
Due to $\bigl[\kappatilde,\ftilde_0\bigr]=0$,
we obtain the following:
\[
 \Delta\log|\kappatilde|_{\htilde}^2
\leq
 -\frac{\bigl|[\ftilde_0^{\dagger},\kappatilde]\bigr|^2_{\htilde}}
 {|z|^2\cdot |\kappatilde|^2_{\htilde}}
\leq
 -K_{35}\frac{|\qtilde|_{\htilde}}{|z|^2}
\]
We put $\xi':=\sum_{i=1}^{k(I)} i^2=|\kappatilde'|^2_{\htilde}$
and $k:=\log \bigl(\xi^{\prime\,-1}
 |\kappatilde|^2_{\htilde}\bigr)$.
Because of 
$k\leq \xi^{\prime\,-1}|\qtilde|_{\htilde}^2$,
we obtain
$\Delta k\leq -K_{36}\cdot|z|^{-2}\cdot k$.
By an argument in \cite{s2}
(see also the proof of Proposition 7.2 of \cite{mochi2}),
we obtain $k\leq K_{37}\cdot |z|^{\epsilon_{38}}$.
Then, we can derive
$|\qtilde|_{\htilde}\leq K_{39}\cdot |z|^{\epsilon_{38}}$
on $X^{\ast}(T_{40})$.
Hence we obtain
$|q|_h\leq K_{39}\cdot |z|^{\epsilon_{38}/n}$
on $X^{\ast}(T_{40}^n)$.

Let us finish the proof of Lemma \ref{lem;06.8.7.10}.
First we show the estimate on $X^{\ast}(T_{40}^n)$.
We have
$ 0=[\kappa,f_0]
=[\kappa',f_0]+[q,\rho']+[q,f_0-\rho']$.
We have the following on $X^{\ast}(T_{40}^n)$:
\[
\bigl|[q,f_0-\rho']\bigr|_h
\leq
 \frac{K_{41}|z|^{\epsilon_{38}/n}}{-\log|z|}
\leq \frac{K_{42}|z|^{\epsilon_{38}/n}}{n}
\]
Recall we have $|\alpha-\beta|\leq (n-1)^{-1}$
for $\alpha\in S^{(n-1)}_{I,i}$ and $\beta\in S^{(n-1)}_{I,j}$.
Hence  we have
$\bigl|\bigl[q,\rho'\bigr]_{n,I,i,j}\bigr|_h
\leq
K_{43}\cdot |z|^{\epsilon_{38}/n}\cdot n^{-1}$.
Therefore, we obtain
$\bigl| \bigl[\kappa',f_0\bigr]_{n,I,i,j}\bigr|_h
\leq K_{44}\cdot |z|^{\epsilon_{38}/n}\cdot n^{-1}$,
which implies
$\bigl|(f_0)_{n,I,i,j}\bigr|_h\leq K'_{44}\cdot
|z|^{\epsilon_{38}/n}\cdot n^{-1}$ $(i\neq j)$.
Then, we obtain the estimate on $X^{\ast}(T_{40}^n)$:
\[
\bigl|
 \bigl[\rho^{\prime\,\dagger},f_0\bigr]_{n,I,i,j}
\bigr|_h
\leq
 K_{45}\cdot |z|^{\epsilon_{38}/n}\cdot n^{-2}
\leq
 K_{46}\cdot \bigl(-\log|z|\bigr)^{-2}.
\]

On the other hand,
$\bigl[\rho^{\prime\,\dagger},f_0\bigr]_{n,I,i,j}$
is dominated by
$K_{47}\cdot (-\log|z|)^{-1}\cdot n^{-1}$
on $X^{\ast}(T_1)$,
which is obtained by the estimate of $f_0-\rho'$ 
(Lemma \ref{lem;06.8.7.15})
and our choice of 
$S^{(n)}_{I,k}$ $(k=i,j)$.
Outside of $X^{\ast}(T_{40}^n)$,
we have 
$K_{47}\cdot (-\log|z|)^{-1}\cdot n^{-1}
\leq K_{48}\cdot (-\log|z|)^{-2}$.
Thus we are done.
\hfill\qed

\vspace{.1in}
Let us finish the proof of Proposition 
\ref{prop;06.8.5.5}.
We have the following:
\[
 R(h)=-\bigl[\theta,\theta^{\dagger}\bigr]
=-\Bigl(
\bigl[\rho',\,(f_0-\rho')^{\dagger}\bigr]
+\bigl[f_0-\rho',\,\rho^{\prime\,\dagger}\bigr]
+\bigl[(f_0-\rho'),\,(f_0-\rho')^{\dagger}
 \bigr]
\Bigr)\cdot\frac{dz\cdot d\bar{z}}{|z|^2}
\]
The second term is estimated by Lemma \ref{lem;06.8.7.10}.
The first term is adjoint of the second term.
The estimate of the third term follows from
Lemma \ref{lem;06.8.7.15}.
Thus we are done.
\hfill\qed

\section{Uhlenbeck Type Theorem on a Punctured Disc}
\label{section;05.8.29.101}

We put $\hyperh:=
\{(x,y)\,|\,y>0\}$,
and let $\varphi:\hyperh\lrarr X^{\ast}$
be the covering given by
$\varphi(x,y)=\exp\bigl(2\pi\sqrt{-1}(x+\sqrt{-1}y)\bigr)$.
We use the Euclidean  metric
$g_0=dx\cdot dx+dy\cdot dy$.
We put as follows: 
\[
 H_n:=\bigl\{(x,y)\in\hyperh\,\big|\,n-1/2<y<n+1,\,-1<x<1/2\bigr\}
\]
\[
 H_{n,1}:=\bigl\{(x,y)\in H_n\,\big|\,
 -1<x<-1/2
 \bigr\},
\quad
 H_{n,2}:=\bigl\{(x,y)\in H_n\,\big|\,
 0<x<1/2
 \bigr\}.
\]
Let $g_0$ denote the Euclidean metric of $\hyperh$.

Let $\harmonicbundle$ be a tame harmonic bundle
on $X^{\ast}$.
We put $\nabla:=\del_E+\delbar_E$.
After taking the pull back via the map
$\phi_{\gamma}:X^{\ast}\lrarr X^{\ast}$
given by $\phi_{\gamma}(z)=\gamma\cdot z$,
we may assume to have
$C^1$-orthonormal frames $\vecv_n$ of 
 $\varphi^{\ast}E_{|H_n}$
satisfying the following conditions,
due to Proposition \ref{prop;06.8.5.5}
and the Uhlenbeck's result (\cite{u1}):
\begin{itemize}
\item
Let $A_n$ denote the connection one form of $\nabla$
with respect to $\vecv_n$ on $H_n$.
Then the norm of $A_n$ with respect to $g_0$
is dominated as $\bigl|A_n\bigr|_{g_0}\leq C\cdot n^{-2}$,
where $C$ is a good constant.
\end{itemize}

\begin{prop}
 \label{prop;05.8.28.2}
There exists a good constant $\gamma_0$ and
a $C^1$-orthonormal frame $\vecw$ of $E$ on $X^{\ast}(\gamma_0)$
for which $\delbar_E$ is expressed as follows:
\[
 \delbar_E \vecw=
 \vecw\cdot\Bigl(A-\frac{1}{2}\Gamma\Bigr)\cdot
 \frac{d\zbar}{\zbar}
\]
\begin{itemize}
\item
 $\Gamma$ is a constant diagonal matrix whose 
$(i,i)$-th components $\alpha_i$ satisfy
$0\leq \alpha_r\leq \alpha_{r-1}\leq\cdots\leq \alpha_1<1$.
\item
 $|A|\leq C_1\cdot \bigl(-\log|z|\bigr)^{-1}$
 for some good constant $C_1$.
\end{itemize}
\end{prop}
\pf
In the following argument,
the constants $C_i$ and $N_i$ depend only on $C$.
Let $U(r)$ and $\gminiu(r)$ denote 
the $r$-th unitary group and its Lie algebra,
respectively.
Let $s_n:H_n\cap H_{n-1}\lrarr U(r)$
be the $C^1$-function determined by
$\vecv_{n-1}=\vecv_n\cdot s_n$.
From the relation $ds_n=s_n\cdot A_n-A_{n-1}\cdot s_n$,
we have
$\bigl|ds_n\bigr|_{C^0}\leq C_3\cdot n^{-2}$.
If we take a sufficiently large constants
$N_1$ and $C_4$,
we have the expression
$s_n=S_n\cdot \exp\bigl(\widetilde{s}_n\bigr)$
for each $n\geq N_1$,
where
$\widetilde{s}_n:H_n\cap H_{n-1}\lrarr \unitary(r)$
satisfies $\bigl|\widetilde{s}_n\bigr|_{C^1}\leq C_4\cdot n^{-2}$,
and $S_n$ denotes an element of $U(r)$.
Let us consider
$\vecv'_n=\vecv_n\cdot S_n\cdot S_{n-1}\cdot \cdots \cdot S_{N_1}$
instead of $\vecv_n$.
Then we have the function
$s_n':H_n\cap H_{n-1}\lrarr \unitary(r)$
satisfying $\vecv_n':=\vecv'_{n-1}\cdot \exp(s_n')$
and $\bigl|s_n'\bigr|_{C^1}\leq C_4\cdot n^{-2}$.
The connection one form $A'_n$ of $\nabla$
with respect to the frame $\vecv_n'$
satisfies $|A'_n|\leq C\cdot n^{-2}$.
Hence we may and will assume $S_n=1$
for any $n$.

We put $\vecv_{n,a}:=\vecv_{n\,|\,H_{n,a}}$ $(a=1,2)$.
Via the identification
$H_{n,1}\simeq \varphi(H_{n,a}) \simeq H_{n,2}$,
both of
$\vecv_{n,a}$ give the frames of $\varphi^{\ast}E_{|H_{n,1}}$,
which are denoted by the same notation.
Let $t_n:H_{n,1}\lrarr U(r)$ be the $C^1$-function
determined by $\vecv_{n,1}=\vecv_{n,2}\cdot t_n$.
We have $\bigl|dt_n\bigr|_{C^0}\leq C_5\cdot n^{-2}$
as in the case of $s_n$.
If we take constants $N_2$ and $C_6$
appropriately,
we have the expression
$t_n=T_n\cdot \exp\bigl(\widetilde{t}_n\bigr)$
for any $n\geq N_2$,
where $\widetilde{t}_n:H_{n,1}\lrarr \unitary(r)$ 
satisfies $\bigl|\widetilde{t}_{n}\bigr|_{C^1}\leq C_6\cdot n^{-2}$,
and $T_n$ denotes an element of $U(r)$.

The distance between the monodromies
with respect to the two loops contained in $\varphi(H_n)$
can be dominated by $C_7'\cdot n^{-2}$.
Hence the distance of the monodromies
with respect to the two loops contained in 
$\varphi(H_n)$ and $\varphi(H_m)$ are 
dominated by $C_7''\cdot \max\bigl\{n^{-1},m^{-1}\bigr\}$.
Hence the distance between $T_n$ and $T_m$
is dominated by 
$C_7\cdot \max\bigl( n^{-1},m^{-1} \bigr)$
for any $n,m\geq N_2$.

We put $T:=\lim_{n\to \infty} T_n$.
We take $B\in U(r)$ such that
$B^{-1}TB=\exp\bigl(2\pi\sqrt{-1}\cdot\Gamma\bigr)$,
where $\Gamma$ is the diagonal matrix
whose $(i,i)$-th components $\alpha_i$
satisfy $0\leq \alpha_r\leq \cdots\leq \alpha_1<1$.
Since we can consider $\vecv_n'=\vecv_n\cdot B$
instead of $\vecv_n$,
we may and will assume $T=\exp\bigl(2\pi\sqrt{-1}\Gamma\bigr)$.
If we take appropriately large 
constants $N_3$ and $C_8$,
we have the functions $\check{t}_n:H_{n,1}\lrarr \unitary(r)$
satisfying
$t_n=T\cdot\exp\bigl(\check{t}_n\bigr)$
and
$\bigl|\check{t}_n\bigr|_{C^1} \leq C_8\cdot n^{-1}$
for any $n\geq N_3$.

We put $\beta(x,y):=\exp\bigl(2\pi\sqrt{-1}\Gamma\cdot x\bigr)$
and
$\widehat{\vecv}_n
 :=\vecv_n\cdot \beta$.
We put $\widehat{\vecv}_{n,a}:=\widehat{\vecv}_{n\,|\,H_{n,a}}$
$(a=1,2)$,
and we regard them as the frames of
$\varphi^{\ast}E_{|H_{n,1}}$.
We put
$\widehat{t}_n:=\beta_{|H_{n,1}}^{-1}\cdot
 \check{t}_n\cdot \beta_{|H_{n,1}}$.
Then we have
$\widehat{\vecv}_{n,1}
=\widehat{\vecv}_{n,2}\cdot \exp\bigl(\widehat{t}_n\bigr)$
and 
$\bigl|\widehat{t}_n\bigr|_{C^1}
\leq C_9\cdot n^{-1}$.
We put $\widehat{s}_n=\beta^{-1}\cdot \widetilde{s}_n\cdot \beta$.
We have
$\widehat{\vecv}_{n-1}=\widehat{\vecv}_n\cdot
 \exp\bigl(\widehat{s}_n\bigr)$
and the estimates
$\bigl|\widehat{s}_n\bigr|_{C^1}\leq C_{10}\cdot n^{-2}$.
We put
$\widehat{A}_n:=
 \beta^{-1}\cdot A_n\cdot \beta+2\pi\sqrt{-1}\Gamma\cdot dx$.
Then we have the relation
$\nabla\widehat{\vecv}_n=\widehat{\vecv}_n\cdot \widehat{A}_n$
and the estimates
$\bigl|\widehat{A}_n-2\pi\sqrt{-1}\Gamma\cdot dx\bigr|_{C^0}
 \leq C_{11}\cdot n^{-2}$.

\vspace{.1in}

Take a $C^{\infty}$-function $\chi:\hyperh\lrarr [0,1]$
satisfying the following:
\[
 \chi(x,y)=
 \left\{
 \begin{array}{ll}
 1 & (x\leq -1/3),\\
 0 & (x\geq -1/6).
 \end{array}
 \right.
\]
We put $\Psi_n:=\exp\bigl(-\chi(x,y)\cdot \widehat{t}_n\bigr)$,
and $\vecu_n:=\widehat{\vecv}_n\cdot \Psi_{n}$.
Under the natural identification $H_{n,1}\simeq H_{n,2}$,
we have
$ \vecu_{n\,|\,H_{n,2}}
=\widehat{\vecv}_{n,2}
=\widehat{\vecv}_{n,1}\cdot\exp\bigl(-\widehat{t}_n\bigr)
=\vecu_{n\,|\,H_{n,1}}$.
Hence $\vecu_n$ gives the orthonormal frame
of $E_{|\varphi(H_n)}$.
If we take appropriately large constants
$N_4$ and $C_{12}$,
 the following holds:
\begin{itemize}
\item
 We have the functions $\sbar_n:H_n\cap H_{n-1}\lrarr\unitary(r)$
 satisfying
 $\vecu_{n-1}=\vecu_n\cdot \exp\bigl(\sbar_n\bigr)$
 and
 $\bigl|\sbar_n\bigr|_{C^1}\leq C_{12}\cdot n^{-1}$
for any $n\geq N_4$.
\item
 The connection forms $B_n$ of $\nabla$
 with respect to the frames $\vecu_n$ satisfy
 $\bigl|B_n-2\pi\sqrt{-1}\cdot\Gamma\cdot dx\bigr|_{g_0}
 \leq C_{12}\cdot n^{-1}$.
\end{itemize}

Take a $C^{\infty}$-function
$\rho_n:\hyperh\lrarr [0,1]$
satisfying the following:
\[
 \rho_n(x,y):=\left\{
 \begin{array}{ll}
 1 & \bigl(y\leq n\bigr),\\
 0 & (y\geq n+1/2).
 \end{array}
 \right.
\]
We put $\vecw_n:=\vecu_n\cdot \exp\bigl(\rho_n\cdot \sbar_n\bigr)$
on $H_n$.
Let $\widehat{B}_n$ denote the connection one form
of $\nabla$ with respect to $\vecw_n$.
Then we have the estimate
$\bigl|
 \widehat{B}_n-2\pi\sqrt{-1}\Gamma dx
 \bigr|_{g_0}\leq C_{14}\cdot n^{-1}$.

On the intersection
$H_n\cap H_{n-1}$
we have the relation
$\vecw_{n-1}=\vecu_{n-1}=\vecu_n\cdot \exp\bigl(\sbar_n\bigr)=\vecw_n$.
Hence $\bigl\{\vecw_n\bigr\}$ gives the orthonormal frame $\vecw$
of $E$ on $X^{\ast}(\gamma_1)$ for some constant $\gamma_1$,
which depends only on $C$.
Let $B$ be a connection one form of $\nabla$ with respect to
the frame $\vecw$.
Let us denote $\varphi^{\ast}B$ as 
$B^{x}\cdot dx+B^{y}\cdot dy$.
If we take appropriately large constants
$C_{15}$ and $N_5$,
the estimates
$\bigl|B^{x}-2\pi\sqrt{-1}\Gamma\bigr|\leq C_{15}\cdot y^{-1}$
and 
$\bigl|B^{y}\bigr|\leq C_{15}\cdot y^{-1}$
hold on $\bigl\{(x,y)\,\big|\,y>N_5\bigr\}$
by our construction.
We have the following formula on $X^{\ast}(\gamma_1)$:
\[
 \delbar \vecw
=\vecw\cdot\left(
 -\frac{1}{2}\Gamma
+\frac{\sqrt{-1}}{2}\Bigl(
 \frac{B^{x}}{2\pi}-\sqrt{-1}\Gamma\Bigr)
-\frac{1}{4\pi} B^{y}
 \right)\cdot\frac{d\zbar}{\zbar}.
\]
Then $\vecw$ gives the desired frame.
Therefore the proof of Proposition \ref{prop;05.8.28.2}
is accomplished.
\hfill\qed

\section{Construction of Local Holomorphic Frames}
\label{section;05.8.29.200}

\subsection{Setting}
\label{subsection;06.8.5.25}

Let $(E,\delbar_E,\theta,h)$ be a harmonic bundle of rank $r$
on $X^{\ast}$
as in Proposition \ref{prop;05.8.28.50}.
We will construct the desired holomorphic sections
in Proposition \ref{prop;05.8.28.50}.
By considering the tensor product of 
the line bundle with the metric $|z|^{-c}$,
we have only to discuss the case $a=0$.
We use the metrics
$g$ and  $\widetilde{g}$ of $X^{\ast}$ given as follows:
\[
 g:=dz\cdot d\zbar,
\quad
\widetilde{g}:=\frac{dz\cdot d\zbar}{|z|^2\cdot\bigl(-\log|z|\bigr)^2}
\]

By considering a pull back via the map
$\phi_{\gamma}:X^{\ast}\lrarr X^{\ast}$
given by $\phi_{\gamma}(z)=\gamma\cdot z$,
we may assume the following,
due to Proposition \ref{prop;06.8.5.5}
and Proposition \ref{prop;05.8.28.2}.
\begin{assumption}
 \label{assumption;05.8.29.1}
The norm of $R(h)$
with respect to $h$ and $g$ is dominated as follows:
\[
\bigl|
 R(h)
\bigr|_{h,g}
\leq C_1\cdot \frac{1}{|z|^{2}\cdot \bigl(-\log|z|+1\bigr)^2}.
\]
 There exists a $C^1$-orthonormal frame $\vecv$ of $E$,
 for which $\delbar_E$ is represented as follows:
 \[
  \delbar_E\vecv=\vecv\cdot\left(
 -\frac{\Gamma}{2}+A
 \right)\cdot \frac{d\zbar}{\zbar}
 \]
Here $\Gamma$ is a constant diagonal matrix
whose $(i,i)$-th components $\alpha_i$ satisfy
$0\leq \alpha_r\leq \cdots\leq\alpha_1<1$,
and $A$ is a matrix-valued continuous function
such that $|A|\leq C_2\cdot (-\log|z|+1)^{-1}$.
The constants $C_1$ and $C_2$ are good.
\hfill\qed
\end{assumption}

\subsection{Preliminary for a proof}

We recall some results on the solvability of the $\delbar$-equation.
For any real numbers $b$ and $M$,
we put $h(b,M):=h\cdot |z|^{2b}\cdot (-\log |z|)^{M}$.
Let $A^{0,1}_{b,M}(E)$ denote the space of
sections of $E\otimes\Omega^{0,1}$,
which are  $L^2$ with respect to $h(b,M)$ and $\widetilde{g}$.
Let $A^{0,0}_{b,M}(E)$ denote the space of
sections $f$ of $E$ such that
$f$ and $\delbar f$ are $L^2$ 
with respect to $h(b,M)$ and $\widetilde{g}$.
The the norm of $A^{p,q}_{b,M}(E)$ is denoted by
$\|\cdot\|_{b,N}$.
On the other hand,
$|\cdot|_{b,M}$ denote the norm at fibers.
In the following argument,
$B_i$ will denote good constants.

Recall some results of 
\cite{cg} and \cite{mochi2}.
(But the signature is slightly changed.)
Take a sufficiently large good constant $N$,
which depends only on $C_1$
in Assumption \ref{assumption;05.8.29.1}.
\begin{lem}
\label{lem;06.8.7.30}
For any $f_1\in A^{0,1}_{b,N}(E)$,
we have $f_2\in A^{0,0}_{b,N}$
satisfying $\delbar f_2=f_1$ and
$\|f_2\|_{b,N}\leq B_{1}\cdot
 \|f_1\|_{b,N}$.
\hfill\qed
\end{lem}

On the other hand,
if $f$ is a holomorphic section of $E$,
we have the subharmonicity
$ \Delta\log |f|_{b,-N}\leq 0$.
Hence, if we have $\|f\|_{b,N}<\infty$,
the following holds around the origin $O$:
\begin{multline}
 \log|f(z)|^2_{b,-N}
 \leq
  \frac{4}{\pi|z|^2}
 \int_{|w-z|\leq |z|/2}
 \log |f(w)|_{b,-N}^2\cdot \dvol_g \\
\leq
 \log\left(
 \frac{4}{\pi|z|^2}
 \int_{|w-z|\leq |z|/2}
 |f(w)|_{b,-N}^2\cdot \dvol_g
 \right) 
\leq
 \log\bigl(B_{2}\cdot\|f\|_{b,N}^2\bigr)
\end{multline}
Hence, we obtain the following lemma.
\begin{lem}
\label{lem;06.8.9.1}
For a holomorphic section $f$ of $E$
such that $\|f\|_{b,N}<\infty$,
we have
$|f|_h
\leq B_{2}\|f\|_{b,N}\cdot|z|^{-b}\cdot(-\log|z|)^{N}$.
\hfill\qed
\end{lem}
We give one more elementary remark.
\begin{lem}
\label{lem;06.8.9.2}
Let $f$ be a holomorphic section of
$\prolongg{b}{E}$ on $X(\gamma')$.
Then the maximum principle holds
for $H(z):=|f(z)|_h^2\cdot |z|^{2b}\cdot (-\log|z|)^{-N}$
on $X(\gamma'')$ for $\gamma''<\gamma'$,
i.e.,
$\sup_{X^{\ast}(\gamma'')}H(z)
=\max_{\del X(\gamma'')}H(z)$.
\end{lem}
\pf
We put $H_{\epsilon}(z):=
 |f(z)|_h^2\cdot |z|^{2b+\epsilon}\cdot(-\log|z|)^{-N}$
for any $\epsilon>0$.
We have $\Delta \log H_{\epsilon}\leq 0$ on 
$X^{\ast}(\gamma')$
and $\lim_{z\to 0}\log H_{\epsilon}(z)=-\infty$.
Therefore, the maximum principle holds
for $\log H_{\epsilon}$ on $X(\gamma'')$.
Then it is easy to derive the maximum principle for $H$.
\hfill\qed

\vspace{.1in}

\subsection{Proof}

Take $0<\eta\leq \eta_0$.
Let $\Gamma$ and $\vecv=(v_1,\ldots,v_r)$ be
as in Assumption \ref{assumption;05.8.29.1}.
We put $S(\Gamma):=\bigl\{\alpha_1,\ldots,\alpha_r\bigr\}$.
Let $T_A$ denote the section of $\End(E)\otimes\Omega^{0,1}$
determined by $\vecv$ and $A\cdot d\zbar/\zbar$,
i.e.,
$T_A(\vecv)=\vecv\cdot A\cdot d\zbar/\zbar$.
We put $\delbar_0:=\delbar-T_A$.
We put $f_i:=|z|^{\alpha_i}\cdot v_i$.
Then we have $\delbar_0f_i=0$
and $|f_i|_h=|z|^{\alpha_i}$.
In particular,
we have $f_i\in A^{0,0}_{-\alpha_i+\eta,N}(E)$.
Take $g_i\in A^{0,0}_{-\alpha_i+\eta,N}$
satisfying $\delbar g_i=T_A(f_i)$
and 
$\bigl\| g_i \bigr\|_{-\alpha_i+\eta,N}
\leq B_{1}\cdot \bigl\|T_A(f_i)\bigr\|_{-\alpha_i+\eta,N}$
as in Lemma \ref{lem;06.8.7.30}.
We put $F_i:=f_i-g_i$.
Then we have $\delbar F_i=0$,
$F_i\in A^{0,0}_{-\alpha_i+\eta,N}(E)$,
and the following estimate:
\[
 \bigl\| F_i  \bigr\|_{-\alpha_i+\eta,N}
\leq
 \bigl\|f_i\bigr\|_{-\alpha_i+\eta,N}
+B_{1}\cdot\bigl\| T_A(f_i)\bigr\|_{-\alpha_i+\eta,N}.
\]
We have the following:
\begin{equation}
\label{eq;04.9.19.45}
 \delbar_0 g_i=-T_A(g_i)+T_A(f_i).
\end{equation}
Hence we obtain $g_i\in L_1^2(H)$
for any compact subset $H\subset X^{\ast}$,
and the $L_1^2$-norm is dominated by
$||T_A(f_i)||_{-\alpha_i+\eta,N}$
multiplied by some constant depending only on $H$.
Hence for some number $p>2$ and some good constant $C'(H,p)$,
we have the following:
\[
 \bigl\|g_i\bigr\|_{L^p(H)}
\leq
 C'(H,p)\cdot \bigl\|T_A(f_i)\bigr\|_{-\alpha_i+\eta,N}
\]
Due to (\ref{eq;04.9.19.45}),
we have the following, for some good constant $C''(H,p)$:
\begin{equation} \label{eq;04.9.19.51}
 \bigl\|  g_i  \bigr\|_{L_1^p(H)}
\leq
 C''(H,p)\cdot\Bigl(
  \bigl\|T_A(f_i)\bigr\|_{-\alpha_i+\eta,N}
+\sup_H \bigl|T_A(f_i)\bigr|_{h,\widetilde{g}}
 \Bigr).
\end{equation}
By a standard boot strapping argument,
$p$ can be arbitrarily large.

We put $\alphatilde:=\tr(\Gamma)$
and $\zerotilde:=
 \sum_{b\in \Par(\prolong{E})}\gminim(b)\cdot b$.
Since we have $\tr(R(h))=\tr(F(h))=0$,
the induced metric $\det(h)$ of $\det(E)$ is flat.
Hence we have a holomorphic section $s$ 
of $\prolongg{\zerotilde}{\det(E)}=\det(\prolong{E})$
such that $|s|_h=|z|^{-\zerotilde}$
and $\del_{\det(E)}s=s\cdot (-\zerotilde)\cdot dz/z$.
It is easy to see
 $n=\alphatilde+\zerotilde$ is an integer
by considering the limit of the monodromy
of $\det(E)$ around the origin.
We put $\stilde:=z^n\cdot s$,
which gives the section of 
$\prolongg{-\alphatilde}{\det(E)}$.
\begin{rem}
We will show that $-\widetilde{\alpha}=\zerotilde$,
i.e. $s=\widetilde{s}$ later
(Lemma {\rm\ref{lem;05.9.4.30}}).
\hfill\qed
\end{rem}
Let us consider the function $\widetilde{F}$
determined by
$\widetilde{F}\cdot \widetilde{s}=F_1\wedge\cdots\wedge F_r$.
We put $H_0:=\bigl\{z\,\big|\,3^{-1}\leq|z|\leq 2\cdot 3^{-1}\bigr\}$.

\begin{lem}
 \label{lem;04.9.19.55}
 There exists a small good constant $B_{15}$ with the following property:
\begin{itemize}
\item
 Assume the following inequalities hold:
\begin{equation}
\label{eq;04.9.19.50}
  \sup_{H_0}\bigl|A\bigr|_{\widetilde{g}}<B_{15},\quad
  \bigl\|T_A\cdot f_i\bigr\|_{-\alpha_i+\eta,N}<B_{15},\quad
 (i=1,\ldots,r).
\end{equation}
Then, there exist $z_0\in\{z\in\cnum\,|\,|z|=1\}$
and a good constant $0<B_{16}<1/2$
such that 
$\Ftilde(H_0)\subset \{z\in\cnum\,|\,|z-z_0|<B_{16}\}$.
\end{itemize}
\end{lem}
\pf
From (\ref{eq;04.9.19.51}) and (\ref{eq;04.9.19.50}),
we obtain
$\bigl|
F_1\wedge\cdots\wedge F_r-f_1\wedge \cdots\wedge f_r
\bigr|< 4^{-1}$ holds on $H_0$,
if $B_{15}$ is sufficiently small.
Since $v_1,\ldots,v_r$ are orthonormal
and $f_i$ are given as $|z|^{\alpha_i}\cdot v_i$, 
we have 
$f_1\wedge \cdots\wedge f_r=
 \exp\bigl(\sqrt{-1}\kappa\bigr)\cdot \widetilde{s}$
for some real valued functions $\kappa$.
If $B_{15}$ is sufficiently small,
$\kappa$ is a sum of a constant $\kappa_0$ and a function 
$\kappa_1$ satisfying
$\sup_{H_0}|\kappa_1(z)|<100^{-1}$
because of (\ref{eq;04.9.19.50}).
Then the claim of the lemma follows.
\hfill\qed

\vspace{.1in}

For any number $0<\gamma<1$,
let us consider the map
$\phi_{\gamma}:X^{\ast}\lrarr X^{\ast}$
given by $z\longmapsto \gamma\cdot z$.
We put
$\bigl(E(\gamma),\delbar_{E(\gamma)},
 \theta(\gamma),h(\gamma)\bigr)
:=\phi_{\gamma}^{\ast}\harmonicbundle$.
We have the orthonormal frame 
$\phi_{\gamma}^{\ast}\vecv$
of $E(\gamma)$
for which we have the following:
\[
 \delbar_{E(\gamma)} \bigl(\phi_{\gamma}^{\ast}\vecv\bigr)
=\phi_{\gamma}^{\ast}\vecv\cdot
 \left(
 -\frac{1}{2}\Gamma+
 \phi_{\gamma}^{\ast}A
 \right)\cdot \frac{d\zbar}{\zbar}.
\]
Note we have the following:
\begin{equation}
\label{eq;04.9.19.56}
 \bigl|\phi_{\gamma}^{\ast}A
 \bigr|_{h(\gamma),\widetilde{g}}
\leq C_3\cdot 
\frac{-\log|z|+1}{-\log|z|-\log|\gamma|+1}
\cdot \bigl(-\log|z|+1\bigr)^{-1}.
\end{equation}
Hence it is easy to check 
Assumption \ref{assumption;05.8.29.1}
for $\bigl(E(\gamma),\delbar_{E(\gamma)},
h(\gamma)\bigr)$.
We put
$f^{(\gamma)}_i:=|z|^{\alpha_i}\cdot \phi_{\gamma}^{\ast}v_i$.
We construct the sections
$g_i^{(\gamma)}$ and $F_i^{(\gamma)}$ as above.
We also take $\stilde^{(\gamma)}$ and $s^{(\gamma)}$.

\begin{lem}
For $\eta>0$,
there exists $\gamma_1=\gamma_1(\eta)>0$
such that the assumptions of Lemma  {\rm\ref{lem;04.9.19.55}}
are satisfied for
$\bigl(E(\gamma_1),\delbar_{E(\gamma_1)},h(\gamma_1)\bigr)$
and $\phi_{\gamma_1}^{\ast}\vecv$.
\end{lem}
\pf
If $\gamma$ is sufficiently small,
then we may assume
$\sup_{H_0}\bigl|\phi_{\gamma}^{\ast}A\bigr|_{\widetilde{g}}\leq B_{15}$
due to (\ref{eq;04.9.19.56}).
We also have the following:
\begin{multline}
 \int \bigl|T_{\phi_{\gamma}^{\ast}A}\cdot
   f_i^{(\gamma)}\bigr|^2_{h(\gamma),\widetilde{g}}
 \cdot |z|^{-2\alpha_i+2\eta}\cdot \bigl(-\log|z|\bigr)^{N}\cdot
 \dvol_{\widetilde{g}} \\
\leq
 B_{18}\cdot
 \int \left|
 \frac{-\log|z|+1}{-\log|z|-\log \gamma +1}
 \right|^2\cdot |z|^{2\eta}
\cdot\bigl(-\log|z|\bigr)^N\cdot\dvol_{\widetilde{g}}.
\end{multline}
Since the right hand side converges to $0$
in $\gamma\lrarr 0$,
we can take $\gamma_1$
such that the inequality
$\bigl\|T_{\phi_{\gamma_1}^{\ast}A}
 f_i^{(\gamma_1)}\bigr\|_{-\alpha_i+\eta,N}
 <B_{15}$ holds.
\hfill\qed

\vspace{.1in}

Now we have the holomorphic sections
$F_1^{(\gamma_1(\eta))},\ldots,F_r^{(\gamma_1(\eta))}$
of $\prolong{E(\gamma_1)}$,
satisfying 
$|F_i^{(\gamma_1(\eta))}|_{h(\gamma_1(\eta))}\leq
 C(\eta)\cdot |z|^{\alpha_i-\eta}\bigl(-\log|z|\bigr)^N$.
We put
$ a_i(\eta):=\max\bigl\{
 b\in \Par\bigl(\prolong{E}\bigr)\,\big|\,
 b\leq -\alpha_i+\eta
 \bigr\}$,
and then 
$F_i^{(\gamma_1(\eta))}$
 are sections of
$\prolongg{a_i(\eta)}{E(\gamma_1(\eta))}$.

\begin{lem} 
 \label{lem;05.9.4.30}
 \label{lem;05.9.4.31}
We have $\Par(\prolong{E})=S(\Gamma)$
which preserves the multiplicity.
Hence, we have 
$-\widetilde{\alpha}=\zerotilde$.
\end{lem}
\pf
If $\eta$ is sufficiently small,
we have $a_i(\eta)\leq-\alpha_i$
and hence $\sum a_i(\eta)\leq -\alphatilde$.
We put $\gamma_2:=\gamma_1(\eta)$.
Hence we obtain
$\bigl|\bigwedge_{i=1}^rF_i^{(\gamma_2)}\bigr|_{h(\gamma_2)}
=O(|z|^{\alphatilde})$,
which implies $\Ftilde$ is holomorphic on $X$,
where $\Ftilde$ is given by
$\bigwedge_{i=1}^r F_i^{(\gamma_2)}
=\Ftilde\cdot \stilde^{(\gamma_2)}$.
Due to Lemma \ref{lem;04.9.19.55}
and the maximum principle,
we obtain
$B_{20}^{-1}\leq |\Ftilde(z)|\leq B_{20}$
for $z\in X(2/3)$.
Hence, we obtain $\sum a_i(\eta)=-\alphatilde$.

We put $S(b):=\{i\,|\,-\alpha_i=b\}$
for $b\in\Par\bigl(\prolong{E(\gamma_2)}\bigr)$.
For $i\in S(b)$,
we have $F_i^{(\gamma_2)}\in\prolongg{b}{E(\gamma_2)}$,
which induces $\Fbar_i^{(\gamma_2)}\in\Gr^{F}_b(E(\gamma_2))$.
From $B_{20}^{-1}\leq |\Ftilde(z)|\leq B_{20}$
for $z\in X(2/3)$,
we have the lower estimate
$|\bigwedge_{i\in S(b)}F_i^{(\gamma_2)}|_{h(\gamma_2)}
 \geq C_{\delta}|z|^{-|S(b)|\cdot b+\delta}$
for any $\delta>0$.
Hence we obtain the  linearly independence of
$\Fbar_i^{(\gamma_2)}$ $(i\in S(b))$.
Then, it is easy to show that 
$F_1^{(\gamma_2)},\ldots,F_r^{(\gamma_2)}$
give the frame of $\prolong{E(\gamma_2)}$
compatible with the parabolic structure,
whose parabolic degrees are 
$-\alpha_1,\ldots,-\alpha_r$, respectively.
\hfill\qed

\vspace{.1in}

Now let us fix $\eta=\eta_0$.
We put $\gamma_3:=\gamma_1(\eta_0)$.
We have the holomorphic sections
$F_i^{(\gamma_3)}$ of
$\prolong{E(\gamma_3)}$
on $X$
satisfying
$|F_i^{(\gamma_3)}|_{h(\gamma_3)}\leq
 B_{30}\cdot|z|^{\alpha_i(\eta_0)-\eta_0}$.
Since we have $s^{(\gamma_3)}=\stilde^{(\gamma_3)}$,
the function $\Ftilde$ determined by
$F_1^{(\gamma_3)}\wedge
 \cdots\wedge F_i^{(\gamma_3)}=
 \Ftilde\cdot \stilde^{(\gamma_3)}$
is holomorphic on $X$.
Thus, we have
$B_{31}^{-1}\leq |\Ftilde(z)|\leq B_{31}$
for $z\in X(2/3)$
due to the maximum principle and
Lemma \ref{lem;04.9.19.55}.

The holomorphic sections
$F_i^{(\gamma_3)}$ of $\prolong{E(\gamma_3)}$
on $X$ naturally give
the holomorphic sections
$\widehat{F}_i$ of $\prolong{E}$ on $X(\gamma_3)$.
We take $\gamma_0<\gamma_3$ appropriately,
and we put
$F_i:=\widehat{F}_{i\,|\,X(\gamma_0)}$.
It is clear that they satisfy the second and third claims
of Proposition \ref{prop;05.8.28.50}.

For each $a_i(\eta_0)$, we have the number $b_i\in U_0$
such that $a_i(\eta_0)\in \nbigp(b_i)$.
We obtain $F_i\in \prolongg{\bbar_i}{E}$.
Then, the first claim of Proposition \ref{prop;05.8.28.50}
follows from Lemma \ref{lem;06.8.9.2} and the third claim.
\hfill\qed

%% file: 8.tex
\section{Convergence of a Sequence of Tame Harmonic Bundles}

Let $X$ be a smooth projective variety of an arbitrary dimension,
and $D$ be a simple normal crossing divisor of $X$.
Let $(E_m,\delbar_{m},\theta_m,h_m)$ $(m=1,2,\ldots,)$
be a sequence of tame harmonic bundles of rank $r$ on $X-D$.
We have the associated parabolic Higgs bundles
$(\prolongg{\vecc}{E}_{m\ast},\theta_m)$
on $(X,D)$.

\begin{thm}
 \label{thm;05.9.10.2}
Assume that
 the sequence of the sections
 $\{\det(t-\theta_m)\}$ 
 of $\Sym^{\cdot}\Omega^{1,0}_X(\log D)[t]$  are convergent.
Then the following claims hold:
\begin{itemize}
\item
There exists a subsequence
 $\bigl\{(E_m,\delbar_m,\theta_m,h_m)\,\big|\,  m\in I\bigr\}$
which converges to a tame harmonic bundle
$(E_{\infty},\delbar_{\infty},\theta_{\infty},h_{\infty})$ on $X-D$,
weakly in $L_2^p$ locally on $X-D$,
in the sense of Section {\rm\ref{section;05.10.8.1}}.
Here $p$ denotes an arbitrarily large number.
\item
If we are given
a parabolic Higgs sheaf $(\prolongg{\vecc}{E}_{\ast},\theta)$
such that
$\bigl\{(\prolongg{\vecc}{E_m}_{\ast},\theta_m)_{|C}\bigr\}$
converges to $(\prolongg{\vecc}{E}_{\ast},\theta)_{|C}$
for any generic curve $C$.
Then we have a non-trivial holomorphic morphism
$f:(\prolongg{\vecc}{E}_{\ast},\theta)
\lrarr
(\prolongg{\vecc}{E_{\infty\ast}},\theta_{\infty})$.

If $(\prolongg{\vecc}{E}_{\ast},\theta)$ is 
a $\mu_L$-stable reflexive saturated parabolic Higgs sheaf,
$f$ is isomorphic.
(See Lemma {\rm\ref{lem;06.8.22.20}}.)
\end{itemize}
\end{thm}
\pf
The first claim is well known.
We recall only an outline.
The sequence of sections $\{\det(t-\theta_m)\}$
of $\Sym^{\cdot}\Omega^{1,0}_X[t]$ converges to $\det(t-\theta)$.
Hence we obtain the estimate of the norms of $\theta_m$
locally on $X-D$ (See Lemma \ref{lem;05.8.25.4}, for example).
We also obtain the estimate of
the curvatures $R(h_m)$ because of the relation
$R(h_m)+[\theta_m,\theta_m^{\dagger}]=0$.
Therefore, we obtain the local convergence result like the first claim.
(See  \cite{s5} in the page 26--28, for example.)
Thus we obtain the harmonic bundle
$(E_{\infty},\delbar_{\infty},\theta_{\infty},h_{\infty})$.

Let us show the second claim in Subsection 
\ref{subsection;06.8.21.20} after some preparation.

\subsection{On a punctured disc}
\label{subsection;06.8.20.30}

Let us explain the setting in this subsection.
Let $X(\gamma)$ and $X^{\ast}(\gamma)$
denote the disc $\{z\in\cnum\,|\,|z|<\gamma\}$
and the punctured disc $X(\gamma)-\{0\}$.
In the case $\gamma=1$,
we use the notation $X$ and $X^{\ast}$.
We put $D:=\{0\}$.
Let $(E_m,\delbar_m,\theta_m,h_m)$ $(m=1,2,\ldots,\infty)$ be
a sequence of tame harmonic bundles of rank $r$
on a punctured disc $X^{\ast}$.
We have the associated parabolic Higgs bundles
$(\prolongg{c}{E}_{m\,\ast},\theta_m)$ on $(X,D)$
for $c\in\real$.
Assume the following:
\begin{itemize}
\item
 $\bigl\{(E_m,\delbar_m,\theta_m,h_m)\,\big|\,m<\infty\bigr\}$
 converges to $(E_{\infty},\delbar_{\infty},\theta_{\infty},h_{\infty})$
   in $C^1$ locally on $\Delta^{\ast}$
 via the isometries  $\Phi_m:(E_m,h_m)\lrarr (E_{\infty},h_{\infty})$.
\item
 Assumption \ref{assumption;05.9.2.4} is satisfied
 for any $m$. The constants are independent of 
 the choice of $m$.
\item
 There exists a finite subset $U_0\subset\openopen{c-1}{c}$
 and a function $\gminimbar:U_0\lrarr \seisuu_{>0}$
 such that
 $\bigl\{\bigl(\Par(\prolongg{c}{E_m}),\gminim\bigr)
 \,\big|\,m<\infty\bigr\}$
 converges to $(U_0,\gminimbar)$
 in the sense of Section \ref{section;05.10.8.1}.
 We put $u:=\sum_{b\in U_0} \gminimbar(b)\cdot b$.
\end{itemize}

\begin{lem}
\label{lem;06.8.21.30}
We have holomorphic isomorphisms
$\Psi_{m'}:\prolongg{c}{E_{m'}}\lrarr \prolongg{c}{E_{\infty}}$
on $X(\gamma)$
for some $\gamma<1$
and some subsequence $\{m'\}\subset\{m\}$,
with the following properties:
\begin{itemize}
\item
  $\Psi_m-\Phi_m\lrarr 0$ weakly in $L_1^p$ locally on 
  $X^{\ast}(\gamma)$.
\item
 $\Psi_m(\theta_m)-\theta_{\infty}\to 0$
 as holomorphic sections of 
 $\End(\prolongg{c}{E_{\infty}})\otimes\Omega^{1,0}(\log D)$
 on $X(\gamma)$.
\item
 Let $F^{(m)}(\prolongg{c}{E_m})$
 denote the parabolic filtrations
 of $\prolongg{c}{E_m}_{|D}$ induced by $h_m$.
 Then the sequence of the filtrations
$\bigl\{\Psi_m \bigl(F^{(m)}(\prolongg{c}{E_m}_{|D})\bigr)
 \bigr\}$
converges to $F^{(\infty)}(\prolongg{c}{E_{\infty}})_{|D}$
in the sense of Definition {\rm\ref{df;05.8.29.300}}.
\end{itemize}
\end{lem}
\pf
After going to a subsequence,
we may assume that
 Assumption \ref{assumption;06.8.10.1} is satisfied
 for  $(E_m,\delbar_{m},\theta_m,h_m)$ $(m<\infty)$
with some $\eta_0>0$.
We take holomorphic sections
$F_1^{(m)},\ldots, F_r^{(m)}$
of $\prolongg{c}{E_m}$ on $X(\gamma)$
with $b_1^{(m)},\ldots, b_r^{(m)}\in U_0$
as in Proposition \ref{prop;05.8.28.50},
with some $\gamma<1$.
We may assume that $b_i^{(m)}$ are independent of $m$,
which are denoted by $b_i$.
There exists a subsequence  $\{m'\}$
such that
$\bigl\{ \Phi_{m'}(F_i^{(m')}) \bigr\}$
are convergent weakly in $L_1^{p}$ locally on $X(\gamma)^{\ast}$.
The limits are denoted by
$F_i^{(\infty)}$.
They are holomorphic with respect to $\delbar_{\infty}$.
We replace $\{m\}$ with the subsequence $\{m'\}$,
and we assume that the above convergence holds
from the beginning.

For each $b\in U_0$,
we put $\bbar(m):=
 \max\{a\in\Par(\prolongg{c}{E_m})\,|\,|a-b|<\eta_0\}$.
Then, we have
$|F_i^{(m)}|_{h_m}<C\cdot |z|^{-\bbar_i(m)}\cdot(-\log|z|)^N$,
where the constants $C$ and $N$ are independent of $m$.
Since we have $\bbar_i(m)\to b_i$ for $m\to\infty$,
we obtain
$|F_i^{(\infty)}|_{h_{\infty}}
 <C\cdot |z|^{-b_i}\cdot (-\log|z|)^N$,
and hence $F_i^{(\infty)}\in \prolongg{b_i}{E_{\infty}}$.

We put $\ctilde(m):=
 \sum_{b\in\Par(\prolongg{c}{E_m})}
 b\cdot \gminim(b)$.
The sequence $\{\ctilde(m)\}$ converges to $u$.
We have
$C_1^{-1}\cdot |z|^{-\ctilde(m)}
\leq \bigl|\bigwedge_{i=1}^r F_i^{(m)}\bigr|_{h_m}
 \leq C_1\cdot |z|^{-\ctilde(m)}$,
and hence
$C_1^{-1}\cdot |z|^{-u}
\leq \bigl|\bigwedge_{i=1}^r F_i^{(\infty)}\bigr|_{h_{\infty}}
 \leq C_1\cdot |z|^{-u}$.
We put $S_b:=\{i\,|\,b_i=b\}$.
For $i\in S_b$,
we have $F_i^{(\infty)}\in \prolongg{b}{E_{\infty}}$,
which induces $\Fbar_i^{(\infty)}\in\Gr^F_b(E_{\infty|D})$.
We have the lower estimate
$\bigl|\bigwedge_{i\in S_b}F_i^{(\infty)}\bigr|_{h_{\infty}}
 \geq C_{\delta}\cdot |z|^{-|S_b|\cdot b+\delta}$
for any $\delta>0$,
from which we obtain the linearly independence of
$\Fbar_i^{(\infty)}$ $(i\in S_b)$
in $\Gr^F_{b}(E_{\infty|D})$.
Then, it can be shown that
the sections  $F_1^{(\infty)},\ldots,F_r^{(\infty)}$ give
a holomorphic frame of $\prolongg{c}{E_{\infty}}$,
which is compatible with the parabolic structure,
and $b_i^{(\infty)}$ are the degrees of $F_i^{(\infty)}$
with respect to the parabolic structure.

\vspace{.1in}

We construct the holomorphic map
$\Psi_m:\prolongg{c}{E_m}\lrarr \prolongg{c}{E_{\infty}}$
on $X(\gamma)$
by the correspondence $\Psi_m(F_i^{(m)})= F_i^{(\infty)}$.
The first and third claims are satisfied from our construction.
We also have that the convergence of
$\Psi_m(\theta_m)-\theta_{\infty}$ to $0$ weakly in $L^p$
on any compact subset of $X^{\ast}(\gamma)$.
Since $\Psi_m(\theta_m)$ and $\theta_{\infty}$ are holomorphic,
the second claim also holds.
\hfill\qed

\subsection{On a curve}
\label{subsection;05.9.6.30}

Let us explain the setting in this subsection.
Let $C$ be a smooth projective curve
with a finite subset $D_C\subset C$.
Let $(E_{m},\delbar_m,h_m,\theta_m)$ $(m=1,2,\ldots,\infty)$
be a sequence of harmonic bundles of rank $r$ on $C-D_C$.
We have the associated Higgs bundles
$(\prolongg{\vecc}{E}_{m\,\ast},\theta_m)$,
where $\vecc=(c(P)\,|\,P\in D)\in\real^D$.
We assume the following:
\begin{itemize}
\item 
 The sequence $\bigl\{(E_m,\delbar_m,h_m,\theta_m)\bigr\}$
converges to
 $(E_{\infty},\delbar_{\infty},h_{\infty},\theta_{\infty})$
in $C^1$ locally on $C-D_C$
via isometries $\Phi_m:(E_m,h_m)\lrarr (E_{\infty},h_{\infty})$.
\item
 For each $i$, a finite subset
 $U(P)\subset \openopen{c(P)-1}{c(P)}$
 and a function $\gminim:U(P)\lrarr\seisuu_{>0}$ are given,
 and $\{\bigl(\Par(E_m,P),\gminim\bigr)\,|\,m<\infty\}$ converges to
 $\bigl(U(P),\gminim\bigr)$.
\end{itemize}
By the first condition, the sequence
$\det(t-\theta_m)\in \Sym^{\cdot}\Omega^{1,0}_C(\log D_C)$
converges to $\det(t-\theta_{\infty})$.
 Around each point $P\in D_C$,
 we can take  a coordinate neighbourhood $V_P$
 such that  Assumption \ref{assumption;05.9.2.4}
 is satisfied  on $V_P$  for any $m<\infty$,
 and that the constants are independent of $m$.

\begin{lem}
\label{lem;06.8.21.35}
$\bigl\{(\prolongg{\vecc}{E}_{m'},\vecF^{(m')},\theta_{m'})
 \,\big|\,m'\in I\bigr\}$
converges to
$(\prolongg{\vecc}{E_{\infty}},\vecF^{(\infty)},\theta_{\infty})$
for an appropriate subsequence $I\subset \{m\}$
in the sense of Definition {\rm\ref{df;05.8.29.300}}.
\end{lem}
\pf
We would like to replace $\Phi_{m'}$
with $\Psi_{m'}:\prolongg{\vecc}{E}_{m'}
 \lrarr \prolongg{\vecc}{E_{\infty}}$
for an appropriate subsequence $\{m'\}\subset\{m\}$.
By shrinking $V_P$ appropriately,
we take the holomorphic maps
$\lefttop{P}\Psi_{m'}:
 \prolongg{c(P)}{E_{m'}}
\lrarr
 \prolongg{c(P)}{E_{\infty}}$ on $V_P$
for some subsequence $\{m'\}\subset\{m\}$
for each point $P\in D_C$,
as in Lemma \ref{lem;06.8.21.30}.
We replace $\{m\}$ with $\{m'\}$.

Let $\chi_P:C\lrarr [0,1]$ denote a $C^{\infty}$-function
which is constantly $1$ around $P$,
and constantly $0$ on $C-V_P$.
Let $\Psi_m:E_m\lrarr E_{\infty}$ 
be the $L_1^p$-map  given as follows:
\begin{equation}
\label{eq;06.8.21.2}
 \Psi_m:=\sum_P\chi_P\cdot \lefttop{P}\Psi_m
+\Bigl(1-\sum_P\chi_P\Bigr)\cdot \Phi_m.
\end{equation}
If $m$ is sufficiently large,
then $\Psi_m$ are isomorphisms.
We have the following:
\begin{multline}
 \Psi_m\circ\delbar_m-\delbar_{\infty}\circ\Psi_m
=\sum\delbar\chi_P\cdot(\lefttop{P}\Psi_m-\Phi_m)\\
+\Bigl(1-\sum\chi_P\Bigr)\cdot
 \Bigl(
  \Phi_m\circ\delbar_m-\delbar_{\infty}\circ\Phi_m
 \Bigr).
\end{multline}
Hence the sequence
$\bigl\{
\Psi_m\circ\delbar_m-\delbar_{\infty}\circ\Psi_m
\bigr\}$ converges to $0$ weakly in $L^p$ on $C$.
By construction,
the sequence of the parabolic filtrations
of $\prolongg{\vecc}{E}_{m\,\ast}$
converges that of $\prolongg{\vecc}{E_{\infty\,\ast}}$.
We also have the convergence of
$\Psi_m(\theta_m)-\theta_{\infty}$ to $0$
weakly in $L^p$ on $C$.
Hence we obtain the convergence of
$\bigl\{(\prolongg{\vecc}{E_m},\vecF^{(m)},\theta_m)\,\big|\,
 m<\infty\bigr\}$ to
$\bigl(\prolongg{\vecc}{E_{\infty}},\vecF^{(\infty)},\theta_{\infty}\bigr)$
weakly in $L_1^p$ on $C$.
\hfill\qed

\subsection{The end of Proof of Theorem \ref{thm;05.9.10.2}}
\label{subsection;06.8.21.20}

Let us return to the setting for Theorem 
\ref{thm;05.9.10.2}.
Let $(E_{\infty},\delbar_{\infty},\theta_{\infty},h_{\infty})$
be a harmonic bundle obtained as a limit.
We obtain the parabolic Higgs bundle
$(\prolongg{\vecc}{E}_{\infty\ast},\theta_{\infty})$.
We would like to show the existence of a non-trivial holomorphic 
homomorphism
$\bigl(\prolongg{\vecc}{E}_{\ast},\theta\bigr)
\lrarr \bigl(
\prolongg{\vecc}{E}_{\infty\ast},\theta_{\infty}\bigr)$.
Due to Lemma \ref{lem;06.8.22.2},
we have only to show 
the existence of a non-trivial map
$f_C:\bigl(\prolongg{\vecc}{E}_{\ast},\theta\bigr)_{|C}
\lrarr
 \bigl(\prolongg{\vecc}{E}_{\infty\,\ast},\theta_{\infty}\bigr)_{|C}$
for some sufficiently ample generic curve $C\subset X$.
We may and will assume that $c_i\not\in\Par(\prolongg{\vecc}{E},i)$.

We have the convergence of the sequence
$\bigl\{
 \bigl(\prolongg{\vecc}{E_m}_{\ast},\theta_m\bigr)_{|C}
 \,\big|\,m
 \bigr\}$
to $\bigl(\prolongg{\vecc}{E}_{\ast},\theta\bigr)_{|C}$ on $C$.
In particular,
we have the convergence
$\bigl\{(\Par(\prolongg{\vecc}{E_m}_{|C},P),\gminim)\,\big|\,m<\infty\bigr\}$
to $(\Par(\prolongg{\vecc}{E}_{|C},P),\gminim)$
for any $P\in C\cap D$.
The sequence
$\bigl\{(E_m,\delbar_m,\theta_m,h_m)_{|C\setminus D}\bigr\}$
is convergent to
$(E_{\infty},\delbar_{\infty},
 \theta_{\infty},h_{\infty})_{|C\setminus D}$
in $C^1$ locally on $C\setminus D$.
After going to a subsequence,
we obtain the convergence of
$\bigl\{
\bigl(\prolongg{\vecc}{E_m}_{\ast},,\theta_m\bigr)_{|C}
 \,\big|\,m
 \bigr\}$
to
$\bigl(\prolongg{\vecc}{E_{\infty\ast}},
 \theta_{\infty}\bigr)_{|C}$
weakly in $L_1^p$  on $C$,
due to Lemma \ref{lem;06.8.21.35}.
Thus we obtain the existence of the desired non-trivial map $f_C$ 
due to Corollary \ref{cor;05.8.29.301}.
Thus the proof of Theorem \ref{thm;05.9.10.2} is finished.
\hfill\qed

\section{Preparation for the Proof of Theorem \ref{thm;05.9.8.150}}

Let $C$ be a smooth projective curve
with a simple effective divisor $D$.
Let $\{(\prolongg{\vecc}{E}_{m\ast},\theta_m)\}$
be a sequence of stable parabolic Higgs bundles
on $(C,D)$ with
$\pardeg(\prolongg{\vecc}{E}_{m\ast})=0$,
which converges to a stable Higgs bundle
$(\prolongg{\vecc}{E}_{\infty\ast},\theta_{\infty})$.
We take pluri-harmonic metrics
$h_{0}^{(m)}$ of $(E_m,\delbar_{E_m},\theta_m)$
adapted to the parabolic structure
$(m=1,2,\ldots,\infty)$
(Proposition \ref{prop;05.9.160}),
where $E_m:=\prolongg{\vecc}{E}_{m|C-D}$.
We put $\nbigd_m:=\delbar_{E_m}+\theta_m$ and
$\nbigd^{\star}_{m}:=
 \del_{E_m,h_0^{(m)}}
 +\theta^{\dagger}_{m,h_0^{(m)}}$
 $(m=1,2,\ldots,\infty)$.

Take a sequence of small positive numbers $\{\epsilon_m\}$.
For each $P\in D$,
let $(V_P,z)$ be a holomorphic coordinate around $P$
such that $z(P)=0$.
Let $N$ be a large positive number,
for example $N>10$.
Let $g_{m}$ be Kahler metrics
of $C-D$ with the following form on $V_P$ for each $P\in D$:
\[
\bigl(
 \epsilon_m^{N+2}|z|^{2\epsilon_m}
+|z|^{2}
\bigr)\frac{dz\cdot d\zbar}{|z|^2}.
\]
We assume that $\{g_{m}\}$ converges
to a smooth Kahler metric $g_{0}$ of $C$
in the $C^{\infty}$-sense locally on $C-D$.

In the following argument,
$\|\rho\|_{h,g}$ will denote the $L^2$-norm
of a section $\rho$ of $E_m\otimes\Omega^{p,q}_{C-D}$
or $\End(E)\otimes\Omega^{p,q}_{C-D}$,
with respect to a metric $g$ of $C-D$
and a metric $h$ of $E_m$.
On the other hand,
$|\rho|_{h,g}$ will denote the norm at fibers.

\begin{prop}
\label{prop;06.8.21.5}
Let $h^{(m)}$ $(m<\infty)$
be hermitian metrics of $E_m$ 
with the following properties:
\begin{enumerate}
\item
 Let $s^{(m)}$ be determined by
 $h^{(m)}=h_0^{(m)}\cdot s^{(m)}$.
 Then $s^{(m)}$ is bounded with respect to
 $h^{(m)}_0$,
 and we have $\det s^{(m)}=1$.
 We also have the finiteness
$ \bigl\|\nbigd_{m} s^{(m)}
 \bigr\|_{h_0^{(m)},g_{m}}<\infty$.
(The estimates may depend on $m$.)
\item
We have
$\|F(h^{(m)})\|_{h^{(m)},g_{m}}
 <\infty$
and $\lim_{m\to\infty}
 \|F(h^{(m)})\|_{h^{(m)},g_{m}}=0$.
\item
 There exists a tame harmonic bundle
 $(E',\delbar_{E'},\theta',h')$
 such that the sequence
 $\bigl\{(E_{m},\delbar_{E_m},\theta_m,h^{(m)})\bigr\}$
 converges to
 $(E',\delbar_{E'},\theta',h')$ in $C^1$
 locally on $C-D$.
\end{enumerate}
Then, after going to a subsequence,
$\{(\prolongg{\vecc}{E_m}_{\ast},\theta_{m})\}$
converges to
$(\prolongg{\vecc}{E'}_{\ast},\theta')$
weakly in $L_1^p$ on $C$.
\end{prop}
\pf
We may and will assume that
$\{(E_m,\delbar_{E_m},\theta_m,h_m)\}$
converges to
$(E_{\infty},\delbar_{E_{\infty}},\theta_{\infty},h_{\infty})$
via the isometries
$\Phi_m:(E_m,h_m)\lrarr (E_{\infty},h_{\infty})$,
due to Theorem \ref{thm;05.9.10.2}.
First, 
let us show that $s^{(m)}$ are bounded
independently of $m$.

\subsection{Uniform boundedness of $s^{(m)}$}

For any point $P\in C-D$,
let $SE(s^{(m)})(P)$ denote the maximal eigenvalue
of $s^{(m)}_{|P}$.
There exists a constant $0<C_1<1$
such that
$ C_1\cdot |s^{(m)}_{|P}|_{h_0^{(m)}}
\leq SE(s^{(m)})(P)
\leq |s^{(m)}_{|P}|_{h^{(m)}_0}$.
Because of
$\det s^{(m)}_{|P}=1$,
we have
$SE(s^{(m)})(P)\geq 1$ for any $P$.

Let us take $b_{m}>0$
satisfying 
$2\leq b_{m}\cdot\sup_P SE(s^{(m)})(P)
\leq 3$.
We put
 $\widetilde{s}^{(m)}=b_{m}\cdot s^{(m)}$
and 
$\widetilde{h}^{(m)}:=
 h^{(m)}_0\cdot \widetilde{s}^{(m)}$.
Then $\widetilde{s}^{(m)}$
are uniformly bounded with respect to $h_0^{(m)}$.
We remark
$F(\widetilde{h}^{(m)})
=F(h^{(m)})$.
We also remark that $h^{(m)}$
and $\widetilde{h}^{(m)}$ induce
the same metric of $\End(E_m)$.

Recall the following equality
(Lemma 3.1 of \cite{s1}):
\begin{equation}
 \label{eq;06.1.15.250}
 \Delta_{g_0,h_0^{(m)}} \widetilde{s}^{(m)}
=\widetilde{s}^{(m)}
 \sqrt{-1}\Lambda_{g_0}
 F(\widetilde{h}^{(m)})
+\sqrt{-1}\Lambda_{g_0}
 \nbigd_m\widetilde{s}^{(m)}
 \bigl(\widetilde{s}^{(m)\,-1}\bigr)
 \nbigdstar_{m} \widetilde{s}^{(m)}.
\end{equation}
Because of
$\|\nbigd_m s^{(m)}\|_{h^{(m)}_0,g_{m}}
=\|\nbigd_m s^{(m)}\|_{h^{(m)}_0,g_0}<\infty$
and the boundedness of $\stilde^{(m)}$,
we have
$\int
 \Delta_{g_0}\tr \widetilde{s}^{(m)}
 \cdot\dvol_{g_0}=0$.
Hence, we obtain the following inequality
from (\ref{eq;06.1.15.250})
and the uniform boundedness of
$\widetilde{s}^{(m)}$:
\begin{multline}
\int\bigl|
 \nbigd_m \widetilde{s}^{(m)}\cdot
 \widetilde{s}^{(m)\,-1/2}
 \bigr|_{g_0,h_0^{(m)}}^2
 \dvol_{g_0}\leq
 A\cdot\int
 \bigl|\tr \Lambda_{g_0} F(\widetilde{h}^{(m)})\bigr|
 \cdot\dvol_{g_0} \\
=A\cdot\int 
 \bigl|\tr \Lambda_{g_{m}}F(\widetilde{h}^{(m)})
 \bigr|\cdot\dvol_{g_{m}}
\leq
 A'\cdot\bigl\|F(\widetilde{h}^{(m)})
 \bigr\|_{\widetilde{h}^{(m)},g_{m}}.
\end{multline}
Here, $A$ and $A'$ denote the constants
which are independent of $m$.
In particular,
we obtain the following inequality
for some $A''$:
\begin{equation}
\label{eq;06.8.21.40}
 \bigl\|\nbigd_m \widetilde{s}^{(m)}
 \bigr\|_{h_0^{(m)},g_0}
\leq
 A''\cdot\bigl\| F(\widetilde{h}^{(m)})
 \bigr\|_{\widetilde{h}^{(m)},g_{m}}
\end{equation}
We put $\ttilde^{(m)}:=\Phi_m(\stilde^{(m)})\in \End(E_{\infty})$.
\begin{lem}
\label{lem;06.8.21.1}
After going to an appropriate subsequence,
$\bigl\{\ttilde^{(m)}\bigr\}$ converges
to a positive constant multiplication
weakly in $L_1^2$ locally on $C-D$.
\end{lem}
\pf
$\bigl\{\widetilde{t}^{(m)}\bigr\}$ is $L^2_1$-bounded
on any compact subset of $C-D$
due to (\ref{eq;06.8.21.40}).
By going to an appropriate subsequence,
it is weakly $L^2_1$-convergent locally on $C-D$.
Let $\widetilde{t}^{(\infty)}$ denote the weak limit.
We obtain $\nbigd_{\infty} \widetilde{t}^{(\infty)}=0$
from (\ref{eq;06.8.21.40}).
By construction,
$\widetilde{t}^{(\infty)}$ is also bounded with respect to
$h^{(\infty)}_0$.
Therefore
$\widetilde{t}^{(\infty)}$ gives an automorphism of
$(\prolongg{\vecc}{E}_{\infty\ast},\theta_{\infty})$.
Due to the stability of
$(\prolongg{\vecc}{E}_{\infty\ast},\theta_{\infty})$,
$\widetilde{t}^{(\infty)}$ is a constant multiplication.

We would like to show $\widetilde{t}^{(\infty)}\neq 0$.
Let us take any point $Q_{m}\in C-D$
satisfying the following:
\[
  SE(s^{(m)})(Q_{m})
\geq
\frac{9}{10}\cdot
 \sup_{P\in C-D} SE(s^{(m)})(P).
\]
Then we have
$ \log \tr \widetilde{t}^{(m)}(Q_{m})
\geq \log (9/5)$.
By taking an appropriate subsequence,
we may assume that the sequence $\{Q_{m}\}$ converges
to a point $Q_{\infty}$.
We have two cases
(i) $Q_{\infty}\in D$
(ii) $Q_{\infty}\not\in D$.
We discuss only the case (i).
The other case is similar and easier.

We have $\tr\stilde^{(m)}=\tr\ttilde^{(m)}$,
which we do not distinguish in the following.
We use the coordinate neighbourhood $(U,z)$
such that $z(Q_{\infty})=0$.
For any point $P\in U$,
we put $\Delta(P,T):=\{Q\in U\,|\,|z(P)-z(Q)|<T\}$.
Let $g=dz\cdot d\zbar$ denote 
the standard metric of $U$.
We have the following inequality on $U-\{Q_{\infty}\}$
(Lemma 3.1 of \cite{s1}):
\[
 \Delta_g\log \tr\widetilde{s}^{(m)}
\leq
 \bigl|\Lambda_g F(\widetilde{h}^{(m)})
 \bigr|_{\widetilde{h}^{(m)}}.
\]
Let $B^{(m)}$ be the endomorphism of $E_m$
determined as follows:
\[
  F(\widetilde{h}^{(m)})
=F(h^{(m)})
=B^{(m)}\cdot \frac{dz \cdot d\zbar}{|z|^2},
\]
Then we have the following estimate,
which is independent of $m$:
\[
  \int \bigl|B^{(m)}\bigr|^2_{\widetilde{h}^{(m)}}
 \bigl(\epsilon_m^{N+1}|z|^{2\epsilon_m}+|z|^2\bigr)^{-1}
 \frac{\dvol_{g}}{|z|^2}
\leq
 A\int
 \bigl| F(\widetilde{h}^{(m)})
 \bigr|^2_{\widetilde{h}^{(m)},g_{m}}
 \cdot\dvol_{g_{m}}.
\]
Here $A$ denotes a constant independent of $m$.
Due to Lemma \ref{lem;05.8.25.10},
there exist $v^{(m)}$ such that
the following inequalities hold
for some positive constant $A'$:
\[
 \delbar\del v^{(m)}
=\bigl|B^{(m)}\bigr|_{\widetilde{h}^{(m)}}
 \frac{dz \cdot d\zbar}{|z|^2},
\quad\quad
 \bigl|v^{(m)}(z)\bigr|
\leq
 A'\cdot \bigl\|F(\widetilde{h}^{(m)})
 \bigr\|_{\widetilde{h}^{(m)},g_m}
\]
Then we have
$ \Delta_g\bigl(
\log \tr \widetilde{t}^{(m)}-v^{(m)}
 \bigr)\leq 0$ on $U-\{Q_{\infty}\}$.
Since $s^{(m)}$ and $s^{(m)\,-1}$ are bounded
on $C-D$,
$\log \tr s^{(m)}$ is bounded on $C-D$.
Hence, $ \Delta_g\bigl(
\log \tr \widetilde{t}^{(m)}-v^{(m)}
 \bigr)\leq 0$ holds 
on $U$ as distributions.
(See Lemma 2.2 of \cite{s2}, for example.)
Therefore, we obtain the following:
\[
 \log \tr\widetilde{t}^{(m)}(Q_{m})
-v^{(m)}(Q_{m})
\leq
A''\cdot
 \int_{\Delta(Q_{m},1/2)}
 \Bigl(
 \log \tr\widetilde{t}^{(m)}-v^{(m)}
\Bigr)\cdot\dvol_g.
\]
Here $A''$ denotes a positive constant
independent of $m$.
Then we obtain the following inequalities,
for some positive constants $C_i$ $(i=1,2)$
which are independent of $m$:
\[
 \log (9/5)\leq
\log \tr\widetilde{t}^{(m)}(Q_{m})
\leq
C_1\cdot\int_{\Delta(Q_{m},1/2)}
 \log\tr\widetilde{t}^{(m)}
 \cdot\dvol_g
+C_2.
\]
Recall that $\log \tr \widetilde{t}^{(m)}$ are
uniformly bounded from above.
Therefore there exists a positive constant $C_3$
such that the following holds for any sufficiently
large $m$:
\[
 \int_{\Delta(Q_{m},1/2)}
 -\min (0,\log\tr\widetilde{t}^{(m)})
 \cdot\dvol_g
\leq C_3.
\]
Due to Fatou's lemma,
we obtain the following:
\[
 \int_{\Delta(Q_{\infty},1/2)}
 -\min\bigl(0,\log \tr \widetilde{t}^{(\infty)}\bigr)
 \cdot\dvol_{g}\leq C_3.
\]
It means $\widetilde{t}^{(\infty)}$
is not constantly $0$ on $\Delta(Q_{\infty},1/2)$.
In all,
we can conclude that
$\widetilde{t}^{(\infty)}$ is a positive constant multiplication.
Thus the proof of Lemma \ref{lem;06.8.21.1}
is finished.
\hfill\qed

\vspace{.1in}

Let 
$\bigl\{\widetilde{t}^{(m')}\bigr\}$
be a subsequence as in Lemma \ref{lem;06.8.21.1}.
It is almost everywhere convergent to some constant multiplication.
Then we obtain the convergence of
$\bigl\{\det \widetilde{t}^{(m')}
=b_{m'}^{\rank E}\cdot\id_{\det(E)}\bigr\}$
to a positive constant multiplication,
i.e.,
$\{b_{m'}\}$ is convergent to a positive constant.
It means the uniform boundedness of
$\{s^{(m')}\}$ with respect to $h_0^{(m')}$.

\subsection{Construction of maps}

By assumption,
we are given $C^{1}$-isometries
$\Phi_m':(E_m,h_m)\lrarr (E',h')$
for which
$\bigl\{(E_m,\delbar_{E_m},\theta_m)\bigr\}$
converges to $(E',\delbar_{E'},\theta')$.
By modifying them,
we would like to construct the maps
$\Psi_m':\prolongg{\vecc}{E_m}\lrarr \prolongg{\vecc}{E'}$
for which a subsequence of
$\bigl\{(\prolongg{\vecc}{E}_{m\ast},\theta_{m})\bigr\}$
converges to 
$(\prolongg{\vecc}{E'}_{\ast},\theta')$.
The argument is essentially same as 
that in Subsections
\ref{subsection;06.8.20.30}--\ref{subsection;05.9.6.30}.

We put $V_P^{\ast}:=V_P-\{P\}$.
We will shrink $V_P$ in the following argument
if it is necessary.
We may assume that 
Assumption \ref{assumption;05.9.2.4}
 is satisfied on $V_P$  for any $m<\infty$,
and that the constants are independent of $m$.
We have the convergence
$\bigl\{(\Par(\prolongg{\vecc}{E}_{m\ast},P),\gminim)\bigr\}$
to $\bigl(\Par(\prolongg{\vecc}{E}_{\infty\ast},P),\gminim\bigr)$.
Take $\eta>0$,
and we may assume that 
Assumption \ref{assumption;06.8.10.1} is satisfied
on $V_P$ for any $m<\infty$, after going to a subsequence.
By applying Proposition \ref{prop;05.8.28.50}
to harmonic bundles 
$(E_m,\delbar_m,\theta_m,h^{(m)}_0)_{|V_P^{\ast}}$,
we obtain holomorphic sections 
$F_1^{(m)},\ldots, F_r^{(m)}$
of $\prolongg{\vecc}{E}_m$ on $V_P$
with numbers $b_1^{(m)},\ldots,b_r^{(m)}$
as in Proposition \ref{prop;05.8.28.50}.
We may assume $b_i^{(m)}$
are independent of the choice of $m$,
which are denoted by $b_i$.
For $b\in \Par(\prolongg{\vecc}{E_{\infty}},P)$,
we put $\bbar(m):=
 \max\{a\in\Par(\prolongg{\vecc}{E_m})\,|\,
 |a-b|<\eta_0\}$.
We put $\widetilde{c}(m):=
 \sum_{a\in \Par(\prolongg{\vecc}{E_m},P)}
 a\cdot \gminim(a)$.
Because of the uniform boundedness of $s^{(m)}$,
we obtain
$|F_i^{(m)}|_{h^{(m)}}\leq
 C\cdot |z|^{-\bbar_i(m)}(-\log|z|)^N$
and
$C_1\cdot |z|^{-\widetilde{c}(m)}
 \leq\bigl|\bigwedge_{i=1}^rF_i^{(m)}
 \bigr|_{h^{(m)}}
 \leq C_2\cdot |z|^{-\widetilde{c}(m)}$,
where the constants are independent of $m$.
After going to a subsequence,
we may assume that
$\bigl\{ \Phi'_m(F_i^{(m')}) \bigr\}$
are convergent weakly in $L_1^{p}$ locally on  $V_P^{\ast}$.
The limits are denoted by $F_i'$,
which are holomorphic with respect to $\delbar_{E'}$.
We have 
$|F_i'|_{h'}\leq C\cdot |z|^{-b_i}(-\log|z|)^N$
and 
$C_1\cdot |z|^{-\widetilde{c}}
 \leq\bigl|\bigwedge_{i=1}^rF_i'\bigr|_{h^{(m)}}
 \leq C_2\cdot |z|^{-\widetilde{c}}$,
where $\widetilde{c}:=
 \sum_{b\in \Par(\prolongg{\vecc}{E}_{\infty},P)}
 \gminim(b)\cdot b$.
By the same argument as the proof of
Lemma \ref{lem;06.8.21.30},
we obtain that
$F_1',\ldots,F_r'$ gives a frame
of $\prolongg{\vecc}{E'}$ around $P$
which is compatible with the parabolic structure.
(In particular, we obtain
 $\Par(\prolongg{\vecc}{E'},P)
=\Par(\prolongg{\vecc}{E_{\infty}},P)$).

We obtain the holomorphic morphism
$\lefttop{P}\Psi_m':\prolongg{\vecc}{E_m}_{|V_P}
\lrarr \prolongg{\vecc}{E'}_{|V_P}$
by the correspondence 
$\lefttop{P}\Psi_m'(F_i^{(m)})=F_i'$.
By our construction,
(i)
$\lefttop{P}\Psi_m'-\Phi'_{m|V_P^{\ast}}$
converges to $0$ weakly in $L_1^p$ locally on $V_P^{\ast}$,
(ii)
$\lefttop{P}\Psi_m'(\theta_m)-\theta'$
 converges to $0$ on $V_P$
 as holomorphic sections 
of $\End(\prolongg{\vecc}{E'})\otimes\Omega^{1,0}(\log P)$,
(iii)
the parabolic filtrations of
$\prolongg{\vecc}{E_m}_{|P}$
converges to the parabolic filtration
of $\prolongg{\vecc}{E'}_{|P}$
via $\lefttop{P}\Psi_m'$.
Then, we construct $\Psi_m'$
similarly to (\ref{eq;06.8.21.2}),
which gives the convergence
of $\bigl\{(\prolongg{\vecc}{E_{m\ast}},\theta_m)\bigr\}$
to $(\prolongg{\vecc}{E'}_{\ast},\theta')$.
\hfill\qed

%% file: 9.tex
\section{The Surface Case}

Let $X$ be a smooth irreducible projective surface,
and $D$ be a simple normal crossing divisor of $X$.
Let $L$ be an ample line bundle,
and $\omega$ be a Kahler form representing $c_1(L)$.

\begin{thm}
\label{thm;05.9.8.150}
Let $\bigl(\prolongg{\vecc}{E},\vecF,\theta\bigr)$
be a $\mu_L$-stable $\vecc$-parabolic Higgs bundle on $(X,D)$.
Assume that the characteristic numbers vanish:
\[
 \pardeg_{L}(\prolongg{\vecc}{E},\vecF)
=\int_X\parch_{2}(\prolongg{\vecc}{E},\vecF)=0.
\]
Then there exists 
a pluri-harmonic metric $h$ of
$(E,\theta)=(\prolongg{\vecc}{E},\theta)_{|X-D}$
which is adapted to the parabolic structure.
\end{thm}
\pf
We may and will assume 
$c_i\not\in\Par(\prolongg{\vecc}{E},\vecF,i)$.
We take a sequence $\{\epsilonbar_m\}$
converging to $0$,
such that $\epsilonbar_m=N_m^{-1}$
for some integers $N_m$
and that $\epsilonbar_m<\gap(\prolongg{\vecc}{E},\vecF)/100\rank(E)$.
We take the perturbation of parabolic structures
$\vecF^{(\epsilonbar_m)}$
as in Section \ref{section;05.7.30.15}.
We put $\epsilon_m=\epsilonbar_m/100$,
and we take the Kahler metrics $\omega_{\epsilon_m}$ of $X-D$
as in Subsection \ref{section;05.8.23.5}.
For simplicity of the notation,
we denote them by 
$\vecF^{(m)}$ and $\omega^{(m)}$, respectively.
We may assume that 
 $\bigl(\prolongg{\vecc}{E},\vecF^{(m)}\bigr)$
 are $\mu_{L}$-stable.

Due to Corollary \ref{cor;06.8.5.30},
we have already known
$\parchern_1(\prolongg{\vecc}{E},\vecF)
=\parchern_1(\prolongg{\vecc}{E},\vecF^{(m)})=0$.
Thus, we can take a pluri-harmonic metric
$h_{\det E}$ of $\det(E)$ adapted to the parabolic structure.
Due to Proposition \ref{prop;05.7.30.10},
we have the Hermitian-Einstein metric $h^{(m)}$ of
$(E,\delbar_E,\theta)$ with respect to $\omega^{(m)}$
such that
 $\Lambda_{\omega^{(m)}} F(h^{(m)})
 =\tr F(h^{(m)})=0$
and $\det(h^{(m)})=h_{\det E}$,
which is adapted to the parabolic structure
$(\prolongg{\vecc}{E},\vecF^{(m)})$.
We remark that
the sequence of the $L^2$-norms
$\|F(h^{(m)})\|_{h^{(m)},\omega^{(m)}}$
of $F(h^{(m)})$ with respect to $h^{(m)}$ and $\omega^{(m)}$
converges to $0$ in $m\to\infty$,
because of the relation $\|F(h_m)\|_{h^{(m)},\omega^{(m)}}^2\!\!=\!\!
 C\cdot \parch_{2,L}(\prolongg{\vecc}{E},\vecF^{(m)})$
for some non-zero constant $C$.
We will show the local convergence 
of the sequence $\bigl\{(E,\delbar_E,\theta,h^{(m)})\bigr\}$
on $X-D$.

\subsection{Local convergence}

In the following argument,
$B_i$ will denote positive constants
which are independent of $m$.
We use the notation $\|\rho\|_{h',\omega'}$
to denote the $L^2$-norm of 
a section $\rho$ of $E'\otimes\Omega^{i,j}$
or $\End(E')\otimes\Omega^{i,j}$,
where $h'$ and $\omega'$
denote metrics of a vector bundle $E'$
and a base space.
On the other hand,
$|\rho|_{h',\omega'}$ denotes the norms
at fibers.

Let $P$ be any point of $X-D$.
We take a holomorphic coordinate $(U,z_1,z_2)$
around $P$ such that $z_i(P)=0$
and that $\omega_{|P}=\sum dz_i\cdot d\zbar_i$
on the tangent space at $P$.
We have the expression
$\theta=\sum f_i\cdot dz_i$.

Let $\eta$ be a positive number.
If $m$ is sufficiently large,
we have
$\|F(h^{(m)})\|_{\omega^{(m)},h^{(m)}} 
 \leq \eta$.
Due to Lemma \ref{lem;05.8.25.4},
there exists a constant $B_1$,
such that $B_1^{-1}\cdot |f_i|_{h^{(m)}}\leq \eta$.
Take a large number $B_2>B_1$,
and we put $w_i:=B_2 \cdot z_i$,
$\Ytilde(T):=\bigl\{(w_1,w_2)\,\big|\,\sum |w_i|^2\leq T\bigr\}$,
$\gtilde:=\sum dw_i\cdot d\wbar_i$
and $\omegatilde^{(m)}:=B_2^2\cdot\omega^{(m)}$.
Then, we obtain the following:
\[
 \bigl\|R(h^{(m)})_{|\Ytilde(1)}\bigr\|_{h^{(m)},\gtilde}
\leq
 \bigl\|F(h^{(m)})_{|\Ytilde(1)}\bigr\|_{h^{(m)},\gtilde}
+\bigl\|[\theta,\theta^{\dagger}_{h^{(m)}}]_{|\Ytilde(1)}
 \bigr\|_{h^{(m)},\gtilde}
\leq B_3\cdot\eta
\]
Let $d^{\ast}$ denote the formal adjoint of
the exterior derivative $d$ on $\Ytilde(1)$
with respect to $\gtilde$.
If $\eta$ is sufficiently small,
we can apply Uhlenbeck's theorem (\cite{u1}).
Namely,
we can take an orthonormal frame $\vecv_m$ of
$\bigl(E,h^{(m)}\bigr)_{|\Ytilde(1)}$
such that the connection form $A_m$ of
$\delbar_E+\del_{E,h^{(m)}}$ with respect to $\vecv_m$
satisfies the conditions:
\begin{description}
\item[(i)] $d^{\ast}A_m=0$,
\item[(ii)] $\|A_m\|_{L_1^p,\gtilde}\leq 
  C(p)\cdot \|dA_m+A_m\wedge A_m
  \|_{L^p,\gtilde}$
 ($p\geq 2$),
 where $C(p)$ denotes the constant
 depending only on $p$.
\end{description}
By our choice of $B_2$,
we also have the following:
\begin{description}
\item[(iii)]
 Let $\Pi^{(m)}$ denote the orthogonal projection of
 $\Omega^2$ onto
 the self-dual part part with respect to
 $\omegatilde_m$.
 Then, 
 $\bigl|
\Pi^{(m)}\bigl(dA_m+A_m\wedge A_m\bigr)
 \bigr|_{\omegatilde^{(m)}}
 \leq B_4\eta$
because of
$\Lambda_{\omegatilde}R(h^{(m)})
=\Lambda_{\omegatilde}
 [\theta,\theta^{\dagger}_{h^{(m)}}]$.
\end{description}
From (i) and (iii),
we have
$\bigl|\bigl(
 d^{\ast}+\Pi^{(m)}\circ d
\bigr)(A_m)
+\Pi^{(m)}(A_m\wedge A_m)
 \bigr|
\leq B_5$.
If $B_2$ and $m$ are  sufficiently large,
$\omegatilde^{(m)}$ and
$\gtilde$ are sufficiently close.
Recall that $d^{\ast}+\Pi\circ d$
is elliptic,
where $\Pi$ denotes the orthogonal projection
of $\Omega^{2}$ onto the self-dual part
with respect to $\gtilde$.
Using the boot strapping argument
of Donaldson for Corollary 23 in \cite{don3},
we obtain that the $L_1^p$-norm of $A_m$
on $\Ytilde(T)$ $(T<1)$
is dominated by a constant $B_{6}$.
Let $\Theta_m$ be determined by
$\theta(\vecv_m)=\vecv_m\cdot \Theta_m$.
The sup norm of $\Theta_m$ with respect to $\gtilde$
is small, due to our choice of $B_2$.
We also obtain the $L_1^p$-bound of $\Theta_m$
because of $\delbar\Theta_m+[A_m^{0,1},\Theta_m]=0$,
where $A_m^{0,1}$ denotes the $(0,1)$-part of $A_m$.

\begin{lem}
\label{lem;06.8.20.10}
After going to a subsequence,
$\bigl\{
 \bigl(E,\delbar_E,h^{(m)},\theta\bigr)\,\big|\,m\in I
 \bigr\}$ converges to 
a tame harmonic bundle
$(E_{\infty},\delbar_{\infty},h_{\infty},\theta_{\infty})$
weakly in $L_2^p$ locally on $X-D$.
\end{lem}
\pf
Due to the above arguments,
we can take a locally finite covering
$\bigl\{(U_{\alpha},z^{(\alpha)}_1,z^{(\alpha)}_2)
 \,\big|\,\alpha\in \Gamma\bigr\}$ 
of $X-D$
and the numbers $\bigl\{m(\alpha)\,\big|\,\alpha\in\Gamma\bigr\}$
with the following property:
\begin{itemize}
\item 
Each $U_{\alpha}$ is relatively compact in $X-D$.
\item
For any $m\geq m(\alpha)$,
we have orthonormal frames $\vecv_{\alpha,m}$ of $(E,h^{(m)})$
on $U_{\alpha}$ such that 
the $L_1^p$-norms of $A_{\alpha,m}$
are sufficiently small
with respect to the metrics $\sum dz_j^{(\alpha)}\cdot d\bar{z}_j^{(\alpha)}$
independently of $m$,
where $A_{\alpha,m}$ denote the connection forms 
of $(\del_{E,h^{(m)}}+\delbar_{E})$
 with respect to $\vecv_{\alpha,m}$.
\item
Let $\Theta_{\alpha,m}$ be the matrix valued $(1,0)$-forms
given by $\theta\cdot\vecv_{\alpha,m}=
 \vecv_{\alpha,m}\cdot\Theta_{\alpha,m}$.
Then the $L_1^p$-norms of $\Theta_{\alpha,m}$
are sufficiently small
with respect to $\sum dz_j^{(\alpha)}\cdot d\bar{z}_j^{(\alpha)}$,
independently of $m$.
\end{itemize}
Let $g_{\beta,\alpha,m}$ be the unitary transformation on
$U_{\alpha}\cap U_{\beta}$ determined by
$\vecv_{\alpha,m}=\vecv_{\beta,m}\cdot g_{\beta,\alpha,m}$.
Once $\alpha$ and $\beta$ are fixed,
the $L_2^p$-norms of $g_{\beta,\alpha,m}$ are bounded
independently of $m$.
By a standard argument,
we can take a subsequence $I\subset \{m\}$
such that
the sequences
 $\bigl\{ A_{\alpha,m}\,\big|\,m\in I\bigr\}$,
 $\bigl\{\Theta_{\alpha,m}\,\big|\,m\in I\bigr\}$
 are weakly $L_1^p$-convergent for each $\alpha$,
and that the sequence 
 $\bigl\{g_{\alpha,\beta,m}\,\big|\,m\in I\bigr\}$
is weakly  $L_2^p$-convergent
 for each $(\alpha,\beta)$.
Then, we obtain the limit Higgs bundle
$(E_{\infty},\delbar_{\infty},\theta_{\infty})$
with the metric $h_{\infty}$ on $X-D$.
From the convergence
$\bigl\|F(h^{(m)})\bigr\|_{L^2,h^{(m)},\omega^{(m)}}\to 0$,
we obtain $\bigl\|F(h_{\infty})\bigr\|_{L^2,h_{\infty},\omega}=0$,
and hence $(E_{\infty},\delbar_{\infty},\theta_{\infty},h_{\infty})$
is a harmonic bundle.
By using the argument of Uhlenbeck \cite{u1},
we obtain locally $L_2^p$-isometries
$\Phi_m:(E,h^{(m)})\lrarr (E_{\infty},h_{\infty})$,
via which $\bigl\{(E,\delbar_E,\theta,h^{(m)})\bigr\}$
converges to $(E_{\infty},\delbar_{\infty},\theta_{\infty},h_{\infty})$
weakly in $L_2^p$ locally on $X-D$.
Since we have $\det(t-\theta)=\det(t-\theta_{\infty})$
by construction,
the tameness of
$(E_{\infty},\delbar_{E_{\infty}},h_{\infty},\theta_{\infty})$
follows.
Thus, Lemma \ref{lem;06.8.20.10} is proved.
\hfill\qed

\vspace{.1in}

We obtain the associated parabolic Higgs bundle
$\bigl(\prolongg{\vecc}{E_{\infty}},\vecF_{\infty},\theta_{\infty}\bigr)$.
We would like to show that it is isomorphic to the given
parabolic Higgs bundle 
$(\prolongg{\vecc}{E},\vecF,\theta)$.
For that purpose,
we have only to show the existence of a non-trivial morphism
$f:(\prolongg{\vecc}{E},\vecF,\theta)
\lrarr
 \bigl(\prolongg{\vecc}{E_{\infty}},\vecF_{\infty},
 \theta_{\infty}\bigr)$,
because of the $\mu_L$-stability of
$(\prolongg{\vecc}{E},\vecF,\theta)$
and the $\mu_L$-polystability of
$(\prolongg{\vecc}{E_{\infty}},\vecF_{\infty},\theta_{\infty})$.
Moreover,
we have only to show 
the existence of a non-trivial map
$f_C:(\prolongg{\vecc}E_{\infty},\vecF_{\infty},\theta_{\infty})_{|C}
\lrarr (\prolongg{\vecc}{E},\vecF,\theta)_{|C}$
for a sufficiently ample generic curve $C\subset X$,
due to Lemma \ref{lem;06.8.22.2}.
So we show that such $f_C$ exists
for almost all $C$, in the next subsections.

\subsection{Selection of a curve}

Let $L^N$ be sufficiently ample.
We put $\nbigv:=H^0(X,L^n)$.
For any $s\in \nbigv$, we put $X_s:=s^{-1}(0)$.
Recall Mehta-Ramanathan type theorem
(Proposition \ref{prop;06.8.12.15}),
and let $\nbigu$ denote the Zariski open subset of $\nbigv$
which consists of the points $s$
with the properties:
(i)
$X_s$ is smooth,
and $X_s\cap D$ is a simple normal crossing divisor,
(ii)
$(\prolongg{\vecc}{E},\vecF,\theta)_{|X_s}$
is $\mu_L$-stable.

We will use the notation $X_s^{\ast}:=X_s\setminus D$
and $D_s:=X_s\cap D$.
We have the metric $\omega^{(m)}_{s}$ of $X_s^{\ast}$,
induced by $\omega^{(m)}$.
The induced volume form of $X_s^{\ast}$
is denoted by $\vol^{(m)}_{s}$.
We put $(\prolongg{\vecc}{E}_{s},\vecF_{s}^{(m)},\theta_s):=
 (\prolongg{\vecc}{E},\vecF^{(m)},\theta)_{|X_s}$.
We have the metric $h^{(m)}_s:=h^{(m)}_{|X_s^{\ast}}$
of $E_s:=E_{|X_s^{\ast}}$.
Since there exists $m_0$ such that
$(\prolongg{\vecc}{E}_s,\vecF^{(m)}_s,\theta_s)$ is
stable for any point $s\in\nbigu$ and for any $m\geq m_0$,
we have the harmonic metric 
$h_{s,0}^{(m)}$ of
$({E}_s,\theta_s)$ adapted to the parabolic structure
$\vecF^{(m)}_s$
with $\det h_{s,0}^{(m)}=h_{\det E\,|\,X_s^{\ast}}$
(Proposition \ref{prop;05.9.160}).
Let $u_s^{(m)}$ be the endomorphism
of $E_{|X_s^{\ast}}$ determined by
$h^{(m)}_{s}
=h_{s,0}^{(m)}\cdot u_s^{(m)}$.
We put $\nbigd_s:=(\delbar_E+\theta)_{|X_s}$.

\begin{lem}
 \label{lem;06.1.18.20}
For almost every $s\in \nbigu$,
the following holds:
\begin{itemize}
\item
We have the following convergence in $m\to\infty$:
\begin{equation}
 \label{eq;06.1.18.7}
  \bigl\|
 F(h^{(m)}_{s})
 \bigr\|_{h^{(m)}_s,\omega^{(m)}_{s}}
\lrarr 0.
\end{equation}
\item
 For each $m$,
 we have the finiteness:
\begin{equation} 
\label{eq;06.1.26.50}
\bigl\|\nbigd_{s}
 u_s^{(m)}
 \bigr\|_{h_{s,0}^{(m)},\omega^{(m)}_{s}}<\infty.
\end{equation}
\end{itemize}
Let $\widetilde{\nbigu}$ denote the set of $s$
for which both of {\rm(\ref{eq;06.1.18.7})} 
and {\rm(\ref{eq;06.1.26.50})} hold.
\end{lem}
\pf
Let us discuss the condition (\ref{eq;06.1.18.7}).
Let us fix $s_1\in \nbigu$.
We take generic $s_i\in\nbigu$ $(i=2,3)$,
i.e.,
$X_{s_1}$ is transversal with $X_{s_i}$ $(i=2,3)$
and $X_{s_1}\cap X_{s_2}\cap X_{s_3}=\emptyset$.
Take open subsets $W_i^{(j)}$ $(j=1,2,\,\,i=2,3)$ such that
(i) $X_{s_1}\cap X_{s_i}
\subset W_i^{(1)} \subset W_i^{(2)}$,
(ii) $W_i^{(1)}$ is relatively compact in $W_i^{(2)}$.
Take an open neighbourhood $U_1$ of $s_1$,
which is relatively compact in $\nbigu$,
such that 
 $X_{s}$ is transversal with $X_{s_i}$ $(i=2,3)$
 and $X_{s}\cap X_{s_i}\subset W_i^{(1)}$
 for any $s\in U_1$.

Take $T>0$, and
we put $\nbigb:=\bigl\{z\in\cnum\,\big|\,|z|\leq T\bigr\}$.
Let $q_i$ denote the projection of
$X\times U_1\times\nbigb$ onto the $i$-th component.
We put
$\nbigz_2:=
 \bigl\{
 (x,s,t)\in X\times U_1\times\proj^1\,\big|\,
 (ts_2+(1-t)s)(x)=0
 \bigr\}$.
The fiber over $s\in U_1$
via $q_{2\,|\,\nbigz_2}$ is the closed region of
the Lefschetz pencil of $s$ and $s_2$.

We fix any Kahler forms $\omega_{U_1}$ and $\omega_{\nbigb}$
of $U_1$ and $\nbigb$.
The induced volume forms are denoted by $\vol_{U_1}$
and $\vol_{\nbigb}$.
Then we have the following convergence
in $m\to\infty$:
\[
 \int_{\nbigz_2}
 q_1^{\ast}\Bigl(
\bigl|F(h^{(m)})\bigr|^2_{h^{(m)},\omega^{(m)}}\cdot
 \dvol_{\omega^{(m)}}
 \Bigr)
\cdot\dvol_{U_1}
\lrarr 0.
\]
We put $\nbigz'_2:=\nbigz_2\setminus
 q_1^{-1}(W_2^{(2)})$.
Then the following convergence is obtained,
in particular:
\begin{equation}
 \label{eq;06.1.18.10}
 \int_{\nbigz_2'}
 q_1^{\ast} \Bigl(
 \bigl|F(h^{(m)})\bigr|^2_{h^{(m)},\omega^{(m)}}
\cdot\dvol_{\omega^{(m)}}
 \Bigr)\cdot\dvol_{U_1}\lrarr 0.
\end{equation}

Let $\psi:\nbigz_2\lrarr U_1\times\nbigb$
denote the projection.
For $(s,t)\in U_1\times\nbigb$,
we put $X_{(s,t)}:=\psi^{-1}(s,t)=
\bigl(ts_2+(1-t)s\bigr)^{-1}(0)=X_{ts_2+(1-t)s}$.
On $X_{(s,t)}$, we have 
the induced Kahler form $\omega^{(m)}_{(s,t)}$,
the induced volume forms $\dvol^{(m)}_{(s,t)}$
and the hermitian metric
$h_{(s,t)}^{(m)}:=h^{(m)}_{|X_{(s,t)}}$.
The family
$\{\dvol^{(m)}_{(s,t)}\,|\,(s,t)\in U_1\times \nbigb\}$
gives the $C^{\infty}$-relative volume form
$\dvol_{\nbigz_2'/U_1\times\nbigb}^{(m)}$
of $\nbigz_2'\lrarr U_1\times\nbigb$.
There exists a constant $A$
such that the following holds on $\nbigz_2'$:
\begin{multline}
 A\cdot q_1^{\ast}\Bigl(
 \bigl|F(h^{(m)})\bigr|^2_{h^{(m)},\omega^{(m)}}
 \dvol_{\omega^{(m)}}
\Bigr)\dvol_{U_1} \\
\geq
 \bigl|
 F(h^{(m)}_{(s,t)})
 \bigr|^2_{h^{(m)}_{(s,t)},\omega^{(m)}_{(s,t)}}
 \dvol_{\nbigz'_2/U_1\times\nbigb}^{(m)}\cdot\dvol_{\nbigb}
 \dvol_{U_1}
\end{multline}
Therefore,
we obtain the following convergence
for almost every $(s,t)\in U_1\times\nbigb$,
from (\ref{eq;06.1.18.10}):
\begin{equation}
 \label{eq;06.1.18.11}
 \int_{X_{(s,t)}^{\ast}\setminus W_2^{(2)}}\bigl|
 F(h^{(m)}_{(s,t)})
 \bigr|^2_{h^{(m)}_{(s,t)},\omega^{(m)}_{(s,t)}}
 \dvol^{(m)}_{(s,t)}
\lrarr 0.
\end{equation}
Let $\nbigs$ denote the set of the points
$(s,t)\in U_1\times \nbigb$
such that the above convergence (\ref{eq;06.1.18.11})
does not hold.
The measure of $\nbigs$ is $0$
with respect to $\dvol_{U_1}\times\dvol_{\nbigb}$.

Let $\nbigj:U_1\times \nbigb\lrarr \nbigv$ denote the map
given by $(s,t)\longmapsto ts_2+(1-t)s$.
We have the open subset
$\nbigj^{-1}(U_1)\subset U_1\times \nbigb$
and the measure of $\nbigs\cap \nbigj^{-1}(U_1)$ is $0$
with respect to $\dvol_{U_1}\cdot\dvol_{\nbigb}$.
We have
$\nbigs\cap \nbigj^{-1}(U_1)
=\nbigj^{-1}\bigl(\nbigj(\nbigs)\cap U_1\bigr)$,
and hence the measure of
$\nbigt(\nbigs)\cap U_1$ is $0$
with respect to $\omega_{U_1}$.
Namely, we have the following convergence
for almost every $s\in U_1$:
\[
 \int_{X_s^{\ast}\setminus W_2^{(2)}}
 \bigl|
 F(h^{(m)}_{s})
 \bigr|^2_{h^{(m)}_s,\omega^{(m)}_s}
 \cdot\dvol^{(m)}_s\lrarr 0.
\]
Similarly,
we can show the following convergence for almost every
$s\in U_1$:
\[
 \int_{X_s^{\ast}\setminus W_3^{(2)}}
 \bigl|
 F(h^{(m)}_{s})
 \bigr|^2_{h^{(m)}_s,\omega^{(m)}_s}
 \cdot\dvol^{(m)}_s\lrarr 0
\]
Then, we obtain that  
the condition (\ref{eq;06.1.18.7}) holds for
almost all $s\in \nbigu$.

The condition (\ref{eq;06.1.26.50}) can be discussed
similarly. We give only an outline.
Let $h_{in}^{(m)}$ be an initial metric
which was used for the construction of
$h^{(m)}$. (See the proof of Proposition \ref{prop;05.7.30.10}.)
We remark that $h_{in}^{(m)}$ and $h^{(m)}$
are mutually bounded.
Let $t^{(m)}$ be determined by
$h^{(m)}=h_{in}^{(m)}\cdot t^{(m)}$.
Then, we have
$\bigl\|\nbigd t^{(m)}\bigr\|_{\omega^{(m)},h^{(m)}}
<\infty$
due to Proposition \ref{prop;05.7.30.35}.
We put 
$h_{s,in}^{(m)}:=h^{(m)}_{in|X_s^{\ast}}$
and $t_s^{(m)}:=t^{(m)}_{|X_s^{\ast}}$
for $s\in\nbigu$.
By an above argument,
we obtain
$\bigl\| \nbigd_s t_s^{(m)}
 \bigr\|_{\omega^{(m)}_s,h_{s,in}^{(m)}}<\infty$
for almost all $s\in \nbigu$.
On the other hand,
let $\ttilde_s^{(m)}$ be determined by
$h_{s,0}^{(m)}=h_{s,in}^{(m)}\cdot\ttilde_s^{(m)}$.
We can use $h^{(m)}_{s,in}$
as the initial metric for the construction of
$h_{s,0}^{(m)}$.
Hence, we have
$\bigl\|\nbigd_s \ttilde^{(m)}_{s}
 \bigr\|_{\omega^{(m)}_s,h^{(m)}_{s,in}}<\infty$.
Since we have
$u_s^{(m)}=\ttilde_s^{(m)\,-1}\cdot t_s^{(m)}$,
the condition (\ref{eq;06.1.26.50}) is satisfied
for almost $s\in\nbigu$.
Thus, the proof of Lemma \ref{lem;06.1.18.20}
is finished.
\hfill\qed

\subsection{End of the proof of Theorem \ref{thm;05.9.8.150}}

Let us finish the proof of Theorem \ref{thm;05.9.8.150}.
Take $s\in \nbigutilde$,
and we put $C=X_s$.
We have the convergence
of $\bigl\{(E,\delbar_{E},\theta,h^{(m)})\bigr\}$
to $(E_{\infty},\delbar_{\infty},\theta_{\infty},h_{\infty})$
weakly in $L_2^p$ locally on $X-D$
via isometries
$\Phi_m:(E,h_m)\lrarr (E_{\infty},h_{\infty})$.
The restriction of $\Phi_m$ to $C\setminus D$
induce the $C^1$-convergence of
$\bigl\{(E,\delbar,\theta,h^{(m)})_{|C\setminus D}\bigr\}$
to $(E_{\infty},\delbar_{\infty},\theta_{\infty},h_{\infty})
 _{|C\setminus D}$.
By using Proposition \ref{prop;06.8.21.5},
we obtain the convergence of
$\bigl\{
(\prolongg{\vecc}{E},\vecF^{(m')},\theta)_{|C}
\bigr\}$
to $(\prolongg{\vecc}{E_{\infty}},
 \vecF_{\infty},\theta_{\infty})_{|C}$
weakly in $L_1^p$ on $C$
for some subsequence.
We also have the convergence of
$\bigl\{(\prolongg{\vecc}{E},
 \vecF^{(m)},\theta)_{|C}\bigr\}$
to $(\prolongg{\vecc}{E},\vecF,\theta)_{|C}$.
Due to Corollary \ref{cor;05.8.29.301},
we obtain the desired non-trivial map
$f_C:(\prolongg{\vecc}E_{\infty},\vecF_{\infty},\theta_{\infty})_{|C}
\lrarr (\prolongg{\vecc}{E},\vecF,\theta)_{|C}$.
Thus we are done.
\hfill\qed

\section{The Higher Dimensional Case}

Now the main existence theorem is given.
\begin{thm} 
 \label{thm;04.10.24.1}
Let $X$ be a smooth 
irreducible projective variety over the complex number field
of dimension $n$.
Let $D=\bigcup_i D_i$ be a simple normal crossing divisor of $X$.
Let $L$ be an ample line bundle on $X$.
Let $\bigl(\vecE_{\ast},\theta\bigr)$ be
a $\mu_L$-stable regular filtered Higgs bundle
with 
$\pardeg_L(\vecE_{\ast})=\int_X\parch_{2,L}(\vecE_{\ast})=0$.
We put $E:=\vecE_{|X-D}$.
Then there exists a pluri-harmonic metric $h$
of $(E,\delbar_E,\theta)$,
which is adapted to the parabolic structure.
Such a metric is unique up to constant multiplication.
\end{thm}
\pf
We may assume that $D$ is ample.
We can also assume that $L$ is sufficiently ample
as in Proposition \ref{prop;06.8.12.15}.
The uniqueness follows from the more general result
(Proposition \ref{prop;05.9.8.300}).
We use an induction on $n=\dim X$.
We have already known the existence for $n=2$
(Theorem \ref{thm;05.9.8.150}).

Let $(\vecE_{\ast},\theta)$ be a regular filtered Higgs bundle on $(X,D)$.
Assume that it is stable with
$\pardeg_L(\vecE_{\ast})=\int_X\parch_{2,L}(\vecE_{\ast})=0$.
For any element $s\in \proj:=\proj\bigl(H^0(X,L)^{\lor}\bigr)$
determines the hypersurface $Y_s=\bigl\{x\in X\,\big|\,s(x)=0\bigr\}$.
The subset $\nbigx_L\subset X\times\proj$
is given by $\nbigx_L:=\{(x,s)\,\big|\,x\in Y_s\}$.
Let $U$ be a Zariski open subset of $\proj$
which consists of $s\in \proj$ such that
$(\vecE_{\ast},\theta)_{|Y_s}$ is $\mu_L$-stable.
Since $L$ is assumed to be sufficiently ample,
$U$ is not empty
(Proposition \ref{prop;06.8.12.15}).
The image $W$ of the naturally defined map $\nbigx_L\times_{\proj}U\lrarr X$
is Zariski open in $X$.
In fact, $X-W$ consists of, at most, finite points of $X$
due to the ampleness of $L$.

Let $s$ be any element of $U$.
We have a pluri-harmonic metric $h_s$ of $(E,\theta)_{|Y_s}$,
which is adapted to the induced parabolic structure,
due to the hypothesis of the induction.

Let $s_i$ $(i=1,2)$ be elements of $U$
such that $Y_{s_1}$ and $Y_{s_2}$ are transversal
and that $Y_{s_1,s_2}:=Y_{s_1}\cap Y_{s_2}$ is transversal to $D$.
We remark that $\dim Y_{s_1}\cap Y_{s_2}\geq 1$.
We may also assume that $(\prolongg{\vecc}{E},\theta)_{|Y_{s_1,s_2}}$
is $\mu_L$-stable (Proposition \ref{prop;06.8.12.15}).
Hence $h_{s_1\,|\,Y_{s_1,s_2}}$
and $h_{s_2\,|\,Y_{s_1,s_2}}$ are same
up to constant multiplication.
Then, we obtain the metric $h$ of $E_{|X-(D\cup W)}$
such that $h_{|Y_s}=h_s$.

Let $P$ be any point of $X-(D\cup W)$.
We can take a coordinate neighbourhood
$(U_P,z_1,\ldots,z_n)$
around $P$
such that each $z_i^{-1}(0)$ is a part of
some $Y_{s}$
and that $U_P\subset X-(D\cup W)$.
In the following, we shrink $U_P$, if necessary.
Since the restriction of $h$ to $z_i^{-1}(0)$
is pluri-harmonic,
we obtain the boundedness of $\theta$
and $\theta^{\dagger}$
with respect to $h$ around $P$.
(See Proposition \ref{prop;05.8.25.2},
 for example.)

For any $Q\in U_P$,
let us take a path $\gamma$
connecting $P$ and $Q$,
which is contained in some $Y_s$.
Then, the parallel transport $\Pi_{P,Q}:E_{|P}\lrarr E_{|Q}$
is induced from the flat connection
associated to the harmonic bundle
$(E,\delbar_E,\theta)_{|Y_s}$ with $h_{|Y_s}$.
The map $\Pi_{P,Q}$ is independent of the choice
of $\gamma$ and $Y_s$.
From the frame of $E_{|P}$,
we obtain the frame $\vecv=(v_1,\ldots,v_r)$
of $E_{|U_P}$.
The trivialization gives the structure of flat bundle
to $E_{|U_P}$.
For the distinction,
we use the notation $(V,\nabla)$
to denote the obtained flat bundle.
The restriction of $h$, $\theta$ and $\theta^{\dagger}$
to $U_P$ are denoted by the same notation.
By the flat structure,
we can regard the metric $h$
as the map $\varphi_h:U_P\lrarr \GL(n)/U(n)$,
and $\theta+\theta^{\dagger}$ can be
regarded as the differential of the map.
Let $d_{\GL(n)/U_n}$ denote the invariant distance
of $\GL(n)/U_n$.
Due to the boundedness of $\theta+\theta^{\dagger}$
with respect to $h$,
there exists a constant $C$
such that 
$d_{\GL(n)/U(n)}\bigl(\varphi_h(\gamma(0)),
 \varphi_h(\gamma(1))\bigr)$
is less than $C$ times the length of $\gamma$
for any path $\gamma$ contained in some $Y_s$.
In particular, $h$ is a continuous metric of $V$.

Let $H$ be the hermitian-matrices valued function
whose $(i,j)$-th component is $h(v_i,v_j)$.
Let $\Theta=(\Theta_{i,j})$ and
$\Theta^{\dagger}=(\Theta^{\dagger}_{i,j})$
be determined by
$\theta v_i=\sum \Theta_{j,i}\cdot v_j$
and $\theta^{\dagger} v_i=\sum \Theta^{\dagger}_{j,i}\cdot v_j$.
We have 
$d\Hbar=\Hbar(\Theta+\Theta^{\dagger})/2$
and $\delbar\Theta+[\Theta^{\dagger},\Theta]=0$
for the point-wise partial derivatives,
which can be shown by considering the restriction
of $(E,\delbar_E,h,\theta)$ to hyper planes $\{z_i=a\}$.
The equality holds as distributions,
which follows from Fubini's theorem
and the boundedness of $\Hbar$,
$\Theta$ and $\Theta^{\dagger}$.
In particular, $H$ and $\Theta$
are locally $L_1^p$,
and hence $\Theta^{\dagger}$
is also locally $L_1^p$.
By a standard boot strapping argument,
we obtain that $H$, $\Theta$
and $\Theta^{\dagger}$ are $C^{\infty}$
functions.
In other words,
$h$ is a $C^{\infty}$-metric of $V$,
and $\theta^{\dagger}$ is a $C^{\infty}$-section
of $\End(V)\otimes\Omega^{0,1}$.
We also obtain that the $C^{\infty}$-structure
of $E$ and $V$ are same
because of $\delbar_E=d_V''-\theta^{\dagger}$,
where $d_V''$ denotes the $(0,1)$-part of $\nabla$.
Thus, we obtain that $h$ is a $C^{\infty}$-metric 
of $E_{|X-(D\cup W)}$.
The pluri-harmonicity of $h$ is easily obtained.
By using a similar elliptic regularity argument,
it can be shown that
$h$ gives the pluri-harmonic metric of $E$ on $X-D$.
Thus we obtain the tame harmonic bundle
$(E,\delbar_E,\theta,h)$ on $X-D$.

We take the prolongment $(\prolongg{\vecc}{E}(h)_{\ast},\theta)$
which is a parabolic Higgs bundle on $(X,D)$.
(See Section \ref{section;05.9.8.110}
 for the prolongment.)
\begin{lem}
There exists a closed subset $W'\subset D$
with the following properties:
\begin{itemize}
\item
The codimension of $W'$  in $X$ is larger than $2$.
\item
The identity of $E$ is extended to 
the holomorphic isomorphism
$\prolongg{\vecc}{E}_{|X-W'}\lrarr \prolongg{\vecc}{E}(h)_{|X-W'}$.
\end{itemize}
\end{lem}
\pf
The restriction of $\prolongg{\vecc}{E}$ and $\prolongg{\vecc}{E}(h)$
to any generic hypersurface $Y$ of $L$ is isomorphic
by the construction.
Then the claim of the lemma easily follows from 
Corollary 2.53 in \cite{mochi2}.
\hfill\qed

\vspace{.1in}
Since both of $\prolongg{\vecc}{E}$ and $\prolongg{\vecc}{E}(h)$
are locally free, they are isomorphic.
In particular, we can conclude that $h$ is adapted to the parabolic structure.
\hfill\qed

%% file: 10.tex
We see that any flat bundle on a smooth irreducible 
quasiprojective variety 
can be deformed to a Variation of Polarized Hodge Structure.
We can derive a result on the fundamental group.

We owe the essential ideas in this chapter to Simpson \cite{s5}.
In fact, our purpose is to show a natural generalization
of his results on smooth projective varieties.
We will use his ideas without mentioning his name.
This section is included for a rather expository purpose.

\section{Torus Action on the Moduli Space of Representations}

\subsection{Notation}

We begin with a general remark.
Let $V$ and $V'$ be vector spaces over $\cnum$,
and $\Phi:V\lrarr V'$ be a linear isomorphism.
Let $\Gamma$ be any group,
and $\rho:\Gamma\lrarr \GL(V)$ be a homomorphism.
Then $\Phi$ and $\rho$ induce the homomorphism
$\Gamma\lrarr \GL(V')$,
which is denoted by $\Phi_{\ast}(\rho)$.
We also use the notation in Subsection
\ref{subsection;05.9.12.1}.

\subsection{Continuity}
\label{subsection;05.9.11.2}

Let $X$ be a smooth irreducible projective variety
with a polarization $L$,
and $D$ be a normal crossing divisor.
Let $x$ be a point of $X-D$.
We put $\Gamma:=\pi_1(X-D,x)$.
Let $(\vecE_{\ast},\theta)$ be a $\mu_L$-polystable regular filtered Higgs bundle
on $(X,D)$
with trivial characteristic numbers.
We put $E:=\vecE_{|X-D}$.
Since $(\vecE_{\ast},t\cdot\theta)$ are also $\mu_L$-polystable,
we have a pluri-harmonic metric $h_{t}$
for $(E,\delbar_{E},t\cdot\theta)$ on $X-D$
adapted to the parabolic structure, due to Theorem \ref{thm;04.10.24.1}.
Therefore, we obtain the family of the representations
$\rho'_t:\Gamma\lrarr \GL(E_{|x})$
$(t\in\cnum^{\ast})$.
We remark that $\rho_t'$ are independent of
the choice of pluri-harmonic metrics $h_t$.

Let $V$ be a $\cnum$-vector space whose rank is same as $\rank E$.
Let $h_V$ be a hermitian vector space of $V$.
For any $t\in \cnum^{\ast}$,
we take isometries $\Phi_t:(E_{|x},h_{t\,|\,x})\lrarr (V,h_V)$,
and then we obtain the family of representations
$\rho_t:=\Phi_{t\,\ast}(\rho_t')\in  R(\Gamma,\GL(V))$.
We remark that $\pi_{\GL(V)}(\rho_t)$ are independent of
choices of $\Phi_t$.
Thus we obtain the map $\nbigp:\cnum^{\ast}\lrarr M(\Gamma,V,h_V)$
by $\nbigp(t)=\pi_{\GL(V)}(\rho_t)$.

\begin{thm}
 \label{thm;05.9.10.5}
The induced map $\nbigp$ is continuous.
\end{thm}
\pf
We may and will assume that $(\vecE_{\ast},\theta)$ is $\mu_L$-stable
for the proof.
Let $\{t_i\in\cnum^{\ast}\,|\,i\in\seisuu_{>0}\}$ be a sequence
converging to $t_0$.
We have only to take a subsequence $\{t_i\,|\,i\in S\}$ 
and a sequence of isometries
$\bigl\{
\Psi_i:\bigl(E_{|x},h_{t_i\,|\,x}\bigr)
\lrarr
 \bigl(E_{|x},h_{t_0\,|\,x}\bigr)\,\big|\,i\in S
\bigr\}$
such that 
$\bigl\{\Psi_{i\,\ast}\bigl(\rho_{t_i}\bigr)
 \,\big|\,i\in S\bigr\}$ converges to $\rho_{t_0}$.
Since the sections
$\det(T-t_i\cdot\theta)$ of $\Sym^{\cdot}\Omega^{1,0}[T]$
converges to $\det(T-t_0\cdot\theta)$,
we may apply Theorem \ref{thm;05.9.10.2}.
Hence there exists a subsequence 
$\bigl\{t_i\,\big|\,i\in S'\bigr\}$
such that 
$\bigl\{
 \bigl(E,\delbar_E,h_{t_i},t_i\!\cdot\!\theta_i\bigr)\,\big|\,
 i\in S'\bigr\}$
converges to a tame harmonic bundle
$(E',\delbar_{E'},h',\theta')$ in $L_2^p$ locally on $X-D$
via some isometries
$\Phi_i:(E,h_{t_i})\lrarr (E',h')$ $(i\in S')$.
It is easy to see that the representations
$\Phi_{i\,|\,x\,\ast}(\rho_{t_i})$ converges to
$\rho'$ in $R(\Gamma,E'_{|x},h'_{|x})$,
where $\rho'$ is associated to the flat connection
$\delbar_{E'}+\del_{E'}+\theta'+\theta^{\prime\,\dagger}$.

We also have the non-trivial holomorphic map
$f:\prolongg{\vecc}{E}'\lrarr \prolongg{\vecc}{E}$
which is compatible with the parabolic structure and the Higgs fields
due to Theorem \ref{thm;05.9.10.2}.
Since $(\prolongg{\vecc}{E}'_{\ast},\theta')$ is $\mu_L$-polystable
and $(\prolongg{\vecc}{E}_{\ast},t_0\cdot\theta)$ is $\mu_L$-stable,
the map $f$ is isomorphic.
Then we have $f_{|x\,\ast}(\rho')=\rho_{t_0}$.
By replacing $f$ appropriately,
we may assume $f:E'\lrarr E$ is isometric
with respect to $h'$ and $h_{t_0}$.
Hence $\Psi_i:=\bigl(f\circ \Phi_i\bigr)_{|x}$ gives the desired isometries.
Thus Theorem \ref{thm;05.9.10.5} is proved.
\hfill\qed

\subsection{Limit}

\begin{lem}
 \label{lem;05.9.11.3}
$\nbigp\bigl(\{t\in\cnum^{\ast}\,|\,|t|<1\}\bigr)$ is
relatively compact in $M(\Gamma,V,h_V)$.
\end{lem}
\pf
The sequence of sections
$\det (T-t\cdot \theta)$  of $\Sym^{\cdot}\Omega^{1,0}[T]$
clearly converges to $T^{\rank E}$ when $t\to 0$.
Hence we may apply the first claim of Theorem \ref{thm;05.9.10.2},
and we obtain a subsequence $\{t_i\}$ converging to $0$
such that
$\bigl\{(E,\delbar_{E},t_i\cdot\theta,h_{t_i})\bigr\}$
converges to a tame harmonic bundle
$(E',\delbar_{E'},\theta',h')$
weakly in $L_2^p$ locally on $X-D$.
Then we easily obtain the convergence
of the sequence $\bigl\{\pi_{\GL(V)}(\rho_{t_i})\bigr\}$
in $M(\Gamma,V,h_V)$.
\hfill\qed

\vspace{.1in}
Ideally, the sequence $\{\nbigp(t)\}$
should converge in $t\to 0$,
and the limit should come from 
a Variation of Polarized Hodge Structure.
We discuss only a partial but useful result about it.

Let us recall relative Higgs sheaves.
In the following,
we put $\cnum_t:=\Spec \cnum[t]$ and $\cnum_t^{\ast}:=\Spec \cnum[t,t^{-1}]$.
For a smooth morphism $Y_1\lrarr Y_2$,
the sheaf of relative holomorphic $(1,0)$-forms
are denoted by $\Omega^{1,0}_{Y_1/Y_2}$.
We put $\gbigx:=X\times \cnum_t$
and $\gbigx^{\ast}:=X\times \cnum_t^{\ast}$.
Similarly, $\gbigd:=D\times\cnum_t$
and $\gbigd^{\ast}:=D\times\cnum_t^{\ast}$.
We put
$\prolongg{\vecc}{\widetilde{E}}_{\ast}
 :=\prolongg{\vecc}{E}_{\ast}\otimes\nbigo_{\cnum_t^{\ast}}$
which is $\vecc$-parabolic bundle on $(\gbigx^{\ast},\gbigd^{\ast})$.
Then, $t\cdot\theta$ gives the relative Higgs field
$\widetilde{\theta}$,
which is a homomorphism
$ \prolongg{\vecc}{\widetilde{E}}_{\ast}\lrarr
 \prolongg{\vecc}{\widetilde{E}}_{\ast}\otimes
 \Omega^{1,0}_{\gbigx^{\ast}/\cnum_t^{\ast}}(\log \gbigd^{\ast})$
such that $\widetilde{\theta}^2=0$.
Using the standard argument of S. Langton \cite{langton},
we obtain the $\vecc$-parabolic sheaf
 $\prolongg{\vecc}{\widetilde{E}'}_{\ast}$ and
relative Higgs field
$\widetilde{\theta}':
  \prolongg{\vecc}{\widetilde{E}'}_{\ast}\lrarr
 \prolongg{\vecc}{\widetilde{E}'}_{\ast}
  \otimes\Omega^{1,0}_{\gbigx/\cnum_t}$
satisfying the following (see \cite{y}):
\begin{itemize}
\item
 $\prolongg{\vecc}{\widetilde{E}'}_{\ast}$ is flat over $\cnum_t$,
and the restriction to $\gbigx^{\ast}$ is $\prolongg{\vecc}{\widetilde{E}}_{\ast}$.
\item
 The restriction of
 $\widetilde{\theta}'$ to $\gbigx^{\ast}$
 is $\widetilde{\theta}$.
\item
 $(\prolongg{\vecc}{\Ehat}'_{\ast},\thetahat'):=
 \bigl(\prolongg{\vecc}{\widetilde{E}'}_{\ast},\widetilde{\theta}'\bigr)_{|X\times\{0\}}$ 
 is $\mu_L$-semistable.
\end{itemize}
Let $(\prolongg{\vecc}{\Ehat}_{\ast},\thetahat)$
denote the reflexive saturated
regular filtered Higgs sheaf
associated to $(\prolongg{\vecc}{\Ehat}',\thetahat')$.
(See Lemma \ref{lem;06.8.5.40}.)
We put $\Ehat:=\prolongg{\vecc}{\Ehat}_{|X-D}$.

\begin{prop}
 \label{prop;04.10.23.101}
Assume that
$(\prolongg{\vecc}{\Ehat}_{\ast},\thetahat)$ is $\mu_L$-stable.
\begin{itemize}
\item
 $(\prolongg{\vecc}{\Ehat}_{\ast},\thetahat)$ is a Hodge bundle,
 i.e.,
  $(\prolongg{\vecc}{\Ehat}_{\ast},\alpha\cdot \thetahat)
 \simeq (\prolongg{\vecc}{\Ehat}_{\ast},\thetahat)$
 for any $\alpha\in\cnum^{\ast}$.
\item
 We have a pluri-harmonic metric $\widehat{h}$
 of a Hodge bundle $(\Ehat,\thetahat)$ on $X-D$,
 which is adapted to the parabolic structure.
 It induces the Variation of Polarized Hodge Structure.
 Thus we obtain the corresponding representation
  $\widehat{\rho}:\pi_1(X-D,x)\lrarr \GL(\Ehat_{|x})$
 which underlies a Variation of Polarized Hodge Structure.
\item
 Take any isometry $G:(\Ehat_{|x},\hhat_{|x})\simeq (V,h_V)$.
 Then the sequence
 $\bigl\{\pi_{\GL(V)}(\rho_t)\bigr\}$ converges to
$\pi_{\GL(V)}\bigl(G_{\ast}(\widehat{\rho})\bigr)$ 
 in $M(\Gamma,V,h_V)$ for $t\to 0$.
\item
 In particular,
 the map $\pi_{\GL(V)}(\rho_t):\cnum^{\ast}\lrarr M(\Gamma,V,h_V)$
 is continuously extended to the map of $\cnum$ to $M(\Gamma,V,h_V)$.
\end{itemize}
\end{prop}
\pf
The argument is essentially due to Simpson \cite{s5}.
The fourth claim follows from the third one.
Let $\{t_i\,|\,i\in\seisuu_{>0}\}$ be a sequence converging to $0$.
Due to Theorem \ref{thm;05.9.10.2},
there exists a subsequence $\{t_i\,|\,i\in S\}$ such that
the sequence $\bigl\{
\bigl(E,\delbar_E,h_{t_i},t_i\cdot\theta\bigr)
\,\big|\,i\in S \bigr\}$
converges to a tame harmonic bundle
$(E',\delbar_{E'},h',\theta')$
weakly in $L_2^p$ locally on $X-D$,
via isometries $\Phi_i:(E,h_{t_i})\lrarr (E',h')$.
Let $\rho':\pi_1(X-D,x)\lrarr \GL(E'_{|x})$ denote 
the representation associated to 
the flat connection
$\delbar_{E'}+\del_{E'}+\theta'+\theta^{\prime\,\dagger}$.
Then we have the convergence of 
$\bigl\{\Phi_{i|x\,\ast}(\rho_{t_i})\,\big|\,i\in S''\bigr\}$
to $\rho'$ in $M(\Gamma,\Ehat_{|x},\hhat_{|x})$.
Due to Theorem \ref{thm;05.9.10.2},
we also have a non-trivial morphism
$f:\prolongg{\vecc}{\Ehat}'\lrarr
\prolongg{\vecc}{E'}$
which is compatible with the parabolic structures and the Higgs fields.
It induces the morphism
$\prolongg{\vecc}{\Ehat}\lrarr \prolongg{\vecc}{E'}$
compatible with the parabolic structures and the Higgs fields.
Then it must be isomorphic
due to $\mu_L$-polystability of
$(\prolongg{\vecc}{E}'_{\ast},\theta')$ and 
$\mu_L$-stability of $(\prolongg{\vecc}{\Ehat}_{\ast},\thetahat)$.
In particular,
$(\prolongg{\vecc}{\Ehat}_{\ast},\thetahat)$ is a
$\mu_L$-stable $\vecc$-parabolic Higgs bundle.
The metric $\widehat{h}$ of $\Ehat$ is given by $h'$ and $f$.
Thus the third claim is obtained.

Let us consider the morphism $\phi_{\alpha}:\cnum_t\lrarr \cnum_t$
given by $t\longmapsto \alpha\cdot t$.
We have the natural isomorphism
$\phi_{\alpha}^{\ast}\bigl(\prolongg{\vecc}{\Etilde}_{\ast},\thetatilde\bigr)
\simeq \bigl(\prolongg{\vecc}{\Etilde}_{\ast},\alpha\cdot\thetatilde\bigr)$
which can be extended to the morphism
$\phi_{\alpha}^{\ast}
 (\prolongg{\vecc}{\widetilde{E}}'_{\ast},
 \widetilde{\theta}')
\lrarr (\prolongg{\vecc}{\widetilde{E}'_{\ast}},
 \alpha\!\cdot\!\widetilde{\theta}')$
such that the specialization 
$(\prolongg{\vecc}{\widehat{E}}_{\ast},\widehat{\theta})\lrarr 
 (\prolongg{\vecc}{\widehat{E}}_{\ast},\alpha\cdot\widehat{\theta})$
at $t=0$ is not trivial.
Since $(\prolongg{\vecc}{\widehat{E}}_{\ast},\widehat{\theta})$ and 
 $(\prolongg{\vecc}{\widehat{E}_{\ast}},\alpha\cdot\widehat{\theta})$
are $\mu_L$-stable,
the map is isomorphic.
Hence $(\prolongg{\vecc}{\Ehat},\thetahat)$ is a Hodge bundle.
Thus the first is proved.

Since $(\widehat{E},\delbar_{\widehat{E}},\widehat{\theta})$
is a Hodge bundle,
we have the the action $\kappa$ 
of $S^1=\{t\in\cnum\,|\,|t|=1\}$ on $\widehat{E}$
such that 
$\kappa(t):
 \bigl(\widehat{E},\delbar_{\widehat{E}},\widehat{\theta}\bigr)\simeq 
 \bigl(\widehat{E},\delbar_{\widehat{E}},t\cdot\widehat{\theta}\bigr)$
for any $t\in S^1$.
The metric $\kappa(t)_{\ast}\widehat{h}$ is determined by 
$\kappa(t)_{\ast}\widehat{h}(u,v)
=\widehat{h}\bigl(\kappa(t)(u),\kappa(t)(v)\bigr)$,
which is also the pluri-harmonic metric of
$(\widehat{E},\delbar_{\widehat{E}},t\!\cdot\!\theta)$.
Since $(\widehat{\vecE}_{\ast},t\cdot\widehat{\theta})$ is $\mu_L$-stable,
the pluri-harmonic metric is
unique up to a positive constant multiplication.
Hence we obtain the map
$\nu:S^1\lrarr \real_{>0}$
such that $\kappa(t)_{\ast}\widehat{h}=\nu(t)\cdot \widehat{h}$.
Let $\Ehat=\bigoplus \Ehat_w$ be the weight decomposition.
For $v_{i}\in \Ehat_{w_i}$ $(w_1\neq w_2)$,
we have 
$ \nu(t)\cdot \hhat(v_{1},v_2)
=\kappa(t)_{\ast}\hhat(v_1,v_2)
=t^{w_1-w_2}\hhat\bigl(v_1,v_2\bigr)$.
Hence, we obtain $\hhat(v_1,v_2)=0$
and $\nu(t)=1$.
Namely, $\widehat{h}$ is $S^1$-invariant,
which means
$\bigl(\widehat{E},\delbar_{\widehat{E}},
  \widehat{\theta},\widehat{h}\bigr)$
gives a Variation of Polarized Hodge Structure.
Thus the second claim is proved.
\hfill\qed

\begin{lem}
 \label{lem;05.9.11.5}
Assume $(\prolongg{\vecc}{\Ehat}_{\ast},\thetahat)$ is not
 $\mu_L$-stable.
Let $\rho_0$ be an element of $R(\Gamma,V)$
such that $\pi_{\GL(V)}(\rho_0)$ is the limit
of a subsequence 
$\bigl\{\pi_{\GL(V)}(\rho_{t_i})\bigr\}$ for $t_i\to 0$.
Then $\rho_0$ is not simple.
\end{lem}
\pf
Let $\{t_i\}$ be a sequence converging to $0$
such that
$\bigl\{(E,\delbar_E,t_i\cdot\theta,h_{t_i})\bigr\}$ converges
to a tame harmonic bundle $(E',\delbar_{E'},\theta',h')$
in $L_2^p$ locally on $X-D$.
We may assume that
$\rho_0$ is the associated representation to $(E',\delbar_{E'},\theta',h')$.
We have a non-trivial map $f:\prolongg{\vecc}{E'}\lrarr
\prolongg{\vecc}{\Ehat}$
compatible with the parabolic structures and the Higgs fields.
If $\rho_0$ is simple,
then $(\prolongg{\vecc}{E}'_{\ast},\theta')$ is $\mu_L$-stable,
and it can be shown that the map $f$ has to be isomorphic.
But it contradicts with the assumption
that $(\prolongg{\vecc}{\Ehat},\thetahat)$ is not $\mu_L$-stable.
\hfill\qed

\subsection{Deformation to a Variation of Polarized Hodge Structure}

Let $Y$ be a smooth irreducible quasiprojective variety over $\cnum$
with a base point $x$.
We may assume $Y=X-D$,
where $X$ and $D$ denote a smooth projective variety
and its simple normal crossing divisor, respectively.
A representation $\rho:\pi_1(Y,x)\lrarr \GL(V)$
induces a flat bundle $(E,\nabla)$.
We say that $\rho$ comes from a Variation of Polarized Hodge Structure,
if $(E,\nabla)$ underlies a Variation of Polarized Hodge Structure.
For simplicity of the notation,
we put $\Gamma:=\pi_1(Y,x)$.

\begin{thm} \label{thm;04.11.8.10}
Let $\rho\in R(\Gamma,V)$ be a representation.
Then it can be deformed to a representation
$\rho'\in R(\Gamma,V)$
which comes from a Variation of Polarized Hodge Structure
on $Y$.
\end{thm}
\pf
We essentially follow the argument of Theorem 3 in \cite{s5}.
Any representation $\rho\in R(\Gamma,V)$
can be deformed to a semisimple representation
$\rho'\in R(\Gamma,V)$.
Therefore we may assume that $\rho$ is semisimple
from the beginning.
Let $(E,\nabla)$ be the corresponding semisimple flat bundle on $X-D$.
We can take a Corlette-Jost-Zuo metric $h$ of $(E,\nabla)$,
and hence we obtain the tame pure imaginary harmonic bundle
$(E,\delbar_E,\theta,h)$.
Let $(\vecE_{\ast},\theta)$ denote the associated
regular filtered Higgs bundle on $(X,D)$.
We have the canonical decomposition (Corollary \ref{cor;05.9.14.5}):
\[
 (\vecE_{\ast},\theta)
=\bigoplus_{j\in \Lambda}(\vecE_{i\ast},\theta_i)\otimes\cnum^{m(j)}.
\]
We put 
$r(\rho):=\sum_{j\in\Lambda} m(j)$.
Note that $r(\rho)\leq \rank E$,
and we have $r(\rho)=\rank E$ if and only if
$(\vecE_{\ast},\theta)$ is a direct sum of Higgs bundles of rank one.
We use a descending induction on $r(\rho)$.

We obtain the family of regular filtered Higgs bundles
$\bigl\{(\vecE_{\ast},t\cdot\theta)\,\big|\,t\in\cnum^{\ast}\bigr\}$
$(t\in\cnum^{\ast})$.
In particular,
we have the associated deformation of representations 
$\{\rho_t\in R(\Gamma,V)\,|\,t\in \real_{>0}\}$
as in Subsection \ref{subsection;05.9.11.2}.
We may assume 
$\rho_1=\rho$.
We have the induced map
$\nbigp:\openclosed{0}{1}\lrarr  M(\Gamma,V,h_V)$
given by $\nbigp(t):=\pi_{\GL(V)}(\rho_t)$,
which is continuous due to Theorem \ref{thm;05.9.10.5}.
The image is relatively compact due to Lemma \ref{lem;05.9.11.3}.
We take a representation $\rho_0\in R(\Gamma,V)$
such that $\pi_{\GL(V)}(\rho_0)$ is the limit
of a subsequence of
$\bigl\{
 \pi_{\GL(V)}(\rho_t)\,\big|\,t\in \openclosed{0}{1}
 \bigr\}$.
We may assume that it comes from a 
tame harmonic bundle 
as in the proof of Lemma \ref{lem;05.9.11.3}.

\vspace{.1in}
\noindent
{\bf The case 1.}
Let $(\vecE_{\ast},\theta)=
 \bigoplus (\vecE_{i\,\ast},\theta_i)^{\oplus m_i}$
be the canonical decomposition.
Assume that each family
$\{(\vecE_{i\,\ast},t\cdot\theta_i)\,|\,t\in\cnum^{\ast}\}$
converges to the $\mu_L$-stable regular filtered Higgs sheaf.
Then $\rho_0$ comes from a Variation of Polarized Hodge Structure
due to Proposition \ref{prop;04.10.23.101}.

We remark that the rank one Higgs bundle is always stable.
Hence the case $r(\rho)=\rank E$ is done,
in particular.

\vspace{.1in}
\noindent
{\bf The case 2.}
Assume that one of the families
$\{(\vecE_{\ast},t\cdot\theta_i)\,|\,t\in\cnum^{\ast}\}$ converges
to the semistable parabolic Higgs sheaf, which is not $\mu_L$-stable.
Then we have $r(\rho)<r(\rho_0)$
due to Lemma \ref{lem;05.9.11.5}.
Hence the induction can proceed.
\hfill\qed

\section{Monodromy Group}
\label{section;05.10.4.2}

We discuss the monodromy group for the Higgs bundles
or flat bundles, by following the ideas in \cite{s5}.

\subsection{The Higgs monodromy group}

Let $X$ be a smooth irreducible projective variety
with an ample line bundle $L$,
and $D$ be a simple normal crossing divisor.
Let $(\vecE_{\ast},\theta)$ be a $\mu_L$-polystable
regular filtered Higgs bundle on $(X,D)$
with trivial characteristic numbers.
For any non-negative integers $a$ and $b$,
we have the regular filtered Higgs bundles
$(T^{a,b}\vecE_{\ast},\theta)$.
(See Subsection \ref{subsection;05.9.29.1}
for the explanation.)
Since we have a pluri-harmonic metric $h$ of $(E,\delbar_E,\theta)$
adapted to the parabolic structure,
the regular filtered Higgs bundles $T^{a,b}(\vecE_{\ast},\theta)$
are also $\mu_L$-polystable.
In particular, we have the canonical decompositions of them.
We recall the definition of the Higgs monodromy group
given in \cite{s5}.
Let $x$ be a point of $X-D$.

\begin{df}
The Higgs monodromy group $M(\vecE_{\ast},\theta,x)$
of $\mu_L$-polystable Higgs bundle $(\vecE_{\ast},\theta)$ is 
the subgroup of $\GL(E_{|x})$ defined as follows:
An element $g\in \GL(E_{|x})$ is contained in $M(\vecE_{\ast},\theta,x)$,
if and only if
$T^{a,b}g$ preserves the subspace
$F_{|x}\subset T^{a,b}E_{|x}$
for any stable component
$(\vecF_{\ast},\theta_F)\subset T^{a,b}(\vecE_{\ast},\theta)$.
\hfill\qed
\end{df}

\begin{rem}
Although
such a Higgs monodromy group should be
defined for semistable parabolic Higgs bundles
as in {\rm\cite{s5}},
we do not need it in this paper.
\hfill\qed
\end{rem}

We have an obvious lemma.
\begin{lem} 
 \label{lem;04.11.8.20}
We have
$M(\vecE_{\ast},\theta,x)=M(\vecE_{\ast},t\cdot\theta,x)$
for any $t\in\cnum^{\ast}$,
i.e.,
the Higgs monodromy group is invariant under the torus action.
\hfill\qed
\end{lem}

Let us take a pluri-harmonic metric $h$ of
the Higgs bundle $(E,\delbar_E,\theta)$ on $X-D$,
which is adapted to the parabolic structure.
Then we obtain the flat connection 
$\DD^1=\delbar_E+\del_E+\theta+\theta^{\dagger}$.
Then we obtain the monodromy group
$M(E,\DD^1,x)\subset \GL(E_{|x})$
of the flat connection.
(See Subsection \ref{subsection;04.11.6.2}.)

\begin{lem} 
 \label{lem;04.11.8.21}
We have $M(E,\DD^1,x)\subset M(\vecE_{\ast},\theta,x)$.
For a tame pure imaginary harmonic bundle,
we have $M(E,\DD^1,x)=M(\vecE_{\ast},\theta,x)$.
\end{lem}
\pf
A stable component
$(\vecF_{\ast},\theta_F)\subset (\vecE_{\ast},\theta)$ 
induces the flat subbundle of 
$F\subset T^{a,b}(E,\DD^1)$.
If $g\in M(E,\DD^1,x)$,
we have $T^{a,b}g(F_{|x})\subset F_{|x}$.
Hence, $M(E,\DD^1,x)\subset M(\vecE_{\ast},\theta,x)$.
In the pure imaginary case,
a flat subbundle $F\subset T^{a,b}(E,\DD^1)$
induces $(\vecF_{\ast},\theta_F)\subset (\vecE_{\ast},\theta)$.
Therefore, we obtain 
$M(E,\DD^1,x)=M(\vecE_{\ast},\theta,x)$.
\hfill\qed

\subsection{The deformation and the monodromy group}
\label{subsection;04.11.8.1}

For simplicity of the description,
we put $\Gamma:=\pi_1(X-D,x)$.
Let $(E,\nabla)$ be a semisimple flat bundle over $X-D$.
We have a Corlette-Jost-Zuo metric $h$
of $(E,\nabla)$,
and thus we obtain a tame pure imaginary harmonic bundle
$(E,\delbar_E,\theta,h)$ on $X-D$.
The associated regular filtered Higgs bundle is denoted by
$(\vecE_{\ast},\theta)$,
which is $\mu_L$-polystable
with trivial characteristic numbers.

As in Subsection \ref{subsection;05.9.11.2},
we have the pluri-harmonic metrics $h_t$
for any $(E,\delbar_E,t\cdot\theta)$ $(t\in\cnum^{\ast})$.
Hence we obtain the flat connections $\DD^{1}_t$ of $E$,
and the representations
$\rho_t:\pi_1(X-D,x)\lrarr \GL(E_{|x})$.

\begin{lem}
We have $M(E,\DD^1_t)\subset M(E,\DD^1_1)$
for $t\in\cnum-\{0\}$,
and 
$M(E,\DD^1_t)=M(E,\DD^1_1)$
for $t\in \real-\{0\}$.
\end{lem}
\pf
It follows from Lemma \ref{lem;04.11.8.20}
and Lemma \ref{lem;04.11.8.21}.
\hfill\qed

\vspace{.1in}

We put $G_0:=M(E,\DD^1_t,x)$ for $t\in\real_{>0}$
which is independent of the choice of $t$.
Let $U(E,h_{t},x)$ denote the unitary group 
for the metrized space $(E_{|x},h_{t\,|\,x})$.
Due to Lemma \ref{lem;04.11.4.1},
$G_0$ is reductive, and 
the intersection $K_{0,t}:=G_0\cap U(E,h_{t},x)$
is a compact real form of $G_0$.

We put $V:=E_{|x}$ and $h_V:=h_{1\,|\,x}$.
We denote $G_0$ and $K_{0,1}$ by $G$ and $K$ respectively,
when we regard it as the subgroup of $\GL(V)$.
Then we can take an isometry
$\nu_t:(E_{|x},h_{t\,|\,x})\simeq (V,h_V)$
such that $\nu_t(G_0)=G$ and $\nu_t(K_{0\,t})=K$
for each $t$.
Such a map is unique up to the adjoint of $N_{G}(h_V)$.
Thus we obtain the family of representations
$\rhotilde_t:=\nu_{t\,\ast}(\rho_t)\in R(\Gamma,G)$
$(t\in \real_{>0})$.

\begin{lem}
The induced map
$\pi_{G}(\rhotilde_t):\real_{>0}\lrarr M(\Gamma,G,h_V)$
is continuous.
\end{lem}
\pf
Let $M'$ denote the subset of $M(\Gamma,G,h_V)$
which consists of the Zariski dense representations.
The natural morphism $M'\lrarr M(\Gamma,V,h_V)$
is injective,
and the image of $\pi_G(\rhotilde_t)$ is contained in $M'$.
Hence the claim of the lemma
follows from Theorem \ref{thm;05.9.10.5}
and the properness of $M(\Gamma,G,h_V)\lrarr M(\Gamma,V,h_V)$.
\hfill\qed

\begin{lem}
The image $\pi_G(\rhotilde_t)\bigl(\openclosed{0}{1}\bigr)$
is relatively compact in $M(\Gamma,G,h_V)$.
\end{lem}
\pf
It follows from Lemma \ref{lem;05.9.11.3}
and the properness of the map $M(\Gamma,G,h_V)\lrarr M(\Gamma,V,h_V)$.
\hfill\qed

\subsection{Non-existence result about fundamental groups}

Let $Y$ be a quasiprojective variety.
We put $\Gamma:=\pi_1(Y,x)$.
Let $V$ be a finite dimensional $\cnum$-vector space.
Let $G$ be a reductive subgroup of $\GL(V)$.
We see the convergence of $\pi_G(\rhotilde_t)$ $(t\to 0)$
in a simple case.

\begin{lem}
 \label{lem;05.9.11.30}
Let $\rho$ be an element of $R(\Gamma,G)$.
We assume that there exists a subgroup $\Gamma_0$
such that $\rho_{|\Gamma_0}:\Gamma_0\lrarr G$ is Zariski dense and rigid.
Then we can take a deformation $\rho'\in R(\Gamma,G)$ of $\rho$
which comes from a Variation of Polarized Hodge Structure on $Y$.
\end{lem}
\pf
We take a tame pure imaginary pluri-harmonic bundle
$(E,\delbar_E,\theta,h)$
whose associated representation gives $\rho$,
and we take the deformation 
$\pi_G(\rhotilde_t)$.
Let us take $\rho_0\in R(\Gamma,\gun)$
such that some sequence $\{\pi_G(\rhotilde_{t_i})\}$
converges to $\pi_G(\rho_0)$.
We remark that $\rho_{0\,|\,\Gamma_0}:\Gamma_0\lrarr G$ is
also Zariski dense and rigid
(Lemma \ref{lem;05.9.10.1}).
If $\rho_0$ comes from a Variation of Polarized Hodge Structure,
we are done.
If $\rho_0$ does not come from a Variation of Polarized Hodge Structure,
we deform $\rho_0$ as above, again.
The process will stop in the finite steps
by Theorem \ref{thm;04.11.8.10}.
\hfill\qed

\vspace{.1in}

The following lemma is a straightforward
generalization of Lemma 4.4 in \cite{s5}.
(See also Lemma \ref{lem;04.11.4.1},
where we will see the argument of
Lemma 4.4 can be generalized
in our situation.)
\begin{lem}
 \label{lem;04.11.8.100}
Let $\rho:\Gamma\lrarr G$ be a Zariski dense homomorphism.
If $\rho$ comes from a Variation of Polarized Hodge Structure,
then the real Zariski closure $W$ of $\rho$ is 
a real form of $G$, and $W$ is 
a group of Hodge type in the sense of Simpson.
(See the page {\rm 46} in {\rm\cite{s5}}.)
\hfill\qed
\end{lem}

The following lemma is essentially same as Corollary 4.6 in \cite{s5}.
\begin{prop}
 \label{prop;04.11.6.60}
Let $G$ be a complex reductive algebraic group,
and $W$ be a real form of $G$.
Let $\rho:\Gamma\lrarr G$ be a representation
such that $\Image\rho\subset W$.
Assume that there exists a subgroup $\Gamma_0\subset\Gamma$
such that
$\rho_{|\Gamma_0}$ is rigid and Zariski dense in $G$.
Then $W$ is a group of Hodge type,
in the sense of Simpson.
\end{prop}
\pf
We reproduce the argument of Simpson.
Since $\rho(\Gamma_0)$ is Zariski dense in $G$,
$W$ is also the real Zariski closure of $\rho(\Gamma_0)$.
We take a deformation $\rho'$ of $\rho$,
which comes from a Variation of Polarized Hodge Structure
as in Lemma \ref{lem;05.9.11.30}.
Then there exists an element $u\in N(G,U)$
such that $\ad(u)\circ \rho_{|\Gamma_0}\simeq \rho'_{|\Gamma_0}$
due to Lemma \ref{lem;05.9.10.1}.
Let $W'$ denote the real Zariski closure of $\rho'(\Gamma_0)$,
which is also the real Zariski closure of $\rho'$.
It is a group of Hodge type
(Lemma \ref{lem;04.11.8.100}).
Since $W$ and $W'$ are isomorphic,
we are done.
\hfill\qed

\begin{cor}
Let $\Gamma_0$ be a rigid discrete subgroup
of a real algebraic group, which is not of Hodge type.
Then $\Gamma_0$ cannot be a split quotient
of the fundamental groups of any smooth irreducible
quasiprojective variety.
\end{cor}
\pf
It follows from Lemma \ref{lem;04.11.8.100}
and Proposition \ref{prop;04.11.6.60}.
(See the pages 52--54 of \cite{s5}).
\hfill\qed

%% file: 11.tex
\section{$G$-Principal Bundles with Flat Structure 
 or Holomorphic Structure}
\label{section;04.10.26.91}

We recall the Tannakian consideration
about harmonic bundles given in \cite{s5} by Simpson.

\subsection{A characterization of algebraic subgroup of $\GL$}

We recall some facts on algebraic groups.
(See also I. Proposition 3.1 in \cite{deligne-milne},
 for example.)
Let $V$ be a vector space over a field $k$ of characteristic $0$.
We put $T^{a,b}V:=Hom(V^{\otimes\,a},V^{\otimes b})$.
Let $G$ be an algebraic subgroup of $\GL(V)$, defined over $k$.
We have the induced $G$-action on $T^{a,b}V$.
Let $\nbigs(V,a,b)$ denote the set of
$G$-subspaces of $T^{a,b}V$,
and we put $\nbigs(V)=\coprod_{a,b}\nbigs(V,a,b)$.

Let $g$ be an element of $\GL(V)$.
We have the induced element
$T^{a,b}(g)\in \GL\bigl(T^{a,b}V\bigr)$.
Then, it is known that
$g\in \GL(V)$ is contained in $G$,
if and only if
$T^{a,b}(g)W\subset W$ holds
for any $(W,a,b)\in \nbigs(V)$.
Suppose $G$ is reductive.
Then there is an element $v$ of $T^{a,b}(V)$ for some $(a,b)$
such that
$g$ is contained in $G$
if and only if $g\cdot v=v$ holds.

We easily obtain a similar characterization
of Lie subalgebras of $\gl(V)$
corresponding to algebraic subgroups of $\GL(V)$.

\subsection{A characterization of connections of principal
 $G$-bundle }

Let $k$ denote the complex number field $\cnum$
or the real number field $\real$.
Let $G$ be an algebraic group over $k$.
Let $P_G$ be a $G$-principal bundle on a manifold $X$
in the $C^{\infty}$-category.
Let $\kappa:G\lrarr \GL(V)$ be a representation
defined over $k$,
such that the induced morphism $d\kappa:\gminig\lrarr End(V)$
is injective.
We put $E:=P_G\times_G V$.
We have
$T^{a,b}E:=Hom(E^{\otimes\,a}, E^{\otimes\,b})
 \simeq P_G\times_G T^{a,b}V$.
We have the subbundle $E_U=P_G\times_G U$ of $T^{a,b}E$
for each $U\in\nbigs(V,a,b)$.
A connection $\nabla$ on $E$ induces the connection
$T^{a,b}\nabla$ on $T^{a,b}E$.
Let $\nbiga_G(E)$ be the set of the connections $\nabla$ of $E$
such that the induced connections $T^{a,b}\nabla$ preserve
the subbundle $E_U$ for any $(U,a,b)\in\nbigs(V)$.

Let $\nbiga(P_G)$ denote the set of the connections of $P_G$.
If we are given a connection of $P_G$,
the connection $\nabla$ of $E$ is naturally induced.
It is clear that 
the connection $T^{a,b}\nabla$ preserves $E_U\subset T^{a,b}E$
for any $(U,a,b)\in \nbigs(V)$.
Hence we have the map
$\varphi:\nbiga(P_G)\lrarr \nbiga_G(E)$.

\begin{lem}
The map $\varphi$ is bijective.
\end{lem}
\pf
Since $d\kappa$ is injective,
the map $\varphi$ is injective.
Let us take a connection $\nabla\in\nbiga_{G}(E)$
and a connection $\nabla_0$ which comes from a connection
of $P_G$.
Then $f=\nabla-\nabla_0$ is a section of
$\End(E)\otimes\Omega^1$.
Since $T^{a,b}f$ preserves $E_U$
for any $(a,b)$ and $U\subset \nbigs(V,a,b)$,
$f$ comes from a section of
$\ad(P_G)\otimes\Omega^1
 \subset \End(E)\otimes\Omega^1$.
\hfill\qed

\subsection{$K$-Reduction of holomorphic $G$-principal bundle
 and the induced connection}
\label{subsection;05.9.12.10}

Let $G$ be a linear reductive group defined over $\cnum$.
Let $P_G$ be a {\em holomorphic} $G$-principal bundle on $X$.
Let $\kappa:G\lrarr \GL(V)$ be a representation
defined over $\cnum$,
such that $d\kappa:\gminig\lrarr \End(V)$ is injective.
We put $E:=P_G\times_G V$.
Let $K$ be a compact real form of $G$.
Let $P_K\subset P_G$ be a $K$-reduction
in the $C^{\infty}$-category,
i.e.,
$P_K\times_K G\simeq P_G$.
Then the connection of $P_K$ is automatically induced.
We have the canonical $G$-decomposition for each $(a,b)$:
\begin{equation} \label{eq;04.9.29.1}
T^{a,b}V=\bigoplus_{\rho\in \Irrep(G)} V^{(a,b)}_{\rho}.
\end{equation}
Here $\Irrep(G)$ denotes the set of the equivalence classes
of irreducible representations of $G$.
Each $V^{(a,b)}_{\rho}$ is isomorphic to
the tensor product of the irreducible representation $\rho$
and the trivial representation $\cnum^{m(a,b,\rho)}$.
The decomposition (\ref{eq;04.9.29.1})
is same as the canonical $K$-decomposition.
Take a $K$-invariant hermitian metric $h$ of $V$.
It induces the hermitian metric $T^{a,b}h$ of $T^{a,b}V$,
for which the decomposition  (\ref{eq;04.9.29.1})
is orthogonal.
The restriction of $T^{a,b}h$ to $V^{(a,b)}_{\rho}$
is isomorphic to a tensor product
of a $K$-invariant hermitian metric on $\rho$
and a hermitian metric on $\cnum^{m(a,b,\rho)}$.
The metric $h$ induces the hermitian metric of $E$,
which is also denoted by $h$.
From the holomorphic structure $\delbar_E$
and the metric $h$,
we obtain the unitary connection $\nabla=\del_E+\delbar_E$.
The induced connection $T^{a,b}\nabla$ on $T^{a,b}E$ 
is the unitary connection determined 
by $T^{a,b}h$ and the holomorphic structure of $T^{a,b}E$.
Then it is easy to see that $T^{a,b}\nabla$
preserves $E_U$ for any $U\in\nbigs(a,b,V)$.
Hence the connection $\nabla$ comes from $P_G$.
Since $\nabla$ also preserves the unitary structure,
we can conclude that $\nabla$ comes from the connection of $P_K$.

\subsection{The monodromy group}
\label{subsection;04.11.6.2}

We recall the monodromy group
of flat bundles (\cite{s5}).
Let $X$ be a connected complex manifold with a base point $x$.
The monodromy group of a flat bundle $(E,\nabla)$ at $x$
is defined to be the Zariski closure of
the induced representation $\pi_1(X,x)\lrarr\GL(E_{|x})$.
It is denoted by $M(E,\nabla,x)$.
Let us recall the case of principal bundles.
Let $G$ be a linear algebraic group over $\real$ or $\cnum$,
and $P_G$ be a $G$-principal bundle on $X$
with a flat connection in the $C^{\infty}$-category.
Take a point $\widetilde{x}\in P_{G|x}$.
Then we obtain the representation $\rho:\pi_1(X,x)\lrarr G$.
Then the monodromy group $M(P_G,\widetilde{x})\subset G$
is defined to be the Zariski closure of the image of $\rho$.
We obtain the canonical reduction of principal bundles
$P_{M(P_G,\widetilde{x})}\subset P_G$.
The monodromy groups of flat vector bundles
and flat principal bundles are related as follows.
Let $\kappa:G\lrarr \GL(V)$ be an injective representation.
Then we have the flat bundle
$E=P_G\times_G V=P_{M(P_G,\widetilde{x})}\times_{M(P_G,\xtilde)}V$.
Via the identification $V=E_{|x}$ given by $\widetilde{x}$,
we are given the inclusion
$M(P_G,\xtilde)\subset\GL(E_{|x})$.
Clearly $M(P_G,\xtilde)$ is same as $M(E,\nabla,x)$
and it is independent of the choice of $\widetilde{x}$.
Hence we can reduce the problems of
the monodromy groups of flat principal $G$-bundles
to those for flat vector bundles.

For a flat bundle $(E,\nabla)$,
let $T^{a,b}E$ denote
the flat bundle $Hom(E^{\otimes\,a},E^{\otimes\,b})$
provided the canonically induced flat connection.
Let $\nbigs(E,a,b)$ denote the set of
flat subbundles $U$ of $T^{a,b}E$,
and we put $\nbigs(E):=\coprod_{(a,b)}\nbigs(E,a,b)$.
Let $g$ be an element of $\GL(E_{|x})$.
Then $g$ is contained in $M(E,\nabla,x)$
if and only if $T^{a,b}g$ preserves $U_x$
for any $(U,a,b)\in\nbigs(E)$.
If $M(E,\nabla,x)$ is reductive,
we can find some $(a,b)$ and $v\in T^{a,b}E_{|x}$
such that 
$g\in M(E,\nabla,x)$ if and only if $g\cdot v=v$.
Hence there exists a flat subbundle
$W\subset T^{a,b}E$
such that $g\in M(E,\nabla,x)$ if and only if
$T^{a,b}g_{|W}=\id_W$.

\section{Definitions}
\label{section;04.10.26.92}

\subsection{A $G$-principal Higgs bundle and
 a pluri-harmonic reduction}
\label{subsection;04.10.23.1}

Let $G$ be a linear reductive group defined over $\cnum$,
and $K$ be a compact real form.
Let $X$ be a complex manifold
and $P_G$ be a holomorphic $G$-principal bundle on $X$.
Let $\ad(P_G)$ be the adjoint bundle of $P_G$,
i.e., $\ad(P_G)=P_G\times_G\gminig$.
Recall that 
a Higgs field of $P_G$ is defined to be 
a holomorphic section $\theta$ of $\ad(P_G)\otimes\Omega^{1,0}$ 
such that $\theta^2=0$. 

Let $P_K\subset P_G$ be a $K$-reduction of $P_G$
in $C^{\infty}$-category,
then we have the natural connection $\nabla$ of $P_K$,
as is seen in Subsection \ref{subsection;05.9.12.10}.
We also have the adjoint $\theta^{\dagger}$ of $\theta$,
which is a $C^{\infty}$-section of $\ad(P_G)\otimes\Omega^{0,1}$.
Then we obtain the connection
$\DD^1:=\nabla+\theta+\theta^{\dagger}$ of the principal bundle $P_G$.

\begin{df} \label{df;04.9.29.11}
If $\DD^1$ is flat,
then the reduction $P_K\subset P_G$ is called 
pluri-harmonic,
and the tuple $(P_K\subset P_G,\theta)$ is called
a $G$-harmonic bundle.
\hfill\qed
\end{df}

Let $V$ be a $\cnum$-vector space.
A representation $\kappa:G\lrarr \GL(V)$ is called 
immersive if $d\kappa$ is injective,
in this paper.
Take an immersive representation
$\kappa:G\lrarr \GL(V)$
and a $K$-invariant metric $h_V$.
From a $G$-principal Higgs bundle $(P_G,\theta)$
with a $K$-reduction $P_K\subset P_G$,
we obtain the Higgs bundle
$(E,\delbar_E,\theta)$ with the hermitian metric $h$.

\begin{lem}
Let $(P_G,\theta)$ be a $G$-principal Higgs bundle,
and $P_K\subset P_G$ be a $K$-reduction.
The following conditions are equivalent.
\begin{enumerate}
\item
 The reduction $P_K\subset P_G$ is pluri-harmonic.
\item
 For any representation $G\lrarr \GL(V)$
 and any $K$-invariant hermitian metric of $\cnum$-vector space $V$,
 the induced Higgs bundle with the hermitian metric 
 is a harmonic bundle.
\item
 There exist an immersive representation $G\lrarr \GL(V)$
 and a $K$-invariant hermitian metric of $\cnum$-vector space $V$,
 such that
 the induced Higgs bundle with the hermitian metric 
 is a harmonic bundle.
\end{enumerate}
\end{lem}
\pf
If $G\lrarr \GL(V)$ is immersive,
then a connection of $P_G$ is flat
if and only if the induced connection on $P_G\times_GV$ is flat.
Therefore the desired equivalence is clear.
\hfill\qed

\subsection{A flat $G$-bundle and a pluri-harmonic reduction}

Let $G$ be a linear reductive group over $\real$ or $\cnum$,
and let $(P_G,\nabla)$ be a flat $G$-bundle over a complex manifold $X$.
If a $K$-reduction $P_K\subset P_G$ is given,
we obtain the connection $\nabla_0$ of $P_K$
and self adjoint section $\varphi\in \ad(P_G)\otimes\Omega^1$
such that $\nabla=\nabla_0+\varphi$
(\cite{corlette4}),
which can be shown by a Tannakian consideration
as in Subsection \ref{subsection;05.9.12.10},
for example.
Let $\nabla_0=\nabla_0'+\nabla_0''$
and $\varphi=\theta+\theta^{\dagger}$
be the decomposition into
the $(1,0)$-part and the $(0,1)$-part.
The connection $\nabla_0$ induces the connection
on $\ad(P_G)$,
which is also denoted by $\nabla_0=\nabla_0'+\nabla_0''$.

\begin{df}
\label{df;04.9.29.10}
A reduction $P_K\subset P_G$ is called pluri-harmonic,
if $\theta^2=0$ and $\nabla_0''(\theta)=0$ hold.
\hfill\qed
\end{df}

Let $V$ be a vector space over $\cnum$.
Let $\kappa:G\lrarr \GL(V)$ be a representation,
which induces the flat bundle $(E,\nabla_E)$.
We take a $K$-invariant metric $h_V$,
which induces the metric $h_E$ of $E$.
We obtain the decomposition
$\nabla_E=\delbar_E+\del_E+\theta_E+\theta^{\dagger}_E$
as in Section 21.4.3 of \cite{mochi2}.
They are induced by $\nabla_0''$, $\nabla_0'$,
$\theta$ and $\theta^{\dagger}$, respectively.
Thus, if $P_K\subset P_G$ is pluri-harmonic,
we have $\theta_E^2=\delbar_E\theta_E=0$.
Recall that they imply $\delbar_E^2=0$.
Hence, $(E,\nabla_E,h)$ is a harmonic bundle.
On the contrary,
if $\kappa$ is immersive and $(E,\nabla_E,h)$ is a harmonic bundle,
we obtain the vanishings $\theta^2=\nabla_0''\theta=0$.
Hence, $P_K\subset P_G$ is pluri-harmonic.
Therefore, we obtain the following lemma.

\begin{lem}
The following conditions are equivalent.
\begin{enumerate}
\item
 The reduction $P_K\subset P_G$ is pluri-harmonic,
 in the sense of Definition {\rm \ref{df;04.9.29.10}}.
\item
 For any  representation $\kappa:G\lrarr \GL(V)$
 and any $K$-invariant metric of a vector space $V$ over $\cnum$,
 the induced flat bundle with the hermitian metric 
 is a harmonic bundle.
\item
 There exist an immersive representation $\kappa:G\lrarr \GL(V)$
 and a $K$-invariant metric of a vector space $V$ over $\cnum$,
 such that the induced flat bundle with the hermitian metric
 is a harmonic bundle.
\hfill\qed
\end{enumerate}
\end{lem}

Let $\pi:\widetilde{X}\lrarr X$ denote a universal covering.
Take base points $x\in X$
and $x_1\in \widetilde{X}$ such that $\pi(x_1)=x$.
Once we pick a point $\widetilde{x}\in P_{G|x}$,
the homomorphism $\pi_1(X,x)\lrarr G$ is given.
If a $K$-reduction $P_K\subset P_G$ is given,
we obtain a $\pi_1(X,x)$-equivariant map $F:\widetilde{X}\lrarr G/K$,
where the $\pi_1(X,x)$-action on $G/K$ is given by 
the homomorphism $\pi_1(X,x)\lrarr G$.
If $P_K\subset P_G$ is pluri-harmonic,
then $F$ is pluri-harmonic (\cite{zuo})
in the sense that any restriction of $F$
to holomorphic curve is harmonic.

\subsection{A tame pure imaginary $G$-harmonic bundle}

Let $G$ be a linear reductive group over $\cnum$.
Let $\gminih$ denote a Cartan subalgebra of $\gminig$,
and let $W$ denote the Weyl group.
We have the natural real structure $\gminih_{\real}\subset\gminih$.
Hence we have the subspace
$\sqrt{-1}\gminih_{\real}\subset\gminih$.
We have the $W$-invariant metric of $\gminih$,
which induces the distance $d$ of $\gminih/W$.
Let $B\bigl(\sqrt{-1}\gminih_{\real},\epsilon\bigr)$
denote the set of the points $x$ of $\gminih/W$
such that there exists a point $y\in\sqrt{-1}\gminih_{\real}/W$
satisfying $d(x,y)< \epsilon$.

Let $(P_K\subset P_G,\theta)$ be a $G$-harmonic bundle
on $\Delta^{\ast}$.
We have the expression $\theta=f\cdot dz/z$,
where $f$ is a holomorphic section of $\ad(P_G)$
on $\Delta^{\ast}$.
It induces the continuous map $[f]:\Delta^{\ast}\lrarr \gminih/W$.

\begin{df} \label{df;04.10.23.10}\mbox{{}}
\begin{itemize}
\item
 A $G$-harmonic bundle $(P_K\subset P_G,\theta)$ is called tame,
 if $[f]$ is bounded.
\item
 A tame $G$-harmonic bundle $(P_K\subset P_G,\theta)$ is called 
 pure imaginary,
 if for any $\epsilon>0$ there exists a positive number $r$
 such that 
 $\bigl[f(z)\bigr]\in B\bigl(\sqrt{-1}\gminih_{\real},\epsilon\bigr)$
 for any $|z|<r$.
\hfill\qed
\end{itemize}
\end{df}

\begin{lem} \label{lem;04.11.9.1}
Let $(P_K\subset P_G,\theta)$ be a harmonic bundle on $\Delta^{\ast}$.
The following conditions are equivalent.
\begin{enumerate}
\item \label{num;04.9.20.5}
 It is tame (pure imaginary).
\item \label{num;04.9.20.6}
 For any $\kappa:G\lrarr \GL(V)$ and any $K$-invariant metric of $V$,
 the induced harmonic bundle is tame (pure imaginary).
\item \label{num;04.9.20.7}
 For some immersive representation $\kappa:G\lrarr \GL(V)$
 and some $K$-invariant metric of $V$,
 the induced harmonic bundle is tame (pure imaginary).
\end{enumerate}
\end{lem}
\pf
The implications 
$\ref{num;04.9.20.5}\Longrightarrow
 \ref{num;04.9.20.6}\Longrightarrow
 \ref{num;04.9.20.7}$ are clear.
The implication
$\ref{num;04.9.20.7}\Longrightarrow\ref{num;04.9.20.5}$
follows from the injectivity of
$d\kappa:\gminig\lrarr \gl(V)$.
\hfill\qed

\vspace{.1in}

Let $X$ be a smooth projective variety,
and $D$ be a normal crossing divisor.

\begin{df}
A harmonic $G$-bundle
$(P_K\subset P_G,\theta)$ on $X-D$
is called tame (pure imaginary),
if the restriction $(P_K\subset P_G,\theta)_{|C\setminus D}$
is tame (pure imaginary)
for any curve $C\subset X$ which is transversal with $D$.
\hfill\qed
\end{df}

\begin{rem}
Tameness and pure imaginary property
are defined for principal $G$-Higgs bundles.
\hfill\qed
\end{rem}

\begin{rem}
Tameness and pure imaginary property are preserved
by pull back.
We also remark the curve test for usual tame harmonic bundles.
\hfill\qed
\end{rem}

Let us consider the case
where $G$ is a linear reductive group defined over $\real$,
with a maximal compact group $K$.
We have the complexification $G_{\cnum}$
with a maximal compact group $K_{\cnum}$
such that $K=K_{\cnum}\cap G$.
\begin{df}
Let $(P_G,\nabla)$ be a flat bundle.
A pluri-harmonic reduction $(P_K\subset P_G,\nabla)$ is called
a tame pure imaginary,
if the induced reduction $(P_{K_{\cnum}}\subset P_{G_{\cnum}},\nabla)$
is a tame pure imaginary.
\hfill\qed
\end{df}

\begin{lem}
Let $(P_K\subset P_G,\theta)$ be a harmonic bundle on $X-D$.
The following conditions are equivalent.
\begin{enumerate}
\item
 It is tame (pure imaginary).
\item
 For any $\kappa:G\lrarr \GL(V)$ and any $K$-invariant metric of $V$,
 the induced harmonic bundle is tame (pure imaginary).
\item
 There exist an immersive representation $\kappa:G\lrarr \GL(V)$
 and a $K$-invariant metric of $V$
 such that the induced harmonic bundle is tame (pure imaginary).
\hfill\qed
\end{enumerate}
\end{lem}

\section{Semisimplicity and
 Pluri-Harmonic Reduction}
\label{section;04.10.26.300}

\subsection{Preliminary}

Let $X$ be a smooth irreducible quasiprojective variety
with a base point $x$.
We put $\Gamma:=\pi_1(X,x)$ for simplicity
of the notation.
Recall the existence and the uniqueness of 
tame pure imaginary pluri-harmonic metric
(\cite{JZ2}, \cite{mochi3}),
which is called the Corlette-Jost-Zuo metric.
Let $(E,\nabla)$ be a semisimple flat bundle,
and let $\rho:\Gamma\lrarr \GL(E_{|x})$ denote
the corresponding representation.
We have the canonical decomposition of $E_{|x}$:
\[
 E_{|x}=\bigoplus_{\chi\in\Irrep(\Gamma)}E_{|x,\chi}.
\]
Here $\Irrep(\Gamma)$ denotes the set of irreducible
representations,
and $E_{|x,\chi}$ denotes a $\Gamma$-subspace of $E_{|x}$
isomorphic to $\chi^{\oplus m(\chi)}$.
Correspondingly,
we have the canonical decomposition of the flat bundle $(E,\nabla)$:
\[
 (E,\nabla)=\bigoplus_{\chi\in\Irrep(\Gamma)} E_{\chi}.
\]
The flat bundle $E_{\chi}$ is isomorphic to
a tensor product of a trivial bundle $\cnum^{m(\chi)}$
and a flat bundle $L_{\chi}$ whose monodromy is given by $\chi$.

\begin{lem} \mbox{{}} \label{lem;04.9.30.2}
\begin{itemize}
\item
There exists a Corlette-Jost -Zuo metric $h_{\chi}$ of $L_{\chi}$,
which is unique up to positive constant multiplication.
\item
Under the isomorphism
$(E,\nabla)\simeq \bigoplus_{\chi}L_{\chi}\otimes\cnum^{m(\chi)}$,
any Corlette-Jost-Zuo metric of $(V,\nabla)$ is of the following form:
\[
 \bigoplus_{\chi} h_{\chi}\otimes g_{\chi}.
\]
Here $g_{\chi}$ denote any hermitian metrics of $\cnum^{m(\chi)}$.
In other words,
the ambiguity of the Corlette-Jost-Zuo metrics
is a choice of hermitian metrics $g_{\chi}$ of $\cnum^{m(\chi)}$,
once we fix $h_{\chi}$.
\item
The decomposition of flat connection
$\nabla=\del+\delbar+\theta+\theta^{\dagger}$
is independent of a choice of $g_{\chi}$.
\end{itemize}
\end{lem}
\pf
The first claim is proved in \cite{JZ2}.
(See also \cite{mochi3}.)
The second claim easily follows from the proof of
the uniqueness result in \cite{mochi3}.
(See the argument of Proposition \ref{prop;05.9.8.100}).
The third claim follows from the second claim.
\hfill\qed

\vspace{.1in}

We also have the following lemma
(see \cite{sabbah2} or \cite{mochi3})
\begin{lem} \label{lem;04.9.30.1}
If there exists a Corlette-Jost-Zuo metric on a flat bundle $(E,\nabla)$,
then the flat bundle is semisimple.
\hfill\qed
\end{lem}

We have the involution $\chi\longmapsto \overline{\chi}$
on $\Irrep\bigl(\Gamma\bigr)$
such that $\chi\otimes_{\real}\cnum=\chi\oplus\overline{\chi}$.
If $\overline{\chi}=\chi$,
we have the real structure of $L_{\chi}$.
If $\overline{\chi}\neq \chi$,
we have the canonical real structure of
$L_{\chi}\otimes\cnum=L_{\chi}\oplus L_{\overline{\chi}}$.

Let us consider the case where 
a semisimple flat bundle
$(E,\nabla)$ has  the flat real structure $E_{\real}$
such that $E=E_{\real}\otimes_{\real}\cnum$.
Let $\iota:E\lrarr E$ denote the conjugate
with respect to $E_{\real}$.
Then $(E,\nabla)$ is isomorphic to the following:
\[
 \bigoplus_{\overline{\chi}=\chi }
 L_{\chi}\otimes\cnum^{m(\chi)}
\oplus
 \bigoplus_{\overline{\chi}\neq\chi }
 \bigl(L_{\chi}\oplus L_{\overline{\chi}}\bigr)
 \otimes \cnum^{m(\chi)}.
\]
The real structure of $(E,\nabla)$
is induced from the real structures of
$L_{\chi}$ $(\overline{\chi}=\chi)$
and $L_{\chi}\otimes\cnum$ $(\overline{\chi}\neq \chi)$.
For a hermitian metric $h$ of $E$,
the hermitian metric $\iota^{\ast}h$
is given by
$\iota^{\ast}h(u,v)=\overline{h(\iota(u),\iota(v))}$.
Then the following lemma is clear.
\begin{lem}
When $(E,\nabla)$ has a real structure,
there exists a Corlette-Jost-Zuo metric
of $(E,\nabla)$ which is invariant under the conjugation.
The ambiguity of the metric is a choice of
the metrics of the vector spaces $\cnum^{m(\chi)}$.
\hfill\qed
\end{lem}

\subsection{Pluri-harmonic reduction of the principal bundle
 associated with the monodromy group}

Let $G_0\subset \GL(E_{|x})$
denote the monodromy group $M(E,\nabla,x)$.
We obtain the principal $G_0$-bundle $P_{G_0}$
with the flat connection.
If the flat bundle $(E,\nabla)$ is semisimple,
we have a Corlette-Jost-Zuo metric $h$
of $(E,\nabla)$.
Let $U=U(E_{|x},h_{|x})$ denote the unitary group
of the metrized vector space $(E_{|x},h_{|x})$,
and we put $K_0:=G_0\cap U$.

\begin{lem}
 \label{lem;04.11.4.1}
$G_0$ is reductive,
and $K_0$ is a compact real form of $G_0$.
\end{lem}
\pf
The argument was given by Simpson (Lemma 4.4 in \cite{s5})
for a different purpose.
We reproduce it here with a minor change for our purpose.
We have the canonical decomposition
$T^{a,b}(E)=\bigoplus_{\chi\in\Irrep(\Gamma)}
 L_{\chi}\otimes\cnum^{m(a,b,\chi)}$.
The decomposition is orthogonal
with respect to the induced Corlette-Jost-Zuo metric
$T^{a,b}(h)$.
Namely, $T^{a,b}(h)$ is of the form
$\bigoplus_{\chi\in\Irrep(\Gamma)} h_{\chi}
 \otimes h(a,b,\chi)$,
where $h_{\chi}$ denotes a Corlette-Jost-Zuo metric
of $L_{\chi}$, and $h(a,b,\chi)$ denotes
hermitian metric of $\cnum^{m(a,b,\chi)}$.

For any $f\in \End(E_{|x})$,
let $f^{\dagger}$ denote the adjoint of $f$
with respect to $h_{|x}$.
For any $g\in \GL(E_{|x})$, we have the unique expression
$g=u\cdot \exp(y)$,
where $u\in U$ and $y=y^{\dagger}$.
The decomposition is compatible with 
tensor products and $g$-invariant orthogonal decompositions.
It follows that $T^{a,b}u$ and $T^{a,b}y$ preserves 
the components $L_{\chi|x}\otimes\cnum^{m(a,b,\chi)}$.
Namely, we have the decomposition
$T^{a,b}g=(\bigoplus T^{a,b}g)_{\chi}$,
$T^{a,b}u=(\bigoplus T^{a,b}u)_{\chi}$
and $T^{a,b}y=(\bigoplus T^{a,b}y)_{\chi}$.

Let $\kappa$ be an isometric automorphism
of $\bigl(\cnum^{m(a,b,\rho)},h(a,b,\chi)\bigr)$.
Then, $(T^{a,b}g)_{\chi}$ and
$\id_{L_{\chi|x}}\otimes \kappa$ are commutative.
Hence, $(T^{a,b}u)_{\chi}$
and $\id_{L_{\chi|x}}\otimes\kappa$ are commutative,
and thus $(T^{a,b}u)_{\chi}$ is induced by
the automorphism of $L_{\chi|x}$.
Similarly, $(T^{a,b}y)_{\chi}$ is induced by the endomorphism
of $L_{\chi|x}$.
Hence, 
$L_{\chi|x}\otimes H_{\chi}$ is preserved by
$(T^{a,b}u)_{\chi}$ and $(T^{a,b}y)_{\chi}$
for any subspace $H_{\chi}\subset \cnum^{m(a,b,\chi)}$.
Since any $G_0$-invariant subspace of
$T^{a,b}E_{|x}$ is of the form
$\bigoplus L_{\chi|x}\otimes H_{\chi}$,
we obtain
$u\in G_0\cap U=K_0$
and $y\in \gminig_0\subset End(E_{|x})$,
where $\gminig_0$ denotes the Lie subalgebra of $\End(E_{|x})$
corresponding to $G_0$.

Let $\tau:\GL(E_{|x})\lrarr \GL(E_{|x})$ be
the anti-holomorphic involution
such that $\tau(g)=(g^{\dagger})^{-1}$.
We obtain that $\tau(g)=u\cdot \exp\bigl(-y\bigr)$
is contained in $G_0$.
Namely, $\tau$ gives the real structure of $G_0$.
Since we have the decomposition
$g=u\cdot \exp(y)$ for any $g\in G_0$,
$K_0$ intersects with any connected components of $G_0$.
Let $G_0^0$ denote the connected component of $G_0$
containing the unit element.
It is easy to see that $K_0\cap G_0^0$ is
maximal compact in $G_0^0$,
and hence $K_0$ is maximal compact of $G_0$.
Since $K_0\cap G_0^0$ is the fixed point set of $\tau_{|G_0^0}$,
we obtain that $K_0^0$ is a compact real form of $G_0^0$.
Thus $K_0$ is a compact real form of $G_0$.
Since $K_0$ is maximal compact, $G_0$ is reductive.
\hfill\qed

\vspace{.1in}

Let us consider the case where
$(E,\nabla)$ has the real structure.
We have the real parts $E_{\real\,|\,x}\subset E_{|x}$
and  $G_{0\real}:=G_0\cap \GL(E_{\real\,|\,x})$.
We take a Corlette-Jost-Zuo metric
of $h$ which is invariant under the conjugation $\iota$.
We put $K_{0\real}=G_{0\real}\cap K_0=G_{0\real}\cap U$.
The map $\iota$ induces the real endomorphism
of $\End(E_{|x})$
given by $\iota(f)=\iota\circ f\circ\iota$.

\begin{lem} 
\label{lem;04.11.4.10}
$K_{0\,\real}$ is maximal compact in $G_{0\,\real}$.
\end{lem}
\pf
We use the notation in the proof of Lemma \ref{lem;04.11.4.1}.
Since $h_{|x}$ is invariant under the conjugation $\iota$,
$U$ is stable under $\iota$,
and $\tau$ and $\iota$ are commutative.
Let $g$ be an element of $G_{0\,\real}$.
We have the decomposition
$g=u\cdot \exp(y)$ as in the proof of Lemma \ref{lem;04.11.4.1},
where $u$ denotes an element of $K_0$
and $y$ denotes an element of $\gminig_0$
such that $y^{\dagger}=y$.
Since $\iota(g)=g$,
we have
$\iota(u)\cdot \exp\bigl(\iota(y)\bigr)=u\cdot \exp(y)$.
Since we have $\iota(u)\in \iota(U)=U$
and $(\iota(y))^{\dagger}=\iota(y^{\dagger})=-\iota(y)$,
we obtain $\iota(u)=u$ and $\iota(y)=y$.
Namely $u\in K_{0\real}$ and $y\in \gminig_{0\real}$.
Then we can show $K_{0\real}$ is maximal compact in $G_{0\real}$,
by an argument similar to the proof of
Lemma \ref{lem;04.11.4.1}.
\hfill\qed

\vspace{.1in}

\begin{prop}
\label{prop;04.9.30.5}
Assume that $(E,\nabla)$ is semisimple.
Then there exists
the unique tame pure imaginary pluri-harmonic
reduction $P_{K_0}\subset P_{G_0}$.
Assume  $(E,\nabla)$ has the flat real structure, moreover.
Then, it is induced from the pluri-harmonic reduction
of $P_{G_{0\real}}$.
\end{prop}
\pf
Let $h$ be a Corlette-Jost-Zuo metric of $(E,\nabla)$.
For any point $z\in X$,
let $M(E,\nabla,z)$ denote the monodromy group at $z$,
and $U(E_{|z},h_{|z})$ denote the unitary group
of $E_{|z}$ with the metric $h_{|z}$.
Then the intersection $M(E,\nabla,z)\cap U(E_{|z},h_{|z})$
is a maximal compact subgroup
of $M(E,\nabla,z)$, due to Lemma \ref{lem;04.11.4.1}.
Hence they give the reduction $P_{K_0}\subset P_{G_0}$,
which is pluri-harmonic.
By using a similar argument and
Lemma \ref{lem;04.11.4.10},
we obtain the compatibility with the real structure,
if $(E,\nabla)$ has the flat real structure.
The uniqueness of the pluri-harmonic reduction
follows from the uniqueness result
in Lemma \ref{lem;04.9.30.2}.
Hence we are done.
\hfill\qed

\subsection{Characterization of the existence of pluri-harmonic reduction}

Let $G$ be a linear reductive algebraic group over $\cnum$ or $\real$.
Let  $\widetilde{X}$ be a universal covering of $X$.
The following corollary immediately follows from
Proposition \ref{prop;04.9.30.5}.
\begin{cor}
\label{cor;04.10.26.1}
Let $P_G$ be a flat $G$-principal bundle over $X$.
Assume that the image of 
the induced representation $\Gamma\lrarr G$ is Zariski dense in $G$.
Then there exists the unique tame pure imaginary pluri-harmonic
reduction of $P_G$.
Correspondingly, we obtain 
the $\Gamma$-equivariant pluri-harmonic map
$\widetilde{X}\lrarr G/K$.
\hfill\qed
\end{cor}

\begin{prop}
\label{prop;04.10.26.3}
Let $P_G$ be a flat $G$-bundle on $X$.
The monodromy group $G_0$ is reductive
if and only if 
there exists a tame pure imaginary pluri-harmonic reduction
$P_K\subset P_G$.
If such a reduction exists,
the decomposition $\nabla=\nabla_K+(\theta+\theta^{\dagger})$
does not depend on a choice of
a pluri-harmonic reduction $P_K\subset P_G$,
and there is the corresponding $\Gamma$-equivariant
pluri-harmonic map $\widetilde{X}\lrarr G/K$.
\end{prop}
\pf
If a pluri-harmonic reduction exists,
the monodromy group is reductive
due to Lemma \ref{lem;04.11.9.1} and
Lemma \ref{lem;04.11.4.1}.
Assume $G_0$ is reductive.
Let $K_0$ be a maximal compact group of $G_0$.
Then we have the unique tame pure imaginary
pluri-harmonic reduction $P_{K_0}\subset P_{G_0}$.
We take $K$ such as $K\cap G_0=K_0$.
Then the pluri-harmonic reduction $P_{K}\subset P_{G}$ is induced,
and thus the first claim is proved.
The second claim is clear.
\hfill\qed